\documentclass[10pt]{amsart}
\usepackage{amsmath}
\usepackage{amsxtra}
\usepackage{amscd}
\usepackage{amsthm}
\usepackage{amsfonts}
\usepackage{amssymb}
\usepackage{eucal}

\newtheorem{cor}[subsection]{Corollary}
\newtheorem{lem}[subsection]{Lemma}
\newtheorem{prop}[subsection]{Proposition}

\newtheorem{conj}[subsection]{Conjecture}
\newtheorem{thm}[subsection]{Theorem}
\newtheorem{defn}[subsection]{Definition}

\theoremstyle{definition}

\theoremstyle{remark}

\newcommand{\thmref}[1]{Theorem~\ref{#1}}
\newcommand{\secref}[1]{Sect.~\ref{#1}}
\newcommand{\lemref}[1]{Lemma~\ref{#1}}
\newcommand{\propref}[1]{Proposition~\ref{#1}}
\newcommand{\corref}[1]{Corollary~\ref{#1}}
\newcommand{\conjref}[1]{Conjecture~\ref{#1}}

\newcommand{\nc}{\newcommand}
\nc{\renc}{\renewcommand}
\nc{\ssec}{\subsection}
\nc{\sssec}{\subsubsection}
\nc{\on}{\operatorname}

\nc\ol{\overline}
\nc\wt{\widetilde}
\nc\tboxtimes{\wt{\boxtimes}}
\nc{\alp}{\alpha}

\nc{\ZZ}{{\mathbb Z}}
\nc{\NN}{{\mathbb N}}
\nc{\CC}{{\mathbb C}}
\nc{\OO}{{\mathbb O}}
\renc{\SS}{{\mathbb S}}
\nc{\DD}{{\mathbb D}}
\nc{\GG}{{\mathbb G}}
\renewcommand{\AA}{{\mathbb A}}
\nc{\Fq}{{\mathbb F}_q}
\nc{\Fqb}{\ol{{\mathbb F}_q}}
\nc{\Ql}{\ol{{\mathbb Q}_\ell}}
\nc{\id}{\text{id}}
\nc\X{\mathcal X}

\nc{\Hom}{\on{Hom}}
\nc{\Lie}{\on{Lie}}
\nc{\Loc}{\on{Loc}}
\nc{\Pic}{\on{Pic}}
\nc{\Bun}{\on{Bun}}
\nc{\IC}{\on{IC}}
\nc{\Aut}{\on{Aut}}
\nc{\rk}{\on{rk}}
\nc{\Sh}{\on{Sh}}
\nc{\Perv}{\on{Perv}}
\nc{\pos}{{\on{pos}}}
\nc{\Conv}{\on{Conv}}
\nc{\Sph}{\on{Sph}}
\nc{\Sym}{\on{Sym}}
\nc{\BunBb}{\overline{\Bun}_B}
\nc{\Buno}{\overset{o}{\Bun}}
\nc{\BunPb}{{\overline{\Bun}_P}}
\nc{\BunBM}{\overline{\Bun}_{B(M)}}
\nc{\BunPbw}{{\widetilde{\Bun}_P}}
\nc{\BunBP}{\widetilde{\Bun}_{B,P}}
\nc{\GUb}{\overline{G/U}}
\nc{\GUPb}{\overline{G/U(P)}}

\nc{\Hhom}{\underline{\on{Hom}}}
\nc\syminfty{\on{Sym}^{\infty}}
\nc\lal{\ol{\lambda}}
\nc\xl{\ol{x}}
\nc\thl{\ol{\theta}}
\nc\nul{\ol{\nu}}
\nc\mul{\ol{\mu}}
\nc\Sum\Sigma
\nc{\oX}{\overset{o}{X}{}}
\nc{\hl}{\overset{\leftarrow}h}
\nc{\hr}{\overset{\rightarrow}h}
\nc{\M}{{\mathcal M}}
\nc{\N}{{\mathcal N}}
\nc{\F}{{\mathcal F}}
\nc{\D}{{\mathcal D}}
\nc{\Q}{{\mathcal Q}}
\nc{\Y}{{\mathcal Y}}
\nc{\G}{{\mathcal G}}
\nc{\E}{{\mathcal E}}
\nc{\CalC}{{\mathcal C}}
\nc\Dh{\widehat{\D}}
\renewcommand{\O}{{\mathcal O}}
\nc{\C}{{\mathcal C}}
\nc{\K}{{\mathcal K}}
\renewcommand{\H}{{\mathcal H}}

\nc{\T}{{\mathcal T}}
\nc{\V}{{\mathcal V}}
\renc{\P}{{\mathcal P}}
\nc{\A}{{\mathcal A}}
\nc{\B}{{\mathcal B}}
\nc{\U}{{\mathcal U}}
\renewcommand{\L}{{\mathcal L}}
\nc{\Gr}{\on{Gr}}

\nc{\frn}{{\check{\mathfrak u}(P)}}
\nc{\p}{\mathfrak p}
\nc{\q}{\mathfrak q}
\nc\f{{\mathfrak f}}

\nc{\qo}{{\mathfrak q}}
\nc{\po}{{\mathfrak p}}
\nc{\s}{{\mathfrak s}}
\nc\w{\text{w}}

\nc\mathi\iota
\nc\Spec{\on{Spec}}
\nc\Mod{\on{Mod}}
\nc{\tw}{\widetilde{\mathfrak t}}
\nc{\pw}{\widetilde{\mathfrak p}}
\nc{\qw}{\widetilde{\mathfrak q}}
\nc{\jw}{\widetilde j}

\nc{\grb}{\overline{\Gr}}
\nc{\I}{\mathcal I}

\nc{\lambdach}{{\check\lambda}}
\nc{\Lambdach}{{\check\Lambda}{}}
\nc{\much}{{\check\mu}}
\nc{\omegach}{{\check\omega}}
\nc{\nuch}{{\check\nu}}
\nc{\etach}{{\check\eta}}
\nc{\alphach}{{\check\alpha}}
\nc{\betach}{{\check\beta}}
\nc{\rhoch}{{\check\rho}}
\nc{\ch}{{\check h}}

\nc{\Hb}{\overline{\H}}

\nc{\BA}{{\mathbb{A}}}
\nc{\BC}{{\mathbb{C}}}
\nc{\BG}{{\mathbb{G}}}
\nc{\BM}{{\mathbb{M}}}
\nc{\BN}{{\mathbb{N}}}
\nc{\BP}{{\mathbb{P}}}
\nc{\BR}{{\mathbb{R}}}
\nc{\BZ}{{\mathbb{Z}}}
\nc{\BS}{{\mathbb{S}}}

\nc{\CA}{{\mathcal{A}}}
\nc{\CB}{{\mathcal{B}}}

\nc{\CE}{{\mathcal{E}}}
\nc{\CF}{{\mathcal{F}}}
\nc{\CG}{{\mathcal{G}}}
\nc{\CL}{{\mathcal{L}}}
\nc{\CM}{{\mathcal{M}}}
\nc{\CN}{{\mathcal{N}}}
\nc{\CK}{{\mathcal{K}}}
\nc{\CO}{{\mathcal{O}}}
\nc{\CP}{{\mathcal{P}}}
\nc{\CQ}{{\mathcal{Q}}}
\nc{\CR}{{\mathcal{R}}}
\nc{\CS}{{\mathcal{S}}}
\nc{\CT}{{\mathcal{T}}}
\nc{\CU}{{\mathcal{U}}}
\nc{\CV}{{\mathcal{V}}}
\nc{\CW}{{\mathcal{W}}}
\nc{\CZ}{{\mathcal{Z}}}
\nc{\CI}{{\mathcal{I}}}

\nc{\cM}{{\check{\mathcal M}}{}}
\nc{\csM}{{\check{\mathcal A}}{}}
\nc{\oM}{{\overset{\circ}{\mathcal M}}{}}
\nc{\obM}{{\overset{\circ}{\mathbf M}}{}}
\nc{\oCA}{{\overset{\circ}{\mathcal A}}{}}
\nc{\obA}{{\overset{\circ}{\mathbf A}}{}}
\nc{\ooM}{{\overset{\circ}{M}}{}}
\nc{\osM}{{\overset{\circ}{\mathsf M}}{}}
\nc{\vM}{{\overset{\bullet}{\mathcal M}}{}}
\nc{\nM}{{\underset{\bullet}{\mathcal M}}{}}
\nc{\oD}{{\overset{\circ}{\mathcal D}}{}}
\nc{\obD}{{\overset{\circ}{\mathbf D}}{}}
\nc{\oA}{{\overset{\circ}{\mathbb A}}{}}
\nc{\op}{{\overset{\bullet}{\mathbf p}}{}}
\nc{\cp}{{\overset{\circ}{\mathbf p}}{}}
\nc{\oU}{{\overset{\bullet}{\mathcal U}}{}}
\nc{\oZ}{{\overset{\circ}{\mathcal Z}}{}}
\nc{\ofZ}{{\overset{\circ}{\mathfrak Z}}{}}
\nc{\oF}{{\overset{\circ}{\fF}}}

\nc{\fa}{{\mathfrak{a}}}
\nc{\fb}{{\mathfrak{b}}}
\nc{\fg}{{\mathfrak{g}}}
\nc{\fgl}{{\mathfrak{gl}}}
\nc{\fh}{{\mathfrak{h}}}
\nc{\fj}{{\mathfrak{j}}}
\nc{\fm}{{\mathfrak{m}}}
\nc{\fn}{{\mathfrak{n}}}
\nc{\fu}{{\mathfrak{u}}}
\nc{\fp}{{\mathfrak{p}}}
\nc{\fr}{{\mathfrak{r}}}
\nc{\fs}{{\mathfrak{s}}}
\nc{\fsl}{{\mathfrak{sl}}}
\nc{\hsl}{{\widehat{\mathfrak{sl}}}}
\nc{\hgl}{{\widehat{\mathfrak{gl}}}}
\nc{\hg}{{\widehat{\mathfrak{g}}}}
\nc{\chg}{{\widehat{\mathfrak{g}}}{}^\vee}
\nc{\hn}{{\widehat{\mathfrak{n}}}}
\nc{\chn}{{\widehat{\mathfrak{n}}}{}^\vee}

\nc{\fA}{{\mathfrak{A}}}
\nc{\fB}{{\mathfrak{B}}}
\nc{\fD}{{\mathfrak{D}}}
\nc{\fE}{{\mathfrak{E}}}
\nc{\fF}{{\mathfrak{F}}}
\nc{\fG}{{\mathfrak{G}}}
\nc{\fK}{{\mathfrak{K}}}
\nc{\fL}{{\mathfrak{L}}}
\nc{\fM}{{\mathfrak{M}}}
\nc{\fN}{{\mathfrak{N}}}
\nc{\fP}{{\mathfrak{P}}}
\nc{\fU}{{\mathfrak{U}}}
\nc{\fV}{{\mathfrak{V}}}
\nc{\fZ}{{\mathfrak{Z}}}

\nc{\bb}{{\mathbf{b}}}
\nc{\bc}{{\mathbf{c}}}
\nc{\bd}{{\mathbf{d}}}
\nc{\be}{{\mathbf{e}}}
\nc{\bj}{{\mathbf{j}}}
\nc{\bn}{{\mathbf{n}}}
\nc{\bp}{{\mathbf{p}}}
\nc{\bq}{{\mathbf{q}}}
\nc{\bu}{{\mathbf{u}}}
\nc{\bv}{{\mathbf{v}}}
\nc{\bx}{{\mathbf{x}}}
\nc{\bs}{{\mathbf{s}}}
\nc{\by}{{\mathbf{y}}}
\nc{\bw}{{\mathbf{w}}}
\nc{\bA}{{\mathbf{A}}}
\nc{\bK}{{\mathbf{K}}}
\nc{\bB}{{\mathbf{B}}}
\nc{\bC}{{\mathbf{C}}}
\nc{\bG}{{\mathbf{G}}}
\nc{\bD}{{\mathbf{D}}}
\nc{\bH}{{\mathbf{H}}}
\nc{\bM}{{\mathbf{M}}}
\nc{\bN}{{\mathbf{N}}}
\nc{\bV}{{\mathbf{V}}}
\nc{\bW}{{\mathbf{W}}}
\nc{\bX}{{\mathbf{X}}}
\nc{\bZ}{{\mathbf{Z}}}
\nc{\bS}{{\mathbf{S}}}

\nc{\sA}{{\mathsf{A}}}
\nc{\sB}{{\mathsf{B}}}
\nc{\sC}{{\mathsf{C}}}
\nc{\sD}{{\mathsf{D}}}
\nc{\sF}{{\mathsf{F}}}
\nc{\sK}{{\mathsf{K}}}
\nc{\sM}{{\mathsf{M}}}
\nc{\sO}{{\mathsf{O}}}
\nc{\sQ}{{\mathsf{Q}}}
\nc{\sP}{{\mathsf{P}}}
\nc{\sZ}{{\mathsf{Z}}}
\nc{\sfp}{{\mathsf{p}}}
\nc{\sr}{{\mathsf{r}}}
\nc{\sg}{{\mathsf{g}}}
\nc{\sff}{{\mathsf{f}}}
\nc{\sfb}{{\mathsf{b}}}
\nc{\sfc}{{\mathsf{c}}}
\nc{\sd}{{\mathsf{d}}}

\nc{\BK}{{\bar{K}}}

\nc{\tA}{{\widetilde{\mathbf{A}}}}
\nc{\tB}{{\widetilde{\mathcal{B}}}}
\nc{\tg}{{\widetilde{\mathfrak{g}}}}
\nc{\tG}{{\widetilde{G}}}
\nc{\TM}{{\widetilde{\mathbb{M}}}{}}
\nc{\tO}{{\widetilde{\mathsf{O}}}{}}
\nc{\tU}{{\widetilde{\mathfrak{U}}}{}}
\nc{\TZ}{{\tilde{Z}}}
\nc{\tx}{{\tilde{x}}}
\nc{\tbv}{{\tilde{\bv}}}
\nc{\tfP}{{\widetilde{\mathfrak{P}}}{}}
\nc{\tz}{{\tilde{\zeta}}}
\nc{\tmu}{{\tilde{\mu}}}

\nc{\urho}{\underline{\rho}}
\nc{\uB}{\underline{B}}
\nc{\uC}{{\underline{\mathbb{C}}}}
\nc{\ui}{\underline{i}}
\nc{\uj}{\underline{j}}
\nc{\ofP}{{\overline{\mathfrak{P}}}}
\nc{\oB}{{\overline{\mathcal{B}}}}
\nc{\og}{{\overline{\mathfrak{g}}}}
\nc{\oI}{{\overline{I}}}

\nc{\eps}{\varepsilon}
\nc{\hrho}{{\hat{\rho}}}

\nc{\one}{{\mathbf{1}}}
\nc{\two}{{\mathbf{t}}}

\nc{\Rep}{{\mathop{\operatorname{\rm Rep}}}}
\nc{\Tot}{{\mathop{\operatorname{\rm Tot}}}}
\nc{\Ker}{{\mathop{\operatorname{\rm Ker}}}}
\nc{\Hilb}{{\mathop{\operatorname{\rm Hilb}}}}
\nc{\End}{{\mathop{\operatorname{\rm End}}}}
\nc{\Ext}{{\mathop{\operatorname{\rm Ext}}}}
\nc{\CHom}{{\mathop{\operatorname{{\mathcal{H}}\it om}}}}
\nc{\GL}{{\mathop{\operatorname{\rm GL}}}}
\nc{\gr}{{\mathop{\operatorname{\rm gr}}}}
\nc{\Id}{{\mathop{\operatorname{\rm Id}}}}
\nc{\de}{{\mathop{\operatorname{\rm def}}}}
\nc{\length}{{\mathop{\operatorname{\rm length}}}}
\nc{\supp}{{\mathop{\operatorname{\rm supp}}}}

\nc{\Cliff}{{\mathsf{Cliff}}}
\nc{\Fl}{\on{Fl}}
\nc{\Fib}{{\mathsf{Fib}}}
\nc{\Coh}{{\mathsf{Coh}}}
\nc{\FCoh}{{\mathsf{FCoh}}}

\nc{\reg}{{\text{\rm reg}}}

\nc{\cplus}{{\mathbf{C}_+}}
\nc{\cminus}{{\mathbf{C}_-}}
\nc{\cthree}{{\mathbf{C}_*}}
\nc{\Qbar}{{\bar{Q}}}

\nc{\bh}{{\bar{h}}}
\nc{\bOmega}{{\overline{\Omega}}}

\nc{\seq}[1]{\stackrel{#1}{\sim}}

\begin{document}

\title{Uhlenbeck spaces via affine Lie algebras}

\author{Alexander Braverman, Michael Finkelberg, and Dennis Gaitsgory}

\address
{\newline
A.B.: Dept. of Math., Brown Univ., Providence, RI 02912, USA;\newline
M.F.: Independent Moscow Univ., 11 Bolshoj Vlasjevskij per.,
Moscow 119002, Russia; \newline
D.G.: Dept. of Math., The Univ. of Chicago, Chicago, IL 60637, USA}

\date{October 2004}

\email{\newline
braval@math.brown.edu; fnklberg@mccme.ru; gaitsgde@math.uchicago.edu}

\begin{abstract}

Let $G$ be an almost simple simply connected group over $\BC$,
and let $\Bun^a_G(\BP^2,\BP^1)$ be the moduli scheme of principal $G$-bundles
on the projective plane $\BP^2$, of second Chern class $a$,
trivialized along a line $\BP^1\subset \BP^2$. 

We define the Uhlenbeck compactification $\fU^a_G$ of
$\Bun^a_G(\BP^2,\BP^1)$, which classifies, roughly, pairs
$(\F_G,D)$, where $D$ is a $0$-cycle on $\BA^2=\BP^2-\BP^1$ of
degree $b$, and $\F_G$ is a point of $\Bun^{a-b}_G(\BP^2,\BP^1)$,
for varying $b$.

In addition, we calculate the stalks of the Intersection Cohomology
sheaf of $\fU^a_G$. To do that we give a geometric realization
of Kashiwara's crystals for affine Kac-Moody algebras.

\end{abstract}

\maketitle

\section*{Introduction}   \label{intr}

\ssec{}

Let $G$ be an almost simple simply connected group over $\BC$, with
Lie algebra $\fg$, and let $\bS$ be a smooth projective surface.

Let us denote by $\Bun^a_G(\bS)$ the moduli space (stack) of principal
$G$-bundles on $\bS$ of second Chern class $a$. It is easy to see
that $\Bun^a_G(\bS)$ cannot be compact, and the source of the 
non-compactness can be explained as follows:

By checking the valuative criterion of properness, we arrive to the 
following situation: we are given a $G$-bundle $\F_G$ on a
3-dimensional variety $\X$ defined away from a point, and we would
like to extend it to the entire $\X$. However, such an extension does
not always exist, and the obstruction is given by a positive integer,
which one can think of as the second Chern class of the
restriction of $\F_G$ to a suitable $4$-sphere corresponding to the 
point $x$.

However, this immediately suggests how a compactification of
$\Bun^a_G(\bS)$ could look like: it should be a union 
\begin{equation} \label{Uhlenbeck?}
\underset{b\in \BN}\bigcup\, \Bun^{a-b}_G(\bS)\times \on{Sym}^b(\bS).
\end{equation}

In the differential-geometric framework of moduli spaces of $K$-instantons
on Riemannian 4-manifolds (where $K$ is the maximal compact subgroup of $G$)
such a compactification was introduced in the pioneering work ~\cite{u}.
Therefore, we shall call its algebro-geometric version the 
Uhlenbeck space, and denote it by $\fU^a_G(\bS)$.

\medskip

Unfortunately, one still does not know how to construct the spaces
$\fU^a_G(\bS)$ for a general group $G$ and an arbitrary surface $\bS$.
More precisely, one would like to formulate a moduli problem, to
which $\fU^a_G(\bS)$ would be the answer, and so far this is not known.
In this formulation the question of constructing the Uhlenbeck spaces
has been posed (to the best of our knowledge) by V.~Ginzburg. He and
V.~Baranovsky (cf. \cite{BaGi}) have made the first attempts to
solve it, as well as indicated the approach adopted in this paper.

\medskip

A significant simplification occurs for $G=SL_n$.
Let us note that when $G=SL_n$, there exists another natural
compactification of the stack $\Bun^a_n(\bS):=\Bun^a_{SL_n}(\bS)$,
by torsion-free sheaves of generic rank $n$ and of second Chern
class $a$, called the Gieseker compactification, which in this paper
we will denote by $\wt{\fN}^a_n(\bS)$. One expects that there exists
a proper map $\sff:\wt{\fN}^{a}_n(\bS)\to \fU^{a}_{SL_n}(\bS)$, 
described as follows:

A torsion-free sheaf $\M$ embeds into a short exact sequence
$$0\to \M\to \M'\to \M_0\to 0,$$
where $\M'$ is a vector bundle (called the saturation of $\M$), and
$\M_0$ is a finite-length sheaf. The map $\sf$ should send a point
of $\wt{\fN}^a_n(\bS)$ corresponding to $\M$ to the pair 
$(\M',\on{cycle}(\M_0))\in \Bun^{a-b}_n(\bS)\times \on{Sym}^b(\bS)$,
where $b$ is the length of $\M_0$, and $\on{cycle}(\M_0)$ is the cycle
of $\M_0$. In other words, the map $\sf$ must ``collapse" the information
of the quotient $\M'\to \M_0$ to just the information of the length of
$\M_0$ at various points of $\bS$.

Since the spaces $\wt{\fN}^a_n(\bS)$, being a solution of a moduli problem,
are easy to construct, one may attempt to construct the Uhlenbeck spaces
$\fU^a_{SL_n}(\bS)$ by constructing an explicit blow down of the 
Gieseker spaces $\wt{\fN}^a_n(\bS)$. This has indeed been performed in 
the works of J.~Li (cf. \cite{Li}) and J.~W.~ Morgan (cf. \cite{Mo}).

The problem simplifies even further, when we put $\bS=\BP^2$, the 
projective plane, and consider bundles trivialized along a fixed
line $\BP^1\subset \BP^2$.
In this case, the sought-for space $\fU^a_n(\bS)$
has been constructed by S.~Donaldson (cf. Chapter 3 of \cite{dk}) and thoroughly
studied by H.~Nakajima (cf. e.g. \cite{n1}) in his works on quiver
varieties.

\medskip

In the present paper, we will consider the case of an arbitrary group
$G$, but the surface equal to $\BP^2$ (and we will be interested in
bundles trivialized along $\BP^1\subset \BP^2$, i.e., we will work 
in the Donaldson-Nakajima set-up.)

We will be able to construct the Uhlenbeck spaces $\fU^a_G$, but only
up to nilpotents. I.e., we will have several definitions, two of
which admit modular descriptions, and which produce the same answer
on the level of reduced schemes. We do not know, whether the resulting
schemes actually coincide when we take the nilpotents into account.
And neither do we know whether the resulting reduced scheme is normal.

We should say that the problem of constructing the Uhlenbeck spaces
can be posed over a base field of any characteristic. However, the
proof of one of the main results of this paper, 
\thmref{equivalence of definitions}, which insures that our spaces
$\fU^a_G$ are invariantly defined, uses the char.=0 assumption.
It is quite possible that in order to treat the char.=p case,
one needs a finer analysis.

\ssec{}

The construction of $\fU^a_G$ used in this paper is a simplification 
of a suggestion of Drinfeld's (the latter potentially  works for
an arbitrary surface $\bS$).

We are trying to express points of $\fU^a_G$ (one may call them quasi-bundles) 
by replacing the original problem for the surface $\BP^2$, or rather
for a rationally equivalent surface $\BP^1\times \BP^1$, by another
problem for the curve $\BP^1$.

As a motivation, let us consider the following simpler situation. Let
$\M^0$ be the trivial rank-$2$ bundle on a curve $\bC$. A flag in $\M^0$
is by definition a line subbundle $\CL\subset \M^0$, or equivalently
a map from $\bC$ to the flag variety of $GL_2$, i.e. $\BP^1$. 

However, there is a natural generalization of a notion of a flag, also suggested
by Drinfeld: instead of line {\it subbundles} we may consider all pairs
$(\CL,\kappa:\CL\to \M^0)$, where $\CL$ is still a line bundle, but
$\kappa$ need not be a bundle map, just an embedding of coherent
sheaves. We define a {\it quasi-map} $\bC\to \BP^1$ (generalizing 
the notion of a {\it map}) to be such a pair $(\CL,\kappa)$.

In fact, one can introduce the notion of a quasi-map from a curve
(or any projective variety) $\bC$ to another projective variety $\CT$.
When $\CT$ is the flag variety of a semi-simple group $G$, the corresponding
quasi-maps spaces have been studied in \cite{ffkm}, \cite{bg1},
\cite{bgfm}.

\medskip

Our construction of the Uhlenbeck space is based on considering quasi-maps
from $\BP^1$ (thought of as a ``horizontal" component of $\BP^1\times \BP^1$)
to various flag varieties associated to the loop group of $G$, among
them the most important are Kashiwara's thick Grassmannian, 
and the Beilinson-Drinfeld Grassmannian $\Gr^{BD}_G$.

\medskip

The spaces of maps and quasi-maps from a projective curve $\bC$ to Kashiwara's
flag schemes are of independent interest and have been another major source
of motivation for us.

In ~\cite{ffkm}, ~\cite{fkmm} it was shown that if one considers the space of
(based) maps, of multi-degree $\mu$, from $\bC$ to the flag variety 
of a finite-dimensional group $G$, one obtains an affine scheme,
which we denote by $\on{Maps}^\mu(\bC,\CB_\fg)$, endowed with
a symplectic structure, and which admits a Lagrangian projection 
to the space of coloured divisors on $\bC$, denoted $\bC^\mu$. 

Moreover, the irreducible components of the {\it central fiber}
$\fF^\mu_\fg$ of this projection (i.e. the fiber over $\mu\cdot \bc\in
\bC^\mu$, for some point $\bc\in \bC$) form in a natural way a basis 
for the $\mu$-weight piece of $U(\check \fn)$, where $\check \fg\supset \check\fn$ 
are the Langlands dual Lie algebra and its maximal nilpotent
subalgebra, respectively.

One may wonder if this picture can be generalized for an arbitrary
Kac-Moody $\fg'$, instead of the finite-dimensional algebra $\fg$. 
We discuss such a generalization in Parts I and IV of this paper.
We formulate \conjref{lagrange}, which is subsequently
proven for $\fg'$ affine (and, of course, finite), which allows to
define on the set of irreducible components of the central fibers 
$\underset{\mu}\cup\, \fF^\mu_{\fg'}$ a structure of Kashiwara's crystal, 
and thereby link it to the combinatorics of the Langlands dual Lie
algebra $\check \fg'$.
 
In particular, when $\fg'$ is affine, the space 
$\on{Maps}^\mu(\bC,\CB_{\fg'})$ turns out to be closely related 
to the space of bundles on $\BP^2$, and the space of quasi-maps 
$\on{QMaps}^\mu(\bC,\CB_{\fg'})$ --- to the corresponding Uhlenbeck
space. The relation between the irreducible components of $\fF^\mu$
and the Lie algebra $\check\fg'$ mentioned above, allows us to 
explicitly compute the Intersection Cohomology sheaf on $\fU^a_G$, 
and express it in terms of $\check\fg'$ (in this case $\fg'=\fg_{aff}$,
the affinization of $\fg$.)

\ssec{}

The two main results of this paper are construction of the scheme $\fU^a_G$,
so that it has the stratification as in \eqref{Uhlenbeck?} (\thmref{stratification}),
and the explicit description of the Intersection Cohomology sheaf of $\fU^a_G$
(\thmref{ic uhlenbeck}).
Let us now explain the logical structure of the paper and the main
points of each of the parts.

\medskip

Part I is mostly devoted to the preliminaries. In Section 1 
we introduce the notion of a quasi-map (in rather general
circumstances) and prove some of its basic properties. The reader
familiar with any of the works \cite{ffkm}, \cite{bg1} or \cite{bgfm}
may skip Section 1 and return for proofs of statements referred to 
in the subsequent sections.

In Section 2 we collect some facts about Kashiwara's flag schemes
$\CG_{\fg',\fp}$ for a general Kac-Moody Lie algebra $\fg'$, and
study the quasi-maps' spaces from a curve $\bC$ to $\CG_{\fg',\fp}$.
In the main body of the paper we will only use the cases when $\fg'$
is the affine algebra $\fg_{aff}$, or the initial finite-dimensional
Lie algebra $\fg$.

In Section 3 we collect some basic facts about $G$-bundles on a
surface $\BP^1\times \BP^1$ trivialized along the divisor of infinity.

\medskip

In Part II we introduce the Uhlenbeck space $\fU^a_G$ and study
its properties. In Section 4 we give three definitions of
$\fU^a_G$, of which two are almost immediately equivalent,
and the equivalence with the (most invariant) third one is
established later.

Section 5 contains a proof of \thmref{Uhlenbeck=Nakajima}, 
(conjectured in \cite{fgk}) about the isomorphism of our 
definition of the Uhlenbeck space $\fU^a_{SL_n}$ with Donaldson's one.
The proof follows the ideas indicated by Drinfeld.

In Section 6 we prove some additional functoriality properties
of $\fU^a_G$, which, combined with the \thmref{Uhlenbeck=Nakajima}
of the previous section, yields \thmref{equivalence of definitions}
about the equivalence of all three definitions from Section 4.
In addition, in Section 6 we establish the factorization property
of $\fU^a_G$ with respect to the projection on the symmetric power
$\BA^{(a)}$ of the ``horizontal" line, which will be one of the
principal technical tools in the study of the geometry of $\fU^a_G$,
and in particular, for the computation of the IC sheaf. 

In Section 7 we prove \thmref{stratification} saying that 
$\fU^a_G$ indeed has a stratification as in \eqref{Uhlenbeck?}.
In addition, we formulate \thmref{ic uhlenbeck} describing the
stalks of the IC sheaf on the various strata.

In Section 8 we present two moduli problems, whose solutions
provide two more variants of the definition of $\fU^a_G$, and which
coincide with the original one on the level of reduced schemes.

\medskip

In Part III we define ``parabolic" version of Uhlenbeck spaces,
$\fU^\theta_{G,P}$ and $\wt{\fU}^\theta_{G,P}$, which are two
different compactifications of the space 
$\Bun_{G;P}(\bS,\bD_\infty;\bD_0)$, classifying $G$-bundles
on $\BP^2$ with a trivialization along a divisor 
$\BP^1\simeq \bD_\infty\subset \BP^2$, and a reduction to a parabolic 
$P$ along another divisor $\BP^1\simeq \bD_0\subset \BP^2$.
Introducing these more general spaces is necessary for our calculation
of stalks of the IC sheaf.

Thus, in Section 9 we give the definition of $\fU^\theta_{G,P}$ and 
$\wt{\fU}^\theta_{G,P}$, and establish the corresponding factorization
properties. 

In Section 10, we prove
\thmref{stratification of parabolic Uhlenbeck} that describes the stratifications
of $\fU^\theta_{G,P}$ and $\wt{\fU}^\theta_{G,P}$ parallel to those of $\fU^a_G$.

In Section 11, we prove an important geometric property of
$\wt{\fU}^\theta_{G,P}$ when $P=B$ saying that its {\it boundary} (in
a natural sense) is a Cartier divisor.

\medskip

In Part IV we make a digression and discuss a construction of crystals
(in the sense of Kashiwara), using the quasi-maps spaces
$\on{QMaps}(\bC,\CG_{\fg',\fp})$ introduced in Section 2, for an
arbitrary Kac-Moody algebra $\fg'$. Unfortunately, to make this construction
work one has to assume a certain geometric property of the
$\on{QMaps}(\bC,\CG_{\fg',\fp})$ spaces, \conjref{lagrange},
which we verify in \secref{golova!} for $\fg'$ of affine type, using some
geometric properties of the parabolic Uhlenbeck's spaces.
As was mentioned above, it is via Kashiwara's crystals, more precisely
using \thmref{canonical basis}, that we relate the IC stalks on 
$\fU^a_G$ and the Lie algebra $\check \fg_{aff}$.

In Section 12 we recollect some general facts about Kashiwara's 
crystals. In particular, we review what properties are necessary to
prove that a given crystal $\sB_{\fg'}$ is isomorphic to the standard
crystal $\sB^\infty_{\fg'}$ of ~\cite{k5}.

In Section 13 we take our Lie algebra $\fg'$ to be finite dimensional
and spell out our ``new" construction of crystals using the affine
Grassmannian of the corresponding group $G$.

In Section 14 we consider the case of a general Kac-Moody algebra, 
and essentially repeat the construction of the previous section 
using the scheme $\on{QMaps}(\bC,\CG_{\fg',\fp})$ instead of the
affine Grassmannian. 

In Section 15 we verify \conjref{lagrange} and \conjref{golova?}
for finite-dimensional and affine Lie algebras, which insures that
in these cases the crystal of the previous section is well-defined
and can be identified with the standard crystal $\sB^\infty_\fg$.

\medskip

In Part V we perform the calculation of the IC sheaf on
the schemes $\fU^a_G$, $\wt{\fU}^\theta_{G,P}$ and $\fU^\theta_{G,P}$.

In Section 16 we formulate four theorems, which describe the behavior
of the IC sheaf, and in Section 17 we prove all the four statements by
an inductive argument borrowed from \cite{bgfm}.

\medskip

Finally, in the Appendix we reproduce a theorem of A.~Joseph, formulated
in Part I, Section 2,
which says that the space of based maps $\bC\to \CG_{\fg',\fp}$
of given degree is a scheme of finite type for any Kac-Moody Lie
algebra $\fg'$. 

\subsection{Acknowledgements}
It is clear from the above discussion that the present paper owes its
existence to V.~Drinfeld's generous explanations. We have also benefited
strongly from discussions with V.~Baranovsky and V.~Ginzburg. In fact, the
idea that the Uhlenbeck space must be related to the space of quasi-maps
into the affine Grassmannian was formulated by V.Ginzburg back in 1997.

In the course of our study of Uhlenbeck spaces, M.F. has enjoyed the
hospitality and support of the Hebrew University of Jerusalem, the University
of Chicago, and the University of Massachusetts at Amherst. His research
was conducted for the Clay Mathematical Institute, and partially supported
by the CRDF award RM1-2545-MO-03.

D.G. is a prize fellow of the Clay Mathematics Institute. He also wishes
to thank the Hebrew University of Jerusalem where the major part of
this paper was written.

We would also like to thank the referee for helpful remarks.

\bigskip

\centerline{{\bf Part I}:  {\Large Preliminaries on quasi-maps}}

\bigskip

\section{Maps and quasi-maps}  \label{general quasi-maps}

\ssec{}   \label{def quasi-maps}

The simplest framework in which one defines the notion
of quasi-map is the following:

Let $\Y$ be a projective scheme, $\CE$ a vector space,
and $\CT\subset \BP(\CE)$ a closed subscheme.

There exists a scheme, which we will
denote $\on{Maps}(\Y,\CT)$ that represents the functor
that assigns to a test scheme $S$ the set of maps
$\Y\times S\to \CT$.

To show the representability, it is enough to assume
that $\CT$ is the entire $\BP(\CE)$ (since in general
$\on{Maps}(\Y,\CT)$ is evidently a closed subfunctor
in $\on{Maps}(\Y,\BP(\CE))$), and in the latter case
our functor can be rewritten as pairs
$(\CL,\kappa)$,
where $\CL$ is a line bundle on $\Y\times S$, and $\kappa$ is
an injective bundle map
$$\kappa:\CL\hookrightarrow  \CO_{\Y\times S}\otimes \CE.$$
Therefore, we are dealing with an open subset of a suitable
Hilbert scheme. The scheme $\on{Maps}(\Y,\CT)$ splits as a
disjoint union of subschemes, denoted $\on{Maps}^a(\Y,\CT)$,
and indexed by set of connected components of the Picard
stack $\on{Pic}(\Y)$ of $\Y$ (our normalization is such that
$\sigma=(\CL,\kappa)\in \on{Maps}^a(\Y,\CT)$ if
$\CL^{-1}\in\on{Pic}^a(\Y)$), and each $\on{Maps}^a(\Y,\CT)$ is a
quasi-projective scheme.

\medskip

We will introduce a bigger scheme, denoted
$\on{QMaps}^a(\Y,\CT;\CE)$, which contains $\on{Maps}^a(\Y,\CT)$ as an
open subscheme. First, we will consider the case of
$\CT=\BP(\CE)$.

By definition, $\on{QMaps}^a(\Y,\BP(\CE);\CE)$ represents the
functor that assigns to a scheme $S$ the set of pairs $(\CL,\kappa)$,
where $\CL$ is a line bundle on the product $\Y\times S$
belonging to the connected component $\on{Pic}^{-a}(\Y)$ of
the Picard stack, and $\kappa$ is an
{\it injective map of coherent sheaves}
$$\kappa:\CL\hookrightarrow  \CO_{\Y\times S}\otimes \CE,$$
such that the quotient is $S$-flat. (The latter condition is
equivalent to the fact that for every geometric point $s\in S$,
the restriction of $\kappa$ to $\Y\times s$ is injective.)
This functor is also representable by a quasi-projective scheme,
for the same reason.

\medskip

Let now $\CT$ be arbitrary, and let $\CI_\CT=\underset{n\geq 0}\oplus\, \CI_\CT^n
\subset \on{Sym}(\CE^*)$
be the corresponding graded ideal. We define the closed subscheme 
$\on{QMaps}^a(\Y,\CT;\CE)\subset \on{QMaps}^a(\Y,\BP(\CE);\CE)$
by the condition that for every $n$ the composition
$$\CO_{\Y\times S}\otimes \CI_\CT^n\hookrightarrow
\CO_{\Y\times S}\otimes \on{Sym}^n(\CE^*)\to
(\CL^{-1}){}^{\otimes n}$$
vanishes.

We will denote by $\on{QMaps}(\Y,\CT;\CE)$ the union of
$\on{QMaps}^a(\Y,\CT;\CE)$ over all connected components of
$\on{Pic}(\Y)$.

\medskip

In the main body of the present paper the scheme $\Y$ will
be a smooth algebraic curve, but for completeness in this section
we will consider the general case. Note that in the case of
curves, the parameter $a$ amounts to an integer, normalized
to be the negative of $\on{deg}(\CL)$.

In the rest of this section we will study various properties and
generalizations of the notion of a quasi-map introduced above.

\ssec{Variant}

Assume for a moment that $\Y$ is integral. 
Observe that if $\sigma=(\CL,\kappa)$ is an $S$-point of
$\on{QMaps}^a(\Y,\BP(\CE);\CE)$, there exists an open dense subset
$U\subset \Y\times S$, over which $\kappa$ is a bundle
map.

Consider the (automatically closed) subfunctor of
$\on{QMaps}^a(\Y,\BP(\CE);\CE)$, corresponding to the condition
that the resulting {\it map} $U\to \BP(\CE)$ factors through $\CT\subset \BP(\CE)$. 
It is easy to see that this subfunctor coincides with $\on{QMaps}^a(\Y,\CT;\CE)$.

The above definition can be also spelled out as follows.
Let $C(\CT;\CE)$ be the affine cone over $\CT$, i.e. the closure in
$\CE$ of the preimage of $\CT$ under the natural map $(\CE-0)\to
\BP(\CE)$. The multiplicative group $\BG_m$ acts naturally on
$C(\CT;\CE)$ and we can form the stack $C(\CT;\CE)/\BG_m$, which
contains $\CT$ as an open substack.

It is easy to see that a map $S\to \on{QMaps}(\Y,\CT;\CE)$ is the
same as a map $\sigma:\Y\times S\to C(\CT;\CE)/\BG_m$, such that for every
geometric $s\in S$, the map $\Y\simeq \Y\times s\to C(\CT;\CE)/\BG_m$
sends the generic point of $\Y$ into $\CT$.

\medskip

We will now introduce a still bigger scheme, $\on{QQMaps}^p(\Y,\CT;\CE)$.
First, $\on{QQMaps}^p(\Y,\BP(\CE);\CE)$ is the scheme, classifying
pairs $(\CL,\kappa)$, where $\CL$ is an $S$-flat coherent sheaf on
$\Y\times S$, whose restriction to every geometric fiber $\Y\times s$
is of generic rank $1$, and $\kappa$ is a map of coherent sheaves,
injective at each $Y\times s$ (as above $s\in S$ denotes a geometric
point). The superscript $p$ signifies that the Hilbert
polynomial $p$ of $\CL$ (with respect to some ample line bundle
on $\Y$) is fixed. Omitting the superscript means that we are
taking the union over all possible Hilbert polynomials.

As in the case of $\on{QMaps}^p(\Y,\CT;\CE)$,
given an $S$-point $(\CL,\kappa)$ of $\on{QQMaps}^p(\Y,\BP(\CE);\CE)$,
there exists an open dense subset $U\subset \Y\times S$, over which 
$\CL$ is a line bundle, and $\kappa$ is a bundle map.

For a closed subscheme $\CT\subset \BP(\CE)$, we defined the closed 
subfunctor $\on{QQMaps}^p(\Y,\CT;\CE)$ of $\on{QQMaps}^p(\Y,\BP(\CE);\CE)$
by the condition that the resulting map $U\to \BP(\CE)$ factors through $\CT$.

\medskip

Since $\on{QQMaps}^p(\Y,\BP(\CE);\CE)$ is actually the entire
Hilbert scheme, we obtain:

\begin{lem} \label{QQ is proper}
The scheme $\on{QQMaps}^p(\Y,\BP(\CE);\CE)$ is proper.
\end{lem}

Note that by assumption, for every geometric point $s\in S$, the
restriction $\CL_s:=\CL|_{\Y\times s}$ embeds into $\CE\otimes \CO_\Y$;
therefore, $\CL$ is torsion-free. In particular, when $\Y$ is a smooth curve,
$\on{QMaps}(\Y,\CT;\CE)$ coincides with $\on{QQMaps}(\Y,\CT;\CE)$.
 But in general we have an open embedding
$\on{QMaps}(\Y,\CT;\CE)\hookrightarrow \on{QQMaps}(\Y,\CT;\CE)$.

\begin{prop}
Suppose $\Y$ is smooth. Then the scheme
$\on{QMaps}^a(\Y,\CT;\CE)$ is proper as well.
\end{prop}

\begin{proof}

Let us check the valuative criterion of properness.
Using \lemref{QQ is proper} we can assume having
the following set-up:

Let $\bX$ be a curve with a marked point $0_\bX\in \bX$,
and let $\CL$ be a torsion-free coherent sheaf on
$\Y\times \bX$, embedded into $\CE\otimes \CO_{\Y\times \bX}$,
and such that $\CL|_{\Y\times(\bX-0_\bX)}$ is a line bundle.

We claim that $\CL$ itself is a line bundle.
Indeed, $\CL$ is locally free outside of codimension $2$,
and since $\Y\times \bX$ is smooth, it admits a saturation,
i.e. there exists a unique line bundle $\CL^0$ containing $\CL$,
such that $\CL^0/\CL$ is supported in codimension $2$.
Moreover, it is easy to see that the map
$\CL\to \CE\otimes \CO_{\Y\times \bX}$ extends to a map
$\CL^0\to \CE\otimes \CO_{\Y\times \bX}$. But
since the cokernel of $\kappa$ was $\bX$-flat, we obtain that
$\CL\to \CL^0$ must be an isomorphism.

\end{proof}

In particular, we see that when $\Y$ is smooth, $\on{QMaps}(\Y,\CT;\CE)$
is the union of certain of the connected components of
$\on{QQMaps}(\Y,\CT;\CE)$.

\ssec{Based quasi-maps}

Let $\Y'\subset \Y$ be a closed subscheme, and
$\sigma':\Y'\to \CT$ a fixed map.
We introduce the scheme $\on{QMaps}^a(\Y,\CT;\CE)_{\Y',\sigma'}$ as
a (locally closed) subfunctor of $\on{QMaps}^a(\Y,\CT;\CE)$ defined
by the following two conditions:

\begin{itemize}
\item
The restriction of the map $\kappa$ to $\Y'\times S$
is a {\it bundle map}. (Equivalently, the quotient
$\CO_{\Y\times S}\otimes \CE/\CL$ is locally free in
a neighborhood of $\Y'\times S$.)

\item
The the resulting map $\Y'\times S\to \CT$ equals
$\Y'\times S\to \Y'\overset{\sigma'}\to \CT$.
\end{itemize}

When $\Y$ is integral, one defines in a similar way the subscheme
$$\on{QQMaps}^a(\Y,\CT;\CE)_{\Y',\sigma'}\subset
\on{QQMaps}^a(\Y,\CT;\CE).$$

\medskip

The following assertion will be needed in the main body of the paper:
Assume that $\Y$ is a smooth curve $\bC$ and $\Y'$ is a point $\bc$,
so that $\sigma':pt\to \CT$ corresponds to some point of $\CT$.

\begin{lem} \label{qmaps affine}
There exists an affine map
$\on{QMaps}^a(\bC,\CT;\CE)_{\bc,\sigma'}\to
(\bC-\bc)^{(a)}$. In particular, 
the scheme $\on{QMaps}^a(\bC,\CT;\CE)_{\bc,\sigma'}$ is affine.
\end{lem}

\begin{proof}

We can assume that $\CT=\BP(\CE)$, so that $\sigma'$ corresponds to
a line $\ell\subset \CE$. Let us choose a splitting $\CE\simeq
\CE'\oplus \ell$.

Then we have a natural map $\on{QMaps}^a(\bC,\CE)_{\bc,\sigma'}\to
(\bC-\bc)^{(a)}$, that assigns to a pair $(\CL,\kappa:\CL\to
\CO\otimes \CE)$ the divisor of zeroes of the composition
$\CL\to \CO\otimes \CE\to \CO\otimes \ell\simeq \CO$.

Since $\bC-\bc$ is affine, the symmetric power
$(\bC-\bc)^{(a)}$ is affine as well. We claim that the above
morphism $\on{QMaps}^a(\bC,\CE)_{\bc,\sigma'}\to
(\bC-\bc)^{(a)}$ is also affine.

Indeed, given a divisor $D\in (\bC-\bc)^{(a)}$, the fiber of
$\on{QMaps}^a(\bC,\CE)_{\bc,\sigma'}$ over it is the vector
space $\on{Hom}(\CO(-D),\CE')$.

\end{proof}

\ssec{The relative version}

Suppose now that $\Y$ itself is a flat family
of projective schemes over some base $\X$, and 
$\CE$ is a vector bundle on $\Y$, with $\CT\subset \BP(\CE)$
a closed subscheme.

The scheme of maps $\on{Maps}(\Y,\CT)$ assigns to a test
scheme $S$ over $\X$ the set of maps
$\sigma:\Y\underset{\X}\times S\to \CT$, such that the composition
$\Y\underset{\X}\times S\to \CT\to \Y$ is the projection on
the first factor. By definition, $\on{Maps}(\Y,\CT)$ is also
a scheme over $\X$.

By essentially repeating the construction of \secref{def quasi-maps},
we obtain a scheme $\on{QMaps}(\Y,\CT;\CE)$,
which is quasi-projective over $\X$.

\medskip

As before, if $\Y'\subset \Y$ is a closed subscheme, flat over
$\X$, and $\sigma':\Y'\to \CT$ is a $\X$-map, we can define a
locally closed subscheme $\on{QMaps}^a(\Y,\CT;\CE)_{\Y',\sigma'}$ of
based maps. When $\Y\to \X$ has integral fibers, one can define
the relative versions of $\on{QQMaps}(\Y,\CT;\CE)$
and $\on{QQMaps}(\Y,\CT;\CE)_{\Y',\sigma'}$ in a similar fashion.

\medskip

In what follows, we will mostly discuss the ``absolute''
case, leaving the (straightforward) modifications required in
the relative and based cases to the reader.

\ssec{}

Next we will study some functorial properties of
the quasi-maps' spaces. From now on we will assume that
$\Y$ is integral.

For $\CT\subset \BP(\CE)$ as above, let $\CP$ be the very ample
line bundle $\CO(1)|_\CT$. Note that we have a map
$H^0(\CT,\CP)^*\to \CE$.

Let now $\CT_1\subset\BP(\CE_1)$ and $\CT_2\subset\BP(\CE_2)$ be two
projective schemes, and $\phi:\CT_1\to \CT_2$ be a map,
such that $\phi^*(\CP_2)\simeq \CP_1$. Assume, moreover,
that the corresponding map $H^0(\CT_2,\CP_2)\to H^0(\CT_1,\CP_1)^*$ is
extended to a commutative diagram
$$
\CD
H^0(\CT_1,\CP_1)^* @>>> H^0(\CT_2,\CP_2)^* \\
@VVV   @VVV  \\
\CE_1 @>>>   \CE_2.
\endCD
$$

\begin{lem}  \label{map of quasi-maps}
Under the above circumstances the morphism
$\on{Maps}(\Y,\CT_1)\to\on{Maps}(\Y,\CT_2)$ extends
to a Cartesian diagram:
$$
\CD
\on{Maps}(\Y,\CT_1) @>>> \on{QMaps}(\Y,\CT_1,\CE_1) @>>>
\on{QQMaps}(\Y,\CT_1,\CE_1) \\
@VVV   @VVV   @VVV  \\
\on{Maps}(\Y,\CT_2) @>>> \on{QMaps}(\Y,\CT_2,\CE_2) @>>>
\on{QQMaps}(\Y,\CT_2,\CE_2).
\endCD
$$
\end{lem}

\begin{proof}

The existence of the map
$\on{QMaps}(\Y,\CT_1,\CE_1)\to \on{QMaps}(\Y,\CT_2,\CE_2)$
is nearly obvious from the interpretation of quasi-maps
as maps to the affine cone:

Note that the condition of the lemma is equivalent
to the fact that the map $\CT_1\to\CT_2$ is extended
to a commutative diagram of the corresponding affine
cones:
$$
\CD
\CT_1  @<<< C(\CT_1;\CE_1)-0 @>>> C(\CT_1;\CE_1) @>>> \CE_1 \\
@VVV  @VVV  @VVV @VVV \\
\CT_2  @<<< C(\CT_2;\CE_2)-0 @>>> C(\CT_2;\CE_2) @>>> \CE_2.
\endCD
$$

Therefore, we can just compose a map
$\Y\times S\to C(\CT_1;\CE_1)/\BG_m$ with the map
$C(\CT_1;\CE_1)/\BG_m\to C(\CT_2;\CE_2)/\BG_m$.

\medskip

Let us construct the map
$\on{QQMaps}(\Y,\CT_1,\CE_1)\to \on{QQMaps}(\Y,\CT_2,\CE_2)$.
Given a scheme $S$, and a coherent sheaf of generic rank $1$
$\CL$ on $\Y\times S$
with a map $\kappa_1:\CL\to \CO_{\Y\times S}\otimes \CE_1$, let
us consider the composition
$$\kappa_2:\CL\to \CO_{\Y\times S}\otimes \CE_1\to
\CO_{\Y\times S}\otimes \CE_2.$$

We claim that the restriction of $\kappa_2$ to every fiber
$\Y\times s$ is still injective. Indeed, we are dealing with
a map $\CL_s\to \CO_\Y\otimes \CE_2$, which is injective
over an open subset of $\Y$ (because on some open subset $U$,
$\kappa_2$ corresponds to a {\it map} $U\to \CT_2$), and
we know that $\CL_s$ is torsion free.

\medskip

The square
$$
\CD
\on{QMaps}(\Y,\CT_1,\CE_1) @>>> \on{QQMaps}(\Y,\CT_1,\CE_1) \\
@VVV   @VVV   \\
\on{QMaps}(\Y,\CT_2,\CE_2) @>>> \on{QQMaps}(\Y,\CT_2,\CE_2).
\endCD
$$
is Cartesian, because being in $\on{QMaps}$ is just the condition
that $\CL$ is locally free.

\medskip

To prove that the square
$$
\CD
\on{Maps}(\Y,\CT_1) @>>> \on{QMaps}(\Y,\CT_1,\CE_1) \\
@VVV   @VVV    \\
\on{Maps}(\Y,\CT_2) @>>> \on{QMaps}(\Y,\CT_2,\CE_2).
\endCD
$$
is Cartesian, we have to show that if $\on{coker}(\kappa_2)$
is locally free, then so is $\on{coker}(\kappa_1)$. But
this follows from the $4$-term exact sequence:
$$0\to \CO_{\Y\times S}\otimes \on{ker}(\CE_1\to \CE_2)\to
\on{coker}(\kappa_1)\to \on{coker}(\kappa_2)\to
\CO_{\Y\times S}\otimes \on{coker}(\CE_1\to \CE_2)\to 0.$$

\end{proof}

This lemma has a straightforward generalization to the based
case:

If we fix $\Y'$ and $\sigma'_1:\Y'\to \CT_1$ with $\sigma'_2:=
\phi\circ \sigma'_1$, we have a Cartesian square
$$
\CD
\on{QMaps}(\Y,\CT_1,\CE_1)_{\Y',\sigma'_1} @>>>
\on{QMaps}(\Y,\CT_1,\CE_1)  \\
@VVV   @VVV \\
\on{QMaps}(\Y,\CT_2,\CE_2)_{\Y',\sigma'_2} @>>>
\on{QMaps}(\Y,\CT_2,\CE_2),
\endCD
$$
and similarly, for $\on{QQMaps}$.

\medskip

Assume now that we are in the situation of
\lemref{map of quasi-maps}, and that the map $\CT_1\to \CT_2$
is a closed embedding.

\begin{prop}  \label{embedding of quasi-maps}
Under the above circumstances, the resulting morphisms
$$\on{QMaps}(\Y,\CT_1,\CE_1)\to \on{QMaps}(\Y,\CT_2,\CE_2) \text{ and }
\on{QQMaps}(\Y,\CT_1,\CE_1)\to \on{QQMaps}(\Y,\CT_2,\CE_2)$$
are closed embeddings.

If $\CT_1\to \CT_2$ is an isomorphism, and $\Y$ is normal
then the morphism
$\on{QMaps}(\Y,\CT_1,\CE_1)\to \on{QMaps}(\Y,\CT_2,\CE_2)$
induces an isomorphism
on the level of reduced schemes.
\end{prop}

(Due to the Cartesian diagram above, the same assertions
will hold in the based situation.)

\begin{proof}

If $\CE_1\to \CE_2$ is an isomorphism, then the assertion
of the lemma is obvious, as both schemes are closed
subschemes in $\on{QMaps}(\Y,\BP(\CE_2),\CE_2)$
(resp., $\on{QQMaps}(\Y,\BP(\CE_2),\CE_2)$).
Therefore, we can assume that $\CT_1\simeq \CT_2=:\CT$.

\medskip

First, we claim that the map
$$\on{Hom}(S,\on{QQMaps}(\Y,\CT,\CE_1))\to
\on{Hom}(S,\on{QQMaps}(\Y,\CT,\CE_2))$$
is an injection for any scheme $S$.

This follows from the fact that
if $\CL$ is a torsion-free coherent sheaf on $\Y\times S$, then any two maps
$\CL\to\CO_{\Y\times S}\otimes \CE_1$ that coincide over
a dense subset $U\subset \Y\times S$, must coincide globally.

\medskip

Since the scheme $\on{QQMaps}(\Y,\CT,\CE_1)$ is proper,
the map $$\on{QQMaps}(\Y,\CT,\CE_1)\to
\on{QQMaps}(\Y,\CT,\CE_2)$$ is proper. Combined with the
previous statement, we obtain that it is a closed embedding.
Using the Cartesian diagrams of \lemref{map of quasi-maps},
we obtain that
$\on{QMaps}(\Y,\CT,\CE_1)\to \on{QMaps}(\Y,\CT,\CE_2)$
is a closed embedding as well.

\medskip

To prove the last assertion, we claim that for any normal
scheme $S$, the map
$$\on{Hom}(S,\on{QMaps}(\Y,\CT,\CE_1))\to
\on{Hom}(S,\on{QMaps}(\Y,\CT,\CE_2))$$
is in fact a bijection.
Indeed, let us consider the map between the cones
$C(\CT;\CE_1)\to C(\CT;\CE_2)$. This map is finite, and
is an isomorphism away from $0\in C(\CT;\CE_2)$.

For a given element in $\on{Hom}(S,\on{QMaps}(\Y,\CT,\CE_2))$,
consider the Cartesian product
$$(\Y\times S)\underset{C(\CT;\CE_2)/\BG_m}\times
C(\CT;\CE_1)/\BG_m.$$
This is a scheme, finite and generically isomorphic to
$\Y\times S$. Since $\Y\times S$ is normal, we obtain
that it admits a section, i.e. we have a map
$\Y\times S\to C(\CT,\CE_1)/\BG_m$.

\end{proof}

\ssec{}  \label{canonical model}

The definition of quasi-maps given above depends
on the projective embedding of $\CT$. In fact, one can
produce another scheme, which depends only on the
corresponding very ample line bundle $\CP=\CO(1)|_\CT$.
Note that for the discussion below it is crucial that
we work with $\on{QMaps}$ and not $\on{QQMaps}$, because
we will be taking tensor products of the corresponding line
bundles.

For a test scheme $S$, we let $\on{Hom}(S,\on{QMaps}(\Y,\CT;\CP))$
be the set consisting of pairs $(\CL,\kappa)$, where $\CL$ is,
as before, a line bundle on $\Y\times S$, and $\kappa$ is a map
$\CL\to \CO_{\Y\times S}\otimes \Gamma(\CT,\CP)^*$, which 
extends to a map of algebras
$$\underset{n}\oplus\, H^0(\CT,\CP^{\otimes n})\to 
\underset{n}\oplus\, (\CL^*){}^{\otimes n},$$
and such that $\kappa$ is 
injective over every geometric point $s\in S$. (This definition makes
sense for an arbitrary, i.e. not necessarily integral scheme $\Y$.)

\medskip

Resuming the assumption that $\Y$ is integral, we can spell out
the above definition as follows.
Let $C(\CT;\CP)$ be the affine closure of the total space
of $(\CP^{-1}-0)$, i.e. $C(\CT;\CP)$ is the spectrum of
the algebra $\underset{n}\oplus\, H^0(\CT,\CP^{\otimes n})$.
We can form the stack $C(\CT;\CP)/\BG_m$, which contains
$\CT$ as an open subset.

As before, we can consider the functor $\on{QMaps}(\Y,\CT;\CP)$
that assigns to a scheme $S$ the set of maps
$\Y\times S\to C(\CT;\CP)/\BG_m$, such that for every geometric
point $s\in S$, the map $\Y\simeq \Y\times s\to C(\CT;\CP)/\BG_m$
sends the generic point of $\Y$ into $\CT$.

\begin{prop}  \label{changing the line bundle}
For $\CT\subset \BP(\CE)$ and $\CP$ as above we have

{\em (a)} There is a canonical morphism
$\on{QMaps}(\Y,\CT;\CP)\to \on{QMaps}(\Y,\CT;\CE)$.
When $\Y$ is smooth, this map is a closed embedding,
which induces an isomorphism on the level of reduced
schemes.

{\em (b)} For $n\in\BN$, we have a morphism
$\on{QMaps}(\Y,\CT;\CP)\to \on{QMaps}(\Y,\CT;\CP^{\otimes n})$.
When $\Y$ is smooth, this map is also a closed embedding,
which induces an isomorphism on the level of reduced
schemes.

{\em (c)}
If $\phi:\CT_1\to \CT_2$ is a map of projective schemes, such that
$\phi^*(\CP_2)\simeq \CP_1$, we have a morphism
$\on{QMaps}(\Y,\CT;\CP_1)\to \on{QMaps}(\Y,\CT;\CP_2)$, which is
a closed embedding, whenever $\phi$ is.

\end{prop}

\begin{proof}

Observe that we have a finite map
$C(\CT;\CP)\to C(\CT;H^0(\CT,\CP)^*)$, and
moreover, for any positive integer $n$, a finite map
$C(\CT;\CP)\to C(\CT;\CP^{\otimes n})$.
This makes the existence of the maps
$\on{QMaps}(\Y,\CT;\CP)\to \on{QMaps}(\Y,\CT;\CE)$
and $\on{QMaps}(\Y,\CT;\CP)\to \on{QMaps}(\Y,\CT;\CP^{\otimes n})$
obvious.

\medskip

For $k$ large enough, the map of
$C(\CT;\CP)$ into the product
$\underset{i=1,...,k}\Pi\, C(\CT;H^0(\CT,\CP^{\otimes i})^*)$
is a closed embedding. Hence $\on{QMaps}(\Y,\CT;\CP)$
is representable by a closed subscheme in the product
$\underset{i=1,...,k}\Pi\,
\on{QMaps}(\Y,\CT;H^0(\CT,\CP^{\otimes i})^*)$.

In particular, when $\Y$ is smooth, we obtain that
$\on{QMaps}(\Y,\CT;\CP)$ is proper.

\medskip

To finish the proof of the proposition, it would suffice
to show that for a test scheme $S$, the maps
$\Hom(S,\on{QMaps}(\Y,\CT;\CP))\to \Hom(S,\on{QMaps}(\Y,\CT;\CE))$ and
$\Hom(S,\on{QMaps}(\Y,\CT;\CP))\to
\Hom(S,\on{QMaps}(\Y,\CT;\CP^{\otimes n}))$
are injective and are, moreover, bijective if $S$ is normal.
This is done exactly as in the proof of \lemref{embedding of quasi-maps}.

Assertion (c) of the Proposition follows from
\lemref{map of quasi-maps} and \lemref{embedding of quasi-maps}.

\end{proof}

\medskip

Point (c) of the above proposition allows to define $\on{QMaps}(\Y,\CT;\CP)$
as a strict ind-scheme, when $\CT$ is a strict ind-projective ind-scheme
endowed with a very ample line bundle.

Indeed, if $\CT=\underset{\to}{lim}\, \CT_i$, where $\CT_i$ are
projective schemes, and $\CT_i\to \CT_j$ are closed embeddings,
and $\CP_i$ is a compatible system of line bundles on $\CT_i$,
we have a system of closed embeddings
$$\on{QMaps}(\Y,\CT_i;\CP_i)\to  \on{QMaps}(\Y,\CT_j;\CP_j).$$

\medskip

In this section we will discuss the quasi-maps'
spaces in their $\Hom(S,\on{QMaps}(\Y,\CT;\CE))$ incarnation.
The $\Hom(S,\on{QMaps}(\Y,\CT;\CP))$-versions of the corresponding results
are obtained similarly.

\ssec{}  \label{multiple quasi-maps}

Suppose now that we have an embedding $\CT\subset \BP(\CE_1)\times...
\times \BP(\CE_k)$. In this situation for a $k$-tuple $\ol{a}$
of parameters $a_i$ we can introduce a scheme of
quasi-maps:

\medskip

We first define
$\on{QMaps}^{\ol{a}}(\Y,\BP(\CE_1)\times...\times
\BP(\CE_k);\CE_1,...,\CE_k)$ simply as the product
$$\underset{i}\Pi\,\on{QMaps}^{a_i}(\Y,\BP(\CE_i)).$$
For an arbitrary $\CT$, which corresponds to a multi-graded
ideal $$\underset{n_1,...,n_k}\oplus\, \CI_\CT^{n_1,...,n_k}\subset \on{Sym}(\CE^*_1)\otimes...\otimes \on{Sym}(\CE^*_n),$$
we set $\on{QMaps}^{\ol{a}}(\Y,\CT;\CE_1,...,\CE_k)$ to be the closed subscheme in
$\underset{i}\Pi\,\on{QMaps}^{a_i}(\Y,\BP(\CE_i))$,
defined by the condition that for each $n$-tuple
$n_1,...,n_k$ of non-negative integers, the composition
$$\CO_{\Y\times S}\otimes (\CI_\CT)^{n_1,...,n_k}\hookrightarrow
\CO_{\Y\times S}\otimes \on{Sym}^{n_1}(\CE^*_1)\otimes...\otimes
\on{Sym}^{n_k}(\CE^*_k)\to
(\CL_1^{-1})^{\otimes n_1}\otimes...\otimes (\CL^{-1})_k^{\otimes n_k}$$
vanishes.

For $\Y$ integral, this condition is equivalent to demanding that 
the generic point of $\Y$ gets mapped to $\CT$, just as in the definition
of $\on{QMaps}^a(\Y,\CT;\CE)$.

\medskip

Alternatively, let $C(\CT;\CE_1,...,\CE_k)\subset \CE_1\times...\times \CE_k$
be the affine cone over $\CT$. We have an action of $\BG_m^k$ on
$C(\CT;\CE_1,...,\CE_k)$, and the stack-quotient
$C(\CT;\CE_1,...,\CE_k)/\BG_m^k$ again contains $\CT$ as an open subset.
By definition,
$\on{Hom}(S,\on{QMaps}^{\ol{a}}(\Y,\CT;\CE_1,...,\CE_k))$
is the set of maps from $\Y\times S$ to the stack
$C(\CT,\CE_1,...,\CE_k)/\BG_m^k$, such that for every geometric
$s\in S$, the map $\Y\times s\to C(\CT;\CE_1,...,\CE_k)/\BG_m^k$ sends
the generic point of $\Y$ to $\CT$.

\medskip

The above definition has an obvious $\on{QQMaps}$ version, and
\lemref{map of quasi-maps}, \propref{embedding of quasi-maps}
with its based analogues, generalize in a straightforward way.
Note, however, that for the next proposition it is essential
that we work with $\on{QMaps}$, and not with $\on{QQMaps}$,
as we will be taking tensor products of $\CL_i$'s.

\medskip

Observe that $\CT$ can be naturally embedded into
$\BP(\CE_1\otimes...\otimes \CE_k)$ via the Segre embedding
$\BP(\CE_1)\times...\times \BP(\CE_k)\hookrightarrow
\BP(\CE_1\otimes...\otimes \CE_k)$.

\begin{lem}  \label{product of quasi-maps}
There exists a map
$$\on{QMaps}^{\ol{a}}(\Y,\CT;\CE_1,...,\CE_k))\to
\on{QMaps}^{a_1+...+a_k}(\Y,\CT;\CE_1\otimes...\otimes \CE_k),$$
extending the identity map on $\on{Maps}(\Y,S)$.
Moreover, when $\Y$ is smooth, this map is finite.
\end{lem}

\begin{proof}

We have a map of stacks
$$C(\CT;\CE_1\otimes...\otimes \CE_k)/\BG_m\to
C(\CT;\CE_1,...,\CE_k)/\BG_m^k,$$
which makes the existence of the map
$\on{QMaps}(\Y,\CT;\CE_1,...,\CE_k))\to
\on{QMaps}(\Y,\CT;\CE_1\otimes...\otimes \CE_k)$
obvious.

We know that when $\Y$ is smooth,
$\on{QMaps}(\Y,\CT;\CE_1,...,\CE_k)$ is proper; therefore
to check that our map is finite, it is enough to show that
the fiber over every geometric point of the scheme
$\on{QMaps}(\Y,\CT;\CE_1\otimes...\otimes \CE_k)$ is finite.

\medskip

Suppose that we have a line bundle $\CL$ on $\Y$ with
an injective map $\kappa:\CL\to
\CO_\Y\otimes (\CE_1\otimes...\otimes \CE_k)$, such that there exists
an open subset $U\subset \Y$ such that $\kappa|_U$ is a bundle map
corresponding to a map $U\to \CT$. In particular, we have the line
bundles $\CL^U_1,...,\CL^U_k$, defined on $U$, such that
$\CL^U_1\otimes...\otimes\CL^U_k\simeq \CL|_U$.

Let $\CL'_i$ be the (unique) extension of $\CL^U_i$ such that
the map $\CL'_i\to \CO_\Y\otimes \CE_i$ is regular and is a bundle
map away from codimension $2$. Clearly, $\CL$ is a subsheaf in
$\CL'_1\otimes...\otimes\CL'_k$.

Then the fiber of $\on{QMaps}(\Y,\CT;\CE_1,...,\CE_k))\to
\on{QMaps}(\Y,\CT;\CE_1\otimes...\otimes \CE_k)$ is the scheme
classifying subsheaves $\CL_i\subset \CL'_i$ such that
$\CL_1\otimes...\otimes\CL_k\simeq \CL$ and the maps
$\CL_i\to \CO_\Y\otimes \CE_i$ continue to be regular. This
scheme is clearly finite.

\end{proof}

\medskip

An example of the situation of \secref{multiple quasi-maps}
is this:

Suppose that $\CT=\CT_1\times...\times \CT_k$, and
each $\CT_i$ embeds into the corresponding $\CE_i$. Then
from \lemref{product of quasi-maps} we obtain that
there exists a finite morphism
$$\underset{i}\Pi\, \on{QMaps}^{a_i}(\Y,\CT_i;\CE_i)\simeq
\on{QMaps}^{\ol{a}}(\Y,\CT;\CE_1,...,\CE_k)\to
\on{QMaps}^{a_1+...+a_k}(\Y,\CT;\CE_1\otimes...\otimes\CE_k).$$

\medskip

In the main body of the paper we will need the following
technical statement. Let $\CT_0\subset \CE_0$ and $\CT'_0\subset \CE'_0$
be two projective varieties.

Let $\CT$ be another projective variety, embedded into
$\BP(\CE_0)\times\BP(\CE_1)\times...\times \BP(\CE_n)$, and endowed with
a map $\CT\to \CT_0$, such that the two maps from $\CT$ to $\BP(\CE_0)$ coincide.

\begin{lem}  \label{product lemma}
The (a priori finite) map
\begin{align*}
&\on{QMaps}(\Y,\CT;\CE_0,\CE_1,...,\CE_k)\times
\on{QMaps}(\Y,\CT'_0;\CE'_0)\to \\
&\on{QMaps}(\Y,\CT\times \CT'_0;\CE_0\otimes \CE'_0,\CE_1,...,\CE_k)
\underset{\on{QMaps}(\Y,\CT_0\times \CT'_0;\CE_0\otimes \CE'_0)}\times
\on{QMaps}(\Y,\CT_0\times \CT'_0;\CE_0,\CE'_0).
\end{align*}
is an isomorphism.
\end{lem}

The proof is straightforward.

\ssec{}   \label{degenerations of quasi-maps}

Let $\overset{\circ}{\on{QQMaps}}(\Y,\CT;\CE)$ denote the open subset
of $\on{QQMaps}(\Y,\CT;\CE)$ corresponding to the condition that
$\on{coker}(\kappa)$ is a torsion-free coherent sheaf on $\Y$
(over every geometric point $s$ in the corresponding test scheme $S$).
Let us denote by $\partial(\on{QQMaps}(\Y,\CT;\CE))$ the complement
to $\overset{\circ}{\on{QQMaps}}(\Y,\CT;\CE)$ in
$\on{QQMaps}(\Y,\CT;\CE)$.

For a Hilbert polynomial $p$,
let $\on{TFree}^p_1(\Y)$ denote the corresponding
connected component of the stack of
torsion-free coherent sheaves of generic rank $1$ on $\Y$.
For two parameters $p,p'$, let us denote by $\Y^{Q;(p,p')}$ the stack
classifying triples $(\CL,\CL',\beta:\CL\hookrightarrow \CL')$, where
$\CL,\CL'$ are points of $\on{TFree}^p_1(\Y)$ and $\on{TFree}_1^{p'}(\Y)$,
respectively, and $\beta$ is a non-zero map, which is automatically
an embedding.
Note that for $p=p'$, $\Y^{Q;(p,p')}$ projects isomorphically
onto $\on{TFree}_1^p(\Y)$.

\begin{prop}
There is a natural proper morphism
$$\ol{\iota}_{p',p}:\on{QQMaps}^{p'}(\Y,\CT;\CE)\underset{\on{TFree}_1^{p'}(\Y)}\times
\Y^{Q;(p,p')}\to \on{QQMaps}^p(\Y,\CT;\CE).$$ The composition
\begin{align*}
&\iota_{p',p}:\overset{\circ}{\on{QQMaps}}{}^{p'}(\Y,\CT;\CE)
\underset{\on{TFree}_1^{p'}(\Y)}\times \Y^{Q;(p,p')}\hookrightarrow
\on{QQMaps}^{p'}(\Y,\CT;\CE)\underset{\on{TFree}_1^{p'}(\Y)}\times
\Y^{Q;(p,p')}\to \\
&\on{QQMaps}^p(\Y,\CT;\CE)
\end{align*}
is a locally-closed embedding. Moreover, every geometric point
of $\on{QQMaps}^p(\Y,\CT;\CE)$ belongs to the image of $\iota_{p',p}$ for
exactly one $p'$.
\end{prop}

\begin{proof}

The map $\ol{\iota}_{p',p}$ is constructed in a most straightforward
way:
given an $S$-point of $\on{QQMaps}^{p'}(\Y,\CT,\CE)$, i.e.
a sheaf $\CL'$ on $\Y\times S$, with an embedding
$\kappa':\CL'\to \CO_{\Y\times S}\otimes \CE$, and $S$-point
of $\Y^{Q,(p,p')}$, i.e. another sheaf $\CL$ on $\Y\times S$
with an embedding $\beta:\CL\to \CL'$, we define a new
point of $\on{QMaps}(\Y,\CT,\CE)$, by setting
$\kappa$ to be the composition $\kappa'\circ\beta$.

This map is proper, because the scheme $\on{QQMaps}^{p'}(\Y,\CT,\CE)$
is proper, and so is the projection $\Y^{Q;(p,p')}\to
\on{TFree}_1^p(\Y)$.

\medskip

Observe now that the map
$$\underset{p'}\bigcup\, \overset{\circ}{\on{QQMaps}}{}^{p'}(\Y,\CT;\CE)
\underset{\on{TFree}_1^{p'}(\Y)}\times\Y^{Q;(p,p')} \to
\on{QQMaps}^p(\Y,\CT;\CE)$$
is injective on the level of $S$-points for any $S$.

Indeed, for $(\CL,\CL',\kappa',\beta)$ as above with
$(\CL',\kappa')\in \overset{\circ}{\on{QQMaps}}{}^{p'}(\Y,\CT;\CE)$,
the sheaf $\CL'$ is reconstructed as the unique subsheaf
in $\CO_{\Y\times S}\otimes\CE$ containing $\CL$, such that the
quotient $\CO_{\Y\times S}\otimes\CE/\CL'$ is torsion-free over
every geometric point of $S$.
(Moreover, for every geometric point $(\CL,\kappa)$ of
$\on{QQMaps}^p(\Y,\CT;\CE)$, such $\CL'$ exists and equals
the preimage in $\CO_\Y\otimes\CE$ of the maximal
torsion subsheaf in $\on{coker}(\kappa)$.)

\medskip

Thus, we obtain that the open subset in
$\on{QQMaps}^{p'}(\Y,\CT;\CE)\underset{\on{TFree}_1^{p'}(\Y)}\times\Y^{Q;(p,p')}$
equal to
$$(\ol{\iota}_{p',p})^{-1}(\on{QQMaps}^p(\Y,\CT;\CE)-
\ol{\iota}_{p',p}(\partial(\on{QQMaps}^{p'}(\Y,\CT;\CE))))$$
coincides with
$\overset{\circ}{\on{QQMaps}}{}^{p'}(\Y,\CT;\CE)
\underset{\on{TFree}_1^{p'}(\Y)}\times\Y^{Q;(p,p')}$.
Hence, the map
\begin{align*}
&\ol{\iota}_{p',p}:\overset{\circ}{\on{QQMaps}}{}^{p'}(\Y,\CT;\CE)
\underset{\on{TFree}_1^{p'}(\Y)}\times\Y^{Q;(p,p')}\to \\
&\left(\on{QQMaps}^p(\Y,\CT;\CE)-
\ol{\iota}_{p',p}(\partial(\on{QQMaps}^{p'}(\Y,\CT;\CE)))\right)
\end{align*}
is also proper, and being an injection on the level of $S$-points,
it is a closed embedding.

\end{proof}

\medskip

Set now
$\overset{\circ}{\on{QMaps}}(\Y,\CT;\CE):=
\overset{\circ}{\on{QQMaps}}(\Y,\CT;\CE)\cap \on{QMaps}(\Y,\CT;\CE)$,
and $\partial(\on{QMaps}(\Y,\CT;\CE)):=
\partial(\on{QQMaps}(\Y,\CT;\CE))\cap \on{QMaps}(\Y,\CT;\CE)$.

For example, when $\Y$ is a smooth curve (in which case there is no
difference between $\on{QMaps}$ and $\on{QQMaps}$), the locus
$\overset{\circ}{\on{QMaps}}(\Y,\CT;\CE)$ coincides with
$\on{Maps}(\Y,\CT;\CE)$. Indeed, the condition that the quotient
is torsion-free is equivalent in this case to its being a vector bundle.

\medskip

Recall that by $\on{Pic}^a(\Y)$ we denote connected components
of the Picard stack of $\Y$. Observe that
if the parameters $a,a'$ are such that
$\on{Pic}^{-a}(\Y)\times \on{Pic}^{-a'}(\Y)\subset
\on{TFree}^p_1(\Y)\times\on{TFree}^{p'}_1(\Y)$, then
$$\Y^{Q;(p,p')}\underset{\on{TFree}^p_1(\Y)\times\on{TFree}^{p'}_1(\Y)}\times
\left(\on{Pic}^{-a}(\Y)\times \on{Pic}^{-a'}(\Y)\right)\simeq
\on{Pic}^{-a'}(\Y)\times \Y^{(a-a')},$$
 where by $\Y^{(b)}$ we denote the scheme that classifies
pairs $(\CL,s)$, where $\CL\in \on{Pic}^b(\Y)$, and
$s:\CO_{\Y}\to \CL$ is a
non-zero section. Note that when $\Y$ is a smooth curve,
$\Y^{(b)}$ is just the corresponding symmetric power.

\medskip

Thus, from the previous proposition, we have a map
\begin{equation} \label{iota}
\ol{\iota}_{a',a}:\on{QMaps}^{a'}(\Y,\CT;\CE)\times \Y^{(a-a')}\to
\on{QMaps}^a(\Y,\CT;\CE),
\end{equation}
and a locally closed embedding $\iota_{a',a}:
\overset{\circ}{\on{QMaps}}{}^{a'}(\Y,\CT;\CE)\times \Y^{(a-a')}\to
\on{QMaps}^a(\Y,\CT;\CE)$.

Assume now that $\Y$ is smooth. We have:

\begin{lem}  
The map
$\ol{\iota}_{a',a}:\on{QMaps}^{a'}(\Y,\CT;\CE)\times \Y^{(a-a')}\to
\on{QMaps}^a(\Y,\CT;\CE)$ above is proper and finite.
Every geometric point of $\on{QMaps}^a(\Y,\CT;\CE)$ belongs to the
image of exactly one $\iota_{a',a}$.
\end{lem}

\begin{proof}

The properness assertion follows from the fact that the schemes
$\on{QMaps}^{a'}(\Y,\CT;\CE)$ and $\Y^{(a-a')}$ are proper when
$\Y$ is smooth. The finiteness follows from the fact that
given two line bundles $\CL\subset \CL'$ on $\Y$, the scheme
classifying line bundles $\CL''$, squeezed between the two, is
finite.

To prove the last assertion, it suffices to observe that
if $\CL$ is a torsion free sheaf of generic
rank $1$ on a smooth variety, embedded into $\CO_\Y\otimes \CE$,
such that the quotient is torsion free, then $\CL$ is a line bundle.

\end{proof}

The above assertions have an obvious based analog. The only
required modification is that instead of $\Y^{Q;(p,p')}$
(resp., $\Y^{(a-a')}$) we need to consider its open substack
corresponding to the condition that $\CL\to \CL'$ is an isomorphism in
a neighborhood of our subscheme $\Y'$ (resp., the section $s$ has
no zeroes on $\Y'$).

\ssec{Meromorphic quasi-maps}  \label{meromorphic}

In the sequel, we will have to consider the following generalization
of the notion of quasi-map. We will formulate it in the relative
situation.

Let $\Y$ be a flat projective scheme over a base $\X$, and let $\Y_1\subset
\Y$ be a relative Cartier divisor. Let $\CE$ be a vector bundle defined 
over $\Y-\Y_1$, and $\CT\subset \BP(\CE)$ a closed subscheme.

We define the functor $_{\infty\cdot \Y_1}\on{QMaps}(\Y,\CT;\CE)$ on the category 
of schemes over $\X$, to assign to a test scheme $S$ the set of pairs 
$(\CL,\kappa)$, where $\CL$ is a line bundle over
$\Y\underset{\X}\times S$, and $\kappa$ is a map of coherent sheaves
on $(\Y-\Y_1)\underset{\X}\times S$ 
$$\kappa:\CL|_{(\Y-\Y_1)\underset{\X}\times S}\to \CE\underset{\CO_\X}\otimes \CO_S,$$
such that for an open dense subset $U\subset (\Y-\Y_1)\underset{\X}\times S$,
over which $\kappa$ is a bundle map, the resulting map $U\to \BP(\CE)$
factors through $\CT$. 
 
\medskip

\begin{prop}  \label{representability of meromorphic}
The functor $_{\infty\cdot \Y_1}\on{QMaps}(\Y,\CT;\CE)$ is
representable by a strict ind-scheme of ind-finite type over $\X$. 
\end{prop}

\begin{proof}

Suppose that $\Y_2$ is another relative Cartier divisor containing
$\Y_1$. We have a closed embedding of functors 
$_{\infty\cdot \Y_1}\on{QMaps}(\Y,\CT;\CE)\hookrightarrow 
{}_{\infty\cdot \Y_2}\on{QMaps}(\Y,\CT;\CE)$.

Therefore, the question of ind-representability being local on $\X$,
by enlarging $\Y_1$, we can assume that $\CE$ is actually trivial;
we will denote by the same character its extension on the entire $\Y$.

Let $\CL_1$ be a line bundle on $\Y$ and $s:\CO\to \Y$ be a section
such that $\Y_1$ is its set of zeroes. If $\CL_1\in \on{Pic}_\X^{a_1}(\Y)$,
we have a natural closed embedding 
$$\on{QMaps}^a(\Y,\CT;\CE)\to \on{QMaps}^{a+a_1}(\Y,\CT;\CE),$$
as in \eqref{iota}, which sends a point $(\CL,\kappa)\in \on{QMaps}^a(\Y,\CT;\CE)$ to
$(\CL\otimes \CL_1^{-1},\CL\otimes \CL_1^{-1}\overset{s^{-1}}\longrightarrow
\CL\overset{\kappa}\longrightarrow \CE$).

It is now easy to see that $_{\infty\cdot \Y_1}\on{QMaps}(\Y,\CT;\CE)$
is representable by the inductive limit 
$\underset{\longrightarrow}{lim}\, \on{QMaps}^{a+n\cdot a_1}(\Y,\CT;\CE)$.

\end{proof}

If $\Y'\subset \Y$ is another closed subscheme, disjoint from $\Y_1$,
endowed with a section $\sigma':\Y'\to \CT$, in the same manner we define
the corresponding based version $_{\infty\cdot \Y_1}\on{QMaps}(\Y,\CT;\CE)_{\Y',\sigma'}$,
which is also a strict ind-scheme of ind-finite type.

\ssec{}        \label{pro-situation}

Let now $\CE$ be a pro-finite-dimensional vector space, i.e
$\CE=\underset{\leftarrow}{lim}\,\CE_i$.
(In the relative situation, we will assume that $\CE$
is an inverse limit of vector bundles on $\Y$.)

Let us recall the definition of the scheme $\BP(\CE)$ in
this case. By definition, for a test scheme $S$, a map
$S\to \BP(\CE)$ is a line bundle $\CL$ on $S$ together with a
map $S\to \CE$ (which by definition means a compatible system
of maps $S\to \CE_i$), such that locally on $S$, starting from
some index $i$, the map $S\to \CE_i$ is an injective bundle map.

By construction, $\BP(\CE)$ is a union of open sub-schemes
``$\BP(\CE)-\BP(\on{ker}(\CE\to \CE_i))$'', where each
such open is the inverse limit of
finite-dimensional schemes $\BP(\CE_j)-\BP(\on{ker}(\CE_j\to \CE_i))$,
$j\geq i$.

\medskip

When $\Y$ is a projective scheme, the scheme
$\on{Maps}(\Y,\BP(\CE))$ makes sense, and it is also
a union of schemes that are projective limits of schemes
of finite type.

The definition of $\on{QQMaps}(\Y,\BP(\CE))$
proceeds in exactly the same way as when $\CE$ is finite-dimensional.
Namely, for a test scheme $\on{Hom}(S,\on{QQMaps}(\Y,\BP(\CE))$
is the set of pairs $\sigma=(\CL,\kappa)$, where $\CL$
a coherent sheaf of generic rank $1$ on $\Y\times S$,
and $\kappa$ is a compatible system of
maps $\CL\to \CO_{\Y\times S}\otimes \CE_i$,
such that for every geometric point $s\in S$ starting from some
index $i$ the map $\CL_s\to \CO_\Y\otimes \CE_i$ is injective.

It is easy to show that $\on{QQMaps}(\Y,\BP(\CE))$
is a union of open
subschemes, each of which is a projective limit of schemes
of finite type. However, $\on{QQMaps}(\Y,\BP(\CE))$ itself
is generally not of finite type.

\medskip

If now $\CT$ is a closed subscheme on $\BP(\CE)$, we define
$\on{QQMaps}(\Y,\CT;\CE)$ as the corresponding closed
subscheme of $\on{QQMaps}(\Y,\BP(\CE))$. As before, the open
subscheme $\on{QMaps}(\Y,\CT;\CE)\subset\on{QQMaps}(\Y,\CT;\CE)$
is defined by the condition that $\CL$ is a line bundle.

\medskip

The following easily follows from \lemref{map of quasi-maps} and
\secref{canonical model}.
Let $\CT$ be a projective scheme of finite type embedded into
$\BP(\CE)$, where $\CE$ is a pro-finite-dimensional vector space.
Set $\CP:=\CO(1)_\CT$.

\begin{lem}  \label{embedding into pro}
We have a natural closed embedding
$\on{QMaps}(\Y,\CT;\CP)\to \on{QMaps}(\Y,\BP(\CE))$.
\end{lem}

\section{Quasi-maps into flag schemes and Zastava spaces}  \label{Zastava}

In this section $\fg$ will be an arbitrary Kac-Moody
algebra. Below we list some facts related to the flag
scheme (and certain partial flag schemes) attached to
$\fg$. These facts are established in \cite{k} when
$\fg$ is symmetrizable, but according to \cite{k1},
they hold in general. In the main body of the paper
we will only need the cases when $\fg$ is either
finite-dimensional, or non-twisted affine.

\ssec{} We will work in the following set-up:

Let $A=(A_{ij})_{i,j\in I}$ be a finite Cartan matrix.
We fix a root datum, that is, two finitely generated free abelian groups
$\Lambda,\Lambdach$ with a perfect bilinear pairing $\langle,\rangle:\
\Lambda\times\Lambdach\to\BZ$, and embeddings $I\subset\Lambdach,\
i\mapsto\alphach_i$ (simple roots), $I\subset\Lambda,\ i\mapsto\alpha_i$
(simple coroots) such that $\langle\alpha_i,\alphach_j\rangle=A_{ij}$.

We also assume that the subsets $I\subset\Lambda,\ I\subset\Lambdach$ are
both linearly independent, and that the subgroup of $\Lambda$ generated by
$\alpha_i$ is saturated (i.e., the quotient is torsion-free). In the
finite type case, the last condition is equivalent to the fact that
our root datum is simply-connected.

We denote by $\fg(A)=\fg$ (for short) the completed Kac-Moody
Lie algebra associated to the above datum (see ~\cite{k} ~3.2).
In particular, the Cartan subalgebra $\fh$ identifies with
$\Lambda\underset{\BZ}\otimes \BC$.
Sometimes we will denote the
corresponding root datum by $\Lambda_\fg,\Lambdach_\fg$
to stress the relation with $\fg$.

The semigroup of dominant coweights (resp. weights) will be denoted by
$\Lambda^+_\fg\subset\Lambda_\fg$
(resp. $\Lambdach^+_\fg\subset\Lambdach_\fg$).
We say that $\lambda\leq\mu$ if $\mu-\lambda$ is an integral
nonnegative linear combination of simple coroots $\alpha_i,\ i\in I$.
The semi-group of coweights $\lambda$, which are $\geq 0$ in this
sense will be denoted $\Lambda^{pos}_\fg$.

To a dominant weight $\lambdach\in\Lambdach^+_\fg$
one attaches the integrable
$\fg$-module, denoted $\CV_\lambdach$, with a fixed highest weight
vector $v_{\lambdach}\in \V_{\lambdach}$. For a pair of weights
$\lambdach_1,\lambdach_2\in\Lambdach^+_\fg$, there is a canonical map
$\V_{\lambdach_1+\lambdach_2}\to \V_{\lambdach_1}\otimes \V_{\lambdach_2}$
that sends
$v_{\lambdach_1+\lambdach_2}$ to $v_{\lambdach_1}\otimes v_{\lambdach_2}$.

\medskip

The Lie algebra $\fg$ has a triangular decomposition
$\fg=\fn\oplus\fh\oplus\fn^-$ (here $\fn$ is a profinite dimensional
vector space), and we have standard Borel subalgebras $\fb=\fn\oplus\fh,\
\fb^-=\fh\oplus\fn^-$.

If $\fp\subset\fg$ is a standard parabolic, we will denote by
$\fn(\fp)\subset \fn$ its unipotent radical, by $\fp^-\subset \fg$ (resp.,
$\fn(\fp^-)\subset \fn^-$) the corresponding opposite parabolic
(resp., its unipotent radical), and by $\fm(\fp):=\fp\cap\fp^-$
(or just $\fm$) the Levi factor. We will write $\fn(\fm)$ (resp.,
$\fn^-(\fm)$) for the intersections $\fm\cap\fn$ and
$\fm\cap\fn^-$, respectively.

If $i$ is an element of $I$, we will denote the corresponding
subminimal parabolic by $\fp_i$, and by $\fp_i^-$, $\fn_i:=\fn(\fp_i)$, 
$\fn^-_i:=\fn(\fp^-_i)$,
$\fm_i$, respectively, the corresponding associated subalgebras.
For a standard Levi subalgebra $\fm$, we will write that $i\in \fm$
if the corresponding simple root $\check \alpha_i$ belongs to $\fn(\fm)$.

From now on we will only work with parabolics corresponding
to subdiagrams of $I$ of {\it finite type}. In particular,
the Levi $\fm$ is finite-dimensional, and there exists
a canonically defined reductive group $M$, such that $\fm$ is its
Lie algebra. (When $\fp=\fb$, the corresponding Levi is the Cartan torus
$T$, whose set of cocharacters is our lattice $\Lambda_\fg$.)
By the assumption on the root datum, the derived group of $M$ is
simply-connected.

We will denote by $\Lambda_{\fg,\fp}$ the quotient of $\Lambda=\Lambda_\fg$
by the subgroup generated by $\on{Span}(\alpha_i,\, i\in \fm)$.
In other words, $\Lambda_{\fg,\fp}$ is the group of cocharacters of
$M/[M,M]$. By $\Lambda^{pos}_{\fg,\fp}$ we will denote the sub-semigroup
of $\Lambda_{\fg,\fp}$ equal to the positive span of the images
of $\alpha_i\in \Lambda_{\fg}\twoheadrightarrow \Lambda_{\fg,\fp}$
for $i\notin \fm$.

\medskip

In addition, there exists a pro-algebraic group $P$, with Lie algebra
$\fp$, with projects onto $M$, and the kernel $N(P)$ is pro-unipotent.

The entire group $G$ associated to the Lie algebra $\fg$, along
with its subgroup $P^-$, exists as a group-indscheme. Of course,
if $\fg$ is finite-dimensional, $G$ is the corresponding reductive
group. In the untwisted affine case, i.e., for 
 $\fg_{aff}=\fg((x))\oplus K\cdot \BC\oplus d\cdot \BC$
being the affinization of a finite-dimensional simple $\fg$ (cf. \secref{parabolic notation}),
the corresponding group-indscheme is $G_{aff}:=\widehat{G}\times \BG_m$, where
$\widehat{G}$ is the canonical central extension 
$1\to \BG_m\to \widehat{G}\to G((x))\to 1$ of the loop group-indscheme $G((x))$,
corresponding to the minimal $ad_\fg$-invariant scalar product on $\fg$,
such that the induced bilinear form on $\Lambda_\fg$ is integral-valued.

\ssec{Kashiwara's flag schemes}   \label{flag schemes}

Let $\fp\subset \fg$ be a standard parabolic, corresponding
to a subdiagram of $I$ of finite type. Following Kashiwara
(cf. ~\cite{k}), one defines the (partial) flag schemes $\CG_{\fg,\fp}$,
equipped with the action of $G$. Each $\CG_{\fg,\fp}$ has a
unit point $1_{\CG_{\fg,\fp}}$, and its stabilizer in $G$ equals
$P^-$, so that $\CG_{\fg,\fp}$ should be thought of as the quotient
``$G/P^-$''. 

\medskip

The scheme $\CG_{\fg,\fp}$ comes equipped with a closed embedding
(called the {\em Pl{\"u}cker embedding})
$$\CG_{\fg,\fp}\hookrightarrow\prod\BP(\CV_\lambdach^*),$$
where the product is being taken over $\Lambdach^+_{\fg,\fp}:=
\Lambdach_\fg^+\cap \Lambdach_{\fg,\fp}\subset \Lambdach_\fg$, and where each
$\CV_\lambdach^*$ is the profinite dimensional vector space dual
to the integrable module $\CV_\lambdach$.
The subscheme $\CG_{\fg,\fp}\hookrightarrow\prod\BP(\CV_\lambdach^*)$
is cut out by the so-called {\em Pl{\"u}cker equations}:

A collection of lines
$(\ell_\lambdach\subset\CV_\lambdach^*)_{\lambdach\in\Lambdach^+_{\fg,\fp}}$
satisfies Pl{\"u}cker equations if

\smallskip

(a) For the canonical morphism
$\V_{\lambdach_1+\lambdach_2}\to \V_{\lambdach_1}\otimes
\V_{\lambdach_2}$, the map $\CV_{\lambdach_1}^*\hat{\otimes}
\CV_{\lambdach_2}^*\to\CV_{\lambdach_1+\lambdach_2}^*$ sends
$\ell_{\lambdach_1}\otimes\ell_{\lambdach_2}$ to
$\ell_{\lambdach_1+\lambdach_2}$;

\smallskip

(b) For any $\fg$-morphism
$\V_{\much}\to \V_{\lambdach_1}\otimes
\V_{\lambdach_2}$ with $\much\neq \lambdach_1+\lambdach_2$, the map
$\CV_{\lambdach_1}^*\hat{\otimes}
\CV_{\lambdach_2}^*\to\CV_\much^*$ sends
$\ell_{\lambdach_1}\otimes\ell_{\lambdach_2}$ to $0$.

\smallskip

The point $1_{\CG_{\fg,\fp}}$ corresponds to the system of lines
$(\ell^0_\lambdach\subset\CV^*_\lambdach)$ where $\ell^0_\lambdach=
(\CV^*_\lambdach)^{\fn_-}$.

\medskip

The inverse image of the line bundle $\CO(1)$ on $\BP(\CV_\lambdach^*)$ is
the line bundle on $\CG_{\fg,\fp}$ denoted by $\CP^\lambdach_{\fg,\fp}$.
We have $\Gamma(\CG,\CP_{\fg,\fp}^\lambdach)=\CV_\lambda$.

\medskip

The orbits of the action of the pro-algebraic group
$N(P)$ on $\CG_{\fg,\fp}$ are parametrized by the double quotient
$\CW_\fm\backslash\CW_\fg/\CW_\fm$,
where $\CW_\fg$ and $\CW_\fm$ are the Weyl groups of $\fg$ and
$\fm$, respectively. For $w\in \CW_\fm\backslash\CW_\fg/\CW_\fm$, we
will denote by $\CG_{\fg,\fp,w}$ the corresponding orbit
(in our normalization, the unit point $1_{\CG_{\fg,\fp}}$ belongs
to $\CG_{\fg,\fp,e}$, where $e\in \CW_\fg$ is the unit element).
By $\ol{\CG}_{\fg,\fp,w}$ we will denote the closure of $\CG_{\fg,\fp,w}$
(e.g., $\ol{\CG}_{\fg,\fp,e}$ equals the entire $\CG_{\fg,\fp}$), and
by $\CG^w_{\fg,\fp}$ the open subscheme equal to union of orbits with
parameters $w'\leq w$ in the sense of the usual Bruhat order.
For each $w\in \CW_\fm\backslash\CW_\fg/\CW_\fm$ there exists
a canonical subgroup $N(P)_w\subset N(P)$, of finite codimension,
such that $N(P)_w$ acts freely on $\CG^w_{\fg,\fp}$ with a
finite-dimensional quotient.

\medskip

For every simple reflection $s_i\in \CW_\fg$ with $\alpha_i\notin
\fm$, we have a codimension $1$ subscheme $\CG_{\fg,\fp,s_i}$, and its
closure $\ol{\CG}_{\fg,\fp,s_i}$ is an effective Cartier divisor.
The union $\underset{i}\bigcup\,\ol{\CG}_{\fg,\fp,s_i}$ is called
{\em the Schubert divisor} and it equals the complement
$\CG_{\fg,\fp}-\CG^e_{\fg,\fp}$.

\medskip

Consider the affine cone $C(\CG_{\fg,\fp})$ over $\CG_{\fg,\fp}$
corresponding to the Pl{\"u}cker embedding. This is a closed subscheme
of
$\underset{\lambdach\in\Lambdach^+_{\fg,\fp}}\prod\,\CV_\lambdach^*$
consisting of collections of vectors
$(u_\lambdach\subset\CV_\lambdach^*)_{\lambdach\in\Lambdach^+_{\fg,\fp}}$
satisfying the Pl{\"u}cker equations:

\smallskip

(a) For the canonical morphism
$\V_{\lambdach_1+\lambdach_2}\to \V_{\lambdach_1}\otimes
\V_{\lambdach_2}$, the map $\CV_{\lambdach_1}^*\hat{\otimes}
\CV_{\lambdach_2}^*\to\CV_{\lambdach_1+\lambdach_2}^*$ sends
$u_{\lambdach_1}\otimes u_{\lambdach_2}$ to
$u_{\lambdach_1+\lambdach_2}$;

\smallskip

(b) For any $\fg$-morphism
$\V_{\much}\to \V_{\lambdach_1}\otimes
\V_{\lambdach_2}$ with $\much\neq \lambdach_1+\lambdach_2$, the map
$\CV_{\lambdach_1}^*\hat{\otimes}
\CV_{\lambdach_2}^*\to\CV_\much^*$ sends
$u_{\lambdach_1}\otimes u_{\lambdach_2}$ to $0$.

\medskip

We have a natural action of the torus $M/[M,M]$ on $C(\CG_{\fg,\fp})$.
Let $\overset{\circ}C(\CG_{\fg,\fp})$ be the open subset corresponding
to the condition that all $u_\lambdach\neq 0$. Then
$\overset{\circ}C(\CG_{\fg,\fp})$ is a principal $M/[M,M]$-bundle over
$\CG_{\fg,\fp}$.

\ssec{Quasi-maps into flag schemes}  \label{qmaps into flags}

Let $\bC$ be a smooth projective curve, with a marked
point $\infty_\bC\in \bC$, called ``infinity''. We will
denote by $\overset{\circ}\bC$ the complement
$\bC-\infty_\bC$.

We will use a short-hand notation $\on{QMaps}(\bC,\CG_{\fg,\fp})$
for the scheme of {\it based quasi-maps}
$\on{QMaps}(\bC,\CG_{\fg,\fp};\CV_\lambdach, \lambdach\in
\Lambdach^+_{\fg,\fp})_{\infty_\bC,\sigma_{const}}$,
where $\sigma_{const}:\infty_\bC\to \CG_{\fg,\fp}$  is the
constant map corresponding to $1_{\CG_{\fg,\fp}}\in \CG_{\fg,\fp}$.
In other words, $\on{QMaps}(\bC,\CG_{\fg,\fp})$ classifies based
maps from $\bC$ to the stack $C(\CG_{\fg,\fp})/(M/[M,M])$.
By $\on{Maps}(\bC,\CG_{\fg,\fp})$ we will denote the
open subscheme of maps.

The scheme $\on{QMaps}(\bC,\CG_{\fg,\fp})$ splits as a disjoint
union according the degree, which in our case is given by
elements $\theta\in \Lambda^{pos}_{\fg,\fp}$.

\medskip

For a fixed parameter $\theta$ as above, let
$\overset{\circ}\bC{}^\theta$ be the corresponding partially
symmetrized power of $\overset{\circ}\bC$. I.e.,
a point of $\overset{\circ}\bC{}^\theta$ assigns to every
$\lambdach\in \Lambdach^+_{\fg,\fp}$ an element of the
usual symmetric power $\overset{\circ}\bC{}^{(n)}$, where
$n=\langle \theta, \lambdach\rangle$, and this assignment
must be linear in $\lambdach$ in a natural sense.
Explicitly,
if $\theta=\underset{i\in I}\Sigma\, n_i\cdot \alpha_i$, for $i\notin \fm$,
$\overset{\circ}\bC{}^\theta\simeq \underset{i\in I}\prod\,
\overset{\circ}\bC{}^{(n_i)}$.

\medskip

From \secref{degenerations of quasi-maps}, we obtain that
for each $\theta',\theta-\theta'\in \Lambda^{pos}_{\fg,\fp}$, we have
a finite map
\begin{equation} \label{strata of quasi-maps}
\ol{\iota}_{\theta'}:\on{QMaps}^{\theta-\theta'}(\bC,\CG_{\fg,\fp})\times
\overset{\circ}\bC{}^{\theta'}\to
\on{QMaps}^\theta(\bC,\CG_{\fg,\fp}),
\end{equation}
such that the corresponding map
$$\iota_{\theta'}:\on{Maps}^{\theta-\theta'}(\bC,\CG_{\fg,\fp})\times
\overset{\circ}\bC{}^{\theta'}\to
\on{QMaps}^\theta(\bC,\CG_{\fg,\fp})$$
is a locally closed embedding. Moreover, the images of
$\iota_{\theta'}$ define a decomposition of
$\on{QMaps}^\theta(\bC,\CG_{\fg,\fp})$ into locally-closed pieces.
We will denote the embedding of the deepest stratum
$\overset{\circ}\bC{}^{\theta}\to\on{QMaps}^\theta(\bC,\CG_{\fg,\fp})$
by $\fs^\theta_\fp$.

Given a quasi-map $\sigma\in \on{QMaps}^\theta(\bC,\CG_{\fg,\fp})$,
which is the image of $$(\sigma'\in
\on{QMaps}^{\theta-\theta'}(\bC,\CG_{\fg,\fp}),
\theta_k\cdot \bc_k\in \overset{\circ}\bC{}^{\theta'}),$$ we will
say that (a) the defect of $\sigma$ is concentrated in $\underset{k}\cup\, \bc_k$,
(b) the defect of $\sigma$ at $\bc_k$ equals $\theta_k$,
(c) the total defect equals $\theta'$, and (d)
$\sigma'$ is the saturation of $\sigma$. Sometimes, we will assemble
statements (a), (b) and (c) into one by saying that the defect of
$\sigma$ is the colored divisor 
$\underset{k}\Sigma\, \theta_k\cdot \bc_k\in \overset{\circ}\bC{}^{\theta'}$.

\ssec{The finite-dimensional case} \label{fd Zastava}

Assume for a moment that $\fg$ is finite-dimensional.
Let $\Bun_G(\bC)$ (resp., $\Bun_G(\bC,\infty_\bC)$)
denote the stack of $G$-bundles on $\bC$
(resp., $G$-bundles trivialized at $\infty_\bC$).

We can consider a relative situation over $\Bun_G(\bC)$
(resp., over $\Bun_G(\bC,\infty_\bC)$),
when one takes maps or quasi-maps (resp., based maps or
quasi-maps) from $\bC$ to $\F_G\overset{G}\times \CG_{\fg,\fp}$,
where $\F_G$ is a point of $\Bun_G(\bC)$
(resp., $\Bun_G(\bC,\infty_\bC)$).

The corresponding stacks of maps identify, respectively,
with the stack $\Bun_{P^-}(\bC)$ of $P^-$-bundles on $\bC$, and
$\Bun_{P^-}(\bC,\infty_\bC)$--the stack $P^-$-bundles on $\bC$
trivialized at $\infty_\bC$.

In the non-based case, the stack of relative quasi-maps
is denoted by $\ol{\Bun}_{P^-}(\bC)$, and it was studied in
\cite{bg1} and \cite{bgfm}. The corresponding stack of relative
based quasi-maps will be denoted by
$\ol{\Bun}_{P^-}(\bC,\infty_\bC)$.

\ssec{Some properties of quasi-maps' spaces}

Let $\nu\in \Lambda_{\fg}$ be a coweight, such that
$\langle \nu,\alphach_i\rangle =0$ for $i\in \fm$, and
$\langle \nu, \alphach_i\rangle>0$ if $i\notin \fm$.
Consider the corresponding $1$-parameter subgroup
$\BG_m\to T$.
Since the point $1_{\CG_{\fg,\fp}}$ is $T$-stable,
we obtain a $\BG_m$-action on the scheme
$\on{QMaps}^\theta(\bC,\CG_{\fg,\fp})$.

In what follows, if a group $\BG_m$ acts on a
scheme $\Y$, we will say that the action contracts
$\Y$ to a subscheme $\Y'\subset \Y$, if:
(a) the action on $\Y'$ is trivial, and (b)
the action map extends to a morphism
$\BA^1\times \Y\to \Y$, such that $0\times \Y$ is
mapped to $\Y'$.

\begin{prop}  \label{contraction of quasi-maps}
The above $\BG_m$-action contracts
$\on{QMaps}^\theta(\bC,\CG_{\fg,\fp})$ to
$\overset{\circ}{\bC}{}^\theta\overset{\fs^\theta_{\fp}}
\subset \on{QMaps}^\theta(\bC,\CG_{\fg,\fp})$.
\end{prop}

\begin{proof}

The above $T$-action on
$\on{QMaps}^\theta(\bC,\CG_{\fg,\fp})$ corresponds
to an action of $T$ on the cone $C(\CG_{\fp,\fg})$
that takes a collection of vectors
$(u_\lambdach\in \CV^*_\lambdach)_{\lambdach\in \Lambdach^+_{\fg,\fp}}$
to $\lambdach(t)\cdot t\cdot u_\lambdach$, where
$t\cdot u$ denotes the $T$-action on the representation
$\CV^*_\lambdach$.

It is clear now that a $1$-parameter subgroup corresponding
to $\nu$ as in the proposition contracts each
$\CV^*_\lambdach$ to the line $\ell^0_\lambdach$.

However, the subscheme of quasi-maps, which map $\bC$
to $(\ell^0_\lambdach)_{\lambdach\in \Lambdach^+_{\fg,\fp}}\subset
C(\CG_{\fp,\fg})$ coincides with $\fs_{\fp}(\overset{\circ}{\bC}{}^\theta)$.

\end{proof}

\begin{prop}  \label{locally of finite type}
The scheme $\on{Maps}^\theta(\bC,\CG_{\fg,\fp})$
is a union of open subschemes of finite type.
\end{prop}

\noindent{\it Remark.}
One can prove that the scheme
$\on{Maps}^\theta(\bC,\CG_{\fg,\fb})$ is in fact globally 
of finite type, at least when $\fg$ is symmetrizable, 
cf. \secref{AppA}, but the proof is more involved.
Of course, when $\fg$ is finite-dimensional, this is obvious.
When $\fg$ is (untwisted) affine, another proof will 
be given in the next section, using a modular
interpretation of $\CG_{\fg,\fb}$ via bundles on
the projective line. 

\begin{proof}

Let us first remark that if $\Y$ and $\CT$ are projective
schemes, and $\CT^0\subset \CT$ is an open subscheme, the
subfunctor of $\on{Maps}(\Y,\CT)$ consisting of maps landing in
$\CT^0$ is in fact an open subscheme.

\medskip

For an element $w\in\CW_\fm\backslash\CW/\CW_\fm$,
let $\on{Maps}^\theta(\bC,\CG^w_{\fg,\fp})$ be the open
subset in $\on{Maps}^\theta(\bC,\CG_{\fg,\fp})$
of maps that land in the open subset
$\CG^w_{\fg,\fp}\subset \CG_{\fg,\fp}$. Evidently,
$\on{Maps}^\theta(\bC,\CG_{\fg,\fp})\simeq \underset{w}\cup\,
\on{Maps}^\theta(\bC,\CG^w_{\fg,\fp})$.
We claim that each
$\on{Maps}^\theta(\bC,\CG^w_{\fg,\fp})$ is of finite type.

Recall that the subgroup $N(P)_w\subset N(P)$ acts freely
on $\CG^w_{\fg,\fp}$ and that the quotient $\CG^w_{\fg,\fp}/N(P)_w$
is a quasi-projective scheme of finite type. Hence,
the scheme (of based maps)
$\on{Maps}(\bC,\CG^w_{\fg,\fp}/N(P)_w)$ is of
finite type.

We claim now that the natural projection map
$$\on{Maps}^\theta(\bC,\CG^w_{\fg,\fp})\to \on{Maps}(\bC,\CG^w_{\fg,\fp})/N(P)_w$$
is a closed embedding. (We will show moreover that
this map is actually an isomorphism if $\bC$ is of genus $0$.)

\medskip

First, let us show that this map is injective on the level
of $S$-points for any $S$.
Let $\sigma_1,\sigma_2$ be two based maps $\bC\times S\to \CG^w_{\fg,\fp}$
which project to the same map to $\CG^w_{\fg,\fp}/N(P)_w$.
Then by definition, we obtain a map
$\wt{\sigma}:\bC\times S\to N(P)_w$ with $\wt{\sigma}|_{\infty_\bC\times S}\equiv 1$.
Since $\bC$ is projective and $N(P)_w$ is an inverse
limit of groups which are extensions of $\BG_a$, we obtain that
$\wt{\sigma}\equiv 1$.

\medskip

Now let us show that the map
$\on{Maps}^\theta(\bC,\CG^w_{\fg,\fp})\to
\on{Maps}(\bC,\CG^w_{\fg,\fp})/N(P)_w$ is proper by checking
the valuative criterion.

Let $\sigma'$ be a based map $\bC\times \bX\to
\CG^w_{\fg,\fp}/N(P)_w$,
where $\bX$ is an affine curve such that the restriction
$\sigma'|_{\bC\times (\bX-0_\bX)}$ lifts to a map
to $\CG^w_{\fg,\fp}$. Let us filter $N(P)_w$ by normal
subgroups $N(P)_w=N_0\supset N_1\supset\ldots$
with associated graded quotients isomorphic to the
additive group $\BG_a$. Set $\sigma'_0=\sigma'$, and
assume that we have found a lifting of $\sigma'$
to a based map $\sigma'_i:\bC\times \bX\to\CG^w_{\fg,\fp}/N_i$.
Then the obstruction to lifting it to the next level lies in
$H^1(\bC\times \bX,\CO(-\infty_\bC)\boxtimes \CO_\bX)$.
Since $\bX$ is affine, the last group is isomorphic to
the space of $H^1(\bC,\CO(-\infty_\bC))$-valued functions
on $\bX$. By assumption, the function corresponding to
the obstruction class vanishes on the open subset
$\bX-0_\bX$, therefore it vanishes.

\medskip

Note that when $\bC$ is of genus $0$, and an arbitrary affine
base $S$, the same argument shows that any based map
$\sigma':\bC\times S\to \CG^w_{\fg,\fp}/N(P)_w$ lifts to a
map $\bC\times \bX\to\CG^w_{\fg,\fp}$, because in this case
$H^1(\bC,\CO(-\infty_\bC))=0$.

\end{proof}

\medskip

Let $\fp\subset \fp'$ be a pair of standard parabolics, and
set $\fp(\fm')=\fm'\cap \fp$ to be the corresponding
parabolic in $\fm'$. Note that we have an exact sequence
$$0\to \Lambda_{\fm',\fp(\fm')}\to
\Lambda_{\fg,\fp}\to \Lambda_{\fg,\fp'}\to 0.$$

\begin{lem} \label{two parabolics}
We have a natural map
$\on{QMaps}^\theta(\bC,\CG_{\fg,\fp})\to
\on{QMaps}^{\theta'}(\bC,\CG_{\fg,\fp'})$, where
$\theta'$ is the image of $\theta$ under the above map of
lattices. If $\theta'=0$, there is an isomorphism
$\on{QMaps}^\theta(\bC,\CG_{\fg,\fp})\simeq
\on{QMaps}^\theta(\bC,\CG_{\fm',\fp(\fm')})$.
\end{lem}

\begin{proof}

The existence of the map $\on{QMaps}^\theta(\bC,\CG_{\fg,\fp})\to
\on{QMaps}^{\theta'}(\bC,\CG_{\fg,\fp'})$ is immediate from
the definitions.

Note that if $\theta'=0$,
$\on{QMaps}^{\theta'}(\bC,\CG_{\fg,\fp'})$ is a point-scheme
that corresponds to the constant map $\bC\to 1_{\CG_{\fg,\fp'}}\in
\CG_{\fg,\fp'}$. Note also that the preimage of $1_{\CG_{\fg,\fp'}}$
in $\CG_{\fg,\fp}$ identifies naturally with
$\CG_{\fm',\fm'(\fp)}$. In particular, for $\theta$ as above,
we have a closed embedding
$\on{QMaps}^\theta(\bC,\CG_{\fm',\fp(\fm')})\to
\on{QMaps}^\theta(\bC,\CG_{\fg,\fp})$.

To see that it is an isomorphism, note that for
any $\sigma\in \on{QMaps}^\theta(\bC,\CG_{\fg,\fp})$ there
exists an open dense subset $U\subset \bC\times S$, such that the
map $\sigma|_U$ projects to the constant map $U\to\CG_{\fg,\fp'}$,
and hence, has its image in $\CG_{\fm',\fp(\fm')}$. But this implies
that $\sigma$ itself is a quasi-map into $\CG_{\fm',\fp(\fm')}$.

\end{proof}

Assume for example that $\fp=\fb$ and that
$\mu\in \Lambda^{pos}_\fg=\Lambda^{pos}_{\fg,\fb}$ equals
$\alpha_i$ for some $i\notin \fm$. By putting
$\fp'=\fp_i$, from the previous lemma we obtain that the space
$\on{QMaps}^\theta(\bC,\CG_{\fg,\fb})$ in this case identifies
with the space of based quasi-maps $\bC\to \BP^1$ of degree $1$.
In particular, it is empty unless $\bC$ is of genus $0$, and
for $\bC\simeq \BP^1$, it is isomorphic to $\BA^2$.

\medskip

Note that we have a natural map
$\on{Maps}(\bC,\CG_{\fg,\fp})\to \Bun_M(\bC,\infty_\bC)$. Indeed,
for a based map $\bC\times S\to \CG_{\fg,\fp}$ we define
an $M$-bundle on $\bC\times S$ by taking the Cartesian product
$$(\bC\times S)\underset{\CG_{\fg,\fp}}\times
\overset{\circ}{\wt{\CG}}_{\fg,\fp},$$
where $\overset{\circ}{\wt{\CG}}_{\fg,\fp}$ is as in \secref{enhanced}.

The corresponding schemes $\on{Maps}(\bC,\CG_{\fg,\fp})$
and $\on{Maps}(\bC,\CG_{\fg,\fp'})$ are related in the following
explicit way:

Let $P^-(M')$ be the parabolic in $M'$ corresponding to $\fp^-(\fm')$. Let
$\Bun_{M'}(\bC,\infty_\bC)$, be the stack as in \secref{fd Zastava}, and let
$\Bun_{P^-(M')}(\bC,\infty_\bC)\subset \ol{\Bun}_{P^-(M')}(\bC,\infty_\bC)$
be the corresponding stacks of $P^-(M')$-bundles and generalized
$P^-(M')$-bundles, respectively.

The following is straightforward from the definitions:

\begin{lem}  \label{relation between parabolics}
There exists an open embedding
$$\on{Maps}(\bC,\CG_{\fg,\fp'})\underset{\Bun_{M'}(\bC,\infty_\bC)}
\times \ol{\Bun}_{P^-(M')}(\bC,\infty_\bC)\hookrightarrow
\on{QMaps}(\bC,\CG_{\fg,\fp}),$$
whose image is the union of $\on{Im}(\iota_{\theta'})$ over
$\theta'\in \Lambda^{pos}_{\fm',\fp(\fm')}$.
Moreover, the map
$$\on{Maps}(\bC,\CG_{\fg,\fp'})\underset{\Bun_{M'}(\bC,\infty_\bC)}
\times \Bun_{P^-(M')}(\bC,\infty_\bC)\hookrightarrow
\on{Maps}(\bC,\CG_{\fg,\fp})$$
is an isomorphism.
\end{lem}

For $\fp=\fb$ and a vertex $i\in I$, let us denote by
$\partial(\on{QMaps}(\bC,\CG_{\fg,\fb}))_i$ the locally closed
subset equal to $\on{Im}(\iota_{\alpha_i})$.

The above lemma combined with \propref{locally of finite type}
yield the following:

\begin{cor}  \label{codimension 1 locus}
The open subscheme
$\on{Maps}(\bC,\CG_{\fg,\fb})\bigcup \left(\underset{i\in I}\cup\,
\partial(\on{QMaps}(\bC,\CG_{\fg,\fb}))_i\right)$ is (locally) of finite
type.
\end{cor}

\ssec{Projection to the configuration space}  \label{projection to configurations}

We claim that there exists a canonical map
$\varrho^\theta_\fp:\on{QMaps}^\theta(\bC,\CG_{\fg,\fp})\to
\overset{\circ}\bC{}^\theta$. Indeed, let
$\sigma$ be an $S$-point of $\on{QMaps}^\theta(\bC,\CG_{\fg,\fp})$,
i.e., for every $\lambdach\in \Lambdach^+_{\fg,\fp}$ we have
a line bundle $\CL^\lambdach$ on $\bC\times S$ and a map
$\kappa^\lambdach:\CL^\lambdach\to \CO_{\bC\times S}\otimes
\CV_\lambdach^*$. Recall that the choice of the standard parabolic
$\fp$ defines a highest weight vector in $\CV_\lambdach$; hence
by composing $\kappa^\lambdach$ with the corresponding map
$\CO_{\bC\times S}\otimes \CV_\lambdach^*\to \CO_{\bC\times S}$,
we obtain a map $\CL^\lambdach\to \CO$, i.e. a point of
$\bC^{(\langle \theta, \lambdach\rangle)}$. The obtained divisor
avoids $\infty_\bC$ because the composition
$\ell^0_\lambdach\to \CV_\lambdach^*\to \BC$ is non-zero.

\medskip

The above construction of $\varrho^\theta_\fp$ can be alternatively
viewed as follows. Recall the cone $C(\CG_{\fg,\fp})$, and observe that
for each vertex $i\in I$ with $i\notin \fm$, we have a function
on $C(\CG_{\fg,\fp})$ corresponding to the map
$\CV_{\omegach_i}\to \BC$, where $\omegach_i$ is the fundamental weight.
Let us denote by $C(\ol{\CG}_{\fg,\fp,s_i})$ the corresponding
Cartier divisor. Note that the intersection
$\overset{\circ}C(\CG_{\fg,\fp})\cap C(\ol{\CG}_{\fg,\fp,s_i})$
is the preimage of the Cartier divisor $\ol{\CG}_{\fg,\fp,s_i}\subset
\CG_{\fg,\fp}$.

An $S$-point of $\on{QMaps}^\theta(\bC,\CG_{\fg,\fp})$ is
a map $\sigma:\bC\times S\to C(\CG_{\fg,\fp})/(M/[M,M])$,
and by pulling back $C(\ol{\CG}_{\fg,\fp,s_i})$, we obtain a
Cartier divisor on $\bC\times S$, which by the conditions on $\sigma$
is in fact a relative Cartier divisor over $S$, i.e. it gives rise
to a map from $S$ to the suitable symmetric power of $\bC$.

\medskip

Note that the composition of the map $\fs^\theta_\fp$ and
$\varrho^\theta_\fp$ is the identity map on
$\overset{\circ}\bC{}^\theta$. More generally, for
$\theta'\in \Lambda^{pos}_{\fg,\fp}$, the composition
\begin{equation} \label{configuration of saturation}
\varrho^\theta_\fp\circ\iota_{\theta'}:\on{QMaps}^{\theta-\theta'}
(\bC,\CG_{\fg,\fp})\times \overset{\circ}\bC{}^{\theta'}\to
\overset{\circ}\bC{}^\theta
\end{equation}
covers the addition map
$\overset{\circ}\bC{}^{\theta-\theta'}\times
\overset{\circ}\bC{}^{\theta'}\to \overset{\circ}\bC{}^\theta$.

\ssec{Zastava spaces}  \label{Zastava spaces}

We will introduce twisted versions of the schemes
$\on{QMaps}^\theta(\bC,\CG_{\fg,\fp})$, called
{\it Zastava spaces} $\CZ^\theta_{\fg,\fp}(\bC)$, cf. \cite{bgfm}.
We will first define $\CZ^\theta_{\fg,\fp}(\bC)$ as a functor,
and later show that it is representable by a scheme.

Let $\bC$ be a smooth curve (not necessarily complete).
An $S$-point of $\CZ^\theta_{\fg,\fp}(\bC)$ is a data of
$(D^\theta,\F_{N(P)},\kappa)$, where

\begin{itemize}

\item
$D^\theta$ is an $S$-point of $\bC^\theta$.
In particular, we obtain a principal $M/[M,M]$-bundle
$\F_{M/[M,M]}$ on $\bC\times S$, such that for every
$\lambdach\in \Lambdach_{\fg,\fp}$ the associated line bundle
denoted $\CL^\lambdach_{\F_{M/[M,M]}}$ equals
$\CO_{\bC\times S}(-\lambdach(D^\theta))$.

\item
$\F_{N(P)}$ is a principal
$N(P)$-bundle on $\bC\times S$. Note that this makes sense,
since $N(P)$ is a pro-algebraic group. In particular,
for every $N(P)$-integrable $\fg$-module $\CV$, we
can form a pro-vector bundle $(\CV^*)_{\F_{N(P)}}$
on $\bC\times S$. If $v\in \CV$ is an $N(P)$-invariant
vector, it gives rise to a map
$(\CV^*)_{\F_{N(P)}}\to \CO_{\bC\times S}$.

\item
$\kappa$ is a system of maps
$\kappa^\lambdach:\CL^\lambdach_{\F_{M/[M,M]}}\to (\CV_\lambdach^*)_{\F_{N(P)}}$, which
satisfy the Pl{\"u}cker relations, and such that the
composition of $\kappa^\lambdach$ and the projection
$(\CV_\lambdach^*)_{\F_{N(P)}}\to \CO_{\bC\times S}$ corresponding
to the h.w. vector $v_\lambdach\in \CV_\lambdach$ is the
tautological map $\CO_{\bC\times S}(-\lambdach(D^\theta))\to \CO_{\bC\times S}$.

\end{itemize}

We will denote by $\overset{\circ}{\CZ}{}^\theta_{\fg,\fp}(\bC)$
the open subfunctor in $\CZ^\theta_{\fg,\fp}(\bC)$ corresponding
to the condition that the $\kappa^\lambdach$'s are injective bundle
maps.

\medskip

Another way to spell out the definition of $\CZ^\theta_{\fg,\fp}(\bC)$
is as follows:
Let us consider the quotient $N(P)\backslash C(\CG_{\fg,\fp})/(M/[M,M])$. This
is a non-algebraic stack (i.e. all the axioms, but the one
about covering by a scheme, hold). It contains as an open substack
$N(P)\backslash\CG^e_{\fg,\fp}\simeq pt$, and the Cartier divisors
$N(P)\backslash C(\ol{\CG}_{\fg,\fp,s_i})/(M/[M,M])$ for $i\in I$.

It is easy to see that an $S$-point of $\CZ^\theta_{\fg,\fp}(\bC)$
is the same as a map $$\bC\times S\to N(P)\backslash
C(\CG_{\fg,\fp})/(M/[M,M]),$$ such that for every geometric point
$s\in S$, the map $\bC\simeq \bC\times s\to
N(P)\backslash C(\CG_{\fg,\fp})/(M/[M,M])$ sends the generic point
of $\bC$ to $N(P)\backslash\CG^e_{\fg,\fp}\simeq pt$,
and the intersection
$$\bC\cap \left(N(P)\backslash C(\ol{\CG}_{\fg,\fp,s_i})/M/[M,M]\right)$$
is a divisor of degree $\langle \theta,\omegach_i\rangle$ on $\bC$.

From this description of $\CZ^\theta_{\fg,\fp}(\bC)$ we obtain the
following:

\begin{lem}   \label{trivial outside graph}
For an $S$-point $(D^\theta,\F_{N(P)},\kappa)$ of
$\CZ^\theta_{\fg,\fp}(\bC)$, the $N(P)$-bundle $\F_{N(P)}$ admits a
canonical trivialization away from the support of $D^\theta$.
Moreover, in terms of this trivialization,
over this open subset the maps $\kappa^\lambdach$ are
constant maps corresponding to $v_\lambdach\in\ell_\lambdach^0
\subset \CV^*_\lambdach$.
\end{lem}

For the proof it suffices to observe that
for a map $\bC\times S\to N(P)\backslash
C(\CG_{\fg,\fp})/(M/[M,M])$, the open subset
$\bC\times S-\on{supp}(D^\theta)$ is exactly the preimage
of $pt=N(P)\backslash\CG^e_{\fg,\fp}$.

\medskip

Let us consider now the following set-up, suggested in this
generality by Drinfeld.
Let $\CT$ be a (non-necessarily algebraic) stack, which contains
an open substack $\CT^0$ isomorphic to $pt$.
Let $S$ be a scheme, embedded as an open subscheme into a scheme
$S_1$, and let $\sigma:S\to \CT$ be a map such that
$\sigma^{-1}(\CT-\CT^0)$ is closed in $S_1$.

\begin{lem} \label{extension}
There is a canonical bijection between the set
of maps $\sigma$ as above and the set of
maps $\sigma_1:S_1\to \CT$ such that
$\sigma_1^{-1}(\CT-\CT^0)$ is contained in $S$.
\end{lem}

\begin{proof}

Of course, starting from $\sigma_1$, we get the corresponding
$\sigma$ by restriction to $S$.

Conversely, for $\sigma:\Y\times S\to \CT$ as above,
consider the two open subsets $S$ and
$S_1-(\sigma^{-1}(\CT-\CT^0))$, which
cover $S_1$.

By setting $\sigma_1|_S=\sigma$,
$\sigma_1|_{S_1-(\sigma^{-1}(\CT-\CT^0))}$
to be the constant map to $pt=\CT^0$, we have
a gluing data for $\sigma_1$.

\end{proof}

We apply this lemma in the following situation:

\begin{cor}  \label{extension of Zastava}
Let $\bC\to \bC_1$ be an open embedding of curves.
We have an isomorphism
$$\CZ^\theta_{\fg,\fp}(\bC)\simeq
\CZ^\theta_{\fg,\fp}(\bC_1)\underset{\bC_1^\theta}\times \bC^\theta.$$
\end{cor}

The proof follows from the fact that for
$S$-point of $\CZ^\theta_{\fg,\fp}(\bC)$, the support of
the colored divisor $D^\theta$, which is the same as the
preimage of $N(P)\backslash (C(\CG_{\fg,\fp}))/(M/[M,M])-pt$, is
finite over $S$.

\ssec{Factorization principle}  \label{factorization principle}

Let $pt=\CT^0\subset \CT$ be an embedding of stacks as before,
and suppose now that the complement $\CT-\CT^0$ is a union of Cartier
divisors $\CT_i$, for $i$ belonging to some set of indices $I$.

For a set of non-negative integers $\ol{n}=n_i, i\in I$, consider the
functor $\on{Maps}^{\ol{n}}(\bC,\CT)$ that assigns to a scheme $S$
the set of maps $$\sigma:\bC\times S\to \CT,$$
such that each $\sigma^{-1}(\CT_i)\subset \bC\times S$ is
a relative (over $S$) Cartier divisor of degree $n_i$.
It is easy to see that $\on{Maps}^{\ol{n}}(\bC,\CT)$ is a
sheaf in the faithfully-flat topology on the category of schemes.
\footnote{According to Drinfeld, one can formulate a general
hypothesis on $\CT$, under which the functor
$\on{Maps}^{\ol{n}}(\bC,\CT)$ is representable by a scheme.}

By construction, we have a map (of functors)
$\on{Maps}^{\ol{n}}(\bC,\CT)\to \underset{i\in I}\prod\, \bC^{(n_i)}$.
The following factorization principle is due to Drinfeld:

\begin{prop}\label{factorization pattern}
Let $n_i=n'_i+n''_i$ be a decomposition, and let
$(\underset{i\in I}\prod\, \bC^{(n'_i)}\times
\underset{i\in I}\prod\, \bC^{(n''_i)})_{disj}$ be the open subset
corresponding to the condition that the divisors $D'_i\in \bC^{(n'_i)}$
and $D''_i\in \bC^{(n''_i)}$ have disjoint supports. We have a
canonical isomorphism:
\begin{align*}
&\on{Maps}^{\ol{n}}(\bC,\CT)\underset{\underset{i\in I}\prod\, \bC^{(n_i)}}\times
(\underset{i\in I}\prod\, \bC^{(n'_i)}\times
\underset{i\in I}\prod\, \bC^{(n''_i)})_{disj}\simeq \\
&(\on{Maps}^{\ol{n}'}(\bC,\CT)\times \on{Maps}^{\ol{n}''}(\bC,\CT))
\underset{\underset{i\in I}\prod\, \bC^{(n'_i)}\times
\underset{i\in I}\prod\, \bC^{(n''_i)}}\times
(\underset{i\in I}\prod\, \bC^{(n'_i)}\times
\underset{i\in I}\prod\, \bC^{(n''_i)})_{disj}.
\end{align*}
\end{prop}

\begin{proof}

Given an $S$-point $(\sigma, D',D'')$ of
$$\on{Maps}^{\ol{n}}(\bC,\CT)\underset{\underset{i\in I}\prod\, \bC^{(n_i)}}\times
(\underset{i\in I}\prod\, \bC^{(n'_i)}\times
\underset{i\in I}\prod\, \bC^{(n''_i)})_{disj}$$
we produce the $S$-points $\sigma'$ and $\sigma''$ of
$\on{Maps}^{\ol{n}'}(\bC,\CT)$ and $\on{Maps}^{\ol{n}''}(\bC,\CT)$
as follows:

As a map $\bC\times S\to \CT$, $\sigma'$ is set to be
equal to $\sigma$ on the open subset $\bC\times S-\on{supp}(D'')$,
and to the constant map $\bC\times S\to \CT^0=pt$ on the open
subset $\bC\times S-\on{supp}(D')$. By assumption, this gives
a well-defined gluing data for $\sigma'$. The map $\sigma''$
is defined by interchanging the roles of $'$ and $''$.

Conversely, given $\sigma'$ and $\sigma''$, with the corresponding
colored divisors $D'$ and $D''$, respectively, we define
$\sigma$ as follows. On the open subset $\bC\times S-\on{supp}(D'')$,
$\sigma$ is set to be equal to $\sigma'$, and on the open subset
$\bC\times S-\on{supp}(D')$, $\sigma$ is set to be equal to
$\sigma''$. Since $\sigma'$ and $\sigma''$ agree on
$\bC\times S-\on{supp}(D'+D'')$, this
gives a well-defined gluing data for $\sigma$.

\end{proof}

The above proposition admits the following generalization. Let
$p:\wt{\bC}\to \bC$ be an \'etale cover. Let $\wt{\bC}{}^{(n)}_{disj}\subset
\wt{\bC}{}^{(n)}$ be the open subset that consists of divisors
$\wt{D}$ on $\wt{\bC}$, such that the divisor $p^*(p_*(\wt{D}))-\wt{D}$ is
disjoint from $\wt{D}$. 

Essentially the same proof gives the following:
\begin{prop} \label{etale factorization}
There is a canonical isomorphism:
$$\on{Maps}^{\ol{n}}(\bC,\CT)\underset{\underset{i\in I}\prod\, \bC^{(n_i)}}
\times \underset{i\in I}\prod\, \wt{\bC}{}^{(n_i)}_{disj}\simeq
\on{Maps}^{\ol{n}}(\wt{\bC},\CT)\underset{\underset{i\in I}\prod\, \wt{\bC}{}^{(n_i)}}
\times \underset{i\in I}\prod\, \wt{\bC}{}^{(n_i)}_{disj}.$$
\end{prop}

As an application, we take $\CT=N(P)\backslash
(C(\CG_{\fg,\fp}))/(M/[M,M])$, and we obtain the following factorization
property of the Zastava spaces:

\begin{prop}  \label{factorization of Zastava}
For an \'etale map $p:\wt{\bC}\to \bC$ we have a canonical isomorphism
$$\CZ^\theta_{\fg,\fp}(\bC)\underset{\bC^\theta}\times
\wt{\bC}{}^\theta_{disj}\simeq
\CZ^\theta_{\fg,\fp}(\wt{\bC})\underset{\wt{\bC}{}^\theta}\times
\wt{\bC}{}^\theta_{disj}.$$
In particular, for $\theta=\theta_1+\theta_2$, we have a canonical isomorphism
$$\CZ^\theta_{\fg,\fp}(\bC)\underset{\bC^\theta}\times
(\bC^{\theta_1}\times \bC^{\theta_2})_{disj}\simeq
(\CZ^{\theta_1}_{\fg,\fp}(\bC)\times\CZ^{\theta_2}_{\fg,\fp}(\bC))
\underset{\bC^{\theta_1}\times \bC^{\theta_2}}\times
(\bC^{\theta_1}\times \bC^{\theta_2})_{disj}.$$

Moreover, the above isomorphisms preserve the loci
$\overset{\circ}{\CZ}{}^\theta_{\fg,\fp}(\bC)
\subset\CZ^\theta_{\fg,\fp}(\bC)$, $\overset{\circ}{\CZ}{}^\theta_{\fg,\fp}(\wt{\bC})
\subset\CZ^\theta_{\fg,\fp}(\wt{\bC})$
\end{prop}

This proposition together with \lemref{extension} expresses
the locality property of the Zastava spaces $\CZ^\theta_{\fg,\fp}(\bC)$
with respect to $\bC$.

\medskip

Another application of \propref{factorization pattern}
is \propref{Nak fac}, where $\CT$ is taken to be the stack
of coherent sheaves of generic rank $n$ on $\BP^1$ with
a trivialization at $\infty$.

\ssec{The case of genus $0$}

Let now $\bC$ be a projective line, and
$\overset{\circ}\bC=\bC-\infty_\bC$ the corresponding affine
line.

The following proposition will play a key role in this paper.

\begin{prop} \label{Zastava and quasi-maps}
The functor represented by
$\on{QMaps}^\theta(\bC,\CG_{\fg,\fp})$ is naturally
isomorphic to $\CZ^\theta_{\fg,\fp}(\overset{\circ}\bC)$.
\end{prop}

\begin{proof}

By the very definition of $\CZ^\theta_{\fg,\fp}(\overset{\circ}\bC)$,
the assertion of the proposition amounts to the following:
For an $S$-point $(D^\theta,\F_{N(P)},\kappa)$
of $\CZ^\theta_{\fg,\fp}(\overset{\circ}\bC)$, the $N(P)$-bundle
$\F_{N(P)}$ on $\overset{\circ}\bC\times S$ can be canonically
trivialized.

Note that according to \corref{extension of Zastava},
$\F_{N(P)}$ extends to a principal $N(P)$-bundle on the
entire $\bC\times S$, with a trivialization on
$\infty_\bC\times S$. But since $N(P)$ is pro-unipotent, and
$\bC$ is of genus $0$, such a trivialization extends uniquely
on $\bC\times S$.

(Note, however, that this trivialization and the one coming from
\lemref{trivial outside graph} DO NOT agree, but they do
coincide on $\infty_\bC\times S$.)

\end{proof}

\noindent{\it Remark.} For a curve of an arbitrary genus,
it is easy to see that the natural map
$\on{QMaps}^\theta(\bC,\CG_{\fg,\fp})\to
\CZ^\theta_{\fg,\fp}(\overset{\circ}\bC)$ is a closed embedding.

Indeed, if $N(P)$-bundle $\F_{N(P)}$ trivialized at $\infty_\bC$
can be trivialized globally, this can be done in a unique
fashion (since $H^0(\bC,\CO_\bC(-\infty_\bC))=0$), and
the property that it admits such a trivialization
is a closed condition.

\medskip

As a corollary of \propref{Zastava and quasi-maps},
we obtain the following statement:

\begin{cor}
The functor $\CZ^\theta_{\fg,\fp}(\bC)$
is representable by a scheme for any curve $\bC$.
The corresponding open subscheme
$\overset{\circ}{\CZ}{}^\theta_{\fg,\fp}(\bC)$ is
locally of finite type.
\end{cor}

\begin{proof}

Using \corref{extension of Zastava} and
\propref{factorization of Zastava}, we reduce
the assertions of the proposition to the case
when the curve in question is the projective line.
The representability follows now from
\propref{Zastava and quasi-maps}. The fact that
$\overset{\circ}{\CZ}{}^\theta_{\fg,\fp}(\bC)$
is of finite type follows from 
\propref{locally of finite type}.

\end{proof}

\ssec{Further properties of maps' spaces}

In this subsection (until \secref{enhanced}) we will assume that $\bC$
is of genus $0$, and establish certain properties
of the scheme $\on{Maps}^\theta(\bC,\CG_{\fg,\fp})$.
Using \propref{Zastava and quasi-maps}, the same assertions
will hold for the space $\overset{\circ}{\CZ}{}^\theta_{\fg,\fp}(\bC)$
on any curve $\bC$.

First, note that \propref{Zastava and quasi-maps}
and \propref{factorization of Zastava} yield the following
factorization property of $\on{QMaps}^\theta_{\fg,\fp}(\bC)$
with respect to the projection $\varrho^\theta_\fp$:
\begin{equation} \label{factorization of quasi-maps}
\on{QMaps}^\theta_{\fg,\fp}(\bC)\underset{\bC^\theta}\times
(\bC^{\theta_1}\times \bC^{\theta_2})_{disj}\simeq
(\on{QMaps}^{\theta_1}_{\fg,\fp}(\bC)\times\on{QMaps}^{\theta_2}_{\fg,\fp}(\bC))
\underset{\bC^{\theta_1}\times \bC^{\theta_2}}\times
(\bC^{\theta_1}\times \bC^{\theta_2})_{disj}.
\end{equation}

\begin{prop} \label{Zastava is smooth}
The scheme
$\on{Maps}^\theta(\bC,\CG_{\fg,\fp})$
is smooth.
\end{prop}

\begin{proof}

Let
$\sigma:\bC\times S\to \CG_{\fg,\fp}$ be a based map,
where $S$ is an Artinian scheme, and
let $S'\supset S$ be a bigger Artinian scheme. We must
show that $\sigma$ extends to a based map
$\sigma':\bC\times S'\to \CG_{\fg,\fp}$.

Since $\CG_{\fg,\fp}$ is a union of open subschemes,
each of which is a projective limit of smooth schemes
under smooth maps, locally on $\bC$ there is no
obstruction for extending $\sigma$. Therefore, by
induction on the length of $S'$, we obtain that the
obstruction to the existence of $\sigma'$ lies in
$H^1(\bC,\sigma^*(T\CG_{\fg,\fp})(-\infty_\bC))$,
where $T\CG_{\fg,\fp}$ is the tangent sheaf on
$\CG_{\fg,\fp}$.

However, the Lie algebra $\fg$ surjects onto the tangent
space to $T\CG_{\fg,\fp}$ at every point. Therefore,
we have a surjection
$\CO_\bC\otimes \fg\twoheadrightarrow
\sigma^*(T\CG_{\fg,\fp})$. Hence, we have a surjection on
the level of $H^1$, but
$H^1(\bC,\CO_\bC\otimes \fg(-\infty_\bC))\simeq
H^1(\bC,\CO_\bC(-\infty_\bC))\otimes \fg=0$.
\end{proof}

\medskip

Let now $\fp=\fb$, and for an element
$\mu\in \Lambda^{pos}_\fg$ equal
to $\mu=\underset{i}\sum n_i\cdot \alpha_i$ let us define the length
of $\mu$ as $|\mu|=\sum n_i$.

\begin{prop} \label{Zastava is connected}
The scheme $\on{Maps}^\mu(\bC,\CG_{\fg,\fb})$ is connected.
\end{prop}

\begin{proof}

The proof proceeds by induction on the length of $\mu$.
If $|\mu|=1$, or more generally, if $\mu$ is the image
of $n\cdot \alpha_i$ for some $i\in I$, the scheme
$\on{Maps}^\mu(\bC,\CG_{\fg,\fb})$ is
isomorphic to the scheme of based maps $\bC\to \BP^1$ of degree $n$,
and hence is connected.

Let $\mu$ be an element of minimal length for which
$\on{Maps}^\mu(\bC,\CG_{\fg,\fb})$ is disconnected.
By what we said above, we can assume that $\mu$
is not a multiple of one simple coroot.

By the factorization property, Equation \eqref{factorization of quasi-maps},
and the minimality assumption, $\on{Maps}^\mu(\bC,\CG_{\fg,\fb})$
contains a connected component $\bK$, which projects under
$\varrho^\mu_\fb$ to the main diagonal
$\Delta(\overset{\circ}\bC)\subset \overset{\circ}\bC{}^\mu$.
We claim that this leads to a contradiction.

Indeed, by the definition of $\varrho^\mu_\fb$, for any map
$\sigma\in \bK$ there exists a unique point $\bc\in \bC$ with
$\sigma(\bc)\in \CG_{\fg,\fb}-\CG^e_{\fg,\fb}$. Let
$w\in \CW$ be minimal
with the property that there exists $\sigma\in \bK$ as above, such
that $\sigma(\bc)\in \CG_{\fg,\fb,w}$.

\medskip

If $w$ is a simple reflection $s_i$, then $\sigma(\bC)$
intersects only $\CG_{\fg,\fb,s_i}$, i.e. by the definition
of $\varrho^\mu_\fb$, $\theta$ is the multiple of $\alpha_i$,
which is impossible.

Hence, we can assume that $\ell(w)>1$, and let $w'$ be such that
$w=w'\cdot s_i$, $\ell(w')=\ell(w)-1$. We claim that we will be
able to find $\sigma'\in \bK$ such that $\sigma'(\bC)\cap
\CG_{\fg,\fb,w'}\neq \emptyset$.

Indeed, the group $SL_2$ corresponding to $i\in I$ acts on
$\CG_{\fg,\fb}$, and the corresponding $N^-_i\subset SL_2$
preserves the point $1_{\CG_{\fg,\fb}}$. Hence, $N^-_i$
acts on $\on{Maps}^\theta(\bC,\CG_{\fg,\fb})$, and, being
connected, it preserves the connected component $\bK$.
However, for any non-trivial element $u\in N^-_i$,
$u(\CG_{\fg,\fb,w})\subset \CG_{\fg,\fb,w'}$. Hence, if
we define $\sigma'$ as $u(\sigma)$, for the same point
$\bc\in \bC$, $\sigma'(\bc)\in \CG_{\fg,\fb,w'}$.

\end{proof}

\medskip

\noindent{\it Remark.}
For finite dimensional $\fg$ the above Proposition is not new,
see ~\cite{t}, ~\cite{ffkm}, ~\cite{kp}, ~\cite{p}. However all the
proofs avoiding factorization use Kleiman's theorem on generic transversality,
unavailable in the infinite dimensional setting.

\medskip

\begin{cor} \label{dimension of Borel Zastava}
The dimension of $\on{Maps}^\mu(\bC,\CG_{\fg,\fb})$
equals $2|\mu|$.
\end{cor}

\begin{proof}
Using \propref{Zastava is connected} and \propref{Zastava is smooth},
we know that $\on{Maps}^\mu(\bC,\CG_{\fg,\fb})$ is irreducible.
Recall the map
$\varrho^\mu_\fb:\on{Maps}^\mu(\bC,\CG_{\fg,\fb})\to
\overset{\circ}{\bC}{}^\mu$, and consider the open subset in
$\on{Maps}^\mu(\bC,\CG_{\fg,\fb})$ equal to the preimage of the
locus of $\overset{\circ}{\bC}{}^\mu$ corresponding to multiplicity-free
divisors.

Using \eqref{factorization of quasi-maps}, this reduces
us to the case, when $\mu$ is a simple
coroot $\alpha_i$, $i\in I$. But as we have seen before, the
scheme $\on{Maps}^\mu(\bC,\CG_{\fg,\fb})$ in this case is
isomorphic to $\BA^2$.

\end{proof}

Note that since the base $\overset{\circ}{\bC}{}^\mu$ is smooth,
from the above corollary it follows
that for any point $\bc\in \overset{\circ}\bC$, every irreducible
component of the preimage of the corresponding point
$\mu\cdot \bc\in \overset{\circ}{\bC}{}^\mu$
in $\on{Maps}^\mu(\bC,\CG_{\fg,\fb})$ is of dimension
$\geq |\mu|$.

The following Conjecture will be established in the case when
$\fg$ of finite and affine type:

\begin{conj} \label{lagrange}
The projection $\varrho^\mu_\fb:\on{Maps}^\mu(\bC,\CG_{\fg,\fb})\to
\overset{\circ}{\bC}{}^\mu$ is flat. Equivalently, the preimage of
the point $\mu\cdot \bc\in \overset{\circ}{\bC}{}^\mu$ (for any point
$\bc\in \overset{\circ}\bC$) is equidimensional of dimension
$|\mu|$.
\end{conj}

\medskip

For a parabolic $\fp$ and $\theta\in \Lambda^{pos}_{\fg,\fp}$ equal to
the projection of $\wt{\theta}=\sum n_i\cdot \alpha_i\in \Lambda^{pos}_\fg$
with $i\notin \fm$, define
$|\theta|$ as $\sum n_i$, and $|\theta|'=|\theta|-
\langle \wt{\theta},\rhoch_M\rangle$, where $\rhoch_M$ is half
the sum of positive roots of $\fm$.

\begin{cor} \label{dimension of Zastava}
The dimension of $\on{Maps}^\theta(\bC,\CG_{\fg,\fp})$
equals $2|\theta|'$.
\end{cor}

\begin{proof}

Pick an element $\much\in \Lambda^{pos}_\fg$, which projects to $\theta$
under $\Lambda_\fg\to \Lambda_{\fg,\fp}$, and which is sufficiently
dominant with respect to $\fm$, so that the map
$\Bun^\mu_{B^-(M)}(\bC,\infty_\bC)\to \Bun_M(\bC,\bc)$ is smooth.

It is sufficient to show that for any such $\mu$,
\begin{align*}
&\on{dim}\left((\on{Maps}^\theta(\bC,\CG_{\fg,\fp})\underset{\Bun_M(\bC,\bc)}\times
\Bun^\mu_{B^-(M)}(\bC,\infty_\bC)\right)= \\
&2|\theta|'+
\on{dim.rel.}\left(\Bun^\mu_{B^-(M)}(\bC,\infty_\bC),\Bun_M(\bC,\bc)\right).
\end{align*}

However, according to \lemref{relation between parabolics} and
\corref{dimension of Borel Zastava}, the
left-hand side of the above expression equals $2|\mu|$, and
$\on{dim.rel.}\left(\Bun^\mu_{B^-(M)}(\bC,\infty_\bC),\Bun_M(\bC,\bc)\right)$
is readily seen to equal $\langle \mu, 2\rhoch_M\rangle$. Together,
this yields the desired result.

\end{proof}

\ssec{Enhanced quasi-maps} \label{enhanced}

For an element $\lambdach\in \Lambdach^+_{\fg}$, let $\CU_\lambdach$
denote the corresponding integrable module over the Levi $M$.
Note that each such $\CU_\lambdach$ can be
realized as $(\CV_\lambdach)_{\fp^-}$. Therefore, every map of $\fg$-modules
$\CV_\nuch\to \CV_\much\otimes \CV_\lambdach$, gives rise to a map
$\CU_\nuch\to \CU_\much\otimes \CU_\lambdach$.

Consider the subscheme $$\wt{C}(\CG_{\fg,\fp})\subset
\underset{\Lambdach^+_{\fg}}\prod\, \on{Hom}(\CV_\lambdach,\CU_\lambdach)$$
given by the following equations:

A system of maps $\varphi_\lambdach\in \on{Hom}(\CV_\lambdach,\CU_\lambdach)$
belongs to $\wt{C}(\CG_{\fg,\fp})$ if for every
$\nuch,\lambdach,\much\in \Lambdach^+_{\fg}$ and a map
$\CV_\nuch\to \CV_\much\otimes \CV_\lambdach$, the diagram
\begin{equation} \label{Plucker}
\CD
\CV_\nuch   @>>> \CV_\much\otimes \CV_\lambdach \\
@V{\varphi_\nuch}VV   @V{\varphi_\much\otimes \varphi_\lambdach}VV   \\
\CU_\nuch   @>>> \CU_\much\otimes \CU_\lambdach
\endCD
\end{equation}
is commutative.

There is a natural map from $\wt{C}(\CG_{\fg,\fp})$ to
$C(\CG_{\fg,\fp})$, which remembers the data of $\varphi_\lambdach$ for
$\lambdach\in \Lambda^+_{\fg,\fp}$.
Let $\overset{\circ}{\wt{C}}(\CG_{\fg,\fp})$ be the open subscheme of
$\wt{C}(\CG_{\fg,\fp})$ corresponding to the condition that all the
$\varphi_\lambdach$ are surjections. We have a natural map
$$\overset{\circ}{\wt{C}}(\CG_{\fg,\fp})\to
\CG_{\fg,\fp},$$
and the former is a principal $M$-bundle over the latter.

\medskip

For a (not necessarily complete) curve $\bC$ we shall define
now a certain scheme $\wt{\CZ}{}^\theta_{\fg,\fp}(\bC)$,
called the enhanced version of the Zastava space.

By definition,
$\wt{\CZ}{}^\theta_{\fg,\fp}(\bC)$ classifies quadruples
$(D^\theta,\F_{N(P)},\F_M,\kappa)$, where

\begin{itemize}

\item
$(D^\theta,\F_{N(P)})$ are as in the definition of
$\CZ^\theta_{\fg,\fp}(\bC)$,

\item
$\F_M$ is a principal $M$-bundle on $\bC$,
such that the induced $M/[M,M]$-bundle $\F_{M/[M,M]}$
is identified with the one coming from $D^\theta$,

\item
$\kappa$ is a system of generically surjective
maps $$\kappa^\lambdach:(\CV_\lambdach)_{\F_{N(P)}}\to
(\CU_\lambdach)_{\F_M},\, \lambdach\in \Lambdach^+_\fg,$$

\end{itemize}

such that

\begin{itemize}

\item The $\kappa^\lambdach$'s
satisfy the Pl{\"u}cker relations (cf. Equation \eqref{Plucker}) and

\item
For $\lambdach\in \Lambdach^+_{\fg,\fp}$, the composition
$$\CO_\bC\to(\CV_\lambdach)_{\F_{N(P)}}\overset{\kappa^\lambdach}
\longrightarrow \CL^\lambdach_{\F_{M/[M,M]}}$$
equals the tautological embedding
$\CO_\bC\to \CO(\lambdach(D^\theta))$.

\end{itemize}

The proof that $\wt{\CZ}{}^\theta_{\fg,\fp}(\bC)$ is indeed
representable by a scheme is given below.

\medskip

One can reformulate the definition of
$\wt{\CZ}{}^\theta_{\fg,\fp}(\bC)$
as follows: it classifies maps
from $\bC$ to the stack $N(P)\backslash \wt{C}(\CG_{\fg,\fp})/M$,
which send the generic point of $\bC$ to
$N(P)\backslash \CG^e_{\fg,\fp}\simeq pt$.

From this it is easy to see that the analogs of \lemref{trivial outside
graph}, \corref{extension} and \propref{factorization of Zastava} hold.
In particular, for an $S$-point of $\wt{\CZ}{}^\theta_{\fg,\fp}(\bC)$, on the open
subset $\bC\times S-D^\theta$, the bundles $\F_{N(P)}$ and $\F_M$ admit
canonical trivializations, such that the maps $\kappa^\lambdach$
become the projections
$\CO_\bC\otimes \CV_\lambdach \to\CO_\bC\otimes \CU_\lambdach$.

\medskip

If $\bC$ is a complete curve with a marked point $\infty_\bC$,
we define the scheme of based enhanced quasi-maps
$\wt{\on{QMaps}}{}^\theta(\bC,\CG_{\fg,\fp})$
as
$$\wt{\CZ}{}^\theta_{\fg,\fp}(\bC)
\underset{\CZ^\theta_{\fg,\fp}(\bC)}\times
\on{QMaps}^\theta(\bC,\CG_{\fg,\fp}).$$
I.e., $\wt{\on{QMaps}}{}^\theta(\bC,\CG_{\fg,\fp})$
classifies maps from $\bC$ to the stack $\wt{C}(\CG_{\fg,\fp})/M$,
which send a neighborhood of $\infty_\bC$ to
$\overset{\circ}{\wt{C}}(\CG_{\fg,\fp})/M\simeq \CG_{\fg,\fp}$,
and such that $\infty_\bC$ gets sent to $1_{\CG_{\fg,\fp}}\in
\CG_{\fg,\fp}$.

Just as in \propref{Zastava and quasi-maps}, when $\bC$ is of genus
$0$, we have an isomorphism between
$\wt{\CZ}{}^\theta_{\fg,\fp}(\overset{\circ}\bC)$
and $\wt{\on{QMaps}}{}^\theta(\bC,\CG_{\fg,\fp})$, and in general
the latter is a closed subscheme of the former.

\medskip

When $\fg$ is finite-dimensional, we introduce the corresponding
(relative over $\Bun_G(\bC,\infty_\bC)$) version of
$\wt{\on{QMaps}}{}^\theta(\bC,\CG_{\fg,\fp})$, denoted
$\wt{\Bun}_{P^-}(\bC,\infty_\bC)$.

\ssec{}  \label{affine Grassmannian}

To formulate the next assertion we
need to recall some notation related to affine
Grassmannians.

If $M$ is a reductive group, $\bC$ a curve,
we will denote by $\Gr_{M,\bC}$ the corresponding affine Grassmannian.
By definition, this is an ind-scheme classifying triples
$(\bc,\F_M,\beta)$, where $\bc$ is a point of $\bC$, $\F_M$ is
an $M$-bundle, and $\beta$ is a trivialization of $\F_M$ on
$\bC-\bc$. When the point $\bc$ is fixed, we will
denote the corresponding subscheme of $\Gr_{M,\bC}$
by $\Gr_{M,\bc}$, and sometimes simply by $\Gr_M$. 


More generally, for $a\in \BN$, we have the Beilinson-Drinfeld
version of the affine Grassmannian, denoted $\Gr_{M,\bC}^{BD,a}$,
which is now an ind-scheme over $\bC^{(a)}$. For a fixed divisor
$D\in \bC^{(a)}$, we will denote by $\Gr_{M,\bC,D}^{BD,a}$ the
fiber of $\Gr_{M,\bC}^{BD,a}$ over $D$. By definition, the latter
scheme classifies pairs $(\F_M,\beta)$, where $\F_M$ is as before
an $M$-bundle on $\bC$, and $\beta$ is its trivialization off the
support of $D$.

\medskip

Assume now that $M$ is realized as the reductive group
corresponding to a Levi subalgebra $\fm\subset\fg$, for a
Kac-Moody algebra $\fg$, and let $\theta$ be an element
of $\Lambda^+_{\fg,\fp}$.

We define the (finite dimensional) scheme
$\on{Mod}^{\theta,+}_{M,\bC}$ to classify triples
$(D^\theta,\F_M,\beta)$, where $D^\theta$ is a point of
$\bC^\theta$, $\F_M$ is a principal $M$-bundle on $\bC$,
and $\beta$ is a trivialization of $\F_M$ off the support
of $D^\theta$, such that the following conditions hold:

\smallskip

\noindent (1)
The trivialization given by $\beta$ of the induced
$M/[M,M]$-bundle $\F_{M/[M,M]}$ is such that
$\CL^\lambdach_{\F_{M/[M,M]}}\simeq
\CO(\lambdach(D^\theta))$, $\lambdach\in \Lambda^+_{\fg,\fp}$,
where $\CL^\lambdach_{\F_{M/[M,M]}}$ is the line bundle associated
with $\F_{M/[M,M]}$ and the character $\lambdach:M/[M,M]\to\BG_m$.

\smallskip

\noindent (2) For an integrable $\fg$-module $\CV$,
and the corresponding $\fm$-module $\CU:=\CV^{n(\fp)}$,
the (a priori meromorphic) map
$\CO_\bC\otimes\CU\to \CU_{\F_M}$ induced by
$\beta$, is regular.

\medskip

Since $M$ admits a faithful representation of the form $\CV^{n(\fp)}$,
where $\CV$ is an integrable $\fg$-module, we obtain that 
$\on{Mod}^{\theta,+}_{M,\bC}$ is indeed finite-dimensional. 

Let us describe more explicitly the fibers of $\on{Mod}^{\theta,+}_{M,\bC}$ over
$\bC^\theta$. For simplicity, let us take a point in $\bC^\theta$ equal
to $\theta\cdot \bc$, where $\bc$ is some point of $\bC$. Recall that the entire 
affine Grassmannian $\Gr_{M,\bc}$ is the union of Schubert cells, denoted 
$\Gr^\mu_{M,\bc}$, where $\mu$ runs over the set of dominant coweights of $M$.
It is easy to see that the fiber of $\on{Mod}^{\theta,+}_{M,\bC}$ over 
$\theta\cdot \bc$, being a closed subscheme in $\Gr_{M,\bc}$, contains
(equivalently, intersects) $\Gr^\mu_{M,\bc}$ if and only if the followinng two conditions
hold:

\noindent (1) The projection of $\mu$ under $\Lambda_\fm=\Lambda_\fg\to 
\Lambda_{\fg,\fp}$ equals $\theta$. 

\smallskip

\noindent (2) $w_0^M(\mu)\in \Lambda_{\fg}^{pos}$, where $w_0^M$ is the longest
element of the Weyl group of $M$.

In particular, for a given $\theta$, there is only a finite many of such $\mu$,
a fact that follows alternatively from the above finite-dimensionality 
statement.

\medskip

Note that there is a canonical map $\fr_\fp:
\wt{\CZ}{}^\theta_{\fg,\fp}(\bC)\to
\CZ^\theta_{\fg,\fp}(\bC)$, which remembers
the data of $\F_{M/[M,M]}$, and $\kappa^\lambdach$
for $\lambdach\in \Lambdach^+_{\fg,\fp}$. It
is easy to see that $\fr_\fp$ is an isomorphism
over the open subset $\overset{\circ}\CZ{}^\theta_{\fg,\fp}(\bC)$.

\begin{lem}  \label{enhanced quasi-maps in the product}
We have a closed embedding (of functors)
$\wt{\CZ}{}^\theta_{\fg,\fp}(\bC)\to
\CZ^\theta_{\fg,\fp}(\bC)\underset{\bC^\theta}\times
\on{Mod}^{\theta,+}_{M,\bC}$.
\end{lem}

This lemma implies, in particular, that both
$\wt{\CZ}{}^\theta_{\fg,\fp}(\bC)$ and
$\wt{\on{QMaps}}{}^\theta(\bC,\CG_{\fg,\fp})$ are representable,
being closed subfunctors of representable functors.

\begin{proof}

For an $S$-point of $\wt{\CZ}{}^\theta_{\fg,\fp}(\bC)$,
we already know that the corresponding $M$-bundle
admits a trivialization on $\bC\times S-D^\theta$.

Moreover, if $\CU^\lambdach$ is an $\fm$-module with
$\lambdach\in \Lambdach^+_\fg$, the corresponding map
$\beta^\lambdach: \CO_\bC\otimes \CU_\lambdach\to (\CU_\lambdach)_{\F_M}$
equals the composition
\begin{equation}  \label{s volnoj}
\CO_\bC\otimes \CU_\lambdach\to \CO_\bC\otimes \CV_\lambdach\to
(\CU_\lambdach)_{\F_M},
\end{equation}
where the first arrow comes from the embedding
$\CU_\lambdach\simeq (\CV_\lambdach)^{\fn(\fp)}\hookrightarrow
\CV_\lambdach$.

This proves that $\F_M$ with its trivialization indeed defines a point
of $\on{Mod}^{\theta,+}_{M,\bC}$.

\medskip

Conversely, given a point of $\CZ^\theta_{\fg,\fp}(\bC)$, and
an $M$-bundle $\F_M$ trivialized on $\bC\times S-D^\theta$,
from \lemref{trivial outside graph} we obtain that there
is a meromorphic map
$\kappa^\lambdach:\CO_\bC\otimes \CV_\lambdach\to (\CU_\lambdach)_{\F_M}$.
Our data defines a point of $\wt{\CZ}{}^\theta_{\fg,\fp}(\bC)$
if and only if $\kappa^\lambdach$ is regular, which is a
closed condition.

\end{proof}

We will denote by $\varrho^\theta_M$ the projection
$\wt{\CZ}{}^\theta_{\fg,\fp}(\bC)\to \on{Mod}^{\theta,+}_{M,\bC}$.

\ssec{} \label{twisted and enhanced}

Let us once again assume that $\bC$ is complete. We will introduce
yet two more versions of $\CZ^\theta_{\fg,\fp}(\bC)$ (denoted
$\sZ^\theta_{\fg,\fp}(\bC)$, and $\wt{\sZ}{}^\theta_{\fg,\fp}(\bC)$,
respectively, and called the twisted Zastava spaces),
which will be fibered over the stack
$\Bun_M(\bC)$ classifying $M$-bundles on $\bC$.

By definition, $\wt{\sZ}^\theta_{\fg,\fp}(\bC)$ classifies the data of
$(D^\theta,\F_P,\F_M,\kappa)$, where $(D^\theta,\F_M)$ are
as before, $\F_P$ is a principal $P$-bundle on $\bC$, and
$\kappa$ is a collection of maps
$$\kappa^\lambdach:(\CV_\lambdach)_{\F_P} \to
(\CU_\lambdach)_{\F_M},\, \lambdach\in \Lambdach^+_\fg,$$
which are generically surjective and satisfy the
Pl{\"u}cker relations, and for $\lambdach\in \Lambdach^+_{\fg,\fp}$
the compositions
$\CL^\lambdach_{\F_P}\to (\CV_\lambdach)_{\F_P}\to \CL^\lambdach_{\F_M}$
induce isomorphisms $\CL^\lambdach_{\F_M}\simeq
\CL^\lambdach_{\F_P}(\lambdach(D^\theta))$.

\medskip

The stack $\sZ^\theta_{\fg,\fp}(\bC)$ classifies triples
$(D^\theta,\F_P,\kappa)$, where $(D^\theta,\F_P)$
are as above, and $\kappa$ are collections of maps
$\kappa^\lambdach:(\CV_\lambdach)_{N(P)}\to \CO(\lambdach(D^\theta))$,
defined for $\lambdach\in \Lambdach^+_{\fg,\fp}$. The forgetful
map $\wt{\sZ}^\theta_{\fg,\fp}(\bC)\to \sZ^\theta_{\fg,\fp}(\bC)$
will be denoted by the same character $\fr_\fp$.

We have an open substack $\overset{\circ}\sZ{}^\theta_{\fg,\fp}(\bC)\subset
\sZ^\theta_{\fg,\fp}(\bC)$ that corresponds to the condition that the maps
$\kappa^\lambdach$ are surjective. Over it, $\fr_\fp$ is an isomorphism.

\medskip

We have a projection, that we will call $\q_{\fp}$,
from both $\sZ^\theta_{\fg,\fp}(\bC)$ and $\wt{\sZ}^\theta_{\fg,\fp}(\bC)$
to the stack $\Bun_M(\bC)$ that ``remembers'' the data of
the $M$-bundle $\F'_M:=N(P)\backslash \F_P$. In the case of
$\wt{\sZ}^\theta_{\fg,\fp}(\bC)$ we have also the other projection
to $\Bun_M(\bC)$, denoted $\q_{\fp^-}$, that ``remembers'' the data of
$\F_M$. Of course, the preimage of the trivial bundle
$\F^0_M\in \Bun_M(\bC)$ under $\q_{\fp}$ identifies
with the schemes $\CZ^\theta_{\fg,\fp}(\bC)$ and
$\wt{\CZ}{}^\theta_{\fg,\fp}(\bC)$, respectively.

\medskip

It is easy to see that \lemref{trivial outside graph}, \corref{extension}
and \propref{factorization of Zastava} generalize to the context of
twisted Zastava spaces.
In particular, we have the following assertion:

\begin{lem}  \label{locality}
Let $\phi_1,\phi_2:S\to \Bun_M(\bC)$ be two arrows, such that
the corresponding $M$-bundles on $\bC\times S$ are identified over an open
subset $U\subset \bC\times S$. In addition let $D^\theta$ be
the graph of a map $S\to \bC^\theta$, such that $D^\theta\subset U$.
Then the two Cartesian products
$S\underset{\Bun_M(\bC)\times \bC^\theta}\times \wt{\sZ}{}^\theta_{\fg,\fp}(\bC)$,
taken with respect to either $\phi_1$ or $\phi_2$ are naturally
isomorphic; and a similar assertion holds for $\sZ^\theta_{\fg,\fp}(\bC)$.
\end{lem}

This lemma implies that $\wt{\sZ}{}^\theta_{\fg,\fp}(\bC)$
and $\sZ^\theta_{\fg,\fp}(\bC)$ are algebraic stacks, such that
the morphism $\q_{\fp}$ is representable.

\medskip

As in the case of $\CZ^\theta_{\fg,\fp}(\bC)$, we have a
stratification of $\sZ^\theta_{\fg,\fp}(\bC)$ by locally-closed substacks
of the form $\overset{\circ}\sZ{}^{\theta-\theta'}_{\fg,\fp}(\bC)\times \bC^{\theta'}$
for $\theta',\theta-\theta'\in \Lambda^{pos}_{\fg,\fp}$.
To describe the preimages of these strata in $\wt{\sZ}^\theta_{\fg,\fp}(\bC)$
we need to introduce some notation:

\medskip

Let $\H_{M,\bC}$ be the Hecke stack, i.e. a relative over
$\Bun_M(\bC)$ version of $\Gr_{M,\bC}$, which classifies
quadruples $(\bc,\F_M,\F'_M,\beta)$, where $\bc$ is a point
of $\bC$, $\F_M,\F'_M$ are principal $M$-bundles on $\bC$,
and $\beta$ is an identification $\F_M\simeq \F'_M|_{\bC-\bc}$.
Let $\H^{BD,a}_{M,\bC}$, $\H^{+,\theta}_{M,\bC}$ be the
corresponding relative versions of $\Gr^{BD,a}_{M,\bC}$
and $\on{Mod}^{+,\theta}_{M,\bC}$, respectively.
We will denote by $\hl$, $\hr$ the two projections from any
of the stacks $\H_{M,\bC}$, $\H^{BD,a}_{M,\bC}$ or
$\on{Mod}^{+,\theta}_{M,\theta}$ to $\Bun_M(\bC)$ that
``remember'' the data of $\F_M$ and $\F'_M$, respectively.

Note that we have a canonical (convolution) map
$\H_{M,\bC}\underset{\Bun_M(\bC)\times \bC}\times \H_{M,\bC}\to
\H_{M,\bC}$,
which sends $(\bc,\F_M,\F'_M,\beta)\times (\bc,\F'_M,\F''_M,\beta')$
to $(\bc,\F_M,\F''_M,\beta'')$, where $\beta''=\beta'\circ\beta$.
Similarly, we have a map
$\H^{+,\theta}_{M,\bC}\underset{\Bun_M(\bC)}\times
\H^{+,\theta'}_{M,\bC}\to \H^{+,\theta+\theta'}_{M,\bC}$, that
covers the addition map $\bC^\theta\times \bC^{\theta'}\to
\bC^{\theta+\theta'}$.

As in the case of $\CZ^\theta_{\fg,\fp}(\bC)$, we have a map
$\varrho_M^\theta:\wt{\sZ}^\theta_{\fg,\fp}(\bC)\to
\H^{+,\theta}_{M,\bC}$, such that the map
$$\wt{\sZ}^\theta_{\fg,\fp}(\bC)\hookrightarrow
\sZ^\theta_{\fg,\fp}\underset{\Bun_M(\bC)}\times\H^{+,\theta}_{M,\bC}$$
is a closed embedding, and such that the maps $\q_\fp$, $\q_{\fp^-}$ 
are the compositions of $\varrho_M^\theta$ with $\hl$ and $\hr$, 
respectively.

\medskip

\begin{prop} \label{strata s volnoj}
The preimage of a stratum $\overset{\circ}\sZ{}^{\theta-\theta'}_{\fg,\fp}(\bC)
\times \bC^{\theta'}$ in $\wt{\sZ}^\theta_{\fg,\fp}(\bC)$ is isomorphic
to
$\overset{\circ}\sZ{}^{\theta-\theta'}_{\fg,\fp}(\bC)\underset{\Bun_M(\bC)}\times
\H^{+,\theta'}_{M,\bC}$,
where the projections from $\overset{\circ}\sZ{}^{\theta-\theta'}_{\fg,\fp}(\bC)$
and $\H^{+,\theta'}_{M,\bC}$ to $\Bun_M(\bC)$
are $\q_{\fp^-}$ and $\hl$, respectively.
\end{prop}

By fixing the $P$-bundle $\F_P$ (resp., the
$M$-bundle $N(P)\backslash \F_P$) to be trivial,
we obtain the description of the corresponding strata
in $\wt{\on{QMaps}}{}^\theta(\bC,\CG_{\fg,\fp})$
and $\wt{\CZ}^\theta_{\fg,\fp}(\bC)$, respectively.

\begin{proof}

Given a point of $\sZ^\theta_{\fg,\fp}(\bC)$ factoring through
$\overset{\circ}\sZ{}^{\theta-\theta'}_{\fg,\fp}(\bC)
\times \bC^{\theta'}$, we have a triple $(\F_P,\F^1_M,\kappa)$:
$$(\CV_\lambdach)_{\F_P}\overset{\kappa^\theta}
\longrightarrow (\CU_\lambdach)_{\F^1_M}$$
where the maps $\kappa^\theta$ are surjections.

Then it is clear that the scheme of possible $\F_M$'s
with $\beta:\F_M\simeq \F^1_M|_{\bC-D^{\theta'}}$ such that the
(a priori meromorphic) maps
$$(\CV_\lambdach)_{\F_P}\overset{\kappa^\theta}
\longrightarrow (\CU_\lambdach)_{\F_M}$$
continue to be regular, identifies with the fiber of
$\H^{+,\theta'}_{M,\bC}$ over $\F^1{}_M\in \Bun_M(\bC)$.

\end{proof}

We will denote by $\wt{\fs}{}^\theta_\fp$ the embedding of the last
stratum $\H^{\theta,+}_{M,\bC}\to \wt{\sZ}^\theta_{\fg,\fp}(\bC)$.
Note that in terms of the above proposition, the composition
$$\overset{\circ}\sZ{}^{\theta-\theta'}_{\fg,\fp}(\bC)\underset{\Bun_M(\bC)}\times
\H^{+,\theta'}_{M,\bC}\to\wt{\sZ}^\theta_{\fg,\fp}(\bC)\overset{\varrho^\theta_M}
\longrightarrow \H^{+,\theta}_{M,\bC}$$ equals
$$\overset{\circ}\sZ{}^{\theta-\theta'}_{\fg,\fp}(\bC)\underset{\Bun_M(\bC)}\times
\H^{+,\theta'}_{M,\bC}\overset{\varrho^{\theta-\theta'}_M\times\on{id}}\longrightarrow
\H^{+,\theta-\theta'}_{M,\bC}\underset{\Bun_M(\bC)}\times
\H^{+,\theta'}_{M,\bC}\to \H^{+,\theta}_{M,\bC}.$$

\ssec{The meromorphic case} \label{meromorphic Zastava}

Let $\X$ be a scheme, and $D_\X\subset \bC\times \X$ a relative Cartier divisor
disjoint from $\infty_\bC\times\X$. 
We will now develop a notion of a meromorphic quasi-map parallel
to \secref{meromorphic}. These will appear in the following two
contexts.

First, let $\fg$ be finite-dimensional, and let $\F_G$ be a principal
$G$-bundle defined on $\bC\times \X-D_\X$, trivialized along 
$\infty_\bC\times\X$. Let 
$_{\infty\cdot D_\X}\on{QMaps}(\bC,\F_G\overset{G}\times \CG_{\fg,\fp})$ be the
functor on the category of schemes over $\X$, that assigns to a test
scheme $S$ the data of $(\F_{M/[M,M]},\kappa^{\lambdach},\lambdach\in
\Lambda^+_{\fg,\fp})$, where $\F_{M/[M,M]}$ is a principal $M/[M,M]$-bundle
defined on $\bC\times S$, and the $\kappa^{\lambdach}$'s are maps
\begin{equation} \label{merom maps}
\CL^\lambdach_{\F_{M/[M,M]}}|_{(\bC\times\X-D_\X)
\underset{\X}\times S}\to (\CV^*_\lambdach)_{\F_G},
\end{equation}
satisfying the Pl{\"u}cker equations, with a prescribed value at
$\infty_\bC\times S$ corresponding to the trivialization of $\F_G$.
The functor $_{\infty\cdot D_\X}\on{QMaps}^\theta(\bC,\CG_{\fg,\fp})$
is representable by a strict ind-scheme of ind-finite type, by 
\propref{representability of meromorphic}.

\medskip

Similarly, we define the ind-scheme 
$_{\infty\cdot D_\X}\wt{\on{QMaps}}(\bC,\F_G\overset{G}\times\CG_{\fg,\fp})$,
where instead of a $M/[M,M]$-bundle $\F_{M/[M,M]}$, we have an
$M$-bundle $\F_M$, and the maps 
$\kappa^{\lambdach}: (\CV_\lambdach)_{\F_G}\to (\CU_\lambdach)_{\F_M}$
are defined for all $\lambdach\in \Lambda_\fg^+$ and satisfy the
Pl{\"u}cker relations in the sense of \eqref{Plucker}.

\medskip

Let now $\fg$ be arbitrary. We define $_{\infty\cdot D_\X}\on{QMaps}(\bC,\CG_{\fg,\fp})$ 
to be the ind-scheme, which is the union over $\nu\in \Lambda^{pos}_{\fg,\fp}$ 
of schemes, classifying the data of $(\F_{M/[M,M]},\kappa^{\lambdach},\lambdach\in
\Lambda^+_{\fg,\fp})$, where $\F_{M/[M,M]}$ is a principal $M/[M,M]$-bundle
on $\bC\times S$, and the $\kappa^{\lambdach}$'s are maps
$$\CL^\lambdach_{\F_{M/[M,M]}}(-\langle \nu,\lambdach\rangle \cdot D_\X|_{\bC\times S})
\to \CV^*_\lambdach,$$
satisfying the Pl\"ucker relations, with the prescribed value at $\infty_\bC\times S$.
The ind-scheme $_{\infty\cdot D_\X}\wt{\on{QMaps}}(\bC,\CG_{\fg,\fp})$ is defined similarly,
but where the union goes over $\nu\in \Lambda^{pos}_\fg$.

\section{Bundles on $\BP^1\times \BP^1$}  \label{bundles}

In this section $\bC$ will be a projective curve of genus $0$
with a marked infinity $\infty_\bC$, and $\infty_\bX\in\bX$
will be another such curve. We will be interested in the
surface $\bS':=\bC\times \bX$, and we will call
$\bD'_\infty:=\bC\times \infty_\bX\cup \infty_\bC\times\bX\subset \bS'$
the divisor at infinity.

Throughout this section $\fg$ will be a finite-dimensional
simple Lie algebra, and $G$ the corresponding
simply-connected group.

\ssec{}

Consider the stack $\Bun_G(\bS',\bD'_\infty)$ that classifies
$G$-bundles on $\bS'$ with a trivialization on $\bD'_\infty$.
We will see shortly that $\Bun_G(\bS',\bD'_\infty)$
is in fact a scheme. We have
$$\Bun_G(\bS',\bD'_\infty)=\underset{a\in \BN}\cup\Bun^a_G(\bS',\bD'_\infty),$$
where $\Bun^a_G(\bS',\bD'_\infty)$ corresponds to $G$-bundles with
second Chern class equal to $a$. (It is easy to see that
$\Bun_G(\bS',\bD'_\infty)$ contains no points with a negative Chern class.)

Recall now that $\Bun_G(\bX,\infty_\bX)$ denotes the stack
classifying $G$-bundles on $\bX$ with a trivialization at
$\infty_\bX$. The stack $\Bun_G(\bX,\infty_\bX)$ contains
an open subset corresponding to the {\it trivial bundle},
which is isomorphic to $pt$. Since $\Bun_G(\bX,\infty_\bX)$ is
smooth, the complement $\Bun_G(\bX,\infty_\bX)-pt$, being
of codimension $1$, is a Cartier
divisor. I.e., we are in the situation of \secref{factorization
principle}. We will denote the corresponding line bundle on
$\Bun_G(\bX,\infty_\bX)$ by $\CP_{\Bun_G(\bX,\infty_\bX)}$.

\medskip

We will use a short-hand notation of
$\on{Maps}(\bC,\Bun_G(\bX,\infty_\bX))$ for the functor
of based maps from $\bC$ to $\Bun_G(\bX,\infty_\bX)$ that send
$\infty_\bC$ to $pt\subset \Bun_G(\bX,\infty_\bX)$.
By definition,
$$\on{Maps}(\bC,\Bun_G(\bX,\infty_\bX))=\underset{a\in \BN}\cup
\on{Maps}^a(\bC,\Bun_G(\bX,\infty_\bX)),$$ where
each $\on{Maps}^a(\bC,\Bun_G(\bX,\infty_\bX)$ corresponds
to maps $\sigma$ such that $\sigma^*(\CP_{\Bun_G(\bX,\infty_\bX)})$
is of degree $a$. Since $\CP_{\Bun_G(\bX,\infty_\bX)}$ comes equipped
with a section, we obtain a map
$\varpi^a_h:\on{Maps}^a(\bC,\Bun_G(\bX,\infty_\bX))\to
\overset{\circ}{\bC}{}^{(a)}$.

\medskip

We have an obvious isomorphism of functors:
\begin{equation} \label{bundles as maps}
\Bun_G(\bS',\bD'_\infty)\simeq\on{Maps}(\bC,\Bun_G(\bX,\infty_\bX)).
\end{equation}

The following assertion is left to the reader:
\begin{lem}
Under the above isomorphism,
$\Bun^a_G(\bS',\bD'_\infty)$ maps to
$\on{Maps}^a(\bC,\Bun_G(\bX,\infty_\bX))$.
\end{lem}

In particular, we obtain a map
$\varpi^a_h:\Bun^a_G(\bS',\bD'_\infty)\to \overset{\circ}\bC{}^{(a)}$.
By interchanging the roles of $\bX$ and $\bC$, we obtain a map
$\varpi^a_v:\Bun^a_G(\bS',\bD'_\infty)\to \overset{\circ}\bX{}^{(a)}$.

\ssec{}

Let $\CG_{G,\bX}$ denote the ``thick'' Grassmannian corresponding
to $G$ and the curve $\bX$. I.e., $\CG_{G,\bX}$ is the scheme
classifying pairs, $(\F_G,\beta)$, where $\F_G$ is a $G$-bundle
on $\bX$, and $\beta$ is a trivialization of $\F_G$ over the
formal neighborhood of ${\mathcal D}_{\infty_\bX}$ of $\infty_\bX$.

As we shall recall later, $\CG_{G,\bX}$ is one of the partial flag
schemes associated to the affine Kac-Moody algebra corresponding to
$\fg$. The unit point $1_{\CG_{G,\bX}}$ is the pair
$(\F^0_G,\beta_{taut})$, where $\F^0_G$ is the trivial bundle
and $\beta_{taut}$ is its tautological trivialization.
In particular, the scheme of based maps
$\on{Maps}^a(\bC,\CG_{G,\bX})$ makes sense. (And we know from
\propref{locally of finite type} that it is locally of 
finite type and smooth.)

Moreover, $\CG_{G,\bX}$ carries
a very ample line bundle that will be denoted by $\CP_{\CG_{G,\bX}}$,
which is in fact the pull-back of the line bundle
$\CP_{\Bun_G(\bX,\infty_\bX)}$ under the natural projection
$\CP_{\CG_{G,\bX}}\to \Bun_G(\bX,\infty_\bX)$. This
enables us to define the scheme of based quasi-maps
$\on{QMaps}^a(\bC,\CG_{G,\bX})$, studied in the previous
section.

\medskip

Recall also that over the symmetric power $\overset{\circ}\bX{}^{(a)}$
we have the ind-scheme $\Gr_{G,\bX}^{BD,a}$, classifying pairs
$(D,\F_G,\beta)$, where $D\in \overset{\circ}\bX{}^{(a)}$,
$\F_G$ is a $G$-bundle on $\bX$, and $\beta$ is its trivialization
off the support of $D$. We have a section
$\overset{\circ}\bX{}^{(a)}\to \Gr_{G,\bX}^{BD,a}$, which also
corresponds to the trivial bundle with a tautological trivialization.
We will denote by $\on{Maps}(\bC,\Gr_{G,\bX}^{BD,a})$ the (ind)-scheme
of based relative maps $\bC\to \Gr_{G,\bX}^{BD,a}$. By definition,
this is an ind-scheme over $\overset{\circ}\bX{}^{(a)}$.

There exists a natural map $\Gr_{G,\bX}^{BD,a}\to \CG_{G,\bX}\times
\overset{\circ}\bX{}^{(a)}$, which corresponds to the restriction of
the trivialization on $\bX-D$ to ${\mathcal D}_{\infty_\bX}$. It is
easy to see that this map is in fact a closed embedding.
Moreover, the restriction $\CP_{\Gr_{G,\bX}^{BD,a}}$
of the line bundle $\CP_{\CG_{G,\bX}}$ to $\Gr_{G,\bX}^{BD,a}$
is relatively very ample. Therefore, we can introduce the ind-scheme
if based quasi-maps $\on{QMaps}(\bC,\Gr_{G,\bX}^{BD,a})$.

\begin{prop} \label{bundles and maps}
The natural morphisms
$$\on{Maps}^a(\bC,\Gr_{G,\bX}^{BD,a})\to \on{Maps}^a(\bC,\CG_{G,\bX})
\to \on{Maps}^a(\bC,\Bun_G(\bX,\infty_\bX))\simeq \Bun^a_G(\bS',\bD'_\infty)$$
are all isomorphisms.
The resulting map $\Bun^a_G(\bS',\bD'_\infty)\to \overset{\circ}\bX{}^{(a)}$
coincides with the map $\varpi_v^a$.
\end{prop}

\begin{proof}

Let us show first that any map $S\to \Bun^a_G(\bS',\bD'_\infty)$
(for any test scheme $S$) lifts to a map
$S\to \on{Maps}(\bC,\Gr_{G,\bX}^{BD,a})$. Indeed, given
a $G$-bundle $\F_G$ on $\bS'\times S$ and using the map $\varpi^a_v$
we obtain a divisor $D^v\subset \overset{\circ}\bX\times S$,
such that $\F_G$ is trivialized on $\bC\times (\bX\times S-D^v)$.
But this by definition means that we are dealing with a based map
$\bC\times S\to \Gr_{G,\bX}^{BD,a}$, which covers the map
$S\to \overset{\circ}\bX{}^{(a)}$ corresponding to $D^v$.

To prove the proposition it remains to show that if $\F_G$
is a $G$-bundle in $S\times \bS'$ equipped with two trivializations
on $\bC\times {\mathcal D}_{\infty_\bX}\times S$, which agree
on $\infty_\bC \times {\mathcal D}_{\infty_\bX}\times S$, then
these two trivializations coincide.

Indeed, the difference between the trivializations is a map
from $\bC\times S$ to the group of automorphisms of the trivial
$G$-bundle on ${\mathcal D}_{\infty_\bX}$. And since $\bC$
is complete and the group in question is pro-affine, any such map is
constant along the $\bC$ factor.

\end{proof}

\ssec{}

Using \propref{Zastava is smooth} we obtain that the scheme
$\Bun^a_G(\bS',\bD'_\infty)$ is smooth. Moreover, we claim
that $\on{dim}(\Bun^a_G(\bS',\bD'_\infty))=2\cdot \check{h}\cdot a$,
where $\check{h}$ is the dual Coxeter number. One way to see
this is via \corref{dimension of Zastava}, and another
way is as follows:

The tangent space to $\F_G\subset \Bun^a_G(\bS',\bD'_\infty)$
at a point corresponding to a $G$-bundle $\F_G$ equals
$H^1(\bS',\fg_{\F_G}(-\bD'))$. Note that the vector bundle
corresponding $\fg_{\F_G}$ has zero first Chern class,
and the second Chern class equals $2\cdot \check{h}\cdot a$.

Since points of $\Bun^a_G(\bS',\bD'_\infty)$ have no
automorphisms, we obtain that $H^0(\bS',\fg_{\F_G}(-\bD'))=0$,
and by Serre duality we obtain that
$H^2(\bS',\fg_{\F_G}(-\bD'))\simeq H^0(\bS',\fg^*_{\F_G}(-\bD'))^*=0$.
Hence, the dimension of $H^1(\bS',\fg_{\F_G}(-\bD'))$ can be calculated
by the Riemann-Roch formula, which yields $2\cdot \check{h}\cdot a$.
Note that this calculation in fact reproves that
$\Bun^a_G(\bS',\bD'_\infty)$ is smooth, since we have shown that all
the tangent spaces have the same dimension.

\medskip

Let $0_\bX\subset \bX$ be another point, and let us denote
by $\bD'_0\subset \bS'$ the divisor $\bC\times 0_\bX$.
By taking restriction of $G$-bundles we obtain a map
from $\Bun^a_G(\bS',\bD'_\infty)$ to $\Bun_G(\bC,\infty_\bC)$.

\begin{lem}  \label{smoothness of reduction map}
The above map
$\Bun^a_G(\bS',\bD'_\infty)\to\Bun_G(\bC,\infty_\bC)$
is smooth.
\end{lem}

\begin{proof}
Since both $\Bun^a_G(\bS',\bD'_\infty)$ and
$\Bun_G(\bC,\infty_\bC)$ are smooth, it suffices
to check the surjectivity of the corresponding
map on the level of tangent spaces:
$$H^1(\bS',\fg_{\F_G}(-\bD'))\to
H^1(\bC,\fg_{\F_G}|_{\bD'_0}(-\infty_\bC)).$$

By the long exact sequence, the cokernel is given by
$H^2(\bS',\fg_{\F_G}(-\bD'-\bD'_0))$, and we claim that
this cohomology group vanishes. Indeed, by
Serre duality, it suffices to show that
$H^0(\bS',\fg_{\F_G}(-\infty_\bC\times \bX))=0$, and we
have an exact sequence
$$0\to H^0(\bS',\fg_{\F_G}(-\bD'))\to
H^0(\bS',\fg_{\F_G}(-\infty_\bC\times \bX))\to
H^0(\bC,\fg_{\F_G}|_{\bC\times \infty_\bX}(-\infty_\bC))...$$
We know already that the first term vanishes
(since $\Bun^a_G(\bS',\bD'_\infty)$ is a scheme) and
the last term vanishes too (since
$\fg_{\F_G}|_{\bC\times \infty_\bX}$ is trivial).
\end{proof}

\medskip

As another corollary of \propref{bundles and maps} we obtain the following
factorization property of $\Bun^a_G(\bS',\bD'_\infty)$ with
respect to the projection $\varpi^a_h$ (and by symmetry,
with respect to $\varpi^a_v$):
\begin{align} \label{factorization of bundles}
&\Bun^a_G(\bS',\bD'_\infty)\underset{\overset{\circ}\bC{}^{(a)}}\times
(\overset{\circ}\bC{}^{(a_1)}\times
  \overset{\circ}\bC{}^{(a_2)})_{disj}\simeq \\
&(\Bun^{a_1}_G(\bS',\bD'_\infty)\times\Bun^{a_2}_G(\bS',\bD'_\infty))
\underset{\overset{\circ}\bC{}^{(a_1)}\times \overset{\circ}\bC{}^{(a_2)}}\times
(\overset{\circ}\bC{}^{(a_1)}\times \overset{\circ}\bC{}^{(a_2)})_{disj}.
\end{align}

\ssec{Relation with affine Lie algebras}  \label{parabolic notation}

Let $0_\bX\in \bX$ be another chosen point, and let $x$ be a
coordinate on $\bX$ with $x(0_\bX)=\infty$, $x(\infty_\bX)=0$.
Let $\fg_{aff}\simeq \fg((x))\oplus K\cdot \BC\oplus d\cdot \BC$
be the corresponding untwisted affine Kac-Moody algebra, such
that the element $d$ acts via the derivation $x\cdot\partial_x$.
We will denote by $\widehat{\fg}$ the derived algebra of $\fg_{aff}$,
i.e., $\fg((x))\oplus K\cdot \BC$.

The lattice $\Lambda_{\fg_{aff}}$ is by definition the direct sum
$\Lambda_\fg\oplus \delta\cdot \BZ \oplus \BZ$, and we will denote
by  $\widehat\Lambda_\fg$ the direct sum of the first two factors,
i.e. the cocharacter lattice of the corresponding derived group.
(The other factor can be largely ignored, since, for example,
the semi-group $\Lambda^{pos}_{\fg_{aff}}$ is contained in
$\widehat\Lambda_\fg$, and we will sometimes write
$\widehat\Lambda{}^{pos}_\fg$ instead of $\Lambda^{pos}_{\fg_{aff}}$.)
We will write an element $\mu\in \widehat\Lambda_\fg$
as $(\ol{\mu},a)$ for $\ol{\mu}\in \Lambda_\fg$, $a\in \BZ$.
By definition, the element $\delta\in \widehat\Lambda_\fg$ equals
$(0,1)$.

Let $\ol{\alpha}_0$ denote the positive coroot of $\fg$ dual to
the long dominant root. The simple affine coroot $\alpha_0$ equals
$(-\ol{\alpha}_0,1)$.
Note that $(\ol{\mu},a)\in \widehat\Lambda{}^{pos}_\fg$
if and only if $a\geq 0$ and 
$\ol{\mu}+a\cdot\ol{\alpha}_0\in \Lambda^{pos}_\fg$.

\medskip

We will denote by $\fg_{aff}^+$ (resp., $\fg_{aff}^-$)
the subalgebra $\fg[[x]]\oplus K\cdot \BC\oplus d\cdot \BC$ (resp.,
$\fg[x^{-1}]\oplus K\cdot \BC\oplus d\cdot \BC$). These algebras
are the maximal parabolic and its opposite corresponding to
$I\subset I_{aff}$. The corresponding partial flag scheme
$\CG_{\fg_{aff},\fg^+_{aff}}$ identifies with
$\CG_{G,\bX}$. Note that the latter does not depend on the
choice of the point $0_\bX\in \bX$.

\medskip

If $P\subset G$ is a parabolic and $P^-$ is the corresponding
opposite parabolic, let us denote by $\CG_{G,P,\bX}$ the scheme
classifying triples $(\F_G,\beta,\gamma)$, where
$(\F_G,\beta)$ are as in the definition of $\CG_{G,\bX}$, and
$\gamma$ is the data of a reduction to $P$ of the fiber of
$\F_G$ at $0_\bX$.

\medskip

We will denote by $\fp_{aff}^+$ the subalgebra
of $\fg_{aff}^+$ consisting of elements whose value
modulo $x$ belongs to $\fp$. Similarly, we will
denote by $\fp_{aff}^-$ the subalgebra of $\fg_{aff}^-$
consisting of elements whose value modulo $x^{-1}$
belongs to $\fp^-$. The corresponding Lie algebras
$\fp_{aff}^+,\fp_{aff}^-$
are a parabolic and its opposite in $\fg_{aff}$, and
$\CG_{\fg_{aff},\fp_{aff}^+}\simeq \CG_{G,P,\bX}$.
The Levi subgroup corresponding to $\fp_{aff}^+$ is
$M\times \BG_m\times \BG_m$; we will denote by $M_{aff}$ the group
$M\times \BG_m$ corresponding to the first $\BG_m$-factor.

The lattice $\Lambda_{\fg_{aff},\fp_{aff}^+}$
is the direct sum $\Lambda_{\fg,\fp}\oplus \delta\cdot \BZ \oplus \BZ$,
and we will denote by $\widehat\Lambda_{\fg,\fp}\subset
\Lambda_{\fg_{aff},\fp_{aff}}$
the direct sum of the first two factors. For
$\theta\in \widehat\Lambda{}^{pos}_{\fg,\fp}:=
\Lambda^{pos}_{\fg_{aff},\fp_{aff}^+}$,
we have the corresponding scheme
$\on{Mod}^{\theta,+}_{M_{aff},\bC}$.

\medskip

Thus, for a projective curve $\bC$, we can consider the schemes
of based maps $\on{Maps}^\theta(\bC,\CG_{G,P,\bX})$ for
$\theta\in \widehat\Lambda{}^{pos}_{\fg,\fp}$.
For $\theta=(\ol{\theta},a)$,
consider the stack $\Bun^\theta_{G;P}(\bS',\bD'_\infty;\bD'_0)$ that
classifies the data of a $G$-bundle $\F_G\in \Bun^a_G(\bS',\bD'_\infty)$,
and equipped with a reduction to $P$ on $\bD'_0$ of weight $\ol{\theta}$,
compatible with the above trivialization on
$\bD'_\infty\cap \bD'_0=\infty_\bC\times\infty_\bX$.

From \propref{bundles and maps} we obtain an isomorphism
\begin{equation} \label{parabolic bundles as maps}
\Bun^\theta_{G;P}(\bS',\bD'_\infty;\bD'_0)\simeq
\on{Maps}^\theta(\bC,\CG_{G,P,\bX}).
\end{equation}

\medskip

We can also consider the scheme of based quasi-maps
$\on{QMaps}^\theta(\bC,\CG_{G,P,\bX})$. We will denote by
$\varrho^\theta_{\fp^+_{aff}}$ the map
$\on{QMaps}^\theta(\bC,\CG_{G,P,\bX})\to \overset{\circ}\bC{}^\theta$.

In addition, we have the scheme of enhanced based quasi-maps
$\wt{\on{QMaps}}{}^\theta(\bC,\CG_{G,P,\bX})$. We have a projection denoted
$\fr_{\fp^+_{aff}}:\wt{\on{QMaps}}{}^\theta(\bC,\CG_{G,P,\bX})\to
\on{QMaps}^\theta(\bC,\CG_{G,P,\bX})$, and a map
$\varrho^\theta_{M_{aff}}:\wt{\on{QMaps}}{}^\theta(\bC,\CG_{G,P,\bX})\to
\on{Mod}^{\theta,+}_{M_{aff},\bX}$, such that
$$\wt{\on{QMaps}}{}^\theta(\bC,\CG_{G,P,\bX})\to
\on{QMaps}^\theta(\bC,\CG_{G,P,\bX})\underset{\overset{\circ}\bC{}^\theta}
\times \on{Mod}^{\theta,+}_{M_{aff},\bC}$$
is a closed embedding.

Note that the map $\fr_{\fp^+_{aff}}$ is not an isomorphism even
when $P=G$. But it is, of course, an isomorphism when $\fp=\fb$.

\ssec{}   \label{upper estimate}

Finally, let us show that the scheme $\on{Maps}^a(\bC,\CG_{G,\bX})=\Bun^a_G(\bS',\bD'_\infty)$
is globally of finite type. Using Equation \eqref{parabolic bundles as maps},
this would imply that the schemes
$\on{Maps}^\theta(\bC,\CG_{\fg_{aff},\fp_{aff}^+})$
are also of finite type.

Using \propref{bundles and maps}, it suffices to show that any map
$\sigma:\bC\to \Gr^{BD,a}_{G,\bX}$ of degree $a$ has its image in a fixed
finite-dimensional subscheme of $\Gr^{BD,a}_{G,\bX}$. To simplify the
notation, we will fix a divisor $D\in \overset{\circ}\bX{}^{(a)}$
and consider based maps $\bC\to \Gr^{BD,a}_{G,\bX,D}$.

Write $D=\underset{k}\Sum\, n_k\cdot \bx_k$ with $\bx_k$ pairwise distinct.
Then $\Gr^{BD,a}_{G,\bX,D}$ is the product of the affine Grassmannians
$\underset{k}\Pi\, \Gr_{G,\bx_k}$ and consider the subscheme
$\underset{k}\Pi\, \ol{\Gr}{}^{n_k\cdot \ol{\alpha}_0}_{G,\bx_k}\subset
\Gr^{BD,a}_{G,\bX,D}$, where for a dominant coweight
$\ol{\lambda}$ of $G$, $\ol{\Gr}{}^{\ol{\lambda}}_{G,\bx}$ denotes the
corresponding finite-dimensional subscheme of $\ol{\Gr}_{G,\bx}$.

\begin{lem}
For $D\in \overset{\circ}\bX{}^{(a)}$ as above, any based map $\sigma:\bC\to
\Gr^{BD,a}_{G,\bX,D}$  of degree $a$ has its image in the subscheme
$\underset{k}\Pi\, \ol{\Gr}{}^{n_k\cdot \ol{\alpha}_0}_{G,\bx_k}$.
\end{lem}

Of course, an analogous statement holds globally, i.e. when $D$ moves along
$\overset{\circ}\bX{}^{(a)}$.

\begin{proof}

For a fixed point $\F_G\in \Bun^a_G(\bS',\bD'_\infty)$, let
$\sigma_\bC$ be the corresponding based map $\bC\to
\Gr^{BD,a}_{G,\bX}$, and let $\sigma_\bX$ be the corresponding
based map $\bX\to \CG_{G,\bC}$.

Let us fix a point $\bc\in\bC$, which we may as well call
$0_\bC$. We must show that the value of $\sigma_\bC$ at $0_\bC$
belongs to $\underset{k}\Pi\, \ol{\Gr}{}^{n_k\cdot\ol{\alpha}_0}_{G,\bx_k}
\subset \Gr^{BD,a}_{G,\bX,D}$.

From \secref{parabolic notation}, we have a map
$\on{Maps}^a(\bX,\CG_{G,\bC})\to \on{Mod}^{a,+}_{G_{aff},\bX}$
covering the map $\varpi^a_v:\on{Maps}^a(\bX,\CG_{G,\bC})\to
\overset{\circ}\bX{}^{(a)}$.

\medskip

Now, it is easy to see that for our map $\sigma_\bX$,
$\varpi^a_v(\sigma_\bX)=D$ and the fiber of $\on{Mod}^{a,+}_{G_{aff},\bX}$
at $D$ is a closed subscheme of 
$\underset{k}\Pi\, \ol{\Gr}{}^{n_k\cdot\ol{\alpha}_0}_{G,\bx_k}$, cf. \cite{bgfm},
Prop. 1.7; moreover this embedding induces an isomorphism on the
level of reduced schemes. The resulting point of
$\underset{k}\Pi\, \ol{\Gr}{}^{n_k\cdot\ol{\alpha}_0}_{G,\bx_k}$ equals
the value of $\sigma_\bC$ at $0_\bC$, which is what we had to show.

\end{proof}

\bigskip

\centerline{{\bf Part II}:  {\Large Uhlenbeck spaces}}

\bigskip

Throughout Part II, $G$ will be a simple simply-connected group
and $\fg$ its Lie algebra. When $G=SL_n$, the subscript ``$SL_n$'' 
will often be replaced by just ``$n$''.

\section{Definition of Uhlenbeck spaces}

\ssec{Rational surfaces}   \label{rational surfaces}

As was explained in the introduction, Uhlenbeck spaces
$\fU_G^a$ are attached to the surface $\bS\simeq\BP^2$ with a distinguished
``infinity'' line $\bD_{\infty}\simeq \BP^1\subset \BP^2$.
However, in order to define $\fU^a_G$, we will need to replace $\bS$ by
all possible rationally equivalent surfaces isomorphic to
$\BP^1\times \BP^1$.

\medskip

Let $\overset{\circ}\bS\subset \bS$ denote the affine plane
$\bS-\bD_{\infty}$. For two distinct points $\bd_v,\bd_h\in \bD_{\infty}$
we obtain a decomposition of $\overset{\circ}\bS$ as a product of two
affine lines (horizontal and vertical):
$$\overset{\circ}\bS\simeq \overset{\circ}\bC\times
\overset{\circ}\bX,$$
where $\bd_v$ corresponds to the class of parallel lines in
$\overset{\circ}\bS$ that project to a single point in $\overset{\circ}\bC$,
and similarly for $\bd_h$. Let $\bC:=\overset{\circ}\bC\cup \infty_{\bC}$,
$\bX:=\overset{\circ}\bX\cup \infty_{\bX}$ be the corresponding
projective lines.

Let us denote by $\bS'$ the surface $\bC\times \bX$, and by
$\pi_v,\pi_h$ the projections from $\bS'$ to $\bC$ and $\bX$, respectively.
Let $\bD'_\infty:=\infty_{\bC}\times\bX\cup\bC\times\infty_{\bX}$ be the
corresponding divisor ``of infinity" in $\bS'$.
The surfaces $\bS$ and $\bS'$ are connected by a flip-flop. Namely, let
$\bS''$ be the blow-up of $\bS$ at the two points $\bd_v,\bd_h$. Then $\bS'$ is
obtained from $\bS''$ by blowing down the proper transform of
$\bD_\infty$.

In particular, it is easy to see that (a family of) $G$-bundles
on $\bS$ trivialized along $\bD_\infty$ is the same as (a family of)
$G$-bundles on $\bS'$ trivialized along $\bD'_\infty$.

\ssec{}

Let $\sO$ denote the variety
$\bD_{\infty}\times \bD_{\infty}-\Delta(\bD_{\infty})$, i.e., the
variety classifying pairs of distinct points $(\bd_v,\bd_h)\in \bD_{\infty}$,
and let us consider the ``relative over $\sO$'' versions of the
varieties discussed above.

In particular, we have the relative affine (resp., projective)
lines $\overset{\circ}\bC_{\sO}$, $\overset{\circ}\bX_{\sO}$
(resp., $\bC_{\sO}$, $\bX_{\sO}$), and the relative surface
$\bS'_{\sO}$. Let $\pi_{v,\sO}$ (resp., $\pi_{h,\sO}$) denote the
projection $\bS'_\sO\to \bC_\sO$ (resp., $\bS'_\sO\to \bX_\sO$.)
We will denote by $\overset{\circ}\bC{}_{\sO}^{(a)}$,
$\overset{\circ}\bX{}_{\sO}^{(a)}$ the corresponding fibrations into
symmetric powers.

\medskip

Let us recall the following general construction. Suppose that
$\Y_1$ is a scheme (of finite type),
and $\Y_2\to \Y_1$ is an affine morphism (also of finite type).
Then the functor on the category of schemes,
that sends a test scheme $S$ to the set of sections (i.e., $\Y_1$-maps)
$S\times \Y_1\to \Y_2$, is representable by an affine
ind-scheme of ind-finite type, which we will denote by $\on{Sect}(\Y_1,\Y_2)$.

(To show the representability, it is enough to assume that
$\Y_2$ is a total space of a vector bundle $\CE$, in which case
$\on{Sect}(\Y_1,\Y_2)$ is representable by the vector space
$\Gamma(\Y_1,\E)$.)

For example, by applying this construction
to $\Y_1=\sO$, we obtain the ind-schemes
$\on{Sect}(\sO,\overset{\circ}\bC{}^{(a)}_{\sO})$,
$\on{Sect}(\sO,\overset{\circ}\bX{}^{(a)}_{\sO})$,
$\on{Sect}(\sO,\overset{\circ}\bC{}^{(a)}_{\sO}\underset{\sO}\times
\overset{\circ}\bX{}^{(a)}_{\sO})$.

\ssec{}

We are now ready to give the first definition of the Uhlenbeck
space $\fU_G^a$.

For a fixed pair of directions $(\bd_v,\bd_h)\in \bD_{\infty}$, i.e.,
a point of $\sO$, and a divisor $D_v\in \overset{\circ}\bX{}^{(a)}$, consider
the scheme of based quasi-maps $\on{QMaps}^a(\bC,\Gr^{BD,a}_{G,\bX,D_v})$,
cf. \secref{canonical model}.
This ind-scheme is ind-affine, of ind-finite type, cf. \lemref{qmaps affine}.

By making $D_v\in \overset{\circ}\bX{}^{(a)}$ a parameter, we obtain an ind-affine
ind-scheme $\on{QMaps}^a(\bC,\Gr^{BD,a}_{G,\bX})$ fibered over
$\overset{\circ}\bX{}^{(a)}$.
Finally, by letting $(\bd_v,\bd_h)\in \sO$ move, we obtain ind-affine
fibrations
$$\on{QMaps}^a(\bC,\Gr^{BD,a}_{G,\bX})_{\sO}\to
\overset{\circ}\bX{}^{(a)}_{\sO}\to\sO.$$

Thus, we can consider the ind-scheme
$\on{Sect}(\sO,\on{QMaps}^a(\bC,\Gr^{BD,a}_{G,\bX})_{\sO})$.
By construction,
we have a natural map $$\varpi^a_{v,\sO}:
\on{Sect}(\sO,\on{QMaps}^a(\bC,\Gr^{BD,a}_{G,\bX})_{\sO})\to
\on{Sect}(\sO,\overset{\circ}\bX{}^{(a)}_{\sO}).$$
When a pair of directions $(\bd_v,\bd_h)$ is fixed, by further
evaluation we obtain the map
$\varpi^a_v:\on{Sect}(\sO,\on{QMaps}^a(\bC,\Gr^{BD,a}_{G,\bX})_{\sO})\to
\overset{\circ}\bX{}^{(a)}$.

\medskip

Now, we claim that we have a natural map
$$\Bun^a_G(\bS,\bD_\infty)\to
\on{Sect}(\sO,\on{QMaps}^a(\bC,\Gr^{BD,a}_{G,\bX})_{\sO}).$$
Indeed, constructing such a map amounts to giving a map
$\Bun^a_G(\bS,\bD_\infty)\simeq \Bun_G(\bS',\bD'_\infty)
\to \on{QMaps}^a(\bC,\Gr^{BD,a}_{G,\bX})$ for every pair of directions
$(\bd_v,\bd_h)\in \sO$, but this has been done in the previous
section, \propref{bundles and maps}.

Since $\Bun^a_G(\bS,\bD_\infty)$ is a scheme, its image
in $\on{Sect}(\sO,\on{QMaps}^a(\bC,\Gr^{BD,a}_{G,\bX})_{\sO})$ is contained
in a closed subscheme of finite type. In fact, it is contained in the
subscheme described in \secref{upper estimate}.

\begin{defn}
We define $\fU_G^a$ to be the closure of the image of $\Bun^a_G(\bS,\bD_\infty)$
in the ind-scheme
$\on{Sect}(\sO,\on{QMaps}^a(\bC,\Gr^{BD,a}_{G,\bX})_{\sO})$.
\end{defn}

By construction, $\fU_G^a$ is an affine scheme of finite type, which
functorially depends on the pair $(\bS,\bD_\infty)$.

\begin{lem}
The map $\Bun^a_G(\bS,\bD_\infty)\to \fU_G^a$ is an open embedding.
\end{lem}

\begin{proof}

Let us fix a pair of directions $(\bd_v,\bd_h)$ and consider
the corresponding evaluation map
$\on{Sect}(\sO,\on{QMaps}^a(\bC,\Gr^{BD,a}_{G,\bX})_{\sO})\to
\on{QMaps}^a(\bC,\Gr^{BD,a}_{G,\bX})$. We have a composition:
\begin{align*} 
&\Bun^a_G(\bS,\bD_\infty)\to \\
&\to\fU_G^a\underset{\on{QMaps}^a(\bC,\Gr^{BD,a}_{G,\bX})}\times
\on{Maps}^a(\bC,\Gr^{BD,a}_{G,\bX})\to \on{Maps}^a(\bC,\Gr^{BD,a}_{G,\bX})\simeq
\Bun^a_G(\bS,\bD_\infty).
\end{align*}

Since $\Bun^a_G(\bS,\bD_\infty)$ is dense in $\fU_G^a$, we obtain that
all arrows in the above formula are isomorphisms. The assertion of the
lemma follows since $\fU_G^a\underset{\on{QMaps}^a(\bC,\Gr^{BD,a}_{G,\bX})}
\times\on{Maps}^a(\bC,\Gr^{BD,a}_{G,\bX})$ is clearly open in $\fU_G^a$.

\end{proof}

\ssec{}

We will now give two more definitions of the space $\fU_G^a$,
which, on the one hand, are more economical, but on the other hand
possess less symmetry. Of course, later we will establish the
equivalence of all the definitions.

\medskip

Let us fix a point $(\bd_v,\bd_h)\in \bD_{\infty}$
and consider the space of quasi-maps $\on{QMaps}^a(\bC,\CG_{G,\bX})$.
As we have seen in the previous section, this is an affine
scheme of infinite type.

We have a natural map
$$\Bun^a_G(\bS,\bD_\infty)\to \on{QMaps}^a(\bC,\CG_{G,\bX})\times
\on{Sect}(\sO,\overset{\circ}\bX{}^{(a)}_{\sO}),$$
constructed as in the previous section.

We set $'\fU_G^a$ to be the closure of the image of
$\Bun^a_G(\bS,\bD_\infty)$ in the above product.
 From this definition, it is not immediately clear that
$'\fU_G^a$ is of finite type.

\medskip

Again, for a fixed pair of directions $(\bd_v,\bd_h)\in \bD_{\infty}$,
consider the fibration $\Gr_{G,\bX}^{BD,a}\to \overset{\circ}\bX{}^{(a)}$,
and consider the corresponding fibration of quasi-maps' spaces
$\on{QMaps}^a(\bC,\Gr^{BD,a}_{G,\bX})\to \overset{\circ}\bX{}^{(a)}$.
This is an affine ind-scheme of ind-finite type.

We have a natural map
$$\Bun^a_G(\bS,\bD_\infty)\to \on{QMaps}^a(\bC,\Gr^{BD,a}_{G,\bX})
\underset{\overset{\circ}\bX{}^{(a)}}
\times\on{Sect}(\sO,\overset{\circ}\bX{}^{(a)}_{\sO}),$$
where the projection $\on{Sect}(\sO,\overset{\circ}\bX{}^{(a)}_{\sO})\to
\overset{\circ}\bX{}^{(a)}$ corresponds to the evaluation at our fixed
point $(\bd_v,\bd_h)$.

We set $''\fU_G^a$ to be the closure of its image.
By construction, $''\fU_G^a$ is an affine scheme of finite type.

\begin{prop}  \label{prime=two primes}
There is a canonical isomorphism $'\fU_G^a\simeq {}''\fU_G^a$.
\end{prop}

\begin{proof}

Recall that we have a natural map $\Gr_{G,\bX}^{BD,a}\to \CG_{G,\bX}$, such that
the canonical line bundle on $\Gr_{G,\bX}^{BD,a}$ is the
restriction of that on $\CG_{G,\bX}$, moreover the map
$\Gr_{G,\bX}^{BD,a}\to \CG_{G,\bX}\times \overset{\circ}\bX{}^{(a)}$ is a closed
embedding.

Therefore, using \lemref{embedding into pro}, we obtain a closed embedding
$$\on{QMaps}^a(\bC,\Gr_{G,\bX}^{BD,a})\to\on{QMaps}^a(\bC,\CG_{G,\bX})\times
\overset{\circ}\bX{}^{(a)}.$$ This defines a map
$$\on{QMaps}^a(\bC,\Gr_{G,\bX}^{BD,a})\underset{\overset{\circ}\bX{}^{(a)}}
\times \on{Sect}(\sO,\overset{\circ}\bX{}^{(a)}_{\sO})\to
\on{QMaps}^a(\bC,\CG_{G,\bX})\times
\on{Sect}(\sO,\overset{\circ}\bX{}^{(a)}_{\sO}),$$
which is also a closed embedding. This proves the proposition.

\end{proof}

We have also a natural map $\fU^a_G\to{}''\fU^a_G$, which
corresponds to the evaluation map
$\on{Sect}(\sO,\on{QMaps}^a(\bC,\Gr^{BD,a}_{G,\bX})_{\sO})\to
\on{QMaps}^a(\bC,\Gr^{BD,a}_{G,\bX})$.

\begin{thm}  \label{equivalence of definitions}
The map $\fU^a_G\to{}''\fU^a_G$ is an isomorphism.
\end{thm}

This theorem will be proved in the next section for $G=SL_n$,
and in \secref{general case} for $G$ general.

\ssec{}

We conclude this section with the following observation:

Let $(\bd_v,\bd_h)$ be a fixed pair of directions. Using
\secref{projection to configurations}, we obtain a map
$$\on{QMaps}^a(\bC,\CG_{G,\bX})\to \overset{\circ}\bC{}^{(a)}.$$
By composing it with $\on{QMaps}^a(\bC,\Gr^{BD,a}_{G,\bX})\to
\on{QMaps}^a(\bC,\CG_{G,\bX})$, we obtain also the map
$\varpi^a_h:\on{QMaps}^a(\bC,\Gr^{BD,a}_{G,\bX})
\to \overset{\circ}\bC{}^{(a)}$.
By making $(\bd_v,\bd_h)$ vary along $\sO$, we have, therefore, a morphism:
$$\varpi^a_{h,\sO}:\on{Sect}(\sO,\on{QMaps}^a(\bC,\Gr^{BD,a}_{G,\bX})_{\sO})\to
\on{Sect}(\sO,\overset{\circ}\bC{}^{(a)}_{\sO}).$$

\medskip

Note now that the space $\sO$ carries a natural involution, which
interchanges the roles of $\bd_v$ and $\bd_h$. In particular, we have
a map
$\tau:\on{Sect}(\sO,\overset{\circ}\bX{}^{(a)}_{\sO})\to
\on{Sect}(\sO,\overset{\circ}\bC{}^{(a)}_{\sO})$.

We will denote by
$\on{Sect}(\sO,\on{QMaps}^a(\bC,\Gr^{BD,a}_{G,\bX})_{\sO})^\tau$ the
equalizer of the two maps
\begin{equation}
\tau\circ\varpi_{v,\sO} \text{ and }
\varpi_{h,\sO}:
\on{Sect}(\sO,\on{QMaps}^a(\bC,\Gr^{BD,a}_{G,\bX})_{\sO})\to
\on{Sect}(\sO,\overset{\circ}\bC{}^{(a)}_{\sO}).
\end{equation}

We have:
\begin{lem}  \label{transposition invariance}
$\fU_G^a\subset \on{Sect}(\sO,\on{QMaps}^a(\bC,\Gr^{BD,a}_{G,\bX})_{\sO})^\tau$.
\end{lem}

The proof follows from the fact that the maps $\tau\circ\varpi_{v,\sO}$ and
$\varpi_{h,\sO}$ do coincide on $\Bun^a_G(\bS,\bD_\infty)$, cf.
\propref{bundles and maps}.

\section{Comparison of the definitions: the case of $SL_n$}

\ssec{}

In order to prove \thmref{equivalence of definitions} for $SL_n$,
we will recall yet one more definition of Uhlenbeck spaces,
essentially due to S.~Donaldson. To simplify the notation, we will write
$\fU^a_n$ instead of $\fU^a_{SL_n}$.

In case $G=SL_n$, the moduli space $\Bun_n(\bS,\bD_\infty)$ admits the
following linear algebraic ADHM description going back to Barth (see a modern
exposition in ~\cite{n1}). We consider vector spaces $V=\BC^a,\ W=\BC^n$,
and consider the affine space
$$\End(V)\oplus\End(V)\oplus\Hom(W,V)\oplus\Hom(V,W),$$
a typical element of which will be denoted $(B_1,B_2,\imath,\jmath)$.
We define a subscheme 
$M^a_n\subset \End(V)\oplus\End(V)\oplus\Hom(W,V)\oplus\Hom(V,W)$
by the equation $[B_1,B_2]+\imath\jmath=0$.

A quadruple $(B_1,B_2,\imath,\jmath)\in M^a_n$ is called
{\em stable} if $V$ has no proper subspace containing the image of
$\imath$ and invariant with respect to $B_1,B_2$.
A quadruple $(B_1,B_2,\imath,\jmath)\in M^a_n$ is called {\em costable} if
Ker$(\jmath)$ contains no nonzero subspace
invariant with respect to $B_1,B_2$. We have the open subscheme $^{s}M^a_n\subset M^a_n$
(resp. $^{c}M^a_n\subset M^a_n$) formed by stable (resp. costable) quadruples.
Their intersection $^{s}M^a_n\cap {}^{c}M^c$ is denoted by $^{sc}M^a_n$.

According to ~\cite{n1} ~2.1, the natural action of $GL(V)$ on
$^{s}M^a_n$ (resp., $^{c}M^a_n$) is free, and the GIT quotient 
$\wt{\fN}^a_n:={}^{s}M^a_n/GL(V)$ is canonically isomorphic to the
fine moduli space of torsion free sheaves of rank $n$ and second Chern
class $a$ on $\bS$ equipped with a trivialization on $\bD_\infty$.
The open subset $^{sc}M^a_n/GL(V)\subset {}^{s}M^a_n/GL(V)$ corresponds, under this
identification, to the locus of vector bundles
$\Bun^a_n(\bS,\bD_\infty)\subset \wt{\fN}^a_n$.

Finally, consider the categorical quotient $\fN^a_n:=M^a_n//GL(V)$.
The natural projective morphism $\wt{\fN}^a_n \to\fN_n^a$
is the affinization of $\wt{\fN}^a_n$, i.e. $\fN_n^a$ is the spectrum
of the algebra of regular functions on $\wt{\fN}^a_n$. Moreover,
$\fN^a_n$ is reduced and irreducible and the natural map
$\Bun_n(\bS,\bD_\infty)\to \fN^a_n$ is an open embedding.

Our present goal is to construct a map $\fN_n^a\to \fU^a_n$ and
show that the composition $\fN_n^a\to \fU^a_n\to {}'\fU^a_n$
is an isomorphism. The proof given below was indicated
by Drinfeld, and the main step is to interpret $\fN_n^a$
as a coarse moduli space of {\it coherent perverse sheaves} (cf. ~\cite{b}) on
$\bS$.

\ssec{}

Let $\bS$ be a smooth surface.

\begin{defn} We will call a complex
$\CM\in \on{DCoh}^{\geq 0,\leq 1}(\bS)$ a coherent perverse sheaf if

\begin{itemize}

\item $h^1(\CM)$ is a finite length sheaf.

\item $h^0(\CM)$ is a torsion free coherent sheaf,

\end{itemize}
\end{defn}

Coherent perverse sheaves obviously form an additive
subcategory of $\on{DCoh}(X)$, denoted $\on{PCoh}(X)$.
It is easy to see that Serre-Grothendieck duality maps
$\on{PCoh}(X)$ to itself.

\begin{lem} \label{no negative ext}
For $\CM_1,\CM_2\in \on{PCoh}(X)$, the inner Hom satisfies
$\underline{\on{RHom}}^i(\CM_1,\CM_2)=0$ for $i< 0$.
\end{lem}

\medskip

If $S$ is a scheme, an $S$-family of coherent perverse sheaves
on $\bS$ is an object $\CM$ of $\on{DCoh}(\bS\times S)$, such that for
every geometric point $s\in S$ the (derived) restriction
$\CM_s\in \on{DCoh}(\bS)$ belongs to $\on{PCoh}(\bS)$.

\begin{lem}  \label{perverse sheaves are a stack}
The functor Schemes $\to$ Groupoids, that assigns to a test scheme $S$
the full subgroupoid of $\on{DCoh}(\bS\times S)$ consisting
of $S$-families of coherent perverse sheaves on $\bS$, is a sheaf
of categories in the faithfully flat topology.
\end{lem}

The lemma is proved in exactly the same manner as the usual
faithfully flat descent theorem for sheaves, using the fact that
for two $S$-families of coherent perverse sheaves $\M_1$ and $\M_2$,
the cone of any arrow $\M_1\to \M_2$ is a {\it canonically} defined
object of $\on{DCoh}(\bS\times S)$, which follows from
\lemref{no negative ext}.

\ssec{}   \label{coherent sheaves on P^2}

For $\bS=\BP^2$, consider the functor Schemes $\to$ Groupoids,
denoted $\on{Perv}_n^a(\bS,\bD_\infty)$, that
assigns to a test scheme $S$ the groupoid whose objects are
$S$-families of coherent perverse sheaves $\CM$ on $\bS$, such that:

\begin{itemize}

\item $\CM$ is of generic rank $n$, $ch_2(\CM)=-a$.

\item In a neighborhood of $\bD_\infty\times S\subset \bS\times S$,
$\CM$ is a vector bundle, and its restriction to the divisor
$\bD_\infty\times S$ is trivialized.

\end{itemize}
Morphisms in this category are isomorphisms between
coherent perverse sheaves (as objects in $\on{DCoh}(\bS\times S)$),
which respect the trivialization at $\bD_\infty$.

\medskip

The following theorem, due to Drinfeld, is a generalization of Donaldson-Nakajima theory:

\begin{thm}   \label{drinfeld}
The functor $\on{Perv}_n^a(\bS,\bD_\infty)$ is representable
by the stack $M^a_n/GL(V)$.
\end{thm}

\begin{proof}

The proof is a modification of Nakajima's argument in Chapter 2 of
~\cite{n1}. Let us choose homegeneous coordiantes $z_0,z_1,z_2$ on $\BP^2$,
so that the line $\bD_\infty\subset\bS$ is given by equation
$z_0=0$. In particular, this defines a pair of directions $(\bd_v,\bd_h)\in \sO$. 

\begin{lem} For a coherent perverse sheaf $\CM$ on $\bS$ trivialized
at $\bD_\infty$ we have $$H^{\pm1}(\bS,\CM(-1)[1])=H^{\pm1}(\bS,\CM(-2)[1])=
H^{\pm1}(\bS,\CM\otimes\Omega^1[1])=0.$$
\end{lem}

The proof for a torsion free sheaf (in cohomological degree 0)
is given in Chapter 2 of ~\cite{n1}.
The statement of the Lemma obviously holds also for a finite length
sheaf in cohomological degree $1$,
and an arbitrary $\CM$ is an extension of such two perverse sheaves. \qed

\medskip

Let $p_1,p_2$ be the two projections $\bS\times \bS\to \bS$, and 
let $C^\bullet$ be the Koszul complex on $\bS\times\bS$, i.e.,
the complex
$$0\to \CO(-2)\boxtimes \Omega^2(2)\to \CO(-1)\boxtimes \Omega^1(1)\to
\CO\boxtimes \CO\to 0,$$
which is known to be quasi-isomorphic to $\CO_{\Delta(\bS)}$. 
Then the complex of perverse coherent sheaves $p_2^*(\CM(-1))\otimes C^\bullet$
looks like
$$0\to\CO(-2)\boxtimes(\CM(-1)\otimes\Omega^2(2))\to
\CO(-1)\boxtimes(\CM(-1)\otimes\Omega^1(1))\to\CO\boxtimes \CM(-1)\to0$$
The first term of the Beilinson spectral sequence for
$Rp_{1*}(p_2^*(\CM(-1))\otimes C^\bullet)\simeq \CM(-1)$ reduces to
$$\CO(-2)\otimes H^1(\bS,\CM(-2))\to\CO(-1)\otimes
H^1(\bS,\CM\otimes\Omega^1)\to\CO\otimes H^1(\bS,\CM(-1))$$
(in degrees $-1,0,1$). Hence $\CM$ is canonically quasi-isomorphic to
the complex (monad)
$$\CO(-1)\otimes H^1(\bS,\CM(-2))\stackrel{d}{\to}\CO\otimes
H^1(\bS,\CM\otimes\Omega^1)\stackrel{b}{\to}\CO(1)\otimes
H^1(\bS,\CM(-1))$$ and $d$ is injective.

\medskip

Now we are able to go from perverse sheaves to the linear algebraic data
and back.
We set $V=H^1(\bS,\CM(-2)),\ W'=H^1(\bS,\CM\otimes\Omega^1),\
V'=H^1(\bS,\CM(-1))$. We have $\dim V=\dim V'=a$, and $\dim W'=2a+n$.
Since $H^0(\bS,\CO(1))$ has a base $\{z_0,z_1,z_2\}$, we may write
in the above monad $d=z_0d_0+z_1d_1+z_2d_2,\ b=z_0b_0+z_1b_1+z_2b_2$
where $d_i\in\Hom(V,W'),\ b_i\in\Hom(W',V')$.
Nakajima checks in {\em loc. cit.} that $b_1d_2=-b_2d_1$ is an isomorphism from
$V$ to $V'$, and identifies $V'$ with $V$ via this isomorphism.
Nakajima defines $W\subset W'$ as $W:=\Ker(b_1)\cap\Ker(b_2)$, and
identifies $W'$ with $V\oplus V\oplus W$ via
$(d_1,d_2):\ V\oplus V\leftrightarrows W'\ :(-b_2,b_1)$. Note that $\dim W=n$.
Under these identifications, we write
$V\stackrel{d_0}{\to}V\oplus V\oplus W\stackrel{b_0}{\to}V$ as
$d_0=(B_1,B_2,\jmath),\ b_0=(-B_2,B_1,\imath)$ where
$B_1,B_2\in\End(V),\ \imath\in\Hom(W,V),\ \jmath\in\Hom(V,W)$
satisfy the relation $[B_1,B_2]+\imath\jmath=0$.

\medskip

Conversely, given $V,W,B_1,B_2,\imath,\jmath$ as above,
we define $\CM$ as a monad
$$V\otimes\CO(-1)\stackrel{d}{\to}(V\oplus V\oplus W)\otimes\CO
\stackrel{b}{\to}V\otimes\CO(1)$$ (in cohomological degrees $-1,0,1$)
where $d=(z_0B_1-z_1,z_0B_2-z_2,z_0\jmath),\
b=(-z_0B_2+z_2,z_0B_1-z_1,z_0\imath)$.
Evidently, $\CM|_{\bD_\infty}=W\otimes\CO_{\bD_\infty}$.

\end{proof}

\medskip

If $\M$ is a coherent perverse sheaf on $\bS=\BP^2$ trivialized along
$\bD_\infty\simeq \BP^1$, we will denote by the same character $\M$ 
the corresponding coherent perverse sheaf on $\bS'=\bC\times \bX$, 
trivialized along the divisor $\bD'_\infty$.

The following assertion can be deduced from the above proof of
\thmref{drinfeld}

\begin{lem} \label{char pol of B}
The maps $\varpi^a_h$, $\varpi^a_v$ from
$\Bun_n^a(\bS,\bD_\infty)$ to $(\BA^1)^{(a)}\simeq \overset{\circ}\bC{}^{(a)}$ and
$(\BA^1)^{(a)}\simeq \overset{\circ}\bX{}^{(a)}$, respectively
(cf. \propref{bundles and maps}), are given
by the characteristic polynomials of $B_1$ and $B_2$.
\end{lem}

For the proof one has to observe that if $\M$ is a point of
$\Bun_n^a(\bS,\bD_\infty)$,
corresponding to a quadruple $(B_1,B_2,\imath,\jmath)$, and $D_h$ (resp., $D_v$) is the
divisor on $\BA^1$ given by the characteristic polynomial of $B_1$ (resp., $B_2$),
then as a bundle on $\bS'$, $\M$ will be trivilialized on
$(\bC-D_h)\times \bX$ (resp., $\bC\times (\bX-D_v))$.

\ssec{}

We will now construct a map $\on{Perv}_n^a(\bS,\bD_\infty)\to \fU^a_n$.
First, for a fixed pair of directions $(\bd_v,\bd_h)$, we will construct
a map
\begin{equation} \label{from Nakajima to quasi-maps}
\on{Perv}_n^a(\bS,\bD_\infty)\to \on{QMaps}^a(\bC,\CG_{SL_n,\bX}).
\end{equation}

\medskip

Thus, let $\CM$ be an $S$-point of $\on{Perv}_n^a(\bS,\bD_\infty)$,
or $\on{Perv}_n^a(\bS',\bD'_\infty)$.
Over the open subscheme $(\bC\times S)_0$ of $\bC\times S$,
over which $\CM$ is a vector bundle, we do obtain a genuine
map $(\bC\times S)_0\to \CG_{SL_n,\bX}$, and we have to show that this map
extends as a quasi-map on the entire $\bC\times S$.

Let us recall the description of the fundamental representation of
$\hsl_n$ as a Clifford module
${\mathsf C}{\mathsf l}{\mathsf i}{\mathsf f}{\mathsf f}_n$, cf. \cite{fgk},
Section 2.2. Thus, we have to attach to $\CM$ a line bundle $\L_{\CM}$ on
$\bC\times S$, and a map
$$\L_{\CM}\to \O_{\bC\times S}\otimes
\Lambda^{\bullet}\left(x^{-d_1}\BC^n[[x]]/x^{d_2}\BC^n[[x]]\right)\otimes
\left(\on{det}(x^{-d_1}\BC^n[[x]]/\BC^n[[x]])\right)^{-1}$$
for each pair of positive integers $d_1$ and $d_2$.

\medskip

Consider the (derived) direct image
$$\CN_{d}:=(\pi_v\times\on{id})_*\Bigl(\CM\bigl(d\cdot (\bC\times
\infty_{\bX}\times S)\bigr)\Bigr).$$

This is a complex on $\bC\times S$,
whose fiber at every  geometric point of $\bC\times S$ lies in the
cohomological degrees $0$ and $1$. Therefore, locally on $\bC\times S$,
$\CN$ can be represented by a length-2 complex of vector bundles
$\CF^0_{d}\to\CF^1_{d}$. We set $\L_{\CM}:=\on{det}(\CF^0_{d})\otimes
\on{det}(\CF^1_{d})^{-1}\otimes \on{det}(x^{-d}\BC^n[[x]]/\BC^n[[x]])$.

We have also a canonical map on $\bS'\times S$:
$$\CM\to \O^{\oplus n}/\O^{\oplus n}\bigl(-d_2\cdot (\bC\times
\infty_{\bX}\times S)\bigr),$$
which comes from a trivialization of $\CM$ around the divisor
$\bD'_h\times S$. Therefore, we obtain a map in the derived category
$\CN_{d_1}\to \O_{\bC\times S}\otimes (x^{-d_1}\BC^n[[x]]/x^{d_2}\BC^n[[x]])$.

Moreover, by replacing $\CF^0_{d_1}\to\CF^1_{d_1}$ by a
quasi-isomorphic complex of vector bundles, we can
assume that the above map
$\CN_{d_1}\to \O_{\bC\times S}\otimes (x^{-d_1}\BC^n[[x]]/x^{d_2}\BC^n[[x]])$
comes from a map
$\CF^0_{d_1}\to \O_{\bC\times S}\otimes (x^{-d_1}\BC^n[[x]]/x^{d_2}\BC^n[[x]])$.

Thus, we have a map $\CF^0_{d_1}\to \CF^1_{d_1}\oplus
\left(\O_{\bC\times S}\otimes
(x^{-d_1}\BC^n[[x]]/x^{d_2}\BC^n[[x]])\right)$, and, hence, a map
\begin{align*}
&\Lambda^{\on{rk}(\CF_{d_1}^0)}\left(\CF_{d_1}^0\right)\to
\Lambda^{\on{rk}(\CF_{d_1}^0)}\left(\CF_{d_1}^1\oplus \O_{\bC\times S}\otimes
(x^{-d_1}\BC^n[[x]]/x^{d_2}\BC^n[[x]])\right)\twoheadrightarrow \\
&\Lambda^{\on{rk}(\CF_{d_1}^1)}\left(\CF_{d_1}^1\right)\otimes
\Lambda^{\on{rk}(\CF_{d_1}^0)-\on{rk}(\CF_{d_1}^1)}\left(
\O_{\bC\times S}\otimes (x^{-d_1}\BC^n[[x]]/x^{d_2}\BC^n[[x]])\right).
\end{align*}

By tensoring both sides by
$\left(\on{det}(x^{-d_1}\BC^n[[x]]/\BC^n[[x]])\right)^{-1}$, we obtain a map
$\L_{\CM}\to \O_{\bC\times S}\otimes
\Lambda^{\bullet}\left(x^{-d_1}\BC^n[[x]]/x^{d_2}\BC^n[[x]]\right)\otimes
\left(\on{det}(x^{-d_1}\BC^n[[x]]/\BC^n[[x]])\right)^{-1}$, as required.

\medskip

It is easy to check that the definition of this morphism does not
depend on a particular choice of a representing complex, and over
$(\bC\times S)_0$, this is the same map as the one defining the map
$(\bC\times S)_0\to \CG_{SL_n,\bX}$.

\ssec{}

Thus, we have a map $\on{Perv}_n^a(\bS,\bD_\infty)\to
\on{QMaps}^a(\bC,\CG_{SL_n,\bX})$ for every fixed pair of directions
$(\bd_v,\bd_h)$. In particular, we have a map
$\on{Perv}_n^a(\bS,\bD_\infty)\to \overset{\circ}\bC{}^{(a)}$,
and by letting $(\bd_v,\bd_h)$ vary, we obtain a map
$\on{Perv}_n^a(\bS,\bD_\infty)\to
\on{Sect}(\sO, \overset{\circ}\bC{}^{(a)}_{\sO})$. By interchanging
the roles of $\bd_h$, and $\bd_v$, we obtain also a map
$\on{Perv}_n^a(\bS,\bD_\infty)\to
\on{Sect}(\sO, \overset{\circ}\bX{}^{(a)}_{\sO})$.

\medskip

We claim now that for any fixed pair of directions $(\bd_v,\bd_h)$,
the map $$\on{Perv}_n^a(\bS,\bD_\infty)\to \overset{\circ}\bX{}^{(a)}
\times \on{QMaps}^a(\bC,\CG_{SL_n,\bX})$$ factors through the closed subscheme
$\on{QMaps}^a(\bC,\Gr_{SL_n,\bX}^{BD,a})\subset \overset{\circ}\bX{}^{(a)}
\times \on{QMaps}^a(\bC,\CG_{SL_n,\bX})$. Indeed, this is so, because
the corresponding fact is true over the open dense substack
$\Bun^a_n(\bS,\bD_\infty)\subset \on{Perv}_n^a(\bS,\bD_\infty)$. (The
density assertion is a corollary of \thmref{drinfeld}).

Therefore, we obtain a map
\begin{equation} \label{from perverse sheaves to quasi-maps}
\on{Perv}_n^a(\bS,\bD_\infty)\to
\on{Sect}(\sO,\on{QMaps}^a(\bC,\Gr^{BD,a}_{SL_n,\bX})_{\sO}).
\end{equation}
The image of this map lies in $\fU^a_n$, because the open dense
substack $\Bun^a_n(\bS,\bD_\infty)$ does map there.
In other words, we obtain a map $\on{Perv}_n^a(\bS,\bD_\infty)\to \fU^a_n$.

In particular, from \thmref{drinfeld}, we have a $GL(V)$-invariant map 
$M^a_n\to \fU^a_n$, and since $\fU^a_n$ is affine, we thus have a map 
$\fN^a_n\to \fU^a_n$. 

Moreover, from \lemref{char pol of B} we obtain that for a fixed pair
of directions $(\bd_v,\bd_h)$, the composition
$M^a_n\to \fU^a_n\overset{\varpi^a_h}\to
\overset{\circ}\bC{}^{(a)}\simeq (\BA^1)^{(a)}$ is given by the map sending
$B_1$ to its characteristic polynomial, and similarly for $\varpi^a_v$
and $B_2$. Indeed, this is so because the corresponding fact holds for
the open part $^{sc}M^a_n$.

\medskip

\noindent{\it Remark.}
Above we have shown that the map $\Bun^a_n(\bS,\bD_\infty)\to \fU^a_n$
extends to a $GL(V)$-invariant map $M^a_n\to \fU^a_n$. However, if
we used the results of \cite{fgk}, the existence of the latter map
could be proved differently:

In {\it loc. cit.} we constructed a map $\wt{\fN}^a_n\to \fU^a_n$,
using the interpretation of $\wt{\fN}^a_n$ as the moduli space
of torsion free coherent sheaves on $\bS$. Since $\fN^a_n$
is the affinization of $\wt{\fN}^a_n$, and $\fU^a_n$ is affine,
we do obtain a map $\fN^a_n\to \fU^a_n$, and hence a
map $M^a_n/GL(V)\to \fN^a_n$.

\begin{thm} \label{Uhlenbeck=Nakajima}
The maps $\fN^a_n\to \fU^a_n\to {}'\fU^a_n$
are isomorphisms.
\end{thm}

The rest of this section is devoted to the proof of this theorem.
Since all the three varieties that appear in \thmref{Uhlenbeck=Nakajima}
have a common dense open piece, namely, $\Bun^a_n(\bS,\bD_\infty)$,
it is sufficient to prove that the map
$\sff:\fN^a_n\to {}'\fU^a_n$ is an isomorphism.

\ssec{}

Consider the map $\wt{\fN}^a_n\to {}'\fU^a_n$.
It was shown in \cite{fgk}, that this map is proper. Since
$\fN^a_n$ is an affinization of $\wt{\fN}^a_n$, we obtain
that $\sff:\fN^a_n\to {}'\fU^a_n$ is finite.

Choose a point $0_{\bS}\in \overset{\circ}\bS$, and let $0_{\bC}$,
$0_{\bX}$ be the corresponding points on $\bC$ and $\bX$, respectively.
Consider the corresponding action of $\BG_m$ by dilations.
By transport of structure, we obtain a $\BG_m$-action on the scheme
$\on{QMaps}(\bC,\Gr^{BD,a}_{\bX,SL_n})$, and on
$\on{Sect}(\sO,\overset{\circ}\bX{}^{(a)}_{\sO})$.

It is easy to see that the $\BG_m$-action contracts both these
varieties to a single point. Namely,
$\on{Sect}(\sO,\overset{\circ}\bX{}^{(a)}_{\sO})$ is contracted to a
section, that assigns to every $(\bd^1_v,\bd^1_h)\in \bD_{\infty}$
the point $a\cdot 0_{\bX}$ in the symmetric power of
$\overset{\circ}\bX$. Using \propref{contraction of quasi-maps},
we obtain that $\on{QMaps}(\bC,\Gr^{BD,a}_{SL_n,\bX})$ is contracted
to the quasi-map, whose saturation is the constant map $\bC\to\Gr^{BD,a}_{SL_n,\bX}$,
corresponding to the trivial bundle and tautological trivialization,
with defect of order $a$ at $0_{\bC}$.
Let us denote by $\sigma^{\fU}$ the attracting point of $'\fU^a_n$ described
above.

\medskip

The above action of $\BG_m$ on $'\fU^a_n$ is covered by a natural
$\BG_m$-action on the stack $\on{Perv}_n^a(\bS,\bD_\infty)$. Moreover,
if we identify $\overset{\circ}\bS$ with $\BA^2$, such that
$0_{\bS}$ corresponds to the origin, the above action on
$\on{Perv}_n^a(\bS,\bD_\infty)$ corresponds to the canonical
$\BG_m$-action on the variety $\BA^2$ by homotheties.

The induced $\BG_m$-action on $\fN^a_n\simeq M^a_n//GL(V)$ contracts 
this variety to a single point, which we will denote by $\sigma^{\fN}$. 

Thus, it would be sufficient to show that the
scheme-theoretic preimage $\sff^{-1}(\sigma^{\fU})$ is in fact a
point-scheme corresponding to $\sigma^{\fN}$.

\ssec{}

Let $\sg^{\fN}$ (resp., $\sg^{\fU}$) denote the canonical map
from $\on{Perv}_n^a(\bS,\bD_\infty)$ to $\fN^a_n$ (resp., $\fU^a_n$).
We claim that it is sufficient to show that the inclusion
$(\sg^{\fN})^{-1}(z^{\fN}_0)\hookrightarrow
(\sg^{\fU})^{-1}(z^{\fU}_0)$ is an equality. This follows
from the next general observation:

\begin{lem}
Let $\Y_1$ be an affine algebraic variety (in char 0) with an action
of a reductive group $G$. Let $\Y_2$ be another affine variety
and $\sg:\Y_1\to \Y_2$ be a $G$-invariant map, and let us denote by
$\sff:\Y_1//G\to \Y_2$, $\sg':\Y_1\to\Y_1//G$ the corresponding maps.
Then if for some $z'\in \Y_1//G$, $z\in \Y_2$ with $\sff(z')=z$
the inclusion $(\sg')^{-1}(z)\subset (\sg)^{-1}(z')$ is an isomorphism,
then $(\sff)^{-1}(z)$ is a point-scheme.
\end{lem}

\medskip

Consider now the following ind-stack $\on{Perv}_n^a(\bS,\bS-0_{\bS})$:
For a scheme $S$, its $S$-points are $S$-families of coherent perverse 
sheaves on $\bS$ (with $c_1=0$, $ch_2=-a$), equppied with a trivialization on
$\bS-0_{\bS}$ and such that for every pair of directions
$(\bd_v,\bd_h)\in \sO$, the composition
$$S\to \on{Perv}_n^a(\bS)\to \on{QMaps}^a(\bC,\CG_{SL_n,\bX})\to \overset{\circ}\bC{}^{(a)}$$
maps to the point $a\cdot 0_\bC\in \overset{\circ}\bC{}^{(a)}$.

\begin{prop}
The composition
$\on{Perv}_n^a(\bS,\bS-0_{\bS})\to \on{Perv}_n^a(\bS,\bD_\infty)\to \fN^a_n$
is the constant map to the point $\sigma^\fN$.
\end{prop}

\begin{proof}

For a triple $(B_1,B_2,\imath,\jmath)$, representing a point 
of $\on{Perv}_n^a(\bS,\bD_\infty)$, let us denote by $T_W$ 
any endomorphism of the vector space $W$ obtained by composing
the maps $B_1,B_2,\imath,\jmath$, and by $T_V$ any similarly
obtained endomorphism of $V$.

It is easy to see that the space of regular functions on $\fN^a_n$
is obtained by taking matrix coefficients of all possible $T_W$'s
and traces of all possible $T_V$'s. 

\medskip

Let $\M$ be an $S$-point of $\on{Perv}_n^a(\bS,\bS-0_{\bS})$. 
For an integer $m$, let $\M'$ be the (constant) $S$-family of
coherent perverse sheaves on $\bS$, corresponding to the torsion-free  
sheaf 
$$\on{ker}\left(\CO^{\oplus n}\to  \left(\CO/{\fm_{0_\bS}^m}\right)^{\oplus n}\right),$$
where $\fm_{0_\bS}$ is the maximal ideal of the point $0_\bS$. Then,
when $m$ is large enough, we can find a map $\M'\to \M$, which respects
the trivializations of both sheaves on $\bS-0_\bS$. The cone of this
map is set-theoretically supported at $0_\bS$ and has cohomologies in
degrees $0$ and $1$.

\medskip

Let $(V,W,B_1,B_2,\imath,\jmath)$ and $(V',W',B'_1,B'_2,\imath',\jmath')$
denote the linear algebra data corresponding to $\M$ and $\M'$, respectively.
By unraveling the proof of \thmref{drinfeld}, we obtain that there are
maps $V'\to V$ and $W'\simeq W$, which commute with all the endomorphisms.

By the definition of $\M'$, $\jmath'=0$. From this we obtain that 
all the matrices $T_W$ vanish, and the only non-zero $T_V$-matrices are
of the form $B_1^{k_1}\circ B_2^{k_2}\circ...\circ B_1^{k_{m-1}}$. It 
remains to show that any such matrix is traceless.

\medskip

Note that for any matrices $T^1_V$ and $T^2_V$ as above, the trace of
$T^1_V\circ i\circ j\circ T^2_V\in \on{End}(V)$ equals the trace
of the corresponding endomorphism of $V'$, and, hence, vanishes.
The relation $[B_1,B_2]+\imath\circ\jmath=0$ implies that the trace of a 
matrix $B_1^{k_1}\circ B_2^{k_2}\circ...\circ B_1^{k_{m-1}}$ does not
depend on the order of the factors.
Therefore, it is sufficient to show that the characteristic polynomial
of a matrix $B_1+c\cdot B_2$ vanishes for all $c\in \BC$. 

However, using \cite{fgk} Lemma 3.5, for any such $c$, we can find a pair
of directions $(\bd_v,\bd_h)\in \sO$, such that this
characteristic polynomial equals to value of 
$M^a_n\to \fU^a_n\overset{\varpi^a_h}\to
\overset{\circ}\bC{}^{(a)}\simeq (\BA^1)^{(a)}$ at our point of
$M^a_n$.

\end{proof}

Using this proposition, we obtain that in order to prove
\thmref{Uhlenbeck=Nakajima} it would be sufficient to
show that any $S$-point of the stack $(\sg^{\fU})^{-1}(\sigma^{\fU})$
factors through an $S$-point of $\on{Perv}_n^a(\bS,\bS-0_{\bS})$.

\ssec{}

Let $\CM$ be an arbitrary $S$-family of coherent perverse
sheaves on $\bS'$ corresponding to an $S$-point of
$\on{Perv}_n^a(\bS,\bD_\infty)$. Let $(\bd_v,\bd_h)$ be a fixed
configuration, and let us denote by
$D_h\subset \overset{\circ}\bC\times S$,
$D_v\subset \overset{\circ}\bX\times S$ the corresponding divisors.

First of all, by unfolding the definition of the map
$$\on{Perv}_n^a(\bS,\bD_\infty)\to \on{QMaps}(\bC,\CG_{SL_n,\bX})\to
\overset{\circ}\bC{}^{(a)},$$
we obtain that $\CM$ is a {\it vector bundle} away from the
divisors $D_h\times \bX$ and $\bC\times D_v$.
(I.e. a quasi-map necessarily acquires a defect at some point of $\bC$
if $\CM$ has a singularity on the vertical divisor over this point.)

Moreover, when we view
$\CM|_{\bS'\times S-\bC\times D_v}$ as an
$(\bX\times S-D_v)$-family of bundles on $\bC$,
this family is canonically trivialized.

Consider now the restriction of $\CM$ to the complement to
$D_h\times \bX\cup \bC\times D_v\subset \bS\times S$. This is
a vector bundle equipped with two trivializations.
One comes from the fact that we are dealing with a {\it map}
$\bC\times S-D_h\to \Gr^{BD,a}_{SL_n,\bX,D_v}$, and another one comes from
the trivialization of $\CM|_{\bS'\times S-\bC\times D_v}$ mentioned
above.

We claim that these two trivializations actually coincide. This
is so, because the corresponding fact is true over the dense
substack $\Bun_n^a(\bS,\bD_\infty)\subset
\on{Perv}_n^a(\bS,\bD_\infty)$.

\medskip

Going back to the proof of the theorem, assume that $\CM$
corresponds to an $S$-point of $(\sg^{\fU})^{-1}(\sigma^{\fU})$.
Then, first of all, the divisors $D_h$, $D_v$ are $(a\cdot 0_{\bC})\times S$
and $(a\cdot 0_{\bX})\times S$, respectively.

Moreover, $\CM$ is a vector bundle away from $0_{\bS}\times S$,
and it is trivialized away from $\bC\times 0_{\bX}\times S\subset \bS'\times S$.
Therefore, we only have to show that this trivialization extends
across the divisor $\bC\times 0_{\bX}\times S$ over $(\bC-0_{\bC})\times S$.

\medskip

However, by the definition of $\sigma^{\fU}$, the {\it map}
$(\bC-0_{\bC})\times S\to \Gr_{SL_n,\bX}\to
\CG_{SL_n}$ is the constant map corresponding to the trivial bundle.
Therefore, the trivialization does extend.

\section{Properties of Uhlenbeck spaces}

\ssec{}

Let $\phi:G_1\to G_2$ be a homomorphism of simple simply connected
groups, and let $\phi_{\BZ}$ denote the corresponding homomorphism
$\phi_{\BZ}:\BZ\simeq H_3(G_1,\BZ)\to H_3(G_2,\BZ)\simeq \BZ$.
Observe that for a curve $\bX$, the pull-back of the canonical line bundle
$\CP_{\Bun_{G_2}(\bX)}$ under the induced map
$\Bun_{G_1}(\bX)\to \Bun_{G_2}(\bX)$, is
$\CP^{\otimes \phi_{\BZ}(1)}_{\Bun_{G_1}(\bX)}$. In particular,
the map $\Bun_{G_1}(\bS,\bD_\infty)\to\Bun_{G_2}(\bS,\bD_\infty)$
sends $\Bun^a_{G_1}(\bS,\bD_\infty)$ to
$\Bun^{\phi_{\BZ}(a)}_{G_2}(\bS,\bD_\infty)$.

\begin{lem}    \label{functoriality}
The map $\Bun^a_{G_1}(\bS,\bD_\infty)\to
\Bun^{\phi_{\BZ}(a)}_{G_2}(\bS,\bD_\infty)$ extends to morphisms
$\fU^a_{G_1}\to \fU^{\phi_{\BZ}(a)}_{G_2}$ and
$''\fU^a_{G_1}\to {}''\fU^{\phi_{\BZ}(a)}_{G_2}$. When
$\phi$ is injective, both these maps are
closed embeddings.
\end{lem}

\begin{proof}

First, for a curve $\bX$, from $\phi_{\BZ}(a)$
we obtain a map $\bX^{(a)}\to \bX^{(\phi_{\BZ}(a))}$, and the corresponding
map $\Gr_{G_1,\bX}^{BD,a}\to \Gr_{G_2,\bX}^{BD,\phi_{\BZ}(a)}$. The last
morphism is closed embedding if $\phi$ is.

Hence, according to \propref{changing the line bundle}(c),
for a curve $\bC$ with a marked point $\bc\in\bC$
we obtain a morphism of the corresponding
based quasi-maps spaces:
$$\phi_{\on{QMaps}}:\on{QMaps}^a(\bC,\on{Gr}^{BD,a}_{G_1,\bX})\to
\on{QMaps}^{\phi_{\BZ}(a)}
(\bC,\on{Gr}^{BD,\phi_{\BZ}(a)}_{G_2,\bX}).$$

The existence of the maps $\fU^a_{G_1}\to \fU^{\phi_{\BZ}(a)}_{G_2}$ and
$''\fU^a_{G_1}\to {}''\fU^{\phi_{\BZ}(a)}_{G_2}$ follows now from the definition of
$\fU^a_G$ and $''\fU^a_G$. .

\end{proof}

\ssec{Comparison of the definitions: the general case} \label{general case}

The above lemma allows us to establish the equivalence of
the two definitions of $\fU^a_G$:

\begin{proof} (of \thmref{equivalence of definitions})

Let us choose a faithful representation $\phi:G\to SL_n$. By the previous
lemma, $\fU^a_G$ and $''\fU^a_G$
are isomorphic to the closures of $\Bun^a_G(\bS,\bD_\infty)$
in $\fU^{\phi_{\BZ}(a)}_n(\bS,\bD_\infty)$ and
$''\fU^{\phi_{\BZ}(a)}_n(\bS,\bD_\infty)$, respectively.

However, by \thmref{Uhlenbeck=Nakajima}, the map
$\fU^{\phi_{\BZ}(a)}_n(\bS,\bD_\infty)\to {}''\fU^{\phi_{\BZ}(a)}_n(\bS,\bD_\infty)$
is an isomorphism. Hence, $\fU^a_G\to {}''\fU^a_G$ is an isomorphism as well.

\end{proof}

\ssec{Factorization property}

Next, we will establish the factorization property of Uhlenbeck
compactifications. Let us fix a pair of
non-parallel lines $(\bd_v,\bd_h)$, and consider the corresponding projection
$\varpi^a_h:\fU^a_G\to \overset{\circ}\bC{}^{(a)}$.

\begin{prop} \label{factorization of Uhlenbeck}
For $a=a_1+a_2$, there is a natural isomorphism
$$(\overset{\circ}\bC{}^{(a_1)}\times
\overset{\circ}\bC{}^{(a_2)})_{disj}\underset{\overset{\circ}\bC{}^{(a)}}\times
\fU^a_G\simeq (\overset{\circ}\bC{}^{(a_1)}\times
\overset{\circ}\bC{}^{(a_2)})_{disj}\underset{\overset{\circ}\bC{}^{(a_1)}\times
\overset{\circ}\bC{}^{(a_2)}}\times (\fU^{a_1}_G\times \fU^{a_2}_G).$$
\end{prop}

To prove this proposition, we will first consider the case of
$SL_n$:

\begin{prop} \label{Nak fac}
For $a=a_1+a_2$ there are natural isomorphisms:
$$(\overset{\circ}\bC{}^{(a_1}\times
\overset{\circ}\bC{}^{(a_2)})_{disj}\underset{\overset{\circ}\bC{}^{(a)}}\times
\wt{\fN}^a_n\simeq (\overset{\circ}\bC{}^{(a_1)}
\times\overset{\circ}\bC{}^{(a_2)})_{disj}\underset
{\overset{\circ}\bC{}^{(a_1)}\times\overset{\circ}\bC{}^{(a_2)}}\times
(\wt{\fN}^{a_1}_n\times\wt{\fN}^{a_2}_n),$$
$$(\overset{\circ}\bC{}^{(a_1}\times\overset{\circ}\bC{}^{(a_2)})_{disj}
\underset{\overset{\circ}\bC{}^{(a)}}\times
\fN^a_n\simeq
(\overset{\circ}\bC{}^{(a_1)}\times\overset{\circ}\bC{}^{(a_2)})_{disj}\underset
{\overset{\circ}\bC{}^{(a_1)}\times\overset{\circ}\bC{}^{(a_2)}}\times
(\fN^{a_1}_n\times\fN^{a_2}_n).$$
\end{prop}

\begin{proof} The second factorization isomorphism follows from the first
one since $\fN^a_n$ is the affinization of $\wt{\fN}_n^a$.

Let $\Coh_n(\bX,\infty_{\bX})$ be the stack of coherent sheaves on $\bX$
with a trivialization at $\infty_{\bX}$. This is a smooth stack, which
contains an open subset isomorphic to the point scheme, that
corresponds to the trivial rank-$n$ vector bundle. The complement
to this open subset is of codimension $1$, and hence is a Cartier
divisor.

Observe now that the scheme $\wt{\fN}_n^a$ represents the functor
of maps,  $\on{Maps}^a(\overset{\circ}\bC,\Coh_n(\bX,\infty_{\bX}))$,
where the latter is as in \secref{factorization principle}.

Indeed, an $S$-point of $\on{Maps}^a(\overset{\circ}\bC,\Coh_n(\bX,\infty_{\bX}))$,
is according to \lemref{extension} the same as coherent sheaf $\CM$ on
$\bX\times \bC\times S$, trivialized over $\infty_\bX\times \bC\times S$
and $\bX\times \infty_\bC\times S$, and which is $\bC\times
S$-flat. But this implies that for every geometric point $s\in S$, the
restriction $\CM|_{\bS\times s}$ is torsion free.

Hence, the assertion about the factorization of $\wt{\fN}_n^a$ follows
from \propref{factorization pattern}.

\end{proof}

Now we can prove \propref{factorization of Uhlenbeck}:

\begin{proof}

Let us choose a faithful representation $\phi:G\to SL_n$ and consider
the corresponding closed embedding
$$\fU^a_G\to \on{QMaps}^a(\bC,\CG_{G,\bX})\times \fU^{\phi_{\BZ}(a)}_n.$$
The image of this map lies in the closed subscheme
$$\on{QMaps}^a(\bC,\CG_{G,\bX})\underset{\overset{\circ}\bC{}^{(\phi_{\BZ}(a))}}\times
\fU^{\phi_{\BZ}(a)}_n\simeq
\on{QMaps}^a(\bC,\CG_{G,\bX})\underset{\overset{\circ}\bC{}^{(a)}}\times
(\overset{\circ}\bC{}^{(a)}\underset{\overset{\circ}\bC{}^{(\phi_{\BZ}(a))}}\times
\fU^{\phi_{\BZ}(a)}_n).$$

Now, $\on{QMaps}^a(\bC,\CG_{G,\bX})$ factorizes over
$\overset{\circ}\bC{}^{(a)}$
according to \propref{factorization of Zastava}, and the fiber product
$\overset{\circ}\bC{}^{(a)}\underset{\overset{\circ}\bC{}^{(\phi_{\BZ}(a))}}\times
\fU^{\phi_{\BZ}(a)}_n$
factorizes over $\overset{\circ}\bC{}^{(a)}$ because $\fU^{\phi_{\BZ}(a)}_n$ does so
over $\overset{\circ}\bC{}^{(\phi_{\BZ}(a))}$.
This implies the proposition in view of isomorphism \eqref{factorization of bundles}.

\end{proof}

\ssec{}  \label{independence of bd_h}

Let us fix a point $\bd_v\in \bD_\infty$ and observe that the curve
$\bC$ is well-defined without the additional choice of $\bd_h\in
\bD_\infty-\bd_v$. Indeed, $\bC$ is canonically identified with the
projectivization of the tangent space $T\bS_{\bd_v}$. Thus, the map
$\varpi^a_{h,\sO}$ gives rise to a map
$$\varpi^a_{h,\bd_v}:
(\bD_\infty-\bd_v)\times \fU^a_G \to \overset{\circ}\bC{}^{(a)}.$$

\begin{prop}  \label{independence of varpi_h of bd_h}
The above map $\varpi^a_{h,\bd_v}$ is independent of the variable
$\bD_\infty-\bd_v$. Moreover, the corresponding factorization
isomorphisms of \propref{factorization of Uhlenbeck} are also
independent of the choice of $\bd_h$.
\end{prop}

\begin{proof}

Since $\Bun^a_G(\bS,\bD_\infty)$ is dense in $\fU^a_G$, it is enough
to show that the corresponding map
$(\bD_\infty-\bd_v)\times \Bun^a_G(\bS,\bD_\infty)\to \overset{\circ}\bC{}^{(a)}$
is independent of the first $\bd_h$, and similarly for the
factorization isomorphisms.

Consider the surface $\bS^{\bd_v}$ obtained by blowing up $\bS$ at
$\bd_v$. Then the exceptional divisor identifies canonically with $\bC$,
and we have a projection $\pi_v^{\bd_v}:\bS^{\bd_v}\to \bC$. We can
consider the corresponding stack $\Bun^a_G(\bS_{\bd_v},\bD^{\bd_v}_{\infty})$,
which classifies $G$-bundles with a trivialization along
$\bD^{\bd_v}_{\infty}:=\left(\bC\cup
(\pi_v^{\bd_v})^{-1}(\infty_\bC)\right)\subset \bS^{\bd_v}$.

In addition, $\bS^{\bd_v}$ can be regarded as a relative curve
$\bX_{\bC}$ with a marked infinity over $\bC$, and we can consider the
relative stack $\Bun_G(\bX_{\bC},\infty_{\bX_{\bC}})$. This stack
contains an open substack (corresponding to the trivial bundle) and
its complement is a Cartier divisor.

The corresponding space of sections $\bC\to
\Bun_G(\bX_{\bC},\infty_{\bX_{\bC}})$
identifies with $\Bun^a_G(\bS_{\bd_v},\bD^{\bd_v}_{\infty})$,
and thus gives rise to a map
$$\Bun^a_G(\bS_{\bd_v},\bD^{\bd_v}_{\infty})\to \overset{\circ}\bC{}^{(a)},$$
and the corresponding factorization isomorphisms,
cf. \propref{factorization pattern}.
This makes the assertion of the proposition manifest.

\end{proof}

\ssec{}

Let us fix now a point $\bd_h\in \bD_\infty$, and consider the
family of curves $\bC_{\bD_\infty-\bd_h}$ corresponding to the moving
point $\bd_v\in \bD_\infty-\bd_h$.

We have a natural forgetful map
$$\on{Sect}(\sO,\overset{\circ}\bC{}_{\sO}^{(a)})\to
\on{Sect}(\bD_\infty-\bd_h,\overset{\circ}\bC{}_{\bD_\infty-\bd_h}^{(a)}),$$
and by composing with $\varpi^a_{h,\sO}$ we obtain a map
$$\fU^a_G\to\on{Sect}(\bD_\infty-\bd_h,
\overset{\circ}\bC{}_{\bD_\infty-\bd_h}^{(a)}).$$
From
\propref{independence of varpi_h of bd_h} we obtain the following
corollary:

\begin{cor}  \label{fewer lines}
The map
$$\fU^a_G\to \on{QMaps}^a(\bC,\CG_{G,\bX})\times
\on{Sect}(\bD_\infty-\bd_h,\overset{\circ}\bC{}_{\bD_\infty-\bd_h}^{(a)})$$
is a closed embedding, where $\bC$ corresponds to some fixed point of
$\bD_\infty-\bd_h$.
\end{cor}

\begin{proof}

We know that
$\fU^a_G\to \on{Maps}^a(\bC,\CG_{G,\bX})\times
\on{Sect}(\sO,\overset{\circ}\bC{}_{\sO}^{(a)})$ is a closed
embedding, and if we consider the open subset
$\overset{\circ}\sO:=\sO-\bd_h\times (\bD_\infty-\bd_h)$, the map
$$\fU^a_G\to \on{Maps}^a(\bC,\CG_{G,\bX})\times
\on{Sect}(\overset{\circ}\sO,\overset{\circ}\bC{}_{\overset{\circ}\sO}^{(a)})$$
would be a closed embedding as well.
However, from \propref{independence of varpi_h of bd_h}, we obtain
that the map
$\fU^a_G\to
\on{Sect}(\overset{\circ}\sO,\overset{\circ}\bC{}_{\overset{\circ}\sO}^{(a)})$
factors as
$$\fU^a_G\to
\on{Sect}(\bD_\infty-\bd_h,\overset{\circ}\bC{}_{\bD_\infty-\bd_h}^{(a)})\to
\on{Sect}(\overset{\circ}\sO,\overset{\circ}\bC{}_{\overset{\circ}\sO}^{(a)}),
$$
where the last arrow comes from the projection on the first factor
$\overset{\circ}\sO\to (\bD_\infty-\bd_h)$.

This establishes the proposition.
\end{proof}

\section{Stratifications and IC stalks}

\ssec{}

In this section we will introduce a stratification of
$\fU^a_G$ and formulate a theorem describing the intersection cohomology
sheaf $\on{IC}_{\fU^a_G}$.

Let $b$ be an integer $0\leq b\leq a$, and set by definition
$\fU^{a;b}_G:=\Bun^{a-b}_G(\bS,\bD_\infty)\times \on{Sym}^b(\overset{\circ}\bS)$.

\begin{thm}  \label{stratification}
There exists a canonical locally closed embedding
$\iota_b:\fU^{a;b}_G\to \fU^a_G$. Moreover, $\fU^a_G=\underset{b}\cup\, \fU^{a;b}_G$.
\end{thm}

\ssec{Proof of \thmref{stratification}}

Let $(\F'_G,D_{\bS})$ be a point of $\Bun^{a-b}_G(\bS,\bD_\infty)\times
\on{Sym}^b(\overset{\circ}\bS)$. For every pair of directions $(\bd_v,\bd_h)$,
the data of $\F'_G$ defines a based {\it map}
$\sigma':\bC\to \Gr^{BD,a-b}_{G,\bX}$ of degree $a-b$.

Consider the embedding $\Gr^{BD,a-b}_{G,\bX}\hookrightarrow \Gr^{BD,a}_{G,\bX}$
obtained by adding to the divisor $D_v\in \overset{\circ}\bX{}^{(a-b)}$
the divisor $\pi_h(D_{\bS})$. We define the sought-for quasi-map
$\sigma:\bC\to \Gr^{BD,a}_{G,\bX}$ by adding to the map
\begin{equation} \label{eq3}
\bC\overset{\sigma'}\longrightarrow \Gr^{BD,a-b}_{G,\bX}\to \Gr^{BD,a}_{G,\bX}
\end{equation}
the defect equal to $\pi_v(D_{\bS})$, cf. \secref{degenerations of quasi-maps} and
\secref{qmaps into flags} \eqref{strata of quasi-maps}.
Since this can be done for any pair of directions $(\bd_v,\bd_h)$,
we obtain a point of
$\on{Sect}(\sO,\on{QMaps}(\bC,\Gr^{BD,a}_{G,\bX})_\sO)$.
This morphism is a locally closed embedding due to the corresponding
property of the quasi-maps' spaces.

Note that the above construction defines in fact a map
$\ol{\iota}_b:\fU^{a-b}_G\times
\on{Sym}^b(\overset{\circ}\bS)\to
\on{Sect}(\sO,\on{QMaps}(\bC,\Gr^{BD,a}_{G,\bX})_\sO)$.

\medskip

Now our goal is to prove that the image of $\iota_b$ belongs to
$\fU^a_G$, and that every geometric point of $\fU^a_G$ belongs to
(exactly) one of the subschemes $\fU^{a;b}_G$. First, we shall do
that for $G=SL_n$.

Let $\wt{\fN}^{a;b}_n$ be the scheme classifying pairs
$\M\subset \M'$, where $\M\in \fN^a_n$, and
$\M'$ is a vector bundle on $\bS$, such that $\M'/\M$ is a
torsion sheaf of length $b$. We have a natural proper and surjective map
$$\wt{\fN}^{a;b}_n\to \Bun^{a-b}_n(\bS,\bD_\infty)\times
\on{Sym}^b(\overset{\circ}\bS),$$
and a locally-closed embedding $\wt{\fN}^{a;b}_n\hookrightarrow \wt{\fN}^a_n$.
Moreover, it is easy to see that the square:
\begin{equation}  \label{eq2}
\CD
\wt{\fN}^{a;b}_n  @>>> \wt{\fN}^a_n  \\
@VVV  @VVV  \\
\Bun^{a-b}_n(\bS,\bD_\infty)\times \on{Sym}^b(\overset{\circ}\bS)
@>\iota_b>>  \on{Sect}(\sO,\on{QMaps}(\bC,\Gr^{BD,a}_{SL_n,\bX})_\sO)
\endCD
\end{equation}
is commutative, where in the above diagram
the right vertical arrow is the composition of
$\wt{\fN}^a_n\to \fN^a_n\simeq \fU^a_n$ and
$\fU^a_n\hookrightarrow
\on{Sect}(\sO,\on{QMaps}(\bC,\Gr^{BD,a}_{SL_n,\bX})_\sO)$.

This readily shows that the image of $\iota_b$ belongs to
$\fU^a_G$ for $G=SL_n$. Moreover, since
$\wt{\fN}^a_n\simeq \underset{b}\cup\, \wt{\fN}^{a;b}_n$, we obtain that
$\fU^a_n\simeq \underset{b}\cup\, \fU^{a;b}_n$.

\medskip

Now let us treat the case of an arbitrary $G$. To show that
the image of $\iota_b$ belongs to $\fU^a_G$
consider the open subset in the product $\Bun^{a-b}_G(\bS,\bD_\infty)\times
\on{Sym}^b(\overset{\circ}\bS)$ corresponding to pairs
$(\F_G,D_{\bS})$ such that $\pi_v(D_{\bS})$ is a multiplicity-free
divisor, disjoint from $\varpi^{a-b}_h(\F_G)$. It would be enough
to show that the image of this open subset in
$\on{Sect}(\sO,\on{QMaps}(\bC,\Gr^{BD,a}_{G,\bX})_\sO)$
is contained in $\fU^a_G$.

Note, that from Equation \eqref{eq2} we obtain the following compatibility
relation of the map $\iota_b$ with the factorization isomorphisms
of \propref{factorization of Uhlenbeck} for $SL_n$. For $b=b_1+b_2$ the map
$$
\CD
(\Bun^{a-b}_n(\bS,\bD_\infty)\times
\on{Sym}^b(\overset{\circ}\bS))\underset{\overset{\circ}\bC{}^{(a)}}\times
\left(\overset{\circ}\bC{}^{(a-b)}\times \overset{\circ}\bC{}^{(b_1)}\times
\overset{\circ}\bC{}^{(b_2)}\right)_{disj}  \\
@V{\iota_b}VV \\
\fU^a_n \underset{\overset{\circ}\bC{}^{(a)}}\times
\left(\overset{\circ}\bC{}^{(a-b)}\times \overset{\circ}\bC{}^{(b_1)}\times
\overset{\circ}\bC{}^{(b_2)}\right)_{disj}
\endCD
$$
coincides with the composition
$$
\CD
(\Bun^{a-b}_n(\bS,\bD_\infty)\times
\on{Sym}^b(\overset{\circ}\bS))\underset{\overset{\circ}\bC{}^{(a)}}\times
\left(\overset{\circ}\bC{}^{(a-b)}\times \overset{\circ}\bC{}^{(b_1)}\times
\overset{\circ}\bC{}^{(b_2)}\right)_{disj}  \\
@V{\sim}VV  \\
(\Bun^{a-b}_n(\bS,\bD_\infty)\times \on{Sym}^{b_1}(\overset{\circ}\bS)\times 
\on{Sym}^{b_2}(\overset{\circ}\bS))
\underset{\overset{\circ}\bC{}^{(a-b)}\times \overset{\circ}\bC{}^{(b_1)}\times
\overset{\circ}\bC{}^{(b_2)}}\times
\left(\overset{\circ}\bC{}^{(a-b)}\times \overset{\circ}\bC{}^{(b_1)}\times
\overset{\circ}\bC{}^{(b_2)}\right)_{disj} \\
@V{\on{id}\times
\iota_{b_1}\times\iota_{b_2}}VV  \\
(\fU^{a-b}_n\times \fU^{b_1}_n \times \fU^{b_2}_n)
\underset{\overset{\circ}\bC{}^{(a-b)}\times \overset{\circ}\bC{}^{(b_1)}\times
\overset{\circ}\bC{}^{(b_2)}}\times
\left(\overset{\circ}\bC{}^{(a-b)}\times \overset{\circ}\bC{}^{(b_1)}\times
\overset{\circ}\bC{}^{(b_2)}\right)_{disj}  \\
@V{\sim}VV  \\
\fU^a_n \underset{\overset{\circ}\bC{}^{(a)}}\times
\left(\overset{\circ}\bC{}^{(a-b)}\times \overset{\circ}\bC{}^{(b_1)}\times
\overset{\circ}\bC{}^{(b_2)}\right)_{disj}.
\endCD
$$

Hence, by embedding $G$ into $SL_n$ and
using \lemref{functoriality}, we reduce our assertion to the case when
$b=a=1$. Using the action of the group of affine-linear transformations,
we reduce the assertion further to the fact that the image of
the point-scheme,
thought of as $\Bun_G^0(\bS,\bD_\infty)\times 0_\bS$, belongs to $\fU^1_G$.

However, as in the proof of \thmref{Uhlenbeck=Nakajima}, the above point-scheme
is the only attractor of the $\BG_m$-action on
$\on{Sect}(\sO,\on{QMaps}(\bC,\CG_{G,\bX})_\sO)$.
In particular, it is contained
in the closure of $\Bun_G^1(\bS,\bD_\infty)$, which is $\fU^1_G$.

\medskip

To finish the proof of the theorem, we have to show that every geometric
point of $\fU^a_G$ belongs to (exactly) one of the subschemes $\fU^{a;b}_G$.
That these subschemes are mutually disjoint is clear from the fact
that $b$ can be recovered as a total defect of a quasi-map.

Let $\sigma_\sO$ be a point of $\fU^a_G$, thought of as a collection
of quasi-maps $\sigma:\bC\to \Gr^{BD,a}_{G,\bX}$ for every pair
of directions $(\bd_v,\bd_h)$. We have to show that there exists a
point $\sigma'_\sO\in \Bun^{a-b}_G(\bS,\bD_\infty)$, and a $0$-cycle
$D_\bS=\Sigma b_i\cdot \bs_i$ on $\bS$, such that for every pair of
directions $(\bd_v,\bd_h)$, the corresponding quasi-map
$\sigma:\bC\to \Gr^{BD,a}_{G,\bX}$ has as its saturation the
composition
$$\bC\overset{\sigma'}\longrightarrow \Gr^{BD,a-b}_{G,\bX}\hookrightarrow 
\Gr^{BD,a}_{G,\bX},$$
(the last arrow is as in \eqref{eq3}), with defect $\pi_v(D_\bS)$.
(In what follows we will say that the point $\sigma_\sO\in \fU^a_G$
has saturation equal to the $G$-bundle $\F'_G$ corresponding to
$\sigma'_\sO$, and the defect given by the $0$-cycle $D_\bS$.)

For a faithful representation $\phi:G\to SL_n$, let $\sigma_{\sO,n}$
be the corresponding point of $\fU^{\phi_{\BZ}(a)}_n$.
We know that $\sigma_{\sO,n}$
can be described as a pair $(\sigma'_{\sO,n},D_{\bS,n})$, where
$\sigma'_{\sO,n}$ is the saturation of $\sigma_{\sO,n}$, i.e.
a point of $\Bun^{\phi_{\BZ}(a-b)}_n(\bS,\bD_\infty)$, and
$D_{\bS,n}=\Sigma\, b_{i,n}\cdot \bs_i$ is a $0$-cycle on $\bS$.

\medskip

For a fixed pair of directions $(\bd_v,\bd_h)$, let $\sigma''$ be the
saturation of $\sigma$. In particular, $\sigma''$ defines a $G$-bundle on $\bS$,
with a trivialization along $\bD_\infty$, such that the induced
$SL_n$-bundle is isomorphic to $\sigma'_{\sO,n}$. Since the morphism
$\Bun_G(\bS,\bD_\infty)\to\Bun_n(\bS,\bD_\infty)$ is an embedding,
we obtain that this $G$-bundle on $\bS$ is well-defined, i.e., is
independent of $(\bd_v,\bd_h)$. We set $\sigma'_\sO$ to be the
corresponding point of $\Bun_G(\bS,\bD_\infty)$.

It remains to show that the integers $b_{i,n}$ are divisible by
$\phi_{\BZ}(1)$. For that we choose a pair of directions
$(\bd_v,\bd_h)$, so that the points $\pi_v(\bs_i)$ are pairwise
distinct. Then the corresponding $b_i$'s are reconstructed
as the defect of the quasi-map $\sigma$.

This completes the proof of the theorem.

\ssec{}

Note that in the course of the proof of \thmref{stratification} we
have established that the morphisms $\ol{\iota}:\fU^{a-b}_G\times 
\on{Sym}(\overset{\circ}\bS)\to \fU^a_G$ are compatible in a natural
sense with the factorization isomorphisms of \propref{factorization of
Uhlenbeck}. This is so because the corresponding property is true
for $SL_n$.

As a corollary of \thmref{stratification}, we obtain, in
particular, the following:

\begin{cor}
The variety $\Bun_G^a(\bS,\bD_\infty)$ is quasi-affine, and its
affine closure is isomorphic to the normalization of $\fU^a_G$.
\end{cor}

\begin{proof}

Indeed, we know that $\dim{\Bun_G^a(\bS,\bD_\infty)}=2\cdot \check{h}\cdot a$
and from \thmref{stratification}, we obtain that the complement
to $\Bun_G^a(\bS,\bD_\infty)$ inside the normalization of $\fU^a_G$
is of dimension $2\cdot \check{h}\cdot (a-1)+2$. I.e., it is of
codimension $2(\check{h}-1)\geq 2$.

\end{proof}

We do not know whether $\fU^a_G$ is in general a normal variety, but
this is true for $G=SL_n$ since $\fU^a_n\simeq \fN^a_n$, cf.
\cite{CB}.

\ssec{}

For a fixed pair of directions $(\bd_v,\bd_h)$ and a point $0_\bX\in \bX$,
consider the $\BG_m$-action on $\bS$, which acts as identity on $\bC$
and by dilations on $\bX$, fixing $0_\bX$.
By transport of structure, we obtain a $\BG_m$-action on $\fU^a_G$.

\begin{prop}  \label{contraction of Uhlenbeck}
The above $\BG_m$-action on $\fU^a_G$ contracts this space to
$$(\overset{\circ}\bC\times 0_\bX)^{(a)}\subset
\overset{\circ}\bS{}^{(a)}\simeq \fU_G^{a;a}\hookrightarrow \fU^a_G.$$
The contraction map $\fU^a_G\to \overset{\circ}\bC{}^{(a)}$ coincides
with $\varpi^a_h$.
\end{prop}

\begin{proof}

We will use the description of $\fU^a_G$ via \corref{fewer lines},
i.e. we realize it as a closed subset in
$$\on{QMaps}(\bC,\CG_{G,\bX})\times
\on{Sect}(\bD_\infty-\bd_h,\overset{\circ}\bC{}_{\bD_\infty-\bd_h}^{(a)}).$$

Note that the family of curves $\bC_{\bD_\infty-\bd_h}$ can be naturally
identified with the constant family $(\bD_\infty-\bd_h)\times \bC$ via
the projection $\bC\times 0_\bX\overset{\pi_v}\to \bC_{\bd_v}$ for
each $\bd_v\in \bD_\infty-\bd_h$. In terms of this trivialization,
the above $\BG_m$-action contracts the space
$\on{Sect}(\bD_\infty-\bd_h,\overset{\circ}\bC{}_{\bD_\infty-\bd_h}^{(a)})$
to the space of constant sections, which can be identified with
$\overset{\circ}\bC{}^{(a)}$.

In addition, from \propref{contraction of quasi-maps} we obtain that
the above $\BG_m$-action
contracts the space $\on{QMaps}^a(\bC,\CG_{G,\bX})$ to
$\overset{\circ}\bC{}^{(a)}\subset \on{QMaps}^a(\bC,\CG_{G,\bX})$.

By comparing it with our description of the map
$i_b:\fU^{a;b}_G\to \fU^a_G$ for $b=a$ we obtain the assertion
of the proposition.

\end{proof}

\ssec{An example}

Let us describe explicitly the Uhlenbeck space $\fU^a_G$ for
$a=1$. (In particular, using \propref{factorization of Uhlenbeck},
this will imply the description of the singularities of 
$\fU^a_G$ for general $a$ along the strata $\fU^{a;b}_G$ for $b=1$.)

\medskip

Choose a pair of directions $(\bd_v,\bd_h)$, and consider the
corresponding surface $\bS'=\bC\times \bX$. Let us choose
points $0_\bC\in \bC-\infty_\bC$, $0_\bX\in \bX-\infty_\bX$,
and consider the affine Grassmannian $\Gr_G:=\Gr_{G,\bX,0_\bX}$.
Let $\ol{\Gr}^{\ol{\alpha}_0}_G\subset \Gr_G$ be the corresponding closed
subscheme. It contains the unit point $\one_{\Gr}\in \Gr_G$, which is
its only singularity.

By definition, we have a map $\ol{\Gr}^{\ol{\alpha}_0}_G\to \Bun_G(\bX)$,
and, as we know, the stack $\Bun_G(\bX)$ carries a line bubdle 
$\CP_{\Bun_G(\bX)}$ with a section, whose set of zeroes is the locus
of non-trivial bundles. Let us denote by $D_{\ol{\Gr}^{\ol{\alpha}_0}_G}$
the corresponding Cartier divisor in $\ol{\Gr}^{\ol{\alpha}_0}_G$. Note
that the point 
$\one_{\Gr}$ belongs to $\ol{\Gr}^{\ol{\alpha}_0}_G-D_{\ol{\Gr}^{\ol{\alpha}_0}_G}$.

We claim that there is a natural isomorphism 
$$\fU^1_G\simeq \overset{\circ}\bC\times
\overset{\circ}\bX\times
(\ol{\Gr}^{\ol{\alpha}_0}_G-D_{\ol{\Gr}^{\ol{\alpha}_0}_G}).$$

We have the projection
$$\varpi^1_h\times \varpi^1_v:\fU^1_G\to \overset{\circ}\bC\times
\overset{\circ}\bX,$$
and the group $\BA^2$ acts on $\fU^1_G$ via its action on
$\overset{\circ}\bS$ by shifts. Therefore, it is enough to show that
\begin{equation} \label{simple case}
(\varpi^1_h)^{-1}(0_\bC)\cap (\varpi^1_v)^{-1}(0_\bX)\simeq 
\ol{\Gr}^{\ol{\alpha}_0}_G-D_{\ol{\Gr}^{\ol{\alpha}_0}_G}.
\end{equation}
Note that the intersection 
$(\varpi^1_h)^{-1}(0_\bC)\cap (\varpi^1_v)^{-1}(0_\bX)$ contains only 
one geometric point, which is not in $\Bun_{G}(\bS',\bD'_\infty)$, namely,
the point $0_\bC\times 0_\bX$ from the stratum $\fU^{1;1}_G$.

Let $\one_\bC\in \bC$ be a point different from $0_\bC$ and $\infty_\bC$.
Given a point in $(\varpi^1_h)^{-1}(0_\bC)\cap (\varpi^1_v)^{-1}(0_\bX)$,
we have a quasi-map $\sigma:\bC\to \ol{\Gr}^{\ol{\alpha}_0}_G$
(cf. \secref{upper estimate}), such that $\sigma$, restricted to
$\bC-0_\bC$, is a {\it map}, whose image belongs to 
$\ol{\Gr}^{\ol{\alpha}_0}_G-D_{\ol{\Gr}^{\ol{\alpha}_0}_G}$. Hence,
$\sigma\mapsto \sigma(\one_\bC)$ defines a map in one direction
$(\varpi^1_h)^{-1}(0_\bC)\cap (\varpi^1_v)^{-1}(0_\bX)\to 
(\ol{\Gr}^{\ol{\alpha}_0}_G-D_{\ol{\Gr}^{\ol{\alpha}_0}_G})$.

Let us show first that this map is one-to-one. From the proof of 
\thmref{stratification}, it is clear that it sends the unique point in 
$(\varpi^1_h)^{-1}(0_\bC)\cap (\varpi^1_v)^{-1}(0_\bX)\cap \fU^{1;1}_G$
to the point $\one_{\Gr}\in \ol{\Gr}^{\ol{\alpha}_0}_G-D_{\ol{\Gr}^{\ol{\alpha}_0}_G}$.

Now, let $\sigma_1,\sigma_2$ be two elements of 
$(\varpi^1_h)^{-1}(0_\bC)\cap (\varpi^1_v)^{-1}(0_\bX)$. 
Since both $\sigma_1$ and $\sigma_2$ are of degree $1$, if we had
$\one_{\Gr}=\sigma_1(\infty_\bC)=\sigma_2(\infty_\bC)$, and
$\sigma_1(\one_\bC)=\sigma_2(\one_\bC)$, this would imply
$\sigma_1=\sigma_2$. For the same reason, no {\it map} $\sigma$ can
send $\one_\bC$ to $\one_{\Gr}$.

Hence, the map $(\varpi^1_h)^{-1}(0_\bC)\cap
(\varpi^1_v)^{-1}(0_\bX)\to 
(\ol{\Gr}^{\ol{\alpha}_0}_G-D_{\ol{\Gr}^{\ol{\alpha}_0}_G})$
is one-to-one on the level of geometric points. Since the left-hand
side is reduced, and the right-hand side is normal, we obtain that it
is an open embedding. Recall now the action of $\BG_m$ by dilations
along the $\bX$-factor. We claim that it contracts 
$(\ol{\Gr}^{\ol{\alpha}_0}_G-D_{\ol{\Gr}^{\ol{\alpha}_0}_G})$ to the point
$\one_{\Gr}$. Indeed, we have a closed emedding 
$\Gr_G\hookrightarrow \CG_{G,\bX}$, such that $(\ol{\Gr}^{\ol{\alpha}_0}_G-D_{\ol{\Gr}^{\ol{\alpha}_0}_G})$
gets mapped to the
corresponding open cell $\CG_{\fg_{aff},\fg^+_{aff},e}$, cf. \secref{flag schemes},
and for Schubert cells the contraction assertion is evident.

Using \propref{contraction of Uhlenbeck}, this implies that 
$(\varpi^1_h)^{-1}(0_\bC)\cap (\varpi^1_v)^{-1}(0_\bX)\to 
(\ol{\Gr}^{\ol{\alpha}_0}_G-D_{\ol{\Gr}^{\ol{\alpha}_0}_G})$ is in fact 
an isomorphism.

\medskip

\noindent{\it Remark.}
It is well-known that
$(\ol{\Gr}^{\ol{\alpha}_0}_G-D_{\ol{\Gr}^{\ol{\alpha}_0}_G})$ is isomorphic
to the minimal nilpotent orbit closure $\ol{\mathbb O}_{\on{min}}$ in $\fg$.
It was noticed by V.~Drinfeld in 1998 that the transversal slice to the
singularity of $\fU^1_G$ should look like $\ol{\mathbb O}_{\on{min}}$.

\ssec{}

For a fixed $b$, $0\leq b\leq a$, let $\fP(b)$ be its partition, i.e.
$b=\underset{k} \Sigma\, n_k\cdot d_k$, with $d_k>0$ and pairwise
distinct. Let us denote by $\on{Sym}^{\fP(b)}(\overset{\circ}\bS)$
the variety $\underset{k}\Pi\, \on{Sym}^{n_k}(\overset{\circ}\bS)$,
and by $\overset{\circ}{\on{Sym}}{}^{\fP(b)}(\overset{\circ}\bS)$ the
open subset of $\on{Sym}^{\fP(b)}(\overset{\circ}\bS)$ obtained by
removing all the diagonals.

The symmetric power $\on{Sym}^b(\overset{\circ}\bS)$ is the union
$\underset{\fP(b)}\cup \overset{\circ}{\on{Sym}}{}^{\fP(b)}(\overset{\circ}\bS)$
over all the partitions of $b$.
Let $\fU^{a;\fP(b)}_G=\Bun^{a-b}_G(\bS,\bD_\infty)\times
\overset{\circ}{\on{Sym}}{}^{\fP(b)}(\overset{\circ}\bS)$
be the corresponding subscheme in $\fU^{a;\fP(b)}_G$.

\medskip

The next theorem, which can be regarded as the main result of this paper,
describes the restriction of the intersection cohomology sheaf of $\fU^a_G$
to the strata $\fU^{a;\fP}_G$.

\medskip

Let $\check\fg_{aff}$ denote the Langlands dual Lie algebra to
$\fg_{aff}$. Consider the (maximal) parabolic $\fg[[t]]\oplus K\cdot
\BC\oplus d\cdot \BC\subset\fg_{aff}$;
then we obtain the corresponding parabolic
inside $\check\fg_{aff}$, whose unipotent radical $\fV$ is an
(integrable) module over the corresponding dual Levi
subalgebra, i.e. $\fV$ is a representation of the group
$\check G_{aff}=\check G\times \BG_m$. In particular,
we can consider the space $\fV^f$, where $f\in \check \fg$ is
a principal nilpotent element. We regard $\fV^f$ as a bi-graded
vector space, $\fV^f=\underset{m,l}\oplus\, (\fV^f)_l^m$:
where the ``first'' grading, corresponding to the upper index $m$,
is the principal grading coming from the Jacobson-Morozov triple
containing $f$, and the ``second'' grading, corresponding to the
lower index $l$, comes from the action of $\BG_m\subset G_{aff}$.

For a partition $\fP(b)$ as above, we will denote
by $\on{Sym}^{\fP(b)}(\fV^f)$ the graded vector space
$$\underset{k}\bigotimes\, \Bigl(\underset{i\geq 0}\oplus\,
\bigl(\on{Sym}^i(\fV^f)\bigr)_{d_k}[2i]\Bigr)^{\otimes n_k},$$
where the subscript $d_k$ means that we are taking the graded subspace
of the corresponding index with respect to the ``second'' grading,
and the notation $[j]$ means the shift with respect to the ``first''
grading. By declaring the ``first'' grading to be cohomological we
obtain a semi-simple complex of vector spaces.

\medskip

\begin{thm}\label{ic uhlenbeck}
The restriction $\on{IC}_{\fU^a_G}|_{\fU^{a;\fP(b)}_G}$ is locally constant
and is isomorphic to the IC sheaf on this scheme tensored by the (constant)
complex $\on{Sym}^{\fP(b)}(\fV^f)$.
\end{thm}

\section{Functorial definitions}

\ssec{}

An obvious disadvantage of our definition of $\fU^a_G$ is that
we were not able to say what functor it represents. In this
section, we will give two more variants of the definition of the Uhlenbeck
space, we will call them $^{Tann}\fU^a_G$ and $^{Drinf}\fU^a_G$,
respectively, and which will be defined as functors. Moreover,
$^{Tann}\fU^a_G$ will be a scheme, whereas $^{Drinf}\fU^a_G$ only
an ind-scheme.

We will have a sequence of
closed embeddings
$$\fU^a_G\to {}^{Tann}\fU^a_G\to {}^{Drinf}\fU^a_G,$$
which induce isomorphisms on the level of the corresponding reduced
schemes. (In other words, up to nilpotents, one can give a functorial
definition of the Uhlenbeck space.)

\ssec{A Tannakian definition}

In this subsection, we will express the Uhlenbeck space for an
arbitrary $G$ in terms of that of $SL_n$. 

For two integers $n_1$ and $n_2$ consider the homomorphisms
$$SL_{n_1}\times SL_{n_2}\to SL_{n_1+n_2} \text{ and }
SL_{n_1}\times SL_{n_2}\to SL_{n_1\cdot n_2},$$
and the corresponding morphisms
$$\Bun_{n_1}(\bX)\times \Bun_{n_2}(\bX)\to
\Bun_{n_1+n_2}(\bX) \text{ and }
\Bun_{n_1}(\bX)\times \Bun_{n_2}(\bX)\to
\Bun_{n_1\cdot n_2}(\bX).$$

Note that the pull-backs of the line bundles
$\CP_{\Bun_{n_1+n_2}(\bX)}$ and $\CP_{\Bun_{n_1\cdot n_2}(\bX)}$
is isomorphic to $\CP_{\Bun_{n_1}(\bX)}\boxtimes
\CP_{\Bun_{n_2}(\bX)}$ and 
$\CP_{\Bun_{n_1}(\bX)}^{\otimes n_2}\boxtimes 
\CP_{\Bun_{n_2}(\bX)}^{\otimes n_1}$, respectively.

As in \lemref{functoriality}, we obtain the closed embeddings
$\fU^{a_1}_{n_1}\times \fU^{a_2}_{n_2}\to
\fU^{a_1+a_2}_{n_1+n_2}$ and $\fU^{a_1}_{n_1}\times \fU^{a_2}_{n_2}\to
\fU^{a_1\cdot n_2+a_2\cdot n_1}_{n_1\cdot n_2}$, which extend the 
natural morphisms between the corresponding moduli spaces
$\Bun_n^a(\bS,\bD_\infty)$.

\medskip

Moreover, it is easy to deduce from the proof of \thmref{stratification} that
if $\sigma_{\sO,1}$ (resp., $\sigma_{\sO,2}$) is a collection of based
quasi-maps $\bC\to \Gr^{BD,a}_{SL_{n_1},\bX}$ describing a point of 
$\fU^{a_1}_{n_1}$, whose saturation is $\sigma'_{\sO,1}$ and defect
$D_{1,\bS}$ (resp., $\sigma'_{\sO,2}$, $D_{2,\bS}$), then the
corresponding point of $\fU^{a_1+a_2}_{n_1+n_2}$ (resp., 
$\fU^{a_1\cdot n_2+a_2\cdot n_1}_{n_1\cdot n_2}$) will have the
saturation $\sigma'_{\sO,1}\oplus\sigma'_{\sO,2}$ (resp., 
$\sigma'_{\sO,1}\otimes\sigma'_{\sO,2}$) and defect
$D_{1,\bS}+D_{2,\bS}$ (resp., $n_2\cdot D_{1,\bS}+n_1\cdot D_{2,\bS}$).

\ssec{}

We set $^{Tann}\fU^a_G$ to represent the functor that assigns to a test
scheme $S$ the following data:

\begin{itemize}

\item An $S$-point of $\on{Sect}\left(\sO,\on{QMaps}^a(\bC,\Gr^{BD,a}_{G,\bX})_\sO\right)$,

\item For every representation $\phi:G\to SL_n$, a point of
$\fU^{\phi_{\BZ}(a)}_n$,

\end{itemize}
such that
\begin{itemize}

\item
The corresponding $S$-points of
$\on{Sect}\left(\sO,\on{QMaps}^{\phi_{\BZ}(a)}
(\bC,\Gr^{BD,\phi_{\BZ}(a)}_{SL_n,\bX})_\sO\right)$, coming from
$$\on{Sect}\left(\sO,\on{QMaps}^a(\bC,\Gr^{BD,a}_{G,\bX})_\sO\right)\to
\on{Sect}\left(\sO,\on{QMaps}^{\phi_{\BZ}(a)}
(\bC,\Gr^{BD,\phi_{\BZ}(a)}_{SL_n,\bX})_\sO\right)\leftarrow
\fU^{\phi_{\BZ}(a)}_n$$
coincide,

\item
For $\phi=\phi_1\otimes\phi_2$, $n=n_1\cdot n_2$ (resp.,
$\phi=\phi_1\oplus\phi_2$, $n=n_1+n_2$), the corresponding $S$-point
of $\fU^{\phi_{\BZ}(a)}_n$ equals the image of the corresponding point
of $\fU^{\phi_1{}_{\BZ}(a)}_{n_1}\times \fU^{\phi_2{}_{\BZ}(a)}_{n_2}$
under $$\fU^{\phi_1{}_{\BZ}(a)}_{n_1}\times
\fU^{\phi_2{}_{\BZ}(a)}_{n_2}\to\fU^{\phi_{\BZ}(a)}_n.$$

\end{itemize}

\medskip

Since for a faithful representation $\phi:G\to SL_n$, the map
$\on{Sect}\left(\sO,\on{QMaps}^a(\bC,\Gr^{BD,a}_{G,\bX})_\sO\right)\to
\on{Sect}\left(\sO,\on{QMaps}^{\phi_{\BZ}(a)}
(\bC,\Gr^{BD,\phi_{\BZ}(a)}_{SL_n,\bX})_\sO\right)$
is a closed embedding, it is easy to see that $^{Tann}\fU^a_G$ is in fact
a closed subfunctor in $\fU^{\phi_{\BZ}(a)}_n$; in particular,
$^{Tann}\fU^a_G$ is indeed representable by an affine scheme of finite type.
Moreover, by construction, we have a closed embedding
$\fU^a_G\to {}^{Tann}\fU^a_G$.

\begin{prop}  \label{Tannaka=Uhlenbeck}
The above map $\fU^a_G\to {}^{Tann}\fU^a_G$ induces an isomorphism
on the level of reduced schemes.
\end{prop}

We do not know whether it is in general true that $^{Tann}\fU^a_G$ is
actually isomorphic to $\fU^a_G$.

\begin{proof}

Since the map in question is a closed embedding, we only have
to check that it defines a surjection on the level of geometric
points.

Thus, let $\sigma_\sO:\bC\to \Gr^{BD,a}_{G,\bX}$ be the collection
of quasi-maps representing a point of 
$\on{Sect}\left(\sO,\on{QMaps}^a(\bC,\Gr^{BD,a}_{G,\bX})_\sO\right)$,
and $\sigma_{\sO,\phi}$ be a compatible system of points of
$\fU^{\phi_\BZ(a)}_n$ for every representation $\phi:G\to SL_n$.

Using \thmref{stratification}, we have to show that there exists a
point $\F'_G\in \Bun_G(\bS,\bD_\infty)$
(corresponding to a system of maps $\sigma'_\sO$) and a $0$-cycle
$D_\bS$, such that for every $(\bd_v,\bd_h)\in \sO$, the corresponding
quasi-map $\sigma$ has as its saturation the map of \eqref{eq3} and
$\pi_v(D_\bS)$ as its defect.

\medskip

Let $\F'_{\phi}$ be the $SL_n$ bundle corresponding to the saturation
of $\sigma_{\sO,\phi}$. From our description of the maps
$\fU^{a_1}_{n_1}\times \fU^{a_2}_{n_2}\to
\fU^{a_1+a_2}_{n_1+n_2}$ and $\fU^{a_1}_{n_1}\times \fU^{a_2}_{n_2}\to
\fU^{a_1\cdot n_2+a_2\cdot n_1}_{n_1\cdot n_2}$ it follows that
$\F'_{\phi_1\oplus\phi_2}\simeq \F'_{\phi_1}\oplus \F'_{\phi_2}$ and
$\F'_{\phi_1\otimes\phi_2}\simeq \F'_{\phi_1}\otimes \F'_{\phi_2}$.

Therefore, the collection $\{\F'_{\phi}\}$ defines a $G$-bundle
on $\bS$, trivialized along $\bD_\infty$, which we set to be our
$\F'_G$. The $0$-cycle $D_\bS$ is reconstructed using just one
faithful representation of $G$, as in the proof of \thmref{stratification}.

\end{proof}

\ssec{Drinfeld's approach}
Consider the following functor
$\{$ Schemes of finite type $\}\to$ \newline
$\{$ Groupoids $\}$; we will call it $^{Drinf}\fU^a_G$, since the
definition below follows a suggestion of Drinfeld. The category
corresponding to a test scheme $S$ has as objects the data of:

\begin{itemize}

\item
A principal $G$-bundle $\F_G$ defined on an open subscheme
of $U\subset \bS\times S$, such that $\bS\times S-U$ is finite over
$S$, and $U\supset \bD_\infty\times S$, and $\F_G$ is equipped with
a trivialization along $\bD_\infty\times S$. (We will denote by 
$U_\sO$ the corresponding open subset in $\bS'_\sO\times S$, and
by $\F_{G,\sO}$ the corresponding $G$-bundle over it.)

\item
A map
$S\to \on{Sect}(\sO,\on{QMaps}^a(\bC,\Gr^{BD,a}_{G,\bX})_{\sO})^\tau$.
(This expression makes use of the superscript $\tau$, which was introduced
in \lemref{transposition invariance}.)
We will denote by $\sigma_\sO$ the corresponding relative quasi-map
$\bC_\sO\times S\to \Gr^{BD,a}_{G,\bX_\sO}$. Let
$V_\sO\subset \bC_\sO\times S$ denote the open subset over which $\sigma_\sO$
is defined as a map, and $D_{v,\sO}\subset \bX_\sO\times S$ the
corresponding relative divisor. Let $U^{\sigma_\sO}\subset
\bS'_\sO\times S$ be the open subset $V_\sO\underset{\sO}\times
\bX_\sO\cup \bC_\sO\underset{\sO}\times(\bX_\sO-D_{v,\sO})$, and
$\F^{\sigma_\sO}_G$ be the corresponding principal $G$-bundle defined
on $U^{\sigma_\sO}$.

\item
An isomorphism $\alpha$ between $\F_{G,\sO}|_{U^{\sigma_\sO}\cap U_\sO}$
and $\F^{\sigma_\sO}_G|_{U^{\sigma_\sO}\cap U_\sO}$,
respecting the trivializations along $\bD'_{\infty,\sO}\times S$.

\end{itemize}

Morphisms between an object $(\F_G,U,\sigma_\sO,\alpha)$ and
another object $(\F^1_G,U^1,\sigma^1_\sO,\alpha^1)$ are $G$-bundles isomorphisms
$\F_G\simeq \F^1_G|_{U^1_\sO\cap U_\sO}$, which commute with other pieces of
data (in particular, the data of $\sigma_\sO$ must be the same).
In particular, we see that objects in our category have no non-trivial
automorphisms, i.e. $^{Drinf}\fU^a_G$ is a functor
$\{$ Schemes of finite type $\}\to\{$ Sets $\}$.

\begin{lem}  \label{indrepresentability}
The functor $^{Drinf}\fU^a_G$ is representable by an affine ind-scheme of
ind-finite type.
\end{lem}

(Later we will show that $^{Drinf}\fU^a_G\to
\on{Sect}(\sO,\on{QMaps}^a(\bC,\Gr^{BD,a}_{G,\bX})_{\sO})$ is in fact
a closed embedding.)

\begin{proof}

It is enough to show that the forgetful map
$^{Drinf}\fU^a_G\to \on{Sect}(\sO,\on{QMaps}^a(\bC,\Gr^{BD,a}_{G,\bX})_{\sO})$
is ind-representable (and ind-affine).
Thus, for a test-scheme $S$ let $\sigma_\sO$ be a relative quasi-map
$\bC_{\sO}\times S\to\Gr^{BD,a}_{G,\bX_\sO}$, and let $U^{\sigma_\sO}$,
$\F^{\sigma_\sO}_G$ be as above.

Let $(\bd_v,\bd_h)\in \sO$ be some fixed pair of directions,
and set $U\subset \bS'\times S$ (resp., $\F_G$) to be the fiber of 
$U^{\sigma_\sO}$ (resp., $\F^{\sigma_\sO}_G$) over the
corresponding point of $\sO$. We will denote by the same character
$U$ (resp., $\F_G$) the corresponding open subset in $\bS\times S$
(resp., the $G$-bundle over it).

Consider the group ind-scheme over $S$, call it $\on{Aut}_U(\F_G)$,
of automorphisms of $\F_G$, respecting the trivialization at the divisor
of infinity. Consider now the group ind-scheme of $S$-maps 
$S\times \sO\to \on{Aut}_U(\F_G)$. Then it is easy to see that the fiber product
$$S\underset{\on{Sect}(\sO,\on{QMaps}^a(\bC,\Gr^a_G)_{\sO})}\times{}^{Drinf}\fU^a_G$$
is (ind)-represented by the above ind-scheme $\on{Maps}_S(S\times \sO,\on{Aut}_U(\F_G))$.

\end{proof}

The rest of this section will be devoted to the proof of the following result:

\begin{thm}  \label{Tannaka and Drinfeld}
We have a closed embedding $^{Tann}\fU^a_G\to {}^{Drinf}\fU^a_G$,
which induces an isomorphism on the level of reduced schemes.
\end{thm}

\ssec{}

Our first task is to construct the map
$^{Tann}\fU^a_G\to {}^{Drinf}\fU^a_G$. The first case
to consider is $G=SL_n$. In this case, by the definition of
$^{Tann}\fU^a_G$, we have $\fU^a_n\simeq ^{Tann}\fU^a_n$.

Since $^{Drinf}\fU^a_n$ is ind-affine, and $\fU^a_n\simeq \fN^a_n$ is the
affinization of the scheme $\wt{\fN}{}^a_n$, to construct a map
$\fU^a_n\to {}^{Drinf}\fU^a_n$, it suffices to construct a map
$\wt{\fN}{}^a_n\to {}^{Drinf}\fU^a_G$. 

Given an $S$-point $\M$ of $\wt{\fN}{}^a_n$, i.e. an
$S$-family of torsion-free sheaves on $\bS$, there exists an open
subset $U\subset \bS\times S$ whose complement is finite over $S$,
such that $\M$ is a vector bundle when restricted to $U$, i.e.
we obtain a principal $SL_n$-bundle $\F_{SL_n}$ on $U$. 

The data of a relative quasi-map $\sigma_\sO$ are obtained from the
composition $$\wt{\fN}{}^a_n\to \fU^a_n\to
\on{Sect}(\sO,\on{QMaps}^a(\bC,\Gr^a_G)_{\sO}).$$
The isomorphisms
$\alpha:\F_{SL_n,\sO}|_{U^{\sigma_\sO}\cap U_\sO}\simeq
\F^{\sigma_\sO}_{SL_n}|_{U^{\sigma_\sO}\cap U_\sO}$ 
follow from the construction of the map \eqref{from Nakajima to quasi-maps}.

\medskip

Thus, we have a morphism $\fU^a_n\to {}^{Drinf}\fU^a_n$, and
moreover, for $n=n_1+n_2$ (resp., $n=n_1\cdot n_2$) and
$a=a_1+a_2$ (resp., $a=n_1\cdot a_2+n_2\cdot a_1$), the natural
morphism $$^{Drinf}\fU^{a_1}_{n_1}\times{}^{Drinf}\fU^{a_2}_{n_2}\to{}^{Drinf}\fU^a_n$$
is compatible with the corresponding morphism
$\fU^{a_1}_{n_1}\times\fU^{a_2}_{n_2}\to\fU^a_n$. Hence, 
we have the morphism $^{Tann}\fU^a_G\to {}^{Drinf}\fU^a_G$ for an
arbitrary $G$.

\medskip

Next, we will show that the composition
$\fU^a_G\to {}^{Tann}\fU^a_G\to {}^{Drinf}\fU^a_G$ is a closed
embedding. Let $\ol{\fU}^a_G$ be the closure of $\Bun_G(\bS,\bD_\infty)$ in
$^{Drinf}\fU^a_G$.

Then under the projection $^{Drinf}\fU^a_G\to
\on{Sect}(\sO,\on{QMaps}^a(\bC,\Gr^{BD,a}_{G,\bX})_{\sO})$, $\ol{\fU}^a_G$ gets
mapped to $\fU^a_G$, and we obtain a pair of arrows
$$\fU^a_G\leftrightarrows \ol{\fU}^a_G,$$ such
that both compositions induce the identity map on
$\Bun_G(\bS,\bD_\infty)$. Since the latter is dense, we obtain that
the map $\fU^a_G\to \ol{\fU}^a_G$ is in fact an isomorphism.

\medskip

Thus, $\fU^a_G\to {}^{Drinf}\fU^a_G$ is a closed embedding, and let us
show now that it induces an isomorphism from $\fU^a_G$ to the reduced
(ind)-scheme underlying $^{Drinf}\fU^a_G$. For that we have to check
that this map is surjective on the level of geometric points.

\bigskip

Let $(\F_G,\sigma_\sO,\alpha)$ is a geometric point of
$^{Drinf}\fU^a_G$, and let $D_{v,\sO}$ (resp., $D_{h,\sO}$) be the 
relative Cartier divisors in $\bX_\sO$ (resp., $\bC_\sO$) obtained from 
$\sigma_\sO$ via $\varpi^a_{v,\sO}$ (resp., $\varpi^a_{h,\sO}$).
Let us denote by $\F'_G$ the (unique) extension of $\F_G$
to the entire $\bS$. Let $\sigma'_\sO$ be the map
$\bC_\sO\times S\to \Gr^{BD,a-b}_{G,\bX_\sO}$ corresponding to $\F'_G$.
Let $D'_{v,\sO}$ (resp., $D'_{h,\sO}$) be the relative divisors in
$\bX_\sO$ (resp., $\bC_\sO$) corresponding to $\sigma'_\sO$.
 
Set $\bS-U=\underset{i}\cup\, \bs_i$. Our task is to show
that there exist integers $b_i$ such that for every pair
of directions $(\bd_v,\bd_h)\in \sO$, we have:

\smallskip

\noindent (1) $D_v=D'_v+b_i\cdot \pi_h(\bs_i)$, $D_h=D'_h+b_i\cdot \pi_v(\bs_i)$

\smallskip

\noindent (2) The quasi-map $\sigma$ has the saturation equal to
$\bC\overset{\sigma'}\longrightarrow \Gr^{BD,a-b}_{G,\bX,D'_v}\to
\Gr^{BD,a}_{G,\bX,D_v}$,

\smallskip

\noindent (3) The defect of $\sigma$ equals 
$\underset{i}\Sigma\, b_i\cdot \pi_v(\bs_i)$.

\medskip

Let $\sigma'':\bC\to \Gr^{BD,a}_{G,\bX,D_v}$ be
a based map of degree $a-b$, which is the saturation of $\sigma$. We know that
$\sigma$ is obtained from $\sigma''$ by adding to it a divisor
of degree $b$ supported on $\underset{i}\cup\,\pi_v(\bs_i)$. Hence, by
letting $(\bd_v,\bd_h)$ move along $\sO$, we obtain a section $\sO\to
\overset{\circ}\bX{}^{(b)}_\sO$ with the above support property. However,
it is easy to see that any such section is of the form
$$(\bd_v,\bd_h)\mapsto \underset{i}\Sigma\,b_i\cdot\pi_v(\bs_i)$$
for some integers $b_i$.

Now, since the maps $\sigma'$ and $\sigma''$ give rise to the
same bundle on $\bS'$, the two compositions
$$\bC\overset{\sigma'}\longrightarrow \Gr^{BD,a-b}_{G,\bX,D'_v}\to \CG_{G,\bX}
\text{ and } \bC\overset{\sigma''}\longrightarrow \Gr^{BD,a}_{G,\bX,D_v}\to \CG_{G,\bX}$$
coincide. Hence, $\sigma'$ and $\sigma''$ have the same degree $a-b$
and $\varpi^{a-b}_h(\sigma')=\varpi^{a-b}_h(\sigma'')$. Since, by the
property of saturations, cf. \eqref{configuration of saturation}, we have
$\varpi^a_h(\sigma)=\varpi^{a-b}_v(\sigma'')+\underset{i}\Sigma\,b_i\cdot\pi_v(\bs_i)$,
we obtain $D_h=D'_h+b_i\cdot \pi_v(\bs_i)$. By combining this with the
condition that $\sigma_\sO\in
\on{Sect}(\sO,\on{QMaps}^a(\bC,\Gr^{BD,a}_{G,\bX})_{\sO})^\tau$,
cf. \lemref{transposition invariance}, we obtain also that 
$D'_v=D_v+ \underset{i}\Sigma\, b_i\cdot \pi_h(\bs_i)$.
Moreover, we obtain that $\sigma''$ is the composition of $\sigma'$
and the embedding $\Gr^{BD,a-b}_{G,\bX,D'_v}\to\Gr^{BD,a}_{G,\bX,D_v}$, which is what
we had to show. 

\bigskip

To complete the proof of the theorem, it would be enough to show that
the projection $^{Drinf}\fU^a_G\to
\on{Sect}(\sO,\on{QMaps}^a(\bC,\Gr^{BD,a}_{G,\bX})_{\sO})$ is a closed embedding.
We know this for the closed subscheme $\fU^a_G$, and hence, we obtain
a priori that the map in question is ind-finite. Therefore, it would
suffice to show that any tangent vector to $^{Drinf}\fU^a_G$, which
is vertical with respect to the map
$^{Drinf}\fU^a_G\to
\on{Sect}(\sO,\on{QMaps}^a(\bC,\Gr^{BD,a}_{G,\bX})_{\sO})$, is in fact zero.

\medskip

Let $\on{Spec}(\BC[\epsilon]/\epsilon^2)\to {}^{Drinf}\fU^a_G$
be such a tangent vector. Then the relative quasi-map
$\sigma_\sO:\bC_\sO\times \on{Spec}(\BC[\epsilon]/\epsilon^2)\to
\Gr^{BD,a}_{G,\bX_\sO}$ is actually independent of $\epsilon$, and as we
have seen in the proof of \lemref{indrepresentability},
we can assume that the generically defined $G$-bundle $\F_G$ is
also independent of $\epsilon$. Thus, we are dealing with an
infinitesimal automorphism of $\F_G$, and we must show that it is zero.

However, $\F_G$ extends canonically to a $G$-bundle on the entire
$\bS$, and so does our infinitesimal automorphism. Since this
automorphism preserves the trivialization along $\bD_\infty$,
our assertion reduces to the fact that points of
$\Bun_G(\bS,\bD_\infty)$ have no automorphisms.

\bigskip

\centerline{{\bf Part III}:  {\Large Parabolic versions of Uhlenbeck spaces}}

\bigskip

Throughout this part, $\fg$ will be a simple finite-dimensional Lie algebra with 
the corresponding simply-connected group denoted $G$. By $\fg_{aff}$
will denote the corresponding untwisted affine Kac-Moody algebra, 
i.e. $\fg_{aff}=\fg((x))\oplus K\cdot\BC\oplus d\cdot \BC$, as in
\secref{parabolic notation}.

\section{Parabolic Uhlenbeck spaces}

\ssec{}

Let $\BP^1\simeq\bD_0\subset \bS\simeq \BP^2$ be a divisor
different from $\bD_\infty$; we will denote by $\overset{\circ}{\bD}_0$
the intersection $\bD_0\cap \overset{\circ}\bS$.

Let us denote by $\Bun_{G;P}(\bS,\bD_\infty;\bD_0)$ the scheme
classifying the data consisting of a $G$-bundle $\F_G$ on $\bS$,
a trivialization of $\F_G|_{\bD_\infty}$, and a reduction to $P$ of
$\F_G|_{\bD_0}$, compatible with the above trivialization at the point
$\bD_\infty\cap \bD_0$. In other words, if we denote by
$\Bun_{P}(\bD_0;\bD_\infty\cap\bD_0)$ the moduli stack classifying
$P$-bundles on the curve $\bD_0$ with a trivialization at
$\bD_\infty\cap \bD_0$,
we have an identification
\begin{equation}  \label{Bun_P}
\Bun_{G;P}(\bS,\bD_\infty;\bD_0)\simeq \Bun_G(\bS,\bD_\infty)
\underset{\Bun_{G}(\bD_0;\bD_\infty\cap\bD_0)}\times
\Bun_{P}(\bD_0;\bD_\infty\cap\bD_0).
\end{equation}
 From this description, we see that $\Bun_{G;P}(\bS,\bD_\infty;\bD_0)$
splits into connected components indexed by
$\theta\in \widehat\Lambda_{\fg,\fp}\simeq \Lambda_{\fg,\fp}\oplus
\delta\cdot\BZ$, where for $\theta=(\ol{\theta},a)$, the index $a$
enumerates the component of $\Bun_G(\bS,\bD_\infty)$, and $\ol{\theta}$
that of $\Bun_{P}(\bD_0;\bD_\infty\cap\bD_0)$.

Recall the scheme $\CG_{G,P,\bX}$, cf. \secref{parabolic notation},
and observe that from \propref{bundles and maps} it follows that the
scheme $\Bun^\theta_{G;P}(\bS,\bD_\infty;\bD_0)\simeq 
\Bun^\theta_{G;P}(\bS',\bD'_\infty;\bD'_0)$ is isomorphic to the
scheme of based maps, $\on{Maps}^\theta(\bC,\CG_{G,P,\bX})$.
This description implies that $\Bun^\theta_{G;P}(\bS,\bD_\infty;\bD_0)$
is empty unless $\theta\in \widehat\Lambda_{\fg,\fp}^{pos}$.

\medskip

Our present goal is to introduce parabolic versions of the Uhlenbeck
spaces, $\fU^\theta_{G;P}$ and $\wt{\fU}{}^\theta_{G;P}$, both of
which contain $\Bun_{G;P}(\bS,\bD_\infty;\bD_0)$ as an open subset.

Let us first choose a pair of directions $(\bd_v,\bd_h)\in \bD_\infty$,
but with the condition that $\bd_h=\bD_\infty\cap\bD_0$. I.e., we have
a rational surface $\bS'\simeq \bC\times \bX$, and in addition to
$\infty_\bC\in \bC$ and $\infty_\bX\in \bX$, we have a distinguished
point $0_\bX\in \bX$, such that $\bD'_0=\bC\times 0_\bX$ is the proper
transform of $\bD_0$. Henceforth, we will identify $\bD_0\simeq \bD'_0\simeq\bC$.

We will first introduce $\fU^\theta_{G;P}$ and
$\wt{\fU}{}^\theta_{G;P}$ using a choice of $\bd_v$, and then show
that the definition is in fact independent of that choice.

\begin{defn}
For $\theta=(\ol{\theta},a)$, the parabolic Uhlenbeck space
$\fU^\theta_{G;P}$ is defined as the closure
of $\Bun^\theta_{G;P}(\bS',\bD'_\infty;\bD'_0)$ in the product
$\on{QMaps}^\theta(\bC,\CG_{G,P,\bX})
\underset{\on{QMaps}^a(\bC,\CG_{G,\bX})}\times \fU^a_G$.
\end{defn}

Note that for $P=G$, the scheme $\fU^\theta_{G;P}$
is nothing but $\fU^a_G$.

\begin{defn}
The enhanced parabolic Uhlenbeck space $\wt{\fU}{}^\theta_{G;P}$
is defined as the closure of $\Bun^\theta_{G;P}(\bS',\bD'_\infty;\bD'_0)$
in the product
$\wt{\on{QMaps}}{}^\theta(\bC,\CG_{G,P,\bX})
\underset{\on{QMaps}^a(\bC,\CG_{G,\bX})}\times \fU^a_G$.
\end{defn}

As we shall see later,
both $\fU^\theta_{G;P}$ and $\wt{\fU}{}^\theta_{G;P}$ are schemes of finite type.
Note also that from the definition of $\fU^a_G$ and
\thmref{equivalence of definitions}, it follows that
$\fU^\theta_{G;P}$ (resp., $\wt{\fU}{}^\theta_{G;P}$) could be
equivalently defined as the closure of $\Bun^\theta_{G;P}(\bS',\bD'_\infty;\bD'_0)$
inside the product $\on{QMaps}^\theta(\bC,\CG_{G,P,\bX})\times
\on{Sect}\left(\sO,\overset{\circ}\bC{}_\sO^{(a)}\right)$
(resp., $\wt{\on{QMaps}}{}^\theta(\bC,\CG_{G,P,\bX})\times
\on{Sect}\left(\sO,\overset{\circ}\bC{}_\sO^{(a)}\right)$). Note also
that $\wt{\fU}{}^\theta_{G;P}$ is a closed subscheme in the product
$\fU^\theta_{G;P}\underset{\overset{\circ}\bC{}^\theta}\times
\Mod_{M_{aff},\bC}^{\theta,+}$ (this is so because the
corresponding property holds on the level of quasi-maps' spaces). We
will denote the natural projection $\wt{\fU}{}^\theta_{G;P}\to
\fU^\theta_{G;P}$ by $\fr_{\fp_{aff}^+}$. Obviously, for $P=B$, the map
$\fr_{\fp_{aff}^+}$ is an isomorphism, but it is NOT an isomorphism
for $P=G$.

We will denote by $\varrho^\theta_{\fp_{aff}^+}$ and 
$\varrho^\theta_{M_{aff}}$ the natural maps
$\fU^\theta_{G;P}\to\overset{\circ}\bC{}^\theta$
and $\wt{\fU}{}^\theta_{G;P}\to \Mod_{M_{aff},\bC}^{\theta,+}$.

\medskip

The following {\it horizontal factorization property} of
$\fU^\theta_{G;P}$ and $\wt{\fU}{}^\theta_{G;P}$ follows from
the corresponding factorization properties of quasi-maps' spaces
(\propref{factorization of Zastava}) and of the Uhlenbeck space
(\propref{factorization of Uhlenbeck}):

\begin{cor} \label{horizontal factorization}
For $\theta=\theta_1+\theta_2$, we have canonical isomorphisms
$$\fU^\theta_{G;P}\underset{\overset{\circ}\bC{}^\theta}\times
(\overset{\circ}\bC{}^{\theta_1}\times \overset{\circ}\bC{}^{\theta_2})_{disj}\simeq
(\fU^{\theta_1}_{G;P}\times \fU^{\theta_2}_{G;P})\underset
{\overset{\circ}\bC{}^{\theta_1}\times \overset{\circ}\bC{}^{\theta_1}}\times
(\overset{\circ}\bC{}^{\theta_1}\times \overset{\circ}\bC{}^{\theta_2})_{disj},
\text{ and }$$
$$\wt{\fU}{}^\theta_{G;P}\underset{\overset{\circ}\bC{}^\theta}\times
(\overset{\circ}\bC{}^{\theta_1}\times \overset{\circ}\bC{}^{\theta_2})_{disj}\simeq
(\wt{\fU}{}^{\theta_1}_{G;P}\times \wt{\fU}{}^{\theta_2}_{G;P})\underset
{\overset{\circ}\bC{}^{\theta_1}\times \overset{\circ}\bC{}^{\theta_2}}\times
(\overset{\circ}\bC{}^{\theta_1}\times \overset{\circ}\bC{}^{\theta_2})_{disj}.$$
\end{cor}

\ssec{}   \label{two projections}

Now our goal is to prove that the schemes
$\fU^\theta_{G;P}$ and $\wt{\fU}{}^\theta_{G;P}$ are in fact
canonically attached to the surface $\bS$ with the divisor $\bD_0$,
i.e. that they do not depend on the choice of the direction $\bd_v$.
We will consider the case of $\wt{\fU}{}^\theta_{G;P}$, since
that of $\fU^\theta_{G;P}$ is analogous (and simpler).

By \thmref{Tannaka and Drinfeld}, for $S=\fU^a_G$,
the product $\bS\times S$
contains an open subset $U$, such that $\bS\times S-U$ is finite over $S$,
and over which we have a well-defined $G$-bundle $\F_G$. Let $D$ be a
divisor on $\bC\times S$, such that $\bC\times S-D\subset U$. We set
$^{Huge}\wt{\fU}{}^\theta_{G;P}$ to be the ind-scheme of meromorphic 
enhanced quasi-maps 
$$_{D\cdot \infty}\wt{\on{QMaps}}(\bC,\F_G|_{\bC\times S-D}\overset{G}\times
\CG_{\fg,\fp}),$$
cf. \secref{meromorphic Zastava} 
(recall that $\CG_{\fg,\fp}$ is the finite-dimensional flag variety $G/P$).

Clearly, $^{Huge}\wt{\fU}{}^\theta_{G;P}$, is defined in a way independent 
of the choice of $\bd_v$. We claim now that for any $\bd_v$ we have a closed
embedding $\wt{\fU}{}^\theta_{G;P}\hookrightarrow {}^{Huge}\wt{\fU}{}^\theta_{G;P}$.

Indeed, the scheme $\wt{\on{QMaps}}(\bC,\CG_{G,P,\bX})$ embeds as a closed 
subscheme into  the corresponding ind-scheme of meromorphic enhanced 
quasi-maps $_{D\cdot\infty}\wt{\on{QMaps}}(\bC,\CG_{G,P,\bX})$. Now,
$$_{D\cdot\infty}\wt{\on{QMaps}}(\bC,\CG_{G,P,\bX})
\underset{_{D\cdot\infty}\on{QMaps}(\bC,\CG_{G,\bX})}\times S \simeq 
{}^{Huge}\wt{\fU}{}^\theta_{G;P}.$$
This is because a data of an enhanced quasi-map
$(\bC\times S-D)\to \CG_{G,P,\bX}$, that projects to a given map
$(\bC\times S-D)\to \CG_{G,\bX}$ is equivalent to a data of an enhanced quasi-map
$$(\bC\times S-D)\to \F_G|_{\bC\times S-D}\overset{G}\times\CG_{\fg,\fp}.$$

\medskip

Note that the composition
$$\Bun^\theta_{G;P}(\bS',\bD'_\infty;\bD'_0)\to
\wt{\fU}{}^\theta_{G;P}\to {}^{Huge}\wt{\fU}{}^\theta_{G;P}$$
is also independent of $\bd_v$. This shows that we can identify $\wt{\fU}{}^\theta_{G;P}$ with
the closure of $\Bun^\theta_{G;P}(\bS',\bD'_\infty;\bD'_0)$ inside 
$^{Huge}\wt{\fU}{}^\theta_{G;P}$, which makes it manifestly independent of $\bd_v$.
Note also, that in the course of the proof we
established also the following:

\begin{cor}
The map $\wt{\fU}{}^\theta_{G;P}\to \Bun_M(\bC,\infty_\bC)$, which is the composition
of $\varrho^\theta_{M_{aff}}:
\wt{\fU}{}^\theta_{G;P}\to\Mod^{\theta,+}_{M_{aff},\bC}$ and
the natural projection $\Mod^{\theta,+}_{M_{aff},\bC}\to\Bun_M(\bC,\infty_\bC)$,
is independent of the choice of $\bd_v$.
\end{cor}

(Of course, the maps $\varrho^\theta_{M_{aff}}$ and
$\varrho^\theta_{\fp_{aff}^+}:\fU^\theta_{G;P}\to
\overset{\circ}\bC{}^\theta$, do depend on the choice of $\bd_v$.)

\begin{cor}
The schemes $\fU^\theta_{G;P}$ and $\wt{\fU}{}^\theta_{G;P}$ are of
finite type.
\end{cor}

(Indeed, the ind-scheme $^{Huge}\fU^\theta_{G;P}$ is of ind-finite
type, and $\wt{\fU}{}^\theta_{G;P}$ is the closure
of $\Bun^a_{G;P}(\bS,\bD_\infty,\bD_0)$ inside $^{Huge}\fU^\theta_{G;P}$.)

\ssec{Beilinson-Drinfeld-Kottwitz flag space}

Consider the ind-scheme $\Gr^{BD,a}_{G,P,\bX,0_\bX}$ fibered over
over $\overset{\circ}\bX{}^{(a)}$, which classifies the data of

\begin{itemize}

\item
A divisor $D\in \overset{\circ}\bX{}^{(a)}$,

\item
A principal $G$-bundle $\F_G$ on $\bX$,

\item
A trivialization $\F_G\simeq \F^0_G|_{\bX-D}$, and

\item
A reduction to $P$ of the $G$-torsor $\F_G|_{0_\bX}$.

\end{itemize}

We have an evident map $\Gr^{BD,a}_{G,P,\bX,0_\bX}\to \CG_{G,P,\bX}$, which
remembers of the trivialization of $\F_G$ on $\bX-D$ its restriction to the
formal neighborhood of $\infty_\bX$. Moreover, it is easy to see that the map
$\Gr^{BD,a}_{G,P,\bX,0_\bX}\to \CG_{G,P,\bX}\times \overset{\circ}\bX{}^{(a)}$
is in fact a closed embedding.

Thus, we can consider the relative based quasi-maps' space
$\on{QMaps}\left(\bC,\Gr^{BD,a}_{G,P,\bX,0_\bX}\right)$, fibered
over $\overset{\circ}\bX{}^{(a)}$, which is a closed subscheme inside
$\on{QMaps}\left(\bC,\CG_{G,P,\bX}\right)\times \overset{\circ}\bX{}^{(a)}$.

As in \propref{prime=two primes}, we obtain that the parabolic Uhlenbeck space
$\fU_{G;P}^\theta$ can be alternatively defined as the closure of
$\Bun^\theta_{G;P}(\bS',\bD'_\infty;\bD'_0)$ in the product
$$\on{QMaps}^\theta\left(\bC,\Gr^{BD,a}_{G,P,\bX,0_\bX}\right)
\underset{\on{QMaps}^a\left(\bC,\Gr^{BD,a}_{G,\bX}\right)}\times
\fU^a_G,$$
where $a$ is the projection of $\theta$ under $\widehat\Lambda_{\fg,\fp}\to \BZ$.

The above realization of $\fU_{G;P}^\theta$ gives an alternative proof that
$\fU_{G;P}^\theta$ is a scheme of finite type, since $\Gr^{BD,a}_{G,P,\bX,0_\bX}$
is of ind-finite type.

\medskip

We will establish another important property of parabolic Uhlenbeck
spaces (going back to the work ~\cite{v} of G.~Valli),
called the {\it vertical factorization property}:

\begin{prop}  \label{vertical factorization}
Let $\theta=(\ol{\theta},a)$ be an element of $\widehat\Lambda_{\fg,\fp}^{pos}$.
If $a=a_1+a_2$ is such that $\theta_1=\theta-a_2\cdot \delta$ also lies in
$\widehat\Lambda^{pos}_{\fg,\fp}$, then there are
canonical isomorphisms
\begin{align*}
&\fU_{G;P}^\theta\underset{\overset{\circ}\bX{}^{(a)}}\times
(\overset{\circ}\bX{}^{(a_1)}\times (\overset{\circ}\bX-0_\bX){}^{(a_2)})_{disj}
\simeq \\
&(\fU_{G;P}^{\theta_1}\times \fU_G^{a_2})
\underset{\overset{\circ}\bX{}^{(a_1)}\times (\overset{\circ}\bX-0_\bX){}^{(a_2)}}
\times
(\overset{\circ}\bX{}^{(a_1)}\times
(\overset{\circ}\bX-0_\bX){}^{(a_2)})_{disj}, \text{ and } \\
&\wt{\fU}{}_{G;P}^\theta\underset{\overset{\circ}\bX{}^{(a)}}\times
(\overset{\circ}\bX{}^{(a_1)}\times (\overset{\circ}\bX-0_\bX){}^{(a_2)})_{disj}
\simeq \\
&(\wt{\fU}{}_{G;P}^{\theta_1}\times \fU{}_G^{a_2})
\underset{\overset{\circ}\bX{}^{(a_1)}\times (\overset{\circ}\bX-0_\bX){}^{(a_2)}}
\times
(\overset{\circ}\bX{}^{(a_1)}\times (\overset{\circ}\bX-0_\bX){}^{(a_2)})_{disj}.
\end{align*}
If $\theta-\delta\notin \Lambda^{pos}_{aff,\fg,\fp}$, then the composition
$\fU_{G;P}^{\theta}\to \fU^a_G\overset{\varpi^a_v}\longrightarrow
\overset{\circ}\bX{}^{(a)}$ maps to the point $a\cdot 0_\bX\in
\overset{\circ}\bX{}^{(a)}$.
\end{prop}

\begin{proof}

First of all, the usual Uhlenbeck space $\fU^a_G$ has the vertical factorization
property
\begin{equation}  \label{vertical factorization of Uhlenbeck}
\fU^a_G\underset{\overset{\circ}\bX{}^{(a)}}\times
(\overset{\circ}\bX{}^{(a_1)}\times \overset{\circ}\bX{}^{(a_2)})_{disj}\simeq
(\fU^{a_1}_G\times \fU^{a_2}_G)\underset{\overset{\circ}\bX{}^{(a_1)}\times
\overset{\circ}\bX{}^{(a_2)}}\times (\overset{\circ}\bX{}^{(a_1)}\times
\overset{\circ}\bX{}^{(a_2)})_{disj},
\end{equation}
because $\bC$ and $\bX$ play symmetric roles in the definition of $\fU^a_G$.

Using \eqref{Bun_P} we obtain that the corresponding factorization
property for the open subscheme $\Bun^\theta_{G;P}(\bS,\bD_\infty;\bD_0)$
follows from that of $\Bun^a_G(\bS,\bD_\infty)$.

\medskip

To establish the factorization property of $\fU_{G;P}^\theta$ we will use the
fact (cf. \cite{g}) that
\begin{align}  \label{factorization of BDK-Grassmannian}
&\Gr^{BD,a}_{G,P,\bX,0_\bX}\underset{\overset{\circ}\bX{}^{(a)}}
\times(\overset{\circ}\bX{}^{(a_1)}\times
(\overset{\circ}\bX-0_\bX)^{(a_2)})_{disj}\simeq \\
&(\Gr^{BD,a_1}_{G,P,\bX,0_\bX} \times \Gr_{G,\bX}^{BD,a_2})
\underset{\overset{\circ}\bX{}^{(a_1)}\times \overset{\circ}\bX{}^{(a_2)}}\times
(\bX{}^{(a_1)}\times (\overset{\circ}\bX-0_\bX){}^{(a_2)})_{disj}.
\end{align}

Moreover, this decomposition is compatible with the corresponding line
bundles. Therefore, \eqref{factorization of BDK-Grassmannian} combined
with \eqref{vertical factorization of Uhlenbeck} and \lemref{product lemma} 
yield an isomorphism:
\begin{align*}
&\left(\on{QMaps}^{\theta} \left(\bC,\Gr^{BD,a}_{G,P,\bX,0_\bX}\right)
\underset{\on{QMaps}^a \left(\bC,\Gr^{BD,a}_{G,\bX}\right)}\times
\fU^a_G\right)\underset{\overset{\circ}\bX{}^{(a)}}\times 
(\overset{\circ}\bX{}^{(a_1)}\times(\overset{\circ}\bX-0_\bX){}^{(a_2)})_{disj}
\simeq \\
&\Bigl(\on{QMaps}^{\theta_1}
\left(\bC,\Gr^{BD,a_1}_{G,P,\bX,0_\bX}\right)\times
\on{QMaps}^{a_2}\left(\bC,\Gr_{G,\bX}^{BD,a_2}\right)\Bigr) 
\underset{\on{QMaps}^{a_1}
\left(\bC,\Gr^{BD,a_1}_{G,\bX}\right)\times
\on{QMaps}^{a_2}\left(\bC,\Gr_{G,\bX}^{BD,a_2}\right)}\times \\
&\times \left(\fU^{a_1}_G\times \fU^{a_2}_G \right)
\underset{\overset{\circ}\bX{}^{(a_1)}\times (\overset{\circ}\bX-0_\bX){}^{(a_2)}}
\times
(\overset{\circ}\bX{}^{(a_1)}\times(\overset{\circ}\bX-0_\bX){}^{(a_2)})_{disj}.
\end{align*}

Since $\fU_{G;P}^\theta$ is the closure of
$\Bun_{G;P}^\theta(\bS',\bD'_\infty;\bD'_0)$ in the product
$$\on{QMaps}^{\theta} \left(\bC,\Gr^{BD,a}_{G,P,\bX,0_\bX}\right)
\underset{\on{QMaps}^a \left(\bC,\Gr^{BD,a}_{G,\bX}\right)}\times
\fU^a_G,$$
we conclude that it factorizes in the required fashion.

\medskip

To prove the assertion for $\wt{\fU}^\theta_{G;P}$, recall that
it is the closure of $\Bun_{G;P}^\theta(\bS',\bD'_\infty;\bD'_0)$
in the product $\fU_{G;P}^\theta\times \Mod^{\theta,+}_{M_{aff}}$.

Note also that we have a natural closed embedding
$\Mod^{\theta_1,+}_{M_{aff}}\to
\Mod^{\theta,+}_{M_{aff}}$, such that we have a commutative square
$$
\CD
(\overset{\circ}\bX{}^{(a_1)}\times
\overset{\circ}\bX{}^{(a_2)})_{disj}\underset{\overset{\circ}\bX{}^{(a_1)}\times
\overset{\circ}\bX{}^{(a_2)}}\times
\left(\Bun_{G;P}^{\theta_1}(\bS',\bD'_\infty;\bD'_0)
\times \Bun_G^{a_2}(\bS',\bD'_\infty)\right)
@>>> \Mod^{\theta_1,+}_{M_{aff}}  \\
@V{\sim}VV   @VVV   \\
(\overset{\circ}\bX{}^{(a_1)}\times
\overset{\circ}\bX{}^{(a_2)})_{disj}\underset{\overset{\circ}\bX{}^{(a)}}\times
\Bun_{G;P}^\theta(\bS',\bD'_\infty;\bD'_0)  @>>> \Mod^{\theta,+}_{M_{aff}}.
\endCD
$$
Therefore, the factorization property for $\wt{\fU}^\theta_{G;P}$ follows
from that for $\fU_{G;P}^\theta$.

\medskip

Suppose now that $\theta-\delta\notin\widehat\Lambda^{pos}_{\fg,\fp}$.
Since $\Bun_{G;P}^\theta(\bS,\bD_\infty;\bD_0)$ is dense in $\fU_{G;P}^\theta$,
it is enough to prove that
$$\Bun_{G;P}^\theta(\bS,\bD_\infty;\bD_0)\to \Bun_G^a(\bS,\bD_\infty)\overset{\varpi^a_v}\longrightarrow
\overset{\circ}\bX{}^{(a)}$$ is the constant map to $a\cdot 0_\bX\in
\overset{\circ}\bX{}^{(a)}$.

If this were not the case, there would exist $a_1,a_2\neq 0$ with
$a_1+a_2=a$ and a based map
$$\bC\to \Gr^{BD,a_1}_{G,P,\bX,0_\bX}\times \Gr_{G,\bX}^{BD,a_2}$$
of total degree $\theta$, such that its projection onto
the second factor, i.e. $\bC\to \Gr_{G,\bX}^{BD,a_2}$ was not the
constant map. However, the space of such maps is the union over $b\in \BZ$ of
$$\on{Maps}^{\theta-b\cdot \delta}(\bC,\Gr^{BD,a_1}_{G,P,\bX,0_\bX})\times
\on{Maps}^b(\bC,\Gr_{G,\bX}^{BD,a_2}),
$$
which is non-empty only when $b\geq 0$ and
$\theta-b\cdot \delta\in \Lambda^{pos}_{\fg,\fp}$. But this forces
$b=0$, i.e., our map $\bC\to \Gr_{G,\bX}^{BD,a_2}$ must actually be
the constant map.

\end{proof}

\section{Stratification of parabolic Uhlenbeck spaces}

\ssec{}

For an element $\theta\in \widehat\Lambda{}^{pos}_{\fg,\fp}$ consider a decomposition
$\theta=\theta_1+\theta_2+b\cdot \delta$ with 
$\theta_1,\theta_2\in \widehat\Lambda{}^{pos}_{\fg,\fp}-0$, $b\in \BN$.
Let us denote by $\overset{\circ\circ}\bS$ the complement to
$\bD_0$ in $\overset{\circ}\bS$, and let us denote by
$\fU^{\theta;\theta_2,b}_{G;P}$ the following scheme:
$$\Bun_{G;P}^{\theta_1}(\bS,\bD_\infty;\bD_0)\times
\overset{\circ}\bC{}^{\theta_2}\times \on{Sym}^b(\overset{\circ\circ}\bS).$$

The following is a generalization of \thmref{stratification}:

\begin{thm}   \label{stratification of parabolic Uhlenbeck}
There exists a finite map
$\ol{\iota}_{\theta_2,b}:\fU^{\theta_1}_{G;P}\times
\overset{\circ}\bC{}^{\theta_2}\times \on{Sym}^b(\overset{\circ}\bS)\to
\fU^\theta_{G;P}$, compatible with the factorization isomorphisms of
\corref{horizontal factorization} and \propref{vertical factorization}.
Moreover, the composition
$$\iota_{\theta_2,b}:\fU^{\theta;\theta_2,b}_{G;P}\hookrightarrow
\fU^{\theta_1}_{G;P}\times
\overset{\circ}\bC{}^{\theta_2}\times \on{Sym}^b(\overset{\circ}\bS)\to
\fU^\theta_{G;P}$$
is a locally closed embedding, and
$\fU^\theta_{G;P}\simeq \underset{\theta_2,b}\cup\, \fU^{\theta;\theta_2,b}_{G;P}$.
\end{thm}

\begin{proof}

We construct a map 
$\fU^{\theta_1}_{G;P}\times
\overset{\circ}\bC{}^{\theta_2}\times \on{Sym}^b(\overset{\circ}\bS)\to
\on{QMaps}^\theta(\bC,\CG_{G,P,\bX})
\underset{\on{QMaps}^a(\bC,\CG_{G,\bX})}\times \fU^a_G$ as follows:

The projection
$\fU^{\theta_1}_{G;P}\times
\overset{\circ}\bC{}^{\theta_2}\times
\on{Sym}^b(\overset{\circ}\bS)\to \fU^a_G$ is the map
$$\fU^{\theta_1}_{G;P}\times
\overset{\circ}\bC{}^{\theta_2}\times
\on{Sym}^d(\overset{\circ}\bS)\to \fU^{a_1}_G\times
\overset{\circ}\bC{}^{(a_2)}\times \on{Sym}^b(\overset{\circ}\bS)\to
\fU^{a_1}_G\times \on{Sym}^{a_2+b}(\overset{\circ}\bS)\overset{\ol{\iota}_{a_2+b}}
\longrightarrow \fU^a_G,$$
where $\ol{\iota}_{a_2+b}$ is as in \thmref{stratification}.

The projection
$\fU^{\theta_1}_{G;P}\times
\overset{\circ}\bC{}^{\theta_2}\times \on{Sym}^b(\overset{\circ}\bS)\to
\on{QMaps}^\theta(\bC,\CG_{G,P,\bX})$ is
\begin{align*}
&\fU^{\theta_1}_{G;P}\times
\overset{\circ}\bC{}^{\theta_2}\times \on{Sym}^b(\overset{\circ}\bS)\to
\on{QMaps}^{\theta_1}(\bC,\CG_{G,P,\bX})\times \overset{\circ}\bC{}^{\theta_2}\times 
\overset{\circ}\bC{}^{(b)}\to \\
&\on{QMaps}^{\theta_1}(\bC,\CG_{G,P,\bX})\times
\overset{\circ}\bC{}^{\theta_2+b\cdot \delta}\to \on{QMaps}^\theta(\bC,\CG_{G,P,\bX}),
\end{align*}
where the first arrow in the above composition uses the map 
$\pi_v:\on{Sym}^b(\overset{\circ}\bS)\to \overset{\circ}\bC{}^{(b)}$,
the second arrow corresponds to the addition of divisors, and the
third arrow comes from \eqref{strata of quasi-maps}. It is easy to
see that the resulting map $\fU^{\theta_1}_{G;P}\times
\overset{\circ}\bC{}^{\theta_2}\times\on{Sym}^b(\overset{\circ}\bS)\to
\on{QMaps}^\theta(\bC,\CG_{G,P,\bX})\times \fU^a_G$ indeed maps to
the fiber product $\on{QMaps}^\theta(\bC,\CG_{G,P,\bX})
\underset{\on{QMaps}^a(\bC,\CG_{G,\bX})}\times \fU^a_G$ and that the
map $\iota_{\theta_2,b}$ is a locally closed embedding.

Moreover, it is easy to see that the image of
$\ol{\iota}_{\theta_2,b}$ belongs in fact to the subscheme 
$$\on{QMaps}^\theta\left(\bC,\Gr^{BD,a}_{G,P,\bX,0_\bX}\right)
\underset{\on{QMaps}^a\left(\bC,\Gr^{BD,a}_{G,\bX}\right)}\times
\fU^a_G,$$
and that this map is compatible with the isomorphisms of 
\propref{vertical factorization}.

The map $\ol{\iota}_{\theta_2,b}$ is also compatible with 
factorization isomorphisms \eqref{horizontal factorization}
(due to the corresponding property of the map $\ol{\iota}_b$
for the Uhlenbeck space $\fU^a_G$). Therefore, to
prove that the image of $\ol{\iota}_{\theta_2,b}$ belongs to 
$\fU^\theta_{G;P}$ (i.e., to the closure of 
$\Bun_{G;P}^\theta(\bS',\bD'_\infty;\bD'_0)$ in the product
$\on{QMaps}^\theta(\bC,\CG_{G,P,\bX})
\underset{\on{QMaps}^a(\bC,\CG_{G,\bX})}\times \fU^a_G$),
it suffices to analyze separately the following cases:
a) $\theta_2=0, b=0$; b) $\theta_1=0,\theta_2=0$; c) $\theta_1=0, b=0$. Of course,
case a) is trivial, since it corresponds to the open stratum, i.e.
$\Bun_{G;P}^\theta(\bS',\bD'_\infty;\bD'_0)$ itself.

\medskip

The second case, when $\theta_1=0,\theta_2=0$, follows immediately from 
\propref{vertical factorization}.
Indeed, the preimage of $(\overset{\circ}\bX-0_\bX)^{(b)}$ under
$\fU^{b\cdot \delta}_{G;P}\to \fU^b_G \overset{\varpi^b_h}\longrightarrow
\overset{\circ}\bX{}^{(b)}$ in this case is isomorphic
to $$\fU^b_G\underset{\overset{\circ}\bX{}^{(b)}}\times
(\overset{\circ}\bX-0_\bX)^{(b)},$$ which contains $\on{Sym}^b(\overset{\circ\circ}\bS)$
as a closed subscheme, according to \thmref{stratification}.

\medskip

Thus it remains 
to treat the case $\theta_1=0, b=0$. We will  use again the horizontal factorization
property, i.e. \corref{horizontal factorization}, which allows us to
assume that $\theta_2$ is the projection of a simple coroot 
of $\fg_{aff}$. 

To show that the image of $\overset{\circ}\bC\simeq 
\Bun_{G;P}^0(\bS,\bD_\infty;\bD_0)\times
\overset{\circ}\bC\times \on{Sym}^0(\overset{\circ\circ}\bS)$ belongs to
$\fU^\theta_{G;P}$, we will consider a certain 
$\BG_m$-action on $\on{QMaps}^\theta(\bC,\CG_{G,P,\bX})
\underset{\on{QMaps}^a(\bC,\CG_{G,\bX})}\times \fU^a_G$.

We choose a direction $\bd_v\in \bD_\infty-\bd_h$ and consider the 
$\BG_m$-action on $\bS$ as in \propref{contraction of Uhlenbeck}.
This action induces a $\BG_m$-action on $\on{QMaps}^\theta(\bC,\CG_{G,P,\bX})
\underset{\on{QMaps}^a(\bC,\CG_{G,\bX})}\times \fU^a_G$, which we
shall call ``the action of the first kind". Note that the corresponding
action on the first factor, i.e. $\on{QMaps}^\theta(\bC,\CG_{G,P,\bX})$,
corresponds to the canonical homomorphism $\BG_m\to M_{aff}\times\BG_m$,
where the latter is thought of as the Levi subgroup corresponding to
the parabolic $\fp^+_{aff}\subset \fg_{aff}$.

\medskip

Another $\BG_m$-action, which we shall call "the action of the second 
kind", corresponds to a dominant central cocharacter $\ol{\nu}:\BG_m\to Z(M)$,
as in \propref{contraction of quasi-maps}.

Note that the actions of the first and the second kind
commute with one another, and the compound $\BG_m$-action, which we 
shall call ``new", corresponds to a cocharacter
$\BG_m\to M_{aff}\times\BG_m$, which satisfies the dominance condition of
\propref{contraction of quasi-maps}. Using this proposition, combined
with \propref{contraction of Uhlenbeck}, we obtain
that the ``new" action contracts the space $\on{QMaps}^\theta(\bC,\CG_{G,P,\bX})
\underset{\on{QMaps}^a(\bC,\CG_{G,\bX})}\times \fU^a_G$ to 
$\iota_{\theta,0}(\overset{\circ}\bC)$. This contraction map is
dominant, due to its equivariance with respect to the action of 
$\BA^1\simeq \overset{\circ}\bC$ by horizontal shifts.

In particular, we obtain that $\iota_{\theta,0}(\overset{\circ}\bC)$
lies in the closure of $\Bun_{G;P}^\theta(\bS,\bD_\infty;\bD_0)$,
which is what we had to show.

\medskip

Now let us show that every geometric point of $\fU^\theta_{G;P}$ belongs
to one of the locally closed subschemes $\fU^{\theta;\theta_2,b}_{G;P}$.
We will argue by induction on the length of $\theta$.

Consider the map $\fU^\theta_{G;P}\to \fU^a_G$, and let $z\in\fU^a_G$
be the image of our geometric point. If $z$ has a singularity at
$\bs\in \overset{\circ\circ}\bS$, then $\varpi_v^a(z)\neq a\cdot 0_\bX$,
and by \propref{vertical factorization}, $\theta-\delta\in \widehat\Lambda{}^{pos}_{\fg,\fp}$.
Hence, our geometric point is contained
in the image of $\ol{i}_{0,1}$, and our assertion follows from the induction
hypothesis.

If $z$ has a singularity at $\bs\in \overset{\circ}\bC\times 0_\bX$,
then our geometric point is contained in the image of $\ol{i}_{\alpha_0,0}$,
and again the assertion follows by the induction hypothesis.

Hence, it remains to analyze the case when $z$ belongs to
$\Bun^a_G(\bS,\bD_\infty)$, i.e. we are dealing with the locus of
$\fU^\theta_{G;P}$ isomorphic to
$$\Bun_G(\bS,\bD_\infty)\underset{\Bun_G(\bC,\infty_\bC)}\times
\ol{\Bun}^{\ol{\theta}}_P(\bC,\infty_\bC),$$
where $\ol{\Bun}^{\ol{\theta}}_P(\bC,\infty_\bC)$ is the corresponding
version of Drinfeld's stack of generalized $P$-bundles on $\bC$ with
a trivialization at $\infty_\bC\in\bC$.

The assertion in this case follows from the corresponding property
of Drinfeld's compactifications, cf. \cite{bg1}.

\end{proof}

\ssec{}

In the course of the proof of \thmref{stratification of parabolic Uhlenbeck}
we have established the following:

\begin{cor}  \label{contraction of parabolic Uhlenbeck}
The ``new'' $\BG_m$-action on $\on{QMaps}^\theta(\bC,\CG_{G,P,\bX})
\underset{\on{QMaps}^a(\bC,\CG_{G,\bX})}\times \fU^a_G$ preserves
$\fU^\theta_{G;P}$ and contracts it
to the subscheme $\overset{\circ}\bC{}^{\theta}\overset{i_{\theta,0}}
\hookrightarrow \fU^\theta_{G;P}$. The contraction map
$\fU^\theta_{G;P}\to \overset{\circ}\bC{}^{\theta}$ coincides
with the map $\varrho^\theta_{\fp^+_{aff}}$.
\end{cor}

In what follows, we will sometimes use the notation $\fs_{\fp^+_{aff}}$
for the map $\overset{\circ}\bC{}^\theta\to \fU^\theta_{G;P}$.

\ssec{}

Now we turn to the space $\wt{\fU}{}^\theta_{G;P}$. For $\theta=
\theta_1+\theta_2+b\cdot\delta$, let us denote by
$\wt{\fU}{}^{\theta;\theta_2,b}_{G;P}$ the preimage in $\wt{\fU}{}^\theta_{G;P}$
of the locally closed subset $\fU^{\theta;\theta_2,b}_{G;P}\subset \fU^\theta_{G;P}$.

\smallskip

\begin{prop}  \label{description of enhanced strata}
We have an isomorphism on the level of reduced schemes:
$$\wt{\fU}{}^{\theta;\theta_2,b}_{G;P}\simeq
\left(\Bun_{G;P}^{\theta_1}(\bS,\bD_\infty;\bD_0)
\underset{\Bun_M(\bC,\infty_\bC)}
\times \H_{M_{aff},\bC}^{\theta_2,+}\right)\times
\on{Sym}^b(\overset{\circ\circ}\bS).$$
\end{prop}

\begin{proof}

Recall the map
$\fr_{\fp^+_{aff}}:\wt{\on{QMaps}}{}^\theta(\bC,\CG_{G,P,\bX})\to
\on{QMaps}^\theta(\bC,\CG_{G,P,\bX})$, and note
that we have a natural closed embedding
$$\wt{\fU}{}^\theta_{G;P}\to \fU^\theta_{G;P}
\underset{\on{QMaps}^\theta(\bC,\CG_{G,P,\bX})}\times
\wt{\on{QMaps}}{}^\theta(\bC,\CG_{G,P,\bX}).$$

Hence, we {\it a priori} have a closed embedding
$$\wt{\fU}{}^{\theta;\theta_2,b}_{G;P}\to
\left(\Bun_{G;P}^{\theta_1}(\bS,\bD_\infty;\bD_0)\underset{\Bun_M(\bC,\infty_\bC)}
\times \H_{M_{aff},\bC}^{\theta_2+b\cdot\delta,+}\right)\times
\on{Sym}^b(\overset{\circ\circ}\bS),$$
by \propref{strata s volnoj}. However, using the vertical
factorization property, \propref{vertical factorization},
we obtain that its image is in fact contained in
$\left(\Bun_{G;P}^{\theta_1}(\bS,\bD_\infty;\bD_0)\underset{\Bun_M
(\bC,\infty_\bC)}\times \H_{M_{aff},\bC}^{\theta_2,+}\right)\times
\on{Sym}^b(\overset{\circ\circ}\bS)$.
Thus, we have to show that the above inclusion is in fact an equality on
the level of reduced schemes.

\medskip

Using again the vertical factorization property, we can assume that $b=0$.
We know (cf. \lemref{smoothness of reduction map}) that the map
$\Bun_G(\bS,\bD_\infty)\simeq\Bun_G(\bS',\bD'_\infty)\to
\Bun_G(\bC,\infty_\bC)$, given by restriction
of $G$-bundles under $\bC\simeq \bD_0\to \bS$, is smooth. Hence, the map
\begin{align*}
&\Bun_{G;P}(\bS,\bD_\infty;\bD_0)\simeq \\
&\simeq \Bun_G(\bS,\bD_\infty)\underset{\Bun_G(\bC,\infty_\bC)}\times
\Bun_P(\bC,\infty_\bC)\to \Bun_P(\bC,\infty_\bC)\to \Bun_M(\bC,\infty_\bC)
\end{align*}
is smooth as well. Therefore, the scheme
$\Bun_G^{\theta_1}(\bS,\bD_\infty;\bD_0)\underset{\Bun_M(\bC,\infty_\bC)}
\times \H_{M_{aff},\bC}^{\theta_2,+}$ is irreducible, and our
surjectivity assertion is enough to check at the generic point.

\medskip

Consider the map
$$\Bun_{G;P}^{\theta_1}(\bS,\bD_\infty;\bD_0)\underset{\Bun_M(\bC,\infty_\bC)}
\times \H_{M_{aff},\bC}^{\theta_2,+}\to \left(\overset{\circ}\bC{}^{\theta_1}\times
\overset{\circ}\bC{}^{\theta_2}\right).$$
It is sufficient to analyze the locus which projects to pairs of simple
and mutually disjoint divisors in $\overset{\circ}\bC{}^{\theta_1}\times
\overset{\circ}\bC{}^{\theta_2}$. Moreover using the horizontal factorization
property, we reduce the assertion to the case when
$\theta_1=0$, and $\theta_2=\theta$ is the projection of a simple
coroot in $\fg_{aff}$.

For $\theta$ of the form $(\ol{\theta},0)$ we have an isomorphism
$\wt{\fU}{}^{\theta}_{G;P}\simeq \wt{\on{Bun}}{}^{\ol{\theta}}_P(\bC,\infty_\bC)$.
Therefore, if $\theta$ is the projection of a coroot belonging to
$\fg$, our assertion follows from the corresponding fact for
$\wt{\on{Bun}}{}^{\ol{\theta}}_P(\bC,\infty_\bC)$, cf. \cite{bgfm}.
Hence, it remains to consider the case when $\theta$ is the projection
of the simple affine coroot $\alpha_0$.

\medskip

Consider the ``new" $\BG_m$-action on $\wt{\fU}{}^\theta_{G;P}$ as in
\corref{contraction of parabolic Uhlenbeck}.
In the same way we obtain a $\BG_m$-action on $\wt{\fU}{}^\theta_{G;P}$,
which contracts this space to
$$\on{Mod}_{M_{aff},\bC}^{\theta,+} \simeq
\overset{\circ}\bC{}^\theta\underset{\overset{\circ}\bC{}^\theta}
\times\on{Mod}_{M_{aff},\bC}^{\theta,+} \subset
\fU^\theta_{G;P}\underset{\overset{\circ}\bC{}^\theta}
\times\on{Mod}_{M_{aff},\bC}^{\theta,+}.$$
Thus, it is enough to show that the map
$\varrho^\theta_{M_{aff}}:\wt{\fU}^\theta_{G;P}\to
\on{Mod}_{M_{aff},\bC}^{\theta,+}$ is dominant. Of course, it is
sufficient to prove that the restriction of $\varrho^\theta_{M_{aff}}$
to $\Bun_{G;P}^\theta(\bS,\bD_\infty;\bD_0)$ is dominant.

\medskip

Let us fix a point $\bc\in \overset{\circ}\bC$ and consider the the preimage of 
$a\cdot \bc\in \bC$ under $$\Bun_{G;P}^\theta(\bS,\bD_\infty;\bD_0)\to 
\Bun_G^a(\bS,\bD_\infty)\overset{\varpi^a_h}\longrightarrow
\overset{\circ}\bC{}^{(a)},$$ which we will denote by 
$\overset{\circ}{\fF}{}^\theta_{\fg_{aff},\fp^+_{aff}}$, 
cf. \secref{central fiber notation}.
The map $\varrho^\theta_{M_{aff}}$ gives rise to a map 
$\overset{\circ}{\fF}{}^\theta_{\fg_{aff},\fp^+_{aff}}\to 
\Gr^{\theta,+}_{M_{aff}}:=\Gr_{M_{aff}}\cap\on{Mod}_{M_{aff},\bC}^{\theta,+}$.

Recall (cf. \secref{convolution}) that the affine Grassmannian
$\Gr_M$ is a union of locally closed subsets $\Gr_M^\nu$, $\nu\in \Lambda^+_\fm$,
and $\ol{\Gr}{}_M^\nu=\underset{\nu'\leq \nu}\cup\,
\Gr_M^{\nu'}$. Note that on the level of reduced schemes, 
$(\Gr^{\theta,+}_{M_{aff}})_{red}$ is isomorphic to $\ol{\Gr}_M^{w_0^M(\ol{\theta})}$, 
where $w_0^M$ is the longest element in the Weyl group of $M$.

Recall now that we have fixed $\theta$ to be $\alpha_0$. Therefore,
the complement to $\Gr_M^{-w_0^M(\ol{\alpha}_0)}$ in $\ol{\Gr}_M^{-w_0^M(\ol{\alpha}_0)}$ 
is a point-orbit corresponding to $\Gr^{0,+}_{M_{aff}}\overset{\delta}\longrightarrow
\Gr^{\alpha_0,+}_{M_{aff}}$ (unless $M=T$). Therefore, using \corref{locally trivial},
we conclude that if the map $\overset{\circ}{\fF}{}^{\alpha_0}_{\fg_{aff},\fp^+_{aff}}\to
\Gr^{\alpha_0,+}_{M_{aff}}$ were not dominant for some (and hence, every)  
point $\bc\in \bC$, we would obtain that for every geometric point of
$\Bun_{G;P}^\theta(\bS,\bD_\infty;\bD_0)$,
the corresponding $P$-bundle on $\bC$ is such that the induced $M$-bundle is trivial.
However, since $\bC$ is of genus $0$, this would mean the $G$-bundle
on $\bC$ is trivial. But this is a contradiction: we know
(cf. \propref{bundles and maps}) that there must exist a horizontal line $\bC\times
\bx$ such that the restriction of the $G$-bundle to it is non-trivial,
but according to the second part of \propref{vertical factorization},
this point $\bx$ must equal $0_\bX$.

\end{proof}

In what follows, we will denote by $\wt{\fs}_{\fp^+_{aff}}$ the map
$\on{Mod}_{M_{aff},\bC}^{\theta,+}\to \wt{\fU}{}^\theta_{G;P}$.

\ssec{Functorial definitions in the parabolic case}

Let us briefly discuss alternative definitions of $\fU_{G;P}^\theta$
$\wt{\fU}^\theta_{G;P}$, which will be functorial in the spirit
of $^{Drinf}\fU^a_G$. We will consider $\wt{\fU}^\theta_{G;P}$,
since the case of $\fU_{G;P}^\theta$ is analogous (and simpler).

Consider the functor $^{Funct}\wt{\fU}{}^\theta_{G;P}:\{$
Schemes $\}\to\{$ Sets $\}$ that assigns to a scheme $S$ the data of:

\begin{itemize}

\item An $S$-point of $\fU^a_G$. (In particular, we obtain a
$G$-bundle $\F_G$ defined on an open subset $U\subset \bS\times S$
such that $\bS\times S-U$ is finite over $S$. Let $D\subset
\bC\times S$ be a relative Cartier divisor such that 
$\bC\times S-D\subset U\cap \bC\times S$.)

\item A meromorphic enhanced quasi-map
$\sigma\in {}_{\infty\cdot D_h}\wt{\on{QMaps}}{}^{\ol{\theta}}(\bC,
\F_G\overset{G}\times \CG_{\fg,\fp})$. (In particular,
we have an $M$-bundle $\F_M$ defined on $\bC\times S$),

\end{itemize}

such that the following condition is satisfied:

\medskip

\noindent For any direction $\bd_v\in \bD_\infty-\bd_h$, and the corresponding
decomposition $\bS'\simeq \bC\times \bX$, the meromorphic enhanced quasi-map
$$\widehat{\sigma}:\bC\to \CG_{G,P,\bX},$$
obtained from $\F_G$ and $\sigma$, is in fact regular.

\medskip

As before one shows that the functor $^{Funct}\wt{\fU}{}^\theta_{G;P}$
is representable by an ind-scheme of ind-finite type.

\begin{thm}
There exists a canonical closed embedding
$\wt{\fU}{}^\theta_{G;P}\to ^{Funct}\wt{\fU}{}^\theta_{G;P}$,
which induces an isomorphism on the level of geometric points.
\end{thm}

\begin{proof}

The closed embedding $\wt{\fU}{}^\theta_{G;P}\to
^{Funct}\wt{\fU}{}^\theta_{G;P}$ has been constructed in
\secref{two projections}.

Consider a geometric point of $^{Funct}\wt{\fU}{}^\theta_{G;P}$,
and let $\bs_1,...,\bs_k$ be the points of $\overset{\circ\circ}\bS$, where the
corresponding point $z\in \fU^a_G$ has a singularity, and let
$b_1,...,b_k$ be the corresponding multiplicities.
Let also $\bc_1,...,\bc_n$ be the points on $\bD_0\simeq\bC$, where the
enhanced meromorphic quasi-map $\sigma$ (as in the definition of
$^{Funct}\wt{\fU}{}^\theta_{G;P}$) is either non-defined or has a
singularity.

Let us choose a direction $\bd_v$ so that the points $\bc_i$ are
disjoint from $\pi_v(\bs_j)$. Then, the enhanced quasi-map $\widehat{\sigma}$
(as in the definition of $^{Funct}\wt{\fU}{}^\theta_{G;P}$)
has singularities only at the above points $\bc_i$
and $\pi_v(\bs_j)$. Let $\theta_i\in \widehat\Lambda{}^{pos}_{\fg,\fp}$
denote the defect of $\widehat{\sigma}$ at $\bc_i$. The defect of $\widehat{\sigma}$ at
$\pi_v(\bs_j)$ is automatically equal to $b_j\cdot \delta$.

\medskip

Set $\theta'=\theta-\underset{i}{\Sigma} \theta_i
-\underset{j}{\Sigma} b_j\cdot \delta$. According to
\propref{strata s volnoj}, the quasi-map $\widehat{\sigma}$ can be described as its
saturation $\widehat{\sigma}'\in \Bun_{G;P}^{\theta'}(\bS,\bD_\infty;\bD_0)$
and the corresponding positive modification, i.e., a point in
$$pt\underset{\Bun_{M_{aff}}(\bC,\infty_{\bC})}\times\left(\underset{j}{\Pi}\,
\H^{b_j\cdot \delta,+}_{M_{aff},\bC,\pi_v(\bs_j)}\times
\underset{i}{\Pi}\, \H^{\theta_i,+}_{M_{aff},\bC,\bc_i}\right),$$
where $pt\in \Bun_{M_{aff}}$ corresponds to the pair $\F_M$ (cf. the
definition of $^{Funct}\wt{\fU}{}^\theta_{G;P}$) and the line bundle
corresponding to the divisor $\varpi^a_h(z)$ (for $z\in \fU^a_G$ as above).

The theorem now follows from the description of points of $\wt{\fU}{}^\theta_G$ given by
\propref{description of enhanced strata}.

\end{proof}

\section{The boundary}

\ssec{}

Let $\Y$ be a scheme of finite type, whose underlying
reduced sub-scheme $\Y_{red}$ is irreducible; let
$\Y_n\to \Y_{red}$ be its normalization. Let $\Y'\subset \Y$ be a closed 
subset, and  let $\Y'_n$ be the preimage of $\Y'$ in $\Y_n$.

\begin{defn}
$\Y'$ will be called a quasi-effective Cartier divisor if there
exists a line bundle $\CP$ on $\Y$, and a meromorphic 
(i.e., defined over a dense open subset) section
$s:\CO_\Y\to \CP$, such that
\begin{itemize}

\item The section $s$ is an isomorphism on
$\Y-\Y'$,

\item The section $s_n:\CO_{\Y_n}\to \CP\underset{\CO_\Y}\otimes \CO_{\Y_n}$
is regular,

\item The set of zeroes of $s_n$
coincides with $\Y'_n$.
\end{itemize}
\end{defn}

The notion of quasi-effective Cartier divisor is useful because of
the following lemma:

\begin{lem}  \label{dimension drop}
Let $\Y'\subset \Y$ be a quasi-effective Cartier divisor, and let
$Z\subset \Y'$ be a non-empty irreducible subvariety, not contained
in $\Y'$. Then $\on{dim}(Z\cap \Y')=\on{dim}(Z)-1$.
\end{lem}

\begin{proof}
Let $Z_n$ be the preimage of $Z$ in $\Y_n$. It is sufficient to show
that $\on{dim}(Z_n\cap \Y'_n)=\on{dim}(Z_n)-1$, but this is true,
because $\Y'_n$ is a Cartier divisor in $\Y_n$.
\end{proof}

\medskip

Here is a simple criterion how to show that a subset $\Y'$ is
a quasi-effective Cartier divisor. Suppose that we are given
a line bundle $\CP$ on $\Y$ with a trivialization
$s_n:\CO\to \CP$ over $\Y-\Y'$. Suppose also that $\Y'$ contains
a dense subset of the form $\underset{i}\cup\, \Y'_i$, where
$\Y'_i\subset \Y$ are subvarieties of codimension $1$, such that
$s$ has a zero at the generic point of each $\Y'_i$.

\begin{lem}  \label{QCartier criterion}
Under the above circumstances, $\Y'$ is a quasi-effective Cartier divisor.
\end{lem}

\ssec{The finite-dimensional case}

Recall the space $\ol{\Bun}_B(\bC)$ (cf. \secref{fd Zastava}) defined
for any projective curve $\bC$.
By definition, the boundary of $\ol{\Bun}_B(\bC)$ is the closed subset
$\partial(\ol{\Bun}_B(\bC)):=\ol{\Bun}_B(\bC)-\Bun_B(\bC)$.

\begin{thm}  \label{fd Cartier}
The boundary $\partial(\ol{\Bun}_B(\bC))$ is a quasi-effective Cartier
divisor in $\ol{\Bun}_B(\bC)$.
\end{thm}

The idea of the proof given below is borrowed from Faltings's
approach to the construction of the determinant line bundle on
the moduli space of bundles on a curve, cf. \cite{Fa}.

\medskip

Recall that the moduli space $\Bun_G(\bC)$ is equipped with the line bundle
$\CP_{\Bun_G(\bC)}$. Its $2\check{h}$-th tensor power can be identified
with the determinant line bundle $\CP_{\Bun_G(\bC),det}$, whose fiber 
at $\F_G\in \Bun_G(\bC)$ is
$\on{det}\left(R\Gamma(\bC,\fg_{\F_G})\right)$.

Let us consider the line bundle $\CP_{\Bun_T(\bC)}$ on the stack $\Bun_T(\bC)$,
whose fiber at $\F_T\in \Bun_T(\bC)$ is
$$\underset{\check\alpha}\otimes\,
\on{det}\left(R\Gamma(\bC,\CL^{\check\alpha}_{\F_T})\right),$$
where $\check\alpha$ runs over the set of all roots of $\fg$.

We have the canonical projections
$p_G:\ol{\Bun}_B(\bC)\to \Bun_G(\bC)$ and $p_T:\ol{\Bun}_B(\bC)\to
\Bun_T(\bC)$, and we define the line bundle $\CP_{\ol{\Bun}_B(\bC)}$ on
$\ol{\Bun}_B(\bC)$ as $p_G^*((\CP_{\Bun_G(\bC),det})^{-1})\otimes p_T^*(\CP_{\Bun_T(\bC)})$.

Observe that over the open part $\Bun_B(\bC)\subset
\ol{\Bun}_B(\bC)$ we have an (almost canonical)
trivialization of $\CP_{\ol{\Bun}_B(\bC)}$.
Indeed for a point $\F_B\in \Bun_B(\bC)$, with the induced $G$- and
$T$- bundle being $\F_G$ and $\F_T$, respectively,
the vector bundle
$\fg_{\F_G}$ has a filtration with the successive quotients
being $\CL^{\check\alpha}_{\F_T}$, and the trivial bundle $\fh\otimes \CO$. Hence,
$$\on{det}\left(R\Gamma(\bC,\fg_{\F_G})\right)\simeq \underset{\check\alpha}\otimes\,
\on{det}\left(R\Gamma(\bC,\CL^{\check\alpha}_{\F_T})\right) \otimes
\on{det}\left(R\Gamma(\bC,\CO)\right)^{\on{dim}(T)}.$$

In order to be able to apply \lemref{QCartier criterion},
we need to analyze the behavior of the section
$s:\CO\to \CP_{\ol{\Bun}_B(\bC)}$ just constructed at the generic point
at each irreducible component of $\partial(\ol{\Bun}_B(\bC))$.

\medskip

It is well-known (cf. \cite{bg1}) that $\partial(\ol{\Bun}_B(\bC))$
contains as a dense subset the following codimension $1$ locus:
it is the union of connected components
enumerated by the vertices of the Dynkin diagram of $\fg$; 
each such component $\partial(\ol{\Bun}_B(\bC))_i$ is isomorphic
to $\bC\times \Bun_B(\bC)$, and we have an identification
$$\Bun_B(\bC)\cup \partial(\ol{\Bun}_B(\bC))_i\simeq \Bun_{P_i}(\bC)
\underset{\Bun_{M_i}(\bC)}\times \ol{\Bun}_{B_i}(\bC),$$
where $P_i$ (resp., $M_i$) is the corresponding subminimal parabolic
(resp., its Levi quotient), and $\ol{\Bun}_{B_i}(\bC)$ is the
corresponding space for the group $M_i$, cf. 
\lemref{relation between parabolics} and \corref{codimension 1 locus}.

\medskip

This reduces our problem to the following calculation.
Let $G$ be a reductive group of semi-simple rank $1$ and
the derived group isomorphic to $SL_2$, and let $\CV$ be an irreducible $G$-module
of dimension $n+1$. Let $\CP_{\Bun_G(\bC),\CV}$ be the line bundle on
$\Bun_G(\bC)$ given by
$$\F_G\mapsto \on{det}\left(R\Gamma(\bC,\CV_{\F_G})\right),$$
and $\CP_{\Bun_T(\bC),\CV}$ be the line bundle on $\Bun_T(\bC)$ given by
$$\F_T\mapsto \underset{\check\mu}\otimes\,
\on{det}\left(R\Gamma(\bC,\CL^{\check\mu}_{\F_T})\right),$$
where $\check\mu$ runs over the set of weights of $\CV$.
Let $\CP_{\ol{\Bun}_B(\bC),\CV}$ be the corresponding line bundle on $\ol{\Bun}_B(\bC)$,
and $s$ its trivialization over $\Bun_B(\bC)$.

\begin{lem}
The section $s$ has a zero of order $\frac{n(n+1)(n+2)}{6}$ along the boundary
$\partial(\ol{\Bun}_B(\bC))$.
\end{lem}

\begin{proof}

Recall that for groups of semi-simple rank $1$, the stack
$\ol{\Bun}_B(\bC)$ is smooth, and the variety
$\partial(\ol{\Bun}_B(\bC))$ is irreducible.

First of all, it is easy to to see that our initial group can be
replaced by $GL_2$, such that $\CV$ is the $n$-th symmetric power of
the standard representation.
We will construct a map of the projective
line $\BP^1$ into $\ol{\Bun}_B(\bC)$, such that
$\BA^1\subset \BP^1$ maps to the open part $\Bun_B(\bC)$, and the
intersection $\BP^1\cap \partial(\ol{\Bun}_B(\bC))$ is transversal,
and such that the pull-back of $\CP_{\ol{\Bun}_B(\bC),\CV}$ to $\BP^1$ has degree
$\frac{n(n+1)(n+2)}{6}$.

\medskip

Fix a point $\bc\in \bC$ and a point $\F_B\in \Bun_B(\bC)$, i.e.
a vector bundle $\CM$ on $\bC$ and a short exact sequence
$$0\to \CL_1\to \CM\to \CL_2\to 0.$$
Consider the projective line of all elementary lower modifications
of $\CM$ at $\bc$. This $\BP^1$ maps to $\ol{\Bun}_B(\bC)$ by
setting the new $T$-bundle to be $(\CL_1(-\bc),\CL_2)$. The
required properties of the embedding $\BP^1\to\ol{\Bun}_B(\bC)$
are easy to verify.
By construction, the composition
$\BP^1\to\ol{\Bun}_B(\bC)\overset{p_T}\to \Bun_T(\bC)$ is the constant
map; therefore, it suffices to calculate the degree of the pull-back
of $\CP_{\Bun_G(\bC),\CV}$
under $\BP^1\to\ol{\Bun}_B(\bC)\overset{p_G}\to \Bun_G(\bC)$.

But for a lower modification $\CM'\subset \CM$,
the quotient $\on{Sym}^n(\CM)/\on{Sym}^n(\CM')$ has a canonical
$n$-step filtration with successive quotients isomorphic to
$$(\CM/\CM')^{\otimes j}\otimes \on{Sym}^{n-j}(\CM/\CM(-x)).$$
Therefore, on the level of determinant lines,
$\on{det}\left(R\Gamma(\bC,\on{Sym}^n(\CM'))\right)$ is isomorphic 
to a fixed $1$-dimensional vector space, times
$(\CM/\CM')^{\otimes \frac{n(n+1)(n+2)}{6}}$.

\end{proof}

This implies the assertion of the theorem.

\ssec{}

We will now formulate and prove an analog of \thmref{fd Cartier}
for Uhlenbeck spaces on $\BP^2$.

For $\mu\in \widehat\Lambda^{pos}_\fg=\widehat\Lambda^{pos}_{\fg,\fb}$,
consider the open subset $\overset{\circ}\fU{}^\mu_{G;B}$
in $\fU^\mu_{G;B}$ equal to the union of the strata $\fU^{\mu;0,b}_{G;B}$.
In other words, $\overset{\circ}\fU{}^\mu_{G;B}$ consists of points,
which have no singularities along $\bD_0$.
(This means that the induced point
of $\fU^a_G$ gives rise to a bundle $\F_G$ defined in a neighborhood of
$\bD_0$, and the corresponding quasi-map $\bC\to
\F_G|_{\bC}\overset{G}\times \CG_{\fg,\fb}$ is in fact a map.)
Of course, $\overset{\circ}\fU{}^\mu_{G;B}$ contains
$\Bun^\mu_{G;B}(\bS,\bD_\infty;\bD_0)$ as an open subset, and from
\thmref{stratification of parabolic Uhlenbeck} we obtain that the
complement is of codimension at least $2$.

\medskip

We define the boundary $\partial(\fU^\mu_{G;B})\subset
\fU^\mu_{G;B}$ as the complement $\fU^\mu_{G;B}-\overset{\circ}\fU{}^\mu_{G;B}$. From 
\thmref{stratification of parabolic Uhlenbeck}, we obtain that
$\partial(\fU^\mu_{G;B})$ contains a dense subset of codimension
$1$ equal to the disjoint union over $i\in I_{aff}$ of
$$\partial(\fU^\mu_{G;B})_i:=\fU^{\mu;\alpha_i,0}_{G;B}\simeq
\Bun^{\mu-\alpha_i}_{G;B}(\bS,\bD_\infty;\bD_0)\times\overset{\circ}\bC.$$

\medskip

\noindent{\it Remark.}
Note also that $\fU^\mu_{G;B}$ is regular at the generic point
of each $\partial(\fU^\mu_{G;B})_i$ (in particular,
$\fU^\mu_{G;B}$ is regular in codimension $1$.)
Indeed, from \thmref{stratification of parabolic Uhlenbeck},
it is easy to deduce that the open subset
$$\Bun^\mu_{G;B}(\bS,\bD_\infty;\bD_0)\bigcup\, \left(\underset{i\in I_{aff}}\cup\,
\fU^{\mu;\alpha_i,0}_{G;B}\right)$$
projects one-to-one on the open subset in
$\on{QMaps}^\mu(\bC,\CG_{G,P,\bX})$ equal to
$$\on{Maps}(\bC,\CG_{\fg_{aff},\fb^+_{aff}})\bigcup \left(\underset{i\in I_{aff}}\cup\,
\partial(\on{QMaps}(\bC,\CG_{G,P,\bX}))_i\right),$$
and the latter is smooth, by \lemref{relation between parabolics}.

\medskip

\begin{thm}   \label{boundary of Uhlenbeck}
The boundary $\partial(\fU^\mu_{G;B})$ is a quasi-effective
Cartier divisor in $\fU^\mu_{G;B}$.
\end{thm}

The rest of this section is devoted to the proof of this theorem.

\ssec{}

Let us choose a pair of directions $(\bd_v,\bd_h)$, and consider the corresponding
decomposition $\bS'\simeq \bC\times \bX$. From the map 
$\varpi^a_v:\fU^a_G\to \overset{\circ}\bX{}^{(a)}$
we obtain a relative divisor on $\fU^a_G\times \bX$, and the
associated line bundle.

We define the line bundle $\CP_{\fU,G}$ on $\fU^a_G$ as the
restriction of the above line bundle to $\fU^a_G\times 0_\bX\subset
\fU^a_G\times \bX$.

Let $\overset{\circ}{\fU}{}^a_G$ be the open subset in $\fU^a_G$
corresponding to points with no singularity along $\bD_0$.
We have a canonical projection $\overset{\circ}{\fU}{}^a_G\to
\Bun_G(\bC,\infty_{\bC})$.

The next assertion follows from the construction of the map
$\varpi^a_v$:

\begin{lem}
The restriction of $\CP_{\fU,G}$ to $\overset{\circ}{\fU}{}^a_G$
can be canonically identified with the pull-back of $\CP_{\Bun_G(\bC)}$ under
$\overset{\circ}{\fU}{}^a_G\to\Bun_G(\bC)$.
\end{lem}

We have a canonical projection $p_{T,aff}$ from
$\fU^\mu_{G;B}$ to $\Bun_{T_{aff}}(\bC)$. We define the line
bundle $\CP_{\fU,B}$ on $\fU^\mu_{G;B}$
as the tensor product of the pull-back of $\CP_{\Bun_T(\bC)}$ under
$$\fU^\mu_{G;B}\overset{p_{T,aff}}\longrightarrow \Bun_{T_{aff}}(\bC)\to \Bun_T(\bC)$$
and the inverse of the pull-back of $\CP^{\otimes 2\cdot h^\vee}_{\fU,G}$ under
$\fU^\mu_{G;B}\to \fU^a_G$.

We claim that over $\overset{\circ}\fU{}^\mu_{G;B}$
there is a canonical trivialization of $s:\CO\to\CP_{\fU,B}$.
Indeed, we have a projection $\overset{\circ}\fU{}^\mu_{G;B}\to
\Bun_B(\bC)$, and the restriction of $\CP_{\fU,B}$ is the pull-back of
$\CP_{\ol{\Bun}_B(\bC)}$ under this map. Therefore, the existence of the trivialization
follows from the corresponding fact for $\Bun_B(\bC)$, using the fact
that $\CP_{\Bun_G(\bC),det}\simeq \CP_{\Bun_G(\bC)}^{\otimes 2\cdot h^\vee}$.

Thus, to prove \thmref{boundary of Uhlenbeck}, it remains to show
that the above meromorphic section $s$ has a zero at the generic
point of each $\partial(\fU^\mu_{G;B})_i$, $i\in I_{aff}$.

When $i\in I\subset I_{aff}$, this readily follows from \thmref{fd Cartier},
since in this case $\partial(\fU^\mu_{G;B})_i$ is contained in the
preimage of $\overset{\circ}{\fU}{}^a_G$ in $\fU^\mu_{G;B}$.

\medskip

Therefore, it remains to analyze the behavior of $s$
on the open subvariety $$\fU^\mu_{G;B}-\ol{\underset{i\in I}\cup\,
\partial(\fU^\mu_{G;B})_i}$$
at the generic point of $\partial(\fU^\mu_{G;B})_{i_0}$, where
$i_0$ corresponds to the affine root.

\begin{lem}
Let $\Y$ be a variety (regular in codimension $1$),
acted on by the group $\BG_m$;
let $\CP$ be an equivariant line bundle, and $s:\CO\to \CP$
be an equivariant nowhere vanishing section defined outside
an irreducible subvariety $\Y'$ of codimension $1$.

Suppose that there exists a point $y\in \Y-\Y'$ such that
the action map $\BG_m\times y\to \Y$ extends to a map
$\BA^1\to \Y$, such that the image of $0$, call it $y'$, belongs
to $\Y'$. Assume moreover, that the $\BG_m$-action on
the fiber of $\CP$ at $y'$ (note that $y'$ is automatically
$\BG_m$-stable) is given by a positive power
of the standard character. Then $s$ vanishes at the
generic point of $\Y'$.
\end{lem}

\begin{proof}

If we lift $y$ to a geometric point in the normalization
of $\Y$, the same conditions will hold; therefore we
can assume that $\Y$ is normal.

To prove the lemma we have
to exclude the possibility that $s$ has a pole
of order $\geq 0$ on $\Y'$. If it did, the same would
be true for the pull-back of the pair $(\CP,s)$ to
$\BA^1$, endowed with the standard $\BG_m$-action.
In the latter case, to have a pole of non-negative order
means that $\CP\simeq \CO_{\BA^1}(-n\cdot 0)$ with $n\geq 0$,
but the action of $\BG_m$ on the fiber of
$\CP\simeq \CO_{\BA^1}(-n\cdot 0)$ at $0$ is given
by the character $-n$, which is a contradiction.

\end{proof}

Thus, to prove the theorem, we need to construct a $\BG_m$-action
as in the lemma and check its properties. Consider the $\BG_m$-action
``of the first kind" on $\fU^\mu_{G;B}$, as in the proof of
\thmref{stratification of parabolic Uhlenbeck}.
The line bundle $\CP_{\fU,B}$ and the section $s$ are
$\BG_m$-equivariant, by construction.

Let us now construct a point $y\in \fU^\mu_{G;B}$ as in the
lemma. Let us write $\mu=(\ol{\mu},a)$, and recall the projection
$$\varrho^\mu_{\fb^+_{aff}}:\fU^\mu_{G;B}\to \overset{\circ}\bC{}^\mu
\simeq \overset{\circ}\bC{}^{\ol{\mu-a\cdot \delta}}\times\overset{\circ}\bC{}^{(a)}.$$

Let $y$ be any point of $\overset{\circ}\fU{}^\mu_{G;B}$ which projects to a
multiplicity-free point in $\overset{\circ}\bC{}^\mu$. We claim
that the $\BG_m$ action will contract this point to a point
$y'$ on the subscheme
$$\fU^{\mu;a\cdot \alpha_0,0}_{G;B}\subset
\ol{\partial(\fU^\mu_{G;B})_{i_0}}\cap (\fU^\mu_{G;B}-\ol{\underset{i\in I}\cup\,
\partial(\fU^\mu_{G;B})_i}).$$

Indeed, using the horizontal factorization property,
\propref{horizontal factorization}, it suffices to analyze separately
the cases when $a=0$, and when and $\mu=\alpha_0$.
In the former case, the assertion is trivial, since in this case
$\fU^\mu_{G;B}\simeq \ol{\Bun}_B^{\ol{\mu}}(\bC,\bc)$, and the
$\BG_m$-action is trivial.

In the latter case, from \thmref{stratification of parabolic
Uhlenbeck}, we obtain that the projection $$\fU^{\alpha_0}_{G;B}\to
\on{QMaps}^{\alpha_0}(\bC,\CG_{G,B,\bX})$$ is one-to-one. From 
the remark following \lemref{two parabolics}, we obtain that
$\on{QMaps}^{\alpha_0}(\bC,\CB_{\fg_{aff}})\simeq
\overset{\circ}\bC\times \BA^1$, so that the $\BG_m$ action comes 
from the standard $\BG_m$-action on $\BA^1$.
Hence, $\fU^{\alpha_0}_{G;B}\to
\on{QMaps}^{\alpha_0}(\bC,\CB_{\fg_{aff}})$ is an isomorphism,
and the contraction statement in this case is manifest as well.

\medskip

Thus, it remains to check that for any point
$y'\in \fU^{\mu;a\cdot \alpha_0,0}_{G;B}$, the $\BG_m$-action
on the fiber of $\CP_{\fU,B}$ at $y'$ is given by a positive character.
However, by construction this action is the same as the
$\BG_m$-action on the fiber of $\CO(a\cdot 0_\bX)^{\otimes 2\cdot h^\vee}$ at
$0_\bX\in \bX$, i.e., corresponds to the $2\cdot a\cdot h^\vee$-th power of the 
standard character.

Thus, the theorem is proved.

\ssec{Remark}

Let us give another interpretation of the line bundle $\CP_{\fU,B}$
for $G=SL_n$ in terms of the scheme $\wt{\fN}{}^a_n$
classifying torsion-free sheaves on $\bS$.

For $\mu=(\ol{\mu},a)$, and $\ol{\mu}=(\mu_1,...,\mu_n)$,
consider the scheme $\wt{\fN}{}^\mu_{n,flag}$ that classifies
the data of
\begin{itemize}

\item
A torsion free sheaf of generic rank $n$, denoted $\CM$, on $\bS$,
which is locally free near $\bD_\infty$, and with $ch_2(\CM)=-a$,

\item
A trivialization $\CM|_{\bD_\infty}\simeq \CO^{\oplus n}$,

\item
A flag of locally free subsheaves
$$\CM=\CM_0\subset \CM_1\subset...\subset \CM_{n-1}\subset
\CM_n=\CM(\bD_0),$$
such that each $\CM_i/\CM_{i-1}$ is a coherent sheaf on $\bD_0$ of
generic rank $1$ and of degree $\mu_i$.
\end{itemize}

On $\wt{\fN}{}^\mu_{n,flag}$ there is a natural line bundle
$\CP_{\wt{\fN},B}$,
whose fiber at a point $\CM\subset \CM_1\subset...\subset
\CM_{n-1}\subset \CM(\bD_0)$ is the line:
$$\underset{i=1,...,n}\otimes\,
\left(\on{det}\left(R\Gamma(\on{det}(\CM_i/\CM_{i-1}))\right)\otimes
\on{det}\left(R\Gamma(\CM_i/\CM_{i-1})\right)^{-1}\right).$$

\medskip

In \cite{fgk} it was shown that there exists a natural proper map
$\wt{\fN}{}^\mu_{n,flag}\to \fU^\mu_{SL_n,B}$, which is in fact a semi-small
resolution of singularities.
One can show that the pull-back of $\CP_{\fU,B}$ to
$\wt{\fN}{}^\mu_{n,flag}$ is isomorphic to
$\CP_{\wt{\fN},B}^{\otimes 2n}$.

\bigskip

\centerline{{\bf Part IV}: {\Large Crystals}}

\bigskip

In this part, $\fg$ will be an arbitrary Kac-Moody Lie algebra, 
with the exception of \secref{via aff gr}, \secref{finite golova} 
and \secref{golova!}.

\section{Crystals in general}

\ssec{}

We refer to \cite{ks}, Section 3.1 for the definition of
crystals associated to a Kac-Moody algebra.
By a slight abuse of language, by crystal over $\fg$ we will actually
mean a crystal over the Langlands dual algebra of $\fg$, in particular,
the weight function will take values in the lattice $\Lambda$.

Crystals form a category, where morphisms from $\sB_1$ to $\sB_2$
are the maps of sets $\sB_1\to \sB_2$, which preserve the $e_i$, $f_i$
operations, and the $\epsilon_i$, $\phi_i$--functions for $i\in I$.

If $\fm$ is a Levi subalgebra in $\fg$, it corresponds
to a subdiagram of the Dynkin graph, therefore it makes
sense to talk about $\fm$-crystals.
If $\sB$ is a $\fg$-crystal, we will denote by the same
character the underlying $\fm$-crystal. Similarly, for
every $\fm$-crystal $\sB$ we define a $\fg$-crystal structure
on the same underlying set, by leaving $e_i,f_i,\epsilon_i,\phi_i$ for
$i\in \fm$ unchanged, and declaring that for $j\notin \fm$,
$e_j,f_j$ send everything to $0$, and $\epsilon_j=\phi_j=-\infty$.

\ssec{}

Let $\sB_\fg$ be a crystal. Our goal in this section is to
review some extra structures and properties of $\sB_\fg$
which allow to identify it with the standard crystal $\sB_{\fg}^\infty$
that parametrizes the canonical basis in $U(\check{\fn})$,
cf. \cite{ks}, Section 3.2.

For $i\in I$, let $\fm_i$ be the corresponding subminimal
Levi subalgebra. Let $\sB_i$ be the $\fg$-crystal obtained
in the above way from the ``standard'' $\fsl_2$-crystal. I.e., as
a set it is $\ZZ^{\geq 0}$ and consists of elements that we will
denote $\sfb_i(n),\, n\geq 0$ with $wt_i(\sfb_i(n))=n\cdot \check\alpha_i$,
$f_i(\sfb_i(n))=\sfb_i(n-1)$, $\phi_i(\sfb_i(n))=n$, and all the
operations for $j\neq i$ are trivial.

\medskip

According to Proposition 3.2.3 of \cite{ks}, we need to understand
what kind of structure on $\sB_\fg$ allows to construct maps
$\Psi_i:\sB_\fg\to \sB_\fg\otimes \sB_i$, satisfying the conditions
(1)--(7) of {\it loc.cit}.

Our crystal $\sB_\fg$ will automatically satisfy conditions (1)--(4),
and it will be $\phi$-normal, i.e.
$$\text{for  } \sfb\in\sB_\fg,
\,\phi_i(\sfb)=\max\{n\,|\,f_i^n(\sfb)\neq 0\}.$$
Then the $\epsilon$-functions are defined by the rule
$$\phi_i(\sfb)=\epsilon_i(\sfb)+\langle wt(\sfb),\check\alpha_i\rangle.$$

\ssec{}

We suppose that the set $\sB_\fg$ carries {\it another} crystal structure
(also assumed $\phi$-normal), given by operations
$e^*_i$ and $f^*_i$, $i\in I$ with the same weight function
$wt$. Assume also the following:

\medskip

\noindent{(a)}
The operations $e_i$, $f_i$ commute with
$e^*_j$ and $f^*_j$, whenever $i\neq j$.

\smallskip

\noindent{(b)}
For every vertex $i\in I$, we have a set decomposition
\begin{equation} \label{decomposition}
\sB_\fg=\underset{n\in \NN}\cup\,
\sC_{\fm_i}^n \times \sB^n_{\fm_i}\times \sB_i,
\end{equation}
where

\begin{itemize}

\item
$\sC_{\fm_i}^n$ is an abstract set, $\sB^n_{\fm_i}$ is the set underlying
the finite $\fsl_2$-crystal of highest weight $n$, and $\sB_i$ is
as above.

\item
This decomposition respects both the $e_i$, $f_i$ and the
$e^*_i$, $f^*_i$-operations.

\item
The $e_i$, $f_i$--action on each
$\sC_{\fm_i}^n\times \sB^n_{\fm_i}\times \sB_i$
is the trivial action on $\sC_{\fm_i}^n$, times the action on
$\sB^n_{\fm_i}\times \sB_i$ as on the tensor product of crystals.

\item
The $e^*_i$, $f^*_i$--action on each
$\sC_{\fm_i}^n\times \sB^n_{\fm_i}\times \sB_i$ is trivial on the first
two factors times the standard action on $\sB_i$.

\end{itemize}

\ssec{}

We claim that under the above circumstances, we do have canonical
maps $\Psi_i:\sB_\fg\to \sB_\fg\times \sB_i$, satisfying condition (6)
of \cite{ks}, Proposition 3.2.3.

Indeed, such maps obviously exist for $\fg=\fsl_2$, i.e. for
each $i$ we have a map of $\fm_i$-crystals
$$\Psi^{can}_i:\sB_i\to \sB_i\otimes \sB_i.$$

We define the map $\Psi_i$ on $\sB_\fg$ in terms of the decomposition
\eqref{decomposition}. Namely, on each
$\sC_{\fm_i}^n\times \sB^n_{\fm_i}\times \sB_i\subset \sB_\fg$,
$\Psi_i$ is the identity map on the first two factors times $\Psi^{can}_i$
on the third one.

In other words, each $\sfb\in \sB_\fg$ can be uniquely written
in the form $(e_i^*)^n(\sfb')$, where $f_i^*(\sfb')=0$, and
\begin{equation} \label{formula for Psi}
\Psi_i(\sfb)=\sfb'\otimes \sfb_i(n).
\end{equation}

\medskip

We claim that $\Psi_i$ indeed respects the operations $e_j$ and $f_j$.
First, when $j\neq i$, the assertion follows immediately from
\eqref{formula for Psi}, since $f^*_i$, $e^*_i$ commute with
$e_j$ and $f_j$.

For $j=i$ the assertion is also clear, since we are dealing with
the map of $\fm_i$-crystals
$$\sB^n_{\fm_i}\otimes \sB_i\overset{\on{id}\otimes \Psi_i^{can}}
\longrightarrow \sB^n_{\fm_i}\otimes \sB_i\otimes \sB_i.$$

\ssec{} \label{highest weight}

Finally, in order to apply the uniqueness theorem of
\cite{ks}, the crystal $\sB_{\fg}$ must satisfy
the ``highest weight property'' (i.e., condition
(7) of Proposition 3.2.3 of {\it loc.cit}. It can be stated in
either of the two ways:

\noindent For $\sfb\in \sB_\fg$ different from the canonical highest weight
vector $\sfb(0)$

\begin{itemize}

\item
there exists $i\in I$ such that $f_i(\sfb)\neq 0$

\item
there exists $i\in I$ such that $f^*_i(\sfb)\neq 0$.

\end{itemize}

\medskip

Let us denote by $\sB_{\fg}^{\fm_i}$ the subset of $\sB_\fg$ consisting
of elements annihilated by $f_i^*$. We claim that it
also has a natural crystal structure.
Indeed, the operations $e_j$ and $f_j$ obtained by restriction
from $\sB_{\fg}$ preserve $\sB_\fg^{\fm_i}$, since they commute with $e^*_i$.

The operations $e_i$ and $f_i$ are defined in terms of the decomposition
\eqref{decomposition}. Namely,
$$\sB_\fg^{\fm_i}\cap (\sC^n_{\fm_i}\times \sB^n_{\fm_i}\times \sB_i)=
(\sC^n_{\fm_i}\times \sB^n_{\fm_i})\times \sfb_i(0),$$
and we set $e_i$, $f_i$ to act along the $\sB^n_{\fm_i}$--factor.

\medskip

Moreover, note that the map $\Psi_i$ constructed above maps
$\sB_\fg$ to $\sB_{\fg}^{\fm_i}\otimes \sB_i$, and we claim that this
is a map of crystals with respect to the crystal structure
on $\sB_\fg^{\fm_i}$ introduced above.  Moreover, one easily checks
that this map is an isomorphism.

\ssec{}

Thus, if $\sB_\fg$ satisfies the highest weight property,
and admits decompositions as in \eqref{decomposition} for every $i$,
it can be identified with the standard crystal $\sB_\fg^\infty$
of \cite{ks}. In what follows, we will need several additional
properties of $\sB^\infty_{\fg}$.

Let $\fm$ be an arbitrary Levi subalgebra of $\fg$. (In our applications
we will only consider Levi subalgebras that correspond to subdiagrams
of finite type.)
The following result is a generalization of what we have
constructed for $\fm=\fm_i$, cf. \cite{k4},\cite{k5}:

\begin{thm} \label{complete decomposition}
For every Levi subalgebra $\fm$ there exists an isomorphism of
$\fg$-crystals (with respect to the $e_i,f_i,\epsilon_i,\phi_i, i\in I$)
$$\Psi_{\fm}:\sB^\infty_{\fg}\simeq \sB^{\fm}_{\fg}\otimes \sB^\infty_{\fm},$$
(here $\sB^\infty_{\fm}$ is viewed as a $\fg$-crystal), such that

\smallskip

\noindent{\em (a)}
The operations $e^*_i,f^*_i$ for $i\in \fm$ on $\sB^\infty_{\fg}$ go over
to the $e^*_i,f^*_i$ operations along the second factor.

\smallskip

\noindent{\em (b)}
The $\fm$-crystal structure on $\sB^{\fm}_{\fg}$ is normal.

\end{thm}

In particular, from this theorem we obtain that, as a set,
$\sB^{\fm}_{\fg}$ can be realized as a subset of $\sB^\infty_{\fg}$
consisting of elements annihilated by the $f^*_i$, $i\in \fm$.
In terms of this set-theoretic embedding
$\sB_{\fg}^{\fm}\to \sB^\infty_{\fg}$, the operations
$e_j,f_j,\epsilon_j,\phi_j$ for $j\notin \fm$
are induced from those on $\sB^\infty_{\fg}$.
For $i\in\fm$, the $e_i$ and $f_i$ operators on an element
$\sfb\in \sB_{\fg}^{\fm}$ can be explicitly described as follows:

The element $f_i(\sfb)$ (for the action that comes from 
$\sB_{\fg}^{\fm}$-crystal structure) equals $f_i(\sfb)$
(for the action that is induced by 
the $\sB^\infty_\fg$-crystal structure). The element 
$e_i(\sfb)$ (for the action coming from the $\sB_{\fg}^{\fm}$-crystal structure)
equals $e_i(\sfb)$ (for the action induced by the $\sB^\infty_{\fg}$-crystal structure)
if the latter belongs to $\sB_{\fg}^{\fm}$, and $0$ otherwise.

In addition, we have the following assertion, cf. \cite{k3}, Theorem 3:

\begin{thm}   \label{canonical basis}
The set $\sB^{\fm}_{\fg}$ parametrizes the canonical basis
for $U(\fn(\check{\fp}))$, where
$\check{\fp}$ is the corresponding parabolic
subalgebra, and $\fn(\check{\fp})$ is its unipotent radical. Moreover,
as a $\fm$-crystal, $\sB^{\fm}_{\fg}$ splits as a disjoint union
$$\sB^{\fm}_{\fg}=\underset{\nu}\cup\, \sC_{\fm}^\nu\otimes \sB^\nu_\fm,$$
where $\nu$ runs over the set of dominant integral weights of
$\fm$, $\sB^\nu_\fm$ is the crystal associated to the integrable
$\check{\fm}$-module $V^\nu_M$, and $\sC_{\fm}^\nu$ parametrizes a basis in
$\Hom_{\check{\fm}}(V_M^\nu,U(\fn(\check{\fp})))$.
\end{thm}

Note that the set $\sC^\nu_\fm$ of the above theorem embeds into
$\sB^\infty_{\fg}$ as the set of those $\sfb\in \sB^\infty_{\fg}(\nu)$,
which are annihilated by $f_i$ and $f_i^*$ for $i\in \fm$.

\section{Crystals via the affine Grassmannian}  \label{via aff gr}

\ssec{}

In this section we will take $\fg$ to be finite-dimensional, and 
give a construction of a crystal $\sB_\fg$ using the affine
Grassmannian $\Gr_G$ of $G$, in the spirit of \cite{bg}. We will
assume that the reader is familiar with the notation of {\it loc.cit.}

In the next section we will generalize this construction for an arbitrary 
Kac-Moody algebra $\fg$, using the space $\on{QMaps}(\bC,\CG_{\fg,\fb})$ 
instead of the (non-existing in the general case) affine Grassmannian.

\medskip

Recall that for a standard parabolic $P$, we denote by 
$N(P)$ (resp., $M$, $P^-$, $N(P^-)$) its unipotent radical
(resp., the Levi subgroup, the corresponding opposite parabolic, etc.)
By $B(M)$ (resp., $N(M)$, $B^-(M)$, $N^-(M)$) will denote the
standard Borel subgroup in $M$ (resp., unipotent radical of $B(M)$, etc.)
Finally, recall that $T$ denotes the Cartan subgroup of $G$.

\medskip

Let $\bC$ be a (smooth, but not necessarily complete) curve, and
$\bc\in \bC$ a point. We will choose a local parameter on $\bC$
at $\bc$ and call it $t$. Recall that if $H$ is an algebraic group,
we can consider the group-scheme $H[[t]]$, the group-indscheme 
$H((t))$ and the ind-scheme of ind-finite type
$\Gr_H:=H((t))/H[[t]]$, called the affine Grassmannian. 
In what follows we will consider the affine Grassmannians corresponding
to the groups $G$, $M$, $P$, $N(P)$, etc.

For a coweight $\lambda$, we have a canonical point in $T((t))$, 
denoted $t^\lambda$. We will denote by the same character its
image in $\Gr_T$ and $\Gr_G$ via the embedding of $\Gr_T\to \Gr_G$.

We have the maps $i_P:\Gr_P\to \Gr_G$,
$i_{P^-}:\Gr_{P^-}\to \Gr_G$, which are locally closed embeddings
on every connected component, and projections $\q_P:\Gr_P\to\Gr_M$,
$\q_{P^-}:\Gr_{P^-}\to\Gr_M$. The projections $\q_P$, $\q_{P^-}$ induce
a bijection on the set of connected components.
For a coweight $\lambda\in \Lambda$ we will denote by $\Gr_T^\lambda$
(resp., $\Gr_B^\lambda$, $\Gr_{B^-}^\lambda$) the corresponding connected
component of $\Gr_T$ (resp., $\Gr_B$, $\Gr_{B^-}$).

\ssec{}

It is known that the intersection
$\Gr^{\lambda_1}_B\cap \Gr^{\lambda_2}_{B^-}$ is of
pure dimension  $\langle \lambda_1-\lambda_2,\check\rho\rangle$.
Let us denote by $\sB_\fg(\lambda)$ the set of irreducible components
of the above intersection with parameters $\lambda_1-\lambda_2=\lambda$.
(It is easy to see that the action of $t^{\lambda'}\in T((t))$
identifies $\Gr^{\lambda_1}_B\cap \Gr^{\lambda_2}_{B^-}$ with
$\Gr^{\lambda_1+\lambda'}_B\cap \Gr^{\lambda_2+\lambda'}_{B^-}$.)

Our goal now is to show that the set
$\sB_{\fg}:=
\underset{\lambda}\cup\, \sB_{\fg}(\lambda)$ has a natural structure of
$\fg$-crystal.
We need to define the weight function $wt:\sB_{\fg}\to \Lambda$;
the functions $\epsilon_i,\phi_i:\sB_{\fg}\to \ZZ$ for $i\in I$, and
the operations $e_i,f_i:\sB_{\fg}\to \sB_{\fg}\cup 0$.

\medskip

Of course, the function $wt$ is defined so that
$wt(\sB_{\fg}(\lambda))=\lambda$. The function $\phi_i$ is
defined in terms of $f_i$ by the normality condition:
$$\text{For  } \sfb\in\sB_{\fg},
\,\phi_i(\sfb)=max\{n\,|\,f_i^n(\sfb)\neq 0\}.$$
The function $\epsilon_i$ will be defined by the rule
$$\epsilon_i(\sfb)=\phi_i(\sfb)-\langle\check\alpha_i,wt(\sfb)\rangle.$$

Thus, we have to define the operations $e_i$ and $f_i$.

\ssec{}   \label{e_i}

First, let us assume that $G$ is of semi-simple rank $1$.
Then we can identify the coroot lattice with $\ZZ$, and it is easy
to see that $\Gr_B^\lambda\cap \Gr_{B^-}^0$ is non-empty if and only
if $\lambda=n\cdot \alpha$ (where $\alpha$ is the unique positive 
coroot), with $n\geq 0$.
In the latter case, this intersection is an irreducible variety, isomorphic
to $\AA^n-\AA^{n-1}$. Thus,
$\sB_{\fg}$ in this case is naturally $\ZZ^{\geq 0}$, and the definition of
$e_i$ and $f_i$ is evident: they are the raising and the lowering
operators, respectively.

In the case of a general $G$, the operations will be defined
by reduction to the rank $1$ case.

\medskip

For a parabolic $P$ and $\lambda,\mu\in \Lambda$, consider the
intersection
$$\q_P^{-1}(\Gr^{\mu}_{B^-(M)})\cap \Gr^\lambda_{B^-}\subset \Gr_G,$$
projecting by means of $\q_P$ on $\Gr^{\mu}_{B^-(M)}$.

\begin{lem}  \label{product}
The above intersection splits as a direct product
$$\Gr^{\mu}_{B^-(M)}\times \bigl((\q_P)^{-1}(g)\cap \Gr^\lambda_{B^-}\bigr)$$
for any $g\in \Gr^{\mu}_{B^-(M)}$.
\end{lem}

\begin{proof}

The group $N^-(M)((t))$ acts transitively on $\Gr^{\mu}_{B^-(M)}$,
and it also acts on $\Gr_G$ preserving the intersection
$(\q_P)^{-1}(\Gr^{\mu}_{B^-(M)})\cap \Gr^\lambda_{B^-}$.

Moreover, for $g'\in \Gr_P\subset \Gr_G$
and $g=\q_P(g')$, the inclusion
$$\on{Stab}_{N^-(M)((t))}(g',\Gr_G)\subset
\on{Stab}_{N^-(M)((t))}(g,\Gr_M)$$
is an equality, since $N^-(M)\cap N(P)=1$.

Therefore, for any $g$ as above, the action of
$N^-(M)((t))/\on{Stab}_{N^-(M)((t))}(g,\Gr_M)$ on the fiber
$(\q_P)^{-1}(g)\cap \Gr^\lambda_{B^-}$ defines an isomorphism
$$(\q_P)^{-1}(\Gr^{\mu}_{B^-(M)})\cap \Gr^\lambda_{B^-}\simeq
\Gr^{\mu}_{B^-(M)}\times \bigl((\q_P)^{-1}(g)\cap \Gr^\lambda_{B^-}\bigr).$$

\end{proof}

We will use this lemma as follows.
For $\lambda_1,\lambda_2,\mu\in \Lambda$, let us note that the
intersection
$$\Gr_B^{\lambda_1}\cap
(\q_P)^{-1}(\Gr^{\mu}_{B^-(M)})\cap \Gr^{\lambda_2}_{B^-}$$ is the same
as
$$(\q_P)^{-1}(\Gr_{B(M)}^{\lambda_1}\cap \Gr^{\mu}_{B^-(M)})\cap
\Gr^{\lambda_2}_{B^-}.$$

In particular, by putting $\lambda_1=\mu$ we obtain that
the variety $(\q_P)^{-1}(g)\cap \Gr^\lambda_{B^-}$
of \lemref{product} is finite dimensional of
$\on{dim}\leq \langle\mu-\lambda,\check\rho\rangle$.

Of course, the variety
$(\q_P)^{-1}(\Gr^{\mu}_{B^-(M)})\cap \Gr^{\lambda}_{B^-}$ depends
up to isomorphism only on the difference $\mu':=\mu-\lambda$.
Let us denote by $\sB^{\fm,*}_{\fg}(\mu')$ the set of
irreducible components of $(\q_P)^{-1}(g)\cap \Gr^\lambda_{B^-}$
of (the maximal possible) dimension $\langle\mu',\check\rho\rangle$.

Since the group $N^-(B)((t))$ is ind-pro-unipotent, we obtain that the
stabilizer $$\on{Stab}_{N^-(M)((t))}(g,\Gr_M)$$ appearing in
the proof of \lemref{product} is connected.
Therefore, the set is $\sB^{\fm,*}_{\fg}(\mu')$
is well-defined, i.e. is independent of the choice of the point $g$.

\medskip

Clearly, for every irreducible component $\bK$ of
$\Gr_B^{\lambda_1}\cap\Gr^{\lambda_2}_{B^-}$ there exists a
unique $\mu\in \Lambda$ such that the intersection
$\bK\cap (\q_P)^{-1}(\Gr^{\mu}_{B^-(M)})$ is dense in $\bK$.

Using \lemref{product}, we obtain that the set $\sB_{\fg}(\lambda)$
can be canonically decomposed as a union
\begin{equation}  \label{decomposition with respect to parabolic}
\sB_{\fg}(\lambda)=\underset{\mu}\cup \, \sB^{\fm,*}_{\fg}(\mu)\times
\sB_{\fm}(\lambda-\mu).
\end{equation}

\medskip

Finally, we are able to define the operations $e_i$ and $f_i$.
For each $i\in I$ take $P$ to be the corresponding subminimal
parabolic $P_i$ and consider the decomposition
$$\sB_{\fg}(\lambda)=\underset{\mu}\cup\, \sB^{\fm_i,*}_{\fg}(\mu)\times
\sB_{\fm_i}(\lambda-\mu).$$
For an element $\sfb_{\fg}\in \sB_{\fg}(\lambda)$ of the form
$$\sfb'\times \sfb_{\fm_i},\, \sfb'\in \sB^{\fm_i,*}_{\fg}(\mu),\,
\sfb_{\fm_i}\in \sB_{\fm_i}(\lambda-\mu),$$
we set $e_i(\sfb_G)$ to be
$$\sfb'\times e_i(\sfb_{\fm_i})\in \sB^{\fm_i,*}_{\fg}(\mu)\times
\sB_{\fm_i}(\lambda+\alpha_i-\mu)\cup 0\subset
\sB_{\fg}(\lambda+\alpha_i)\cup 0,$$
and similarly for $f_i$.

This is well-defined, because $M_i$ is of semi-simple rank $1$,
and we know how the operators $e_i$ and $f_i$ act on $\sB_{\fm_i}$.

\ssec{}  \label{e_i*}

Let us now define the operations $e^*_i$ and $f^*_i$.
First, when $G$ is of semi-simple rank $1$, we set $e_i^*=e_i$ and
$f_i^*=f_i$. To treat the general case, we simply interchange
the roles of the projections $\q_P$ and $\q_{P^-}$.

In more detail, for any $\lambda$ and $\mu$
as before, let us consider the intersection
$$\Gr_B^{\lambda}\cap(\q_{P^-})^{-1}(\Gr^{\mu}_{B(M)}).$$ As
in \lemref{product}, this intersection can be represented
as a product
$\Gr^{\mu}_{B(M)}\times \bigl((\q_{P^-})^{-1}(g)\cap \Gr_B^{\lambda}\bigr)$
for any $g\in \Gr^{\mu}_{B(M)}$. For $\mu'=\lambda-\mu$,
let us denote by $\sB_{\fg}^{\fm}(\mu')$ the set of
irreducible components of the top dimension of
$(\q_{P^-})^{-1}(g)\cap \Gr_B^{\lambda}$.

By looking at the intersections of irreducible components
of $\Gr_B^{\lambda_1}\cap \Gr_{B^-}^{\lambda_2}$ with the various
$(\q_{P^-})^{-1}(\Gr_{B(M)}^{\mu})$, we obtain a decomposition

\begin{equation}  \label{star decomposition}
\sB_{\fg}(\lambda)=\underset{\mu}\cup\,
\sB_{\fg}^{\fm}(\mu)\times \sB_{\fm}(\lambda-\mu).
\end{equation}

By taking $M=M_i$, for $\sfb_{\fg}\in \sB_{\fg}(\lambda)$ of the form
$$\sfb'\times \sfb_{\fm_i},\,
\sfb'\in \sB_{\fg}^{\fm_i}(\mu),\,
\sfb_{\fm_i}\in \sB_{\fm_i}(\lambda-\mu),$$
we set $e^*_i(\sfb_{\fg})$ to be
$$\sfb'\times e^*_i(\sfb_{\fm_i})\in \sB_{\fg}^{\fm_i}(\mu)\times
\sB_{\fm_i}(\lambda+\alpha_i-\mu) \cup 0
\subset \sB_{\fg}(\lambda+\alpha_i)\cup 0,$$
and similarly for $f_i^*$.

\medskip

\noindent{\it Remark.}
Let us explain the consistence of the notation $\sB_{\fg}^{\fm}(\mu)$.
As we shall see later, the crystals $\sB_{\fg}$ (and in particular
$\sB_{\fm}$) have the highest weight property. If we assume that,
we will obtain that the set $\sB_{\fg}^{\fm}(\mu)$ is precisely
the subset of $\sB_{\fg}(\mu)$ consisting of elements annihilated
by $f^*_i$, $i\in \fm$. Indeed, the latter elements correspond
to the irreducible components of
$$\Gr_B^{0}\cap (\q_{P^-})^{-1}(\Gr_{B(M)}^{\mu}\cap \Gr_{B^-(M)}^{\mu})$$
of dimension $\langle \mu,\check\rho\rangle$, which is the same
as the set $\sB_{\fg}^{\fm}(\mu)$.

\medskip

\begin{prop}
For $i\neq j$, the operations $e_i,f_i$ commute with the operations
$e_j^*$, $f_j^*$.
\end{prop}

\begin{proof}

Consider the intersection
$(\q_{P_i})^{-1}(\Gr_{B^-(M_i)}^{\mu_1})\cap
(\q_{P^-_j})^{-1}(\Gr_{B(M_j)}^{\mu_2})\subset \Gr_G$.
We claim that it decomposes as
$$(\Gr_{B^-(M_i)}^{\mu_1}\times \Gr_{B(M_j)}^{\mu_2})\times
\bigl((\q_{P_i})^{-1}(g_1)\cap (\q_{P^-_j})^{-1}(g_2)\bigr)$$ for any
$g_1\times g_2\in \Gr_{B^-(M_i)}^{\mu_1}\times \Gr_{B(M_j)}^{\mu_2}$.

Indeed, as in the proof of \lemref{product}, the group
$N^-(M_i)((t))\times N(M_j)((t))$ acts transitively on
the base $\Gr_{B^-(M_i)}^{\mu_1}\times \Gr_{B(M_j)}^{\mu_2}$, and
since the subgroups $N^-(M_i),N(M_j)\subset G$ commute with one another,
we obtain the desired decomposition as in \lemref{product}.

For $\mu=\mu_1-\mu_2$, let us denote by $\sC_{i,j}(\mu)$
the set of irreducible components of the top
($=\langle \mu,\check\rho\rangle$) dimension of
$(q_{P_i})^{-1}(g_1)\cap (q_{P^-_j})^{-1}(g_2)$.

\medskip

By intersecting $\Gr_B^{\lambda_1}\cap \Gr_{B^-}^{\lambda_2}$ with
$(\q_{P_i})^{-1}(\Gr_{B^-(M_i)}^{\mu_1})\cap
(\q_{P^-_j})^{-1}(\Gr_{B(M_j)}^{\mu_2})$ for the various $\mu_1$ and $\mu_2$,
we obtain a decomposition
$$\sB_{\fg}(\lambda)\simeq
\underset{\mu_1,\mu_2}\cup\sC_{i,j}(\mu_1-\mu_2)\times
\sB_{\fm_i}(\lambda_1-\mu_1)\times \sB_{\fm_j}(\mu_2-\lambda_2).$$

This decomposition is a refinement of \eqref{decomposition with
respect to parabolic} and \eqref{star decomposition}.
Therefore, for an element $\sfb_G\in \sB_{\fg}(\lambda)$ of the form
$\bc\times \sfb_i\times \sfb_j$ with
$$\sfc\in\sC_{i,j}(\mu_1-\mu_2),\,
\sfb_i\in \sB_{\fm_i}(\lambda-\mu_1),\,
\sfb_j\in \sB_{\fm_j}(\lambda-\mu_2),$$
the $e_i$ and $f_i$ operations
act via their action on $\sfb_i$, and the $e_j^*$ and $f_j^*$ act
via $\sfb_j$, respectively.
This makes the required commutativity property manifest.

\end{proof}

\ssec{}   \label{convolution}

Finally, we are going to define the decompositions of
$\sB_{\fg}$ as in \eqref{decomposition}. In fact, we will
define the decompositions as in \thmref{complete decomposition}
for any Levi subalgebra $\fm$.
First, we need to discuss certain $\fm$-crystals associated with
the convolution diagram of $M$.

Let us recall some notation related to the affine Grassmannian
$\Gr_M$. The group-scheme $M[[t]]$ (being a subgroup of the
group-indscheme $M((t))$) acts naturally on $\Gr_M$, and its
orbits are parametrized by the set $\Lambda_\fm^+$
of dominant coweights of $M$; for $\nu\in \Lambda_\fm^+$ we
will denote the corresponding orbit by $\Gr_M^\nu$.

Let $\Conv_M$ denote the convolution diagram for $M$. It is by definition
the ind-scheme parametrizing the data of
$(\F_M,\F'_M,\beta,\wt{\beta}{}')$, where $\F_M,\F'_M$ are principal 
$M$-bundles on $\bC$, $\beta$ is a trivialization of $\F_M$ on
$\bC-\bc$, and $\wt{\beta}^1)$ is an isomorphism 
$\F_M|_{\bC-\bc}\to \F'_M|_{\bC-\bc}$. We will think of $\Conv_M$
as a fibration over $\Gr_M$ with the typical fiber $\Gr_M$; more 
precisely
$$\Conv_M\simeq M((t))\underset{M[[t]]}\times \Gr_M.$$
Sometimes, we will write $\Conv_M\simeq \Gr_M\star \Gr_M$, so that
the first factor is perceived as a base, and the second one as the
fiber. In particular, if $\Y\subset \Gr_M$ is a sub-scheme, and
$\Y'\subset \Gr_M$ is a $M[[t]]$-invariant subscheme, it makes sense
to consider the subscheme $\Y\star \Y'\subset \Gr_M\star 
\Gr_M\simeq \Conv_M$. For $\nu\in \Lambda_\fm^+$, let 
$\Conv^\nu_M\subset \Conv_M$ be equal to $\Gr_M\star \Gr^\nu_M$.

Note that for a point $(\F_M,\F'_M,\beta,\wt{\beta}{}')\in \Conv_M$,
by taking the composition $\beta^1:=\wt{\beta}^1\circ \beta$ we obtain a
trivialization of $\F'_M$ on $\bC-\bc$, and thus a map 
$\Conv_M\to \Gr_M\times \Gr_M$, which is easily seen to be an isomorphism.
We will denote by $p$ and $p'$ the two projections 
$\Conv_M\to \Gr_M$, which remember the data of $(\F_M,\beta)$ and 
$(\F'_M,\beta')$, respectively. Note that the locally closed subsets
$\Conv^\nu_M$ introduced above are exactly the orbits of the diagonal
$M((t))$ action on $\Conv_M\simeq \Gr_M\times \Gr_M$.

\ssec{}

For another pair of coweights $\lambda_1$ and $\lambda_2$, consider
the intersection
$$(\Gr^{\lambda_1}_{B(M)}\times \Gr^{\lambda_2}_{B^-(M)})^\nu:=
(\Gr^{\lambda_1}_{B(M)}\times \Gr^{\lambda_2}_{B^-(M)})\cap
\Conv_M^\nu.$$

Again, up to isomorphism, the above variety depends only on the difference
$\lambda=\lambda_1-\lambda_2$. Let us denote by
$\sD_{\fm}^\nu(\lambda)$ the set of irreducible components
of the above intersection of the (maximal possible) dimension
$\langle \nu+\lambda,\check\rho\rangle$.

We will define on the set $\sD_{\fm}^\nu:=\underset{\lambda}\cup\,
\sD_{\fm}^\nu(\lambda)$ two structures of a ($\phi$-normal)
$\fm$-crystal.

\medskip

First, let us define the operations $e_i$, $f_i$. For that let us fix
$\lambda_1,\lambda_2\in \Lambda$ and consider the ind-scheme
$$(\Gr_{B^-(M)}^{\lambda_1}\times \Gr^{\lambda_2}_{B^-(M)})^\nu:=
(\Gr_{B^-(M)}^{\lambda_1}\times \Gr^{\lambda_2}_{B^-(M)})\cap
\Conv_M^\nu$$
projecting by means of $p$ onto $\Gr_{B^-(M)}^{\lambda_1}$.
As in \lemref{product}, we obtain that there is an isomorphism
$$(\Gr_{B^-(M)}^{\lambda_1}\times \Gr^{\lambda_2}_{B^-(M)})^\nu\simeq
\Gr_{B^-(M)}^{\lambda_1}\times
\bigl((g\times \Gr^{\lambda_2}_{B^-(M)})\cap
\Conv_M^\nu\bigr)$$ for any $g\in
\Gr_{B^-(M)}^{\lambda_1}$.

For $\lambda=\lambda_1-\lambda_2$ and $g\in \Gr_{B^-(M)}^{\lambda_1}$,
let $\sB_{\fm}^{\nu,*}(\lambda)$ denote the set of irreducible
components of the intersection
$(g\times \Gr^{\lambda_2}_{B^-(M)})\cap
\Conv_M^\nu$
of (the maximal possible) dimension $\langle\nu+\lambda,\check\rho\rangle$.
(As in \lemref{product}, this set does not depend on the choice
of $g\in \Gr_M$).

Thus, we have a set-theoretic decomposition
$$\sD^\nu_{\fm}(\lambda)\simeq \underset{\lambda'}\cup\,
\sB_{\fm}^{\nu,*}(\lambda')\times \sB_{\fm}(\lambda-\lambda').$$

By taking $M=M_i$, we define the operations $e_i$ and $f_i$ to act along the
second multiple of this decomposition. It is easy to see that this
definition agrees with the one of \secref{e_i}.

\medskip

The operations $e^*_i$, $f^*_i$ are defined in a similar fashion.
We fix $\lambda_1$ and $\lambda_2$ and consider the ind-scheme
$$(\Gr_{B(M)}^{\lambda_1}\times \Gr^{\lambda_2}_{B(M)})^\nu:=
(\Gr_{B(M)}^{\lambda_1}\times \Gr^{\lambda_2}_{B(M)})\cap
\Conv_M^\nu$$ projecting by means
of $p'$ onto $\Gr^{\lambda_2}_{B(M)}$. It also splits as a product
$$(\Gr_{B(M)}^{\lambda_1}\times \Gr^{\lambda_2}_{B(M)})^\nu\simeq
\Gr_{B(M)}^{\lambda_2}\times
\bigl((\Gr^{\lambda_1}_{B(M)}\times g)\cap
\Conv_M^\nu\bigr)$$ for any $g\in
\Gr_{B(M)}^{\lambda_2}$.

Let $\sB_{\fm}^{\nu}(\lambda)$ denote the set of irreducible
components of (the maximal possible)
dimension $\langle \nu+\lambda,\check\rho\rangle$
of $(\Gr^{\lambda_1}_{B(M)}\times g)\cap (\Gr_M\times\Gr_M)^\nu$,
and as before, we obtain a set-theoretic decomposition
\begin{equation}
\sD^\nu_{\fm}(\lambda)\simeq \underset{\lambda'}\cup\,
\sB_{\fm}^{\nu}(\lambda')\times \sB_{\fm}(\lambda-\lambda').
\end{equation}

We define the operations $e^*_i$, $f^*_i$ to act along the second
multiple of this decomposition for $M=M_i$, and again this definition
is easily seen to coincide with the one of \secref{e_i*}.

\medskip

Note that by taking $g$ to be the unit element in $\Gr_M$,
we obtain that $\sB_{\fm}^{\nu}(\lambda)$ are exactly the elements
of the crystal associated to the integrable $\check M$-module
with highest weight $\nu$ of weight $\lambda$. Set
$\sB_{\fm}^\nu:=\underset{\lambda}\cup\, \sB_{\fm}^\nu(\lambda)$,
and recall from \cite{bg}, that this set has a canonical structure
of a (normal) $\fm$-crystal.

The following is a generalization of Theorem 3.2 of \cite{bg},
when one of the finite crystals is replaced by $\sB_{\fm}$:

\begin{prop}  \label{BG}
In terms of the above set-theoretic decomposition
$\sD^\nu_{\fm}\simeq \sB_{\fm}^\nu\times \sB_{\fm}$,
the crystal structure on $\sD^\nu_{\fm}$, given by $e_i$ and $f_i$,
corresponds to the tensor product crystal structure on the right-hand
side.
\end{prop}

We omit the proof, since the argument is completely analogous
to the corresponding proof in {\it loc.cit.}

\ssec{}

Finally, we are ready to define the decompositions of
\eqref{decomposition}.

Consider the projection
$$(\q_P\times \q_{P^-}):\Gr_P\underset{\Gr_G}\times \Gr_{P^-}\to
\Gr_M\times \Gr_M\simeq \Conv_M.$$

This projection is $M((t))$-equivariant. For a fixed $M$-dominant
coweight $\nu$ let us denote by
$(\Gr_P\underset{\Gr_G}\times \Gr_{P^-})^\nu$ the preimage of
$\Conv_M^\nu$ under this map. For a point
$g\in \Conv_M^\nu$ consider the scheme
$(\q_P\times \q_{P^-})^{-1}(g)\subset
(\Gr_P\underset{\Gr_G}\times \Gr_{P^-})^\nu$. The set of its
irreducible components of (the maximal possible) dimension
$\langle \nu, \check\rho\rangle$ will be denoted $\sC_{\fm}^\nu$.

Thus, we obtain a decomposition
$$\sB_{\fg}=\underset{\nu}\cup\, \sC_{\fm}^\nu\times \sD_{\fm}^\nu.$$

Moreover, by unraveling the construction of the operations, we obtain
that the above decomposition respects the action of both $e_i,f_i$, and
$e^*_i,f^*_i$, such that all of the operations act via the second
multiple (i.e. $\sD_{\fm}^\nu$) in the way specified above.

This defines the required decomposition in view of \propref{BG}.

\medskip

Finally, by combining what we said above with \thmref{canonical basis},
(and assuming the highest weight property of the crystal $\sB_\fg$),
we obtain that the set $\sC_{\fm}^\nu$, which enumerates
irreducible components of $(\q_P\times \q_{P^-})^{-1}(g)$ for
$g\in \Conv_M^\nu$, parametrizes
a basis of $\Hom_{M}(V_M^\nu,U(\fn(\check \fp)))$.

\section{Crystals via quasi-maps' spaces}

\ssec{} 

In this section $\fg$ will be again an arbitrary Kac-Moody algebra. Recall
that $\CG_{\fg,\fb}$ denotes Kashiwara's flag scheme. We will fix
$\bC$ to be an arbitrary smooth (but not necessarily complete) curve, 
with a distinguished point $\bc\in \bC$. As in the previous section, we 
will denote by $t$ a local coordinate on $\bC$ defined near $\bc$.

Recall that for $\lambda\in \Lambda^{pos}_\fg$
in \secref{Zastava spaces} we introduced the scheme $\CZ^\lambda_{\fg,\fb}(\bC)$ and its
open subset $\oZ^\lambda_{\fg,\fb}(\bC)$. We showed in \propref{Zastava and quasi-maps}
that when $\bC=\BA^1$, then
$\CZ^\lambda_{\fg,\fb}(\bC)$ can be identified with the scheme of based quasi-maps
$\on{QMaps}^\lambda(\BP^1,\CG_{\fg,\fb})$, such that $\CZ^\lambda_{\fg,\fb}(\bC)$
corresponds to the locus of {\it maps} inside quasi-maps.

We have a natural projection
$\varrho^\lambda_\fb:\CZ^\lambda_{\fg,\fb}(\bC)\to \bC^\lambda$, and we set
$\fF^\lambda_{\fg,\fb}:=(\varrho^{\lambda}_\fb)^{-1}(\lambda\cdot 0_{\bC})$ (resp.,
$\overset{\circ}{\fF}{}^\lambda_{\fg,\fb}:=\fF^\lambda_{\fg,\fb}\cap
\overset{\circ}\CZ{}_\fg^\lambda(\bC)$).

It was shown in \cite{ffkm} (cf. also \cite{bgfm}) that when $\fg$ is
finite-dimensional, the scheme
$\overset{\circ}{\fF}{}^\lambda_{\fg,\fb}$ can be identified with the intersection
$\Gr_B^0\cap \Gr_{B^-}^{-\lambda}$ considered in the previous section.
Thus, the contents of the previous section amount to defining a crystal
structure on the union of the sets of irreducible components of the union of
$\overset{\circ}{\fF}{}^\lambda_{\fg,\fb}$ (over all $\lambda\in \Lambda^{pos}_\fg$).

For an arbitrary Kac-Moody algebra $\fg$, set $\sB_\fg(\lambda)$ to be the set
of irreducible components of $\overset{\circ}{\fF}{}^\lambda_{\fg,\fb}$, and
set $\sB_\fg=\underset{\lambda}\cup\, \sB_\fg(\lambda)$.
In this section we will generalize the construction of the previous
section to define on $\sB_\fg$ a structure of a crystal. However,
since in the general case we do not have the affine Grassmannian picture, and we
will have to spell out the definitions using quasi-maps' spaces.

To define the crystal structure on $\sB_\fg$ we will need to
assume Conjecture ~\ref{lagrange}. To identify it with standard crystal
of \cite{ks}, we will need to assume one more conjecture, \conjref{golova?}
(cf. \secref{golova}). Both these conjectures
are verified when $\fg$ is of affine (and, of course, finite) type,
cf. \secref{finite golova} and \secref{golova!}.

\ssec{} \label{central fiber notation}

For a standard parabolic $\fp\subset \fg$, and an element $\theta\in \Lambda^{pos}_{\fg,\fp}$,
let $\fF^\theta_{\fg,\fp}$ (resp., $\wt{\fF}{}^\theta_{\fg,\fp}$, 
$\overset{\circ}{\fF}{}^\theta_{\fg,\fp}$)
denote the preimage of $\theta\cdot \bc$
under $\varrho^\theta_\fb:\CZ^\theta_{\fg,\fp}(\bC)\to \bC^\theta$ (resp., in
$\wt{\CZ}{}^\theta_{\fg,\fp}(\bC)$, $\overset{\circ}{\CZ}{}^\theta_{\fg,\fp}(\bC)$).

Assuming for a moment that $\bC$ is complete (cf. the proof of 
\propref{action on Zastava}, where we get rid of this assumption), 
recall the stacks $\sZ^\theta_{\fg,\fp}(\bC)$ (resp., 
$\wt{\sZ}{}^\theta_{\fg,\fp}(\bC)$, $\overset{\circ}{\sZ}{}^\theta_{\fg,\fp}(\bC)$), 
introduced in \secref{twisted and enhanced}. Recall also the forgetful
map $\Gr_M\to \Bun_M(\bC)$, where $\Gr_M$ is $\Gr_{M,\bc}$.
We define now three ind-schemes
$$\BS^\theta_{\fg,\fp}:=\sZ^\theta_{\fg,\fp}(\bC)
\underset{\Bun_M(\bC)\times \bC^\theta}\times (\Gr_M\times pt),\,\,\,
\wt{\BS}{}^\theta_{\fg,\fp}:=\wt{\sZ}{}^\theta_{\fg,\fp}(\bC)\underset{\Bun_M(\bC)\times \bC^\theta}\times (\Gr_M\times pt) \text{ and }$$
$$\overset{\circ}{\BS}{}^\theta_{\fg,\fp}:=\overset{\circ}{\sZ}{}^\theta_{\fg,\fp}(\bC)
\underset{\Bun_M(\bC)\times \bC^\theta}\times (\Gr_M\times pt),$$
where we are using the projection $\q_\fp:\sZ^\theta_{\fg,\fp}(\bC)\to \Bun_M(\bC)$
(and similarly for $\wt{\sZ}{}^\theta_{\fg,\fp}(\bC)$,
$\overset{\circ}{\sZ}{}^\theta_{\fg,\fp}(\bC)$) to define the fiber
product, and $pt\to \bC^\theta$ corresponds to the point $\theta\cdot \bc$.
By definition, $\fF^\theta_{\fg,\fp}\simeq
\BS^\theta_{\fg,\fp}\underset{\Gr_M}\times pt$,
where $pt\hookrightarrow \Gr_M$ is the unit point (and similarly for 
$\wt{\fF}{}^\theta_{\fg,\fp}$, $\overset{\circ}{\fF}{}^\theta_{\fg,\fp}$).

Using the map
$\varrho^\theta_M:\wt{\sZ}{}^\theta_{\fg,\fp}(\bC)\to \H_{M,\bC}$, we can rewrite
$$\wt{\BS}{}^\theta_{\fg,\fp}\simeq \wt{\sZ}{}^\theta_{\fg,\fp}(\bC)
\underset{\H_{M,\bC}}\times \Conv_M 
\text{ and } \overset{\circ}{\BS}{}^\theta_{\fg,\fp}(\bC)\simeq
\overset{\circ}{\sZ}{}^\theta_{\fg,\fp}\underset{\H_{M,\bC}}\times
\Conv_M.$$

\begin{prop}  \label{action on Zastava}
There exists a canonical action of $M((t))$ on 
$\wt{\BS}{}^\theta_{\fg,\fp}$ (resp., $\BS^\theta_{\fg,\fp}$, 
$\overset{\circ}{\BS}{}^\theta_{\fg,\fp}$) compatible with 
its action on $\Conv_M$ (resp., $\Gr_M$, $\Conv_M$).
\end{prop}

\begin{proof}

We will give a proof for $\wt{\BS}{}^\theta_{\fg,\fp}$, since the
corresponding facts for $\BS^\theta_{\fg,\fp}$ and  
$\overset{\circ}{\BS}{}^\theta_{\fg,\fp}$ are analogous (and simpler).

Let $\D_\bc$ (resp., $\D_\bc^*$) be the formal (resp., formal punctured) disc
in $\bC$ around $\bc$. Using \lemref{trivial outside graph}, the data
defining a point of $\wt{\BS}{}^\theta_{\fg,\fp}$ can be rewritten
using the formal disc $\D_\bc$ instead of the curve $\bC$ as follows:

It consists of a principal $P$-bundle $\F_P$ on $\D_\bc$, a 
principal $M$-bundle $\F_M$ on $\D_\bc$ with a trivialization
$\beta:\F_M|_{\D_\bc^*}\to \F^0_M|_{\D_\bc^*}$ on the formal punctured 
disc $\D_\bc^*$, and a collection of maps
$$\kappa^\lambdach:(\CV_\lambdach)_{\F_P}\to (\CU_\lambdach)_{\F_M},$$
satisfying the Pl{\"u}cker equations, such that for 
$\lambdach\in \Lambdach^+_{\fg,\fp}$ the composed map of line bundles
$$\CL^\lambdach_{\F_P}\to (\CV_\lambdach)_{\F_P}\to \CL^\lambdach_{\F_M}$$
has zero of order $\langle \theta,\lambdach\rangle$ at $\bc$.

(In particular, this description makes it clear that the schemes
$\wt{\BS}{}^\theta_{\fg,\fp}$ etc., do not depend on the global
curve $\bC$, but rather on the formal neighborhood $\D_\bc$ of $\bc$.)

In these terms, the action of $M((t))$ leaves the data of
$(\F_P,\F_M,\kappa^\lambdach)$ intact, and only acts on the data
of the trivialization $\beta$.

\end{proof}

For $\nu\in \Lambda^+_\fm$, recall the subscheme $\Conv^\nu_M\subset
\Conv_M$. Let us denote by
$\BS^\nu_{\fg,\fp}\subset\overset{\circ}{\BS}{}^\theta_{\fg,\fp}$ 
the preimage of $\Conv^\nu_M$ under the map $\varrho^\theta_M$. 
The resulting map $\BS^\nu_{\fg,\fp}\to \Conv_M^\nu$ will be denoted by 
$\varrho^\nu_M$. From \propref{action on Zastava} we obtain the following:

\begin{cor}  \label{locally trivial}
The map $\varrho^\nu_M:\BS^\nu_{\fg,\fp}\to \Conv_M^\nu$ is
a fibration, locally trivial in the smooth topology.
\end{cor}

The proof is immediate from the fact that the $M((t))$-action on
$\Conv_M^\nu$ is transitive.

\ssec{}  \label{crystal B}

For $\lambda_1,\lambda_2\in \Lambda_\fg$, $\lambda=\lambda_1-\lambda_2$,
and $\theta$ being the projection of $\lambda$ under $\Lambda_\fg\to
\Lambda_{\fg,\fp}$, recall the subscheme 
$(\Gr_{B(M)}^{\lambda_1}\times \Gr_{B^-(M)}^{\lambda_2})^\nu\subset \Conv_M$.

By unfolding the definitions of
$\overset{\circ}{\BS}{}^\theta_{\fg,\fp}$ and 
$\overset{\circ}{\BS}{}^\lambda_{\fg,\fp}$, we obtain an isomorphism:

\begin{equation} \label{Parabolic vs Borel Zastava local}
\overset{\circ}{\BS}{}^\theta_{\fg,\fp}
\underset{\Conv_M}\times
(\Gr^{\lambda_1}_{B(M)}\times \Gr^{\lambda_2}_{B^-(M)})\simeq
\overset{\circ}{\BS}{}^\lambda_{\fg,\fb}\underset{\Conv_T}\times
(t^{\lambda_1}\times t^{\lambda_2}).
\end{equation}

(Here we denote by $t^\lambda$ the point-scheme, identified
with the reduced sub-scheme of the corresponding connected component
of $\Gr_T$.) 

\medskip

The isomorphism \eqref{Parabolic vs Borel Zastava local}, combined with
\conjref{lagrange} which we assume, implies the following dimension estimate:

\begin{cor}
The fibers of the projection
$\varrho^\nu_M:\BS^\nu_{\fg,\fp}\to \Conv^\nu_M$ are of
dimension at most $\langle \nu,\check\rho-2\check\rho_M\rangle$.
\end{cor}

\begin{proof}

Pick
$\lambda^1=\lambda=\nu$ and $\lambda^2=0$ so that
$(\Gr^{\lambda_1}_{B(M)}\times \Gr^{\lambda_2}_{B^-(M)})\cap \Conv_M^\nu$
is non-empty, and is in fact of pure dimension
$\langle \nu, 2\check\rho_M\rangle$. Therefore, it would suffice to show
that $\overset{\circ}{\BS}{}^\theta_{\fg,\fp} \underset{\Conv_M}\times
(\Gr^{\lambda_1}_{B(M)}\times \Gr^{\lambda_2}_{B^-(M)})$
is of dimension $\leq \langle \nu,\check\rho\rangle=\langle
\lambda_1-\lambda_2,\check\rho\rangle$.

However, from \conjref{lagrange} (which we assume to hold for our $\fg$)
we obtain that
$$\on{dim}\Bigl(\overset{\circ}{\BS}{}^\theta_{\fg,\fp} \underset{\Conv_M}\times
(\Gr^{\lambda_1}_{B(M)}\times \Gr^{\lambda_2}_{B^-(M)})\Bigr)=
\on{dim}\Bigl(\overset{\circ}{\BS}{}^\lambda_{\fg,\fb}
\underset{\Conv_T}\times (t^{\lambda_1}\times t^{\lambda_2})\Bigr)=
\on{dim}\Bigl(\overset{\circ}{\fF}{}^\lambda_{\fg,\fb}\Bigr)=
\langle \lambda_1-\lambda_2,\check\rho\rangle$$
which implies the desired dimension estimate.

\end{proof}

Since the stabilizers of the $M((t))$-action on $\Conv_M^\nu$ are
connected, we obtain that irreducible components of 
$\overset{\circ}{\BS}{}^\nu_{\fg,\fp}$ are
in one-to-one correspondence with irreducible components of any fiber
of $\varrho^\nu_M:\BS^\nu_{\fg,\fp}\to \Conv^\nu_M$.
Let us denote by $\sC_{\fm}^\nu$ the set of irreducible components
of any such fiber of (the maximal possible) dimension
$\langle \nu,\check\rho-2\check\rho_M\rangle$.

Thus, our set $\sB_\fg(\lambda)$ can be identified  with the set
of irreducible components of the fibers of 
$$\overset{\circ}\BS{}^\lambda_{\fg,\fb}\to 
\Conv^\lambda_T=\underset{\lambda_1-\lambda_2=\lambda}\bigcup\,
t^{\lambda_1}\times t^{\lambda_2}.$$

Using the isomorphism \eqref{Parabolic vs Borel Zastava local}, we obtain a decomposition
\begin{equation}
\sB_\fg\simeq \underset{\nu}\cup\, \sC_{\fm}^\nu\times \sD^\nu_{\fm},
\end{equation}
where $\sD^\nu_{\fm}$ is as in the previous section.

\medskip

This already allows to introduce the operations $e_i$, $f_i$, $e^*_i$,
$f^*_i$.
Indeed, by choosing $\fm$ to be the subminimal Levi subalgebra
$\fm_i$, we let our operations act along the second factor in the
decomposition
$$\sB_\fg\simeq \underset{n}\cup\,
\sC_{\fm_i}^n\times \sD^n_{\fm_i},$$
as in \secref{convolution}.
Moreover, the properties of \eqref{decomposition} hold
due to \propref{BG}.

To be able to apply the uniqueness theorem of \cite{ks}, it remains
to do two things: to check that the operations $e_i$, $f_i$ commute
with $e^*_j$, $f^*_j$ whenever $i\neq j$, and to establish the
highest weight property of $\sB_\fg$. Then the uniqueness theorem of
{\em loc. cit.} would guarantee the isomorphism $\sB_\fg\simeq\sB_\fg^\infty$.

\ssec{}

The above definition of the operations $e_i$, etc. mimics
the definition of \secref{convolution}. We will now give another
(of course, equivalent) construction in the spirit of
\secref{e_i}, which would enable us to prove the commutation relation.

Again, for a standard parabolic $\fp$, let us consider the ind-scheme
$\BS_{\fg,\fp,\fb}^\mu$ equal to $$\overset{\circ}{\BS}{}^\theta_{\fg,\fp}
\underset{\Conv_M}\times (\Gr^\mu_{B^-(M)}\times \Gr^0_{B^-(M)}),$$
where $\theta$ is the image of $\mu$ under $\Lambda_\fg\to \Lambda_{\fg,\fp}$.

In other words, $\BS_{\fg,\fp,\fb}^\mu$ classifies the data of 

\begin{itemize}

\item
A principal $P$-bundle
$\F_P$ (as usual, the induced bundle will be denoted by $\F'_M$),

\item
Regular bundle maps
$(\CV_{\check\lambda})_{\F_P}\to\L^{\check\lambda}_{\F^0_T}$, and

\item
Meromorphic maps
$\L^{\check\nu}_{\F^0_T}\to (\U_{\check\nu})_{\F'_M}$,
\end{itemize}
such that
\begin{itemize}

\item
$\L^{\check\lambda}_{\F^0_T}\to (\U_{\check\lambda})_{\F'_M}\to
(\V_{\check\lambda})_{\F_P}\to
\L^{\check\lambda}_{\F^0_T}$ are the identity maps, and

\item
the (a priori meromorphic) compositions
$(\U_{\check\lambda})_{\F'_M}\to
\L^{\check\lambda}_{\F^0_T}(\langle -\mu,\check\lambda\rangle\cdot \bc)$
are regular and surjective.

\end{itemize}

\medskip

Using \propref{action on Zastava}, we obtain that the natural
projection $$\BS_{\fg,\fp,\fb}^\mu\to (\Gr^\mu_{B^-(M)}\times \Gr^0_{B^-(M)})\to
\Gr^\mu_{B^-(M)}$$ splits as a direct product, as in \lemref{product}.

\medskip

For another element $\lambda\in \Lambda_\fg$, consider the preimage
(call it $\BS_{\fg,\fp,\fb}^{\lambda,\mu}$) of
$\Gr^{\lambda}_{B(M)}\cap \Gr^\mu_{B^-(M)}\subset \Gr^\mu_{B^-(M)}$  
in $\BS_{\fg,\fp,\fb}^{\mu}$. Using \eqref{Parabolic vs Borel Zastava local},
we conclude that $\BS_{\fg,\fp,\fb}^{\lambda,\mu}$
is isomorphic to a locally closed subset inside
$\overset{\circ}\BS{}^\lambda_{\fg,\fb}\underset{\Conv_T}\times
(t^\lambda\times t^0)$. 

Therefore, as in \secref{e_i}, we obtain that the fibers of
$\BS_{\fg,\fp,\fb}^\mu\to \Gr^{\mu}_{B^-(M)}$ are of dimensions
$\leq \langle \mu,\check\rho\rangle$, and if we denote by
$\sB^{\fm,*}_{\fg}(\mu)$ the set of irreducible components of
any such fiber, we obtain a decomposition
\begin{equation} \label{Zastava decomposition}
\sB_\fg(\lambda)=\underset{\mu}\cup\,
\sB^{\fm,*}_{\fg}(\mu)\times \sB_{\fm}(\lambda-\mu).
\end{equation}
Moreover, by unraveling our definition of the $e_i$'s and $f_i$'s
we obtain that they act along the second multiple of the above
decomposition when we choose $\fp=\fp_i$.

\ssec{}

We will now perform a similar procedure for the $e^*_i$ and $f^*_i$.
Let $\BS_{\fg,\fb,\fp}^\mu$ be the ind-scheme equal to
$$\overset{\circ}\BS{}_{\fg,\fp}^\theta\underset{\Conv_M}\times
(\Gr^0_{B(M)}\times \Gr^\mu_{B(M)}).$$

In other words, $\BS_{\fg,\fb,\fp}^\mu$ classifies the data of

\begin{itemize}

\item
A $B$-bundle $\F_B$, such that the induced $T$-bundle $\F_T$ is
trivial,

\item
An $M$-bundle $\F_M$,

\item
Regular bundle maps
$(\V_\lambdach)_{\F_B}\to (\CU_\lambdach)_{\F_M}$, and

\item
Meromorphic maps
$(\U_{\check\nu})_{\F_M}\to \L^{\check\nu}_{\F^0_T},$
\end{itemize}
such that

\begin{itemize}

\item
The compositions
$\L^{\check\lambda}_{\F^0_T}\to (\V_{\check\lambda})_{\F_B}\to
(\U_{\check\lambda})_{\F_M}\to
\L^{\check\lambda}_{\F^0_T}$ are the identity maps, and

\item
The $\kappa^{\lambdach,-}_\fp$'s induce regular bundle maps
$\L^{\check\lambda}_{\F^0_T}\to
(\U_{\check\lambda})_{\F_M}(\langle \mu,\check\lambda\rangle\cdot \bc)$.

\end{itemize}

As above, we have a projection
$\BS_{\fg,\fb,\fp}^\mu\to \Gr^\mu_{B(M)}$, which
splits as a direct product and defines a decomposition

\begin{equation}  \label{Zastava star}
\sB_\fg(\lambda)=\underset{\mu}\cup\,
\sB^{\fm}_{\fg}(\mu)\times \sB_{\fm}(\lambda-\mu),
\end{equation}
where $\sB^{\fm}_{\fg}(\mu)$ is the set of irreducible components
of dimension $\langle \mu,\check\rho\rangle$
of any fiber of $\varrho_{\fb,M}^\mu$.

And as before, the operations $e^*_i$, $f^*_i$ introduced earlier
coincide with those defined in terms of the above decomposition for $M=M_i$.

\ssec{}

At last, we are ready to check the required commutation property of
$e_i$ and $f_i$ with $e^*_j$ and $f^*_j$ for $i\neq j\in I$.

Consider the scheme $\BS_{\fg,i,j}^{\mu_1,\mu_2}$ classifying
the data of

\begin{itemize}

\item
A $P_i$-bundle $\F_{P_i}$ (with the induced
$M_i$-bundle $\F'_{M_i}$),

\item
An $M_j$-bundle $\F_{M_j}$,

\item
Regular maps $(\CV_\nuch)_{\F_{P_i}}\to (\CU_\nuch)_{\F_{M_j}}$,

\item
Meromorphic maps $:\L^{\check\nu}_{\F^0_T}\to
(\U_{\check\nu})_{\F'_{M_i}}$, and

\item
Meromorphic maps $(\U_{\check\nu})_{\F_{M_j}}\to
\L^{\check\nu}_{\F^0_T}$,
\end{itemize}
such that

\begin{itemize}

\item
The compositions
$\L^{\check\lambda}_{\F^0_T}\to (\U_{\check\lambda})_{\F'_{M_i}}\to
(\V_{\check\lambda})_{\F_{P_i}}\to
(\U^{\check\lambda})_{\F_{M_j}}\to
\L^{\check\lambda}_{\F^0_T}$ are the identity maps,

\item
The (a priori meromorphic) compositions
$$(\U_{\check\lambda})_{\F'_{M_i}}\to (\V_{\check\lambda})_{\F_{P_i}}\to
(\CU_\lambdach)_{\F_{M_j}}\to
\L^{\check\lambda}_{\F^0_T}(-\langle \mu_1,\check\lambda\rangle\cdot \bc)$$
are regular bundle maps, and

\item
The (a priori meromorphic) compositions
$$\L^{\check\lambda}_{\F^0_T}(-\langle \mu_2,\check\lambda\rangle\cdot \bc)\to
(\U_{\check\lambda})_{\F'_{M_i}}\to (\V_{\check\lambda})_{\F_{P_i}}\to
(\U_{\check\lambda})_{\F_{M_j}}$$
are regular bundle maps as well.

\end{itemize}

As before, we have a map
$$\BS_{\fg,i,j}^{\mu_1,\mu_2}\to
\Gr^{\mu_1}_{B^-(M_i)}\times \Gr^{\mu_2}_{B(M_j)}.$$
We claim that $\BS_{\fg,i,j}^{\mu_1,\mu_2}$ over
$\Gr^{\mu_1}_{B^-(M_i)}\times \Gr^{\mu_2}_{B(M_j)}$
also splits as a direct product. This follows from the fact that
we can make the group $(N^-(M_i)\times N(M_j))((t))$ act on
$\BS_{\fg,i,j}^{\mu_1,\mu_2}$ lifting its action on
$\Gr^{\mu_1}_{B^-(M_i)}\times \Gr^{\mu_2}_{B(M_j)}$.

\medskip

The isomorphism class of $\BS_{i,j}^{\fg,\mu_1,\mu_2}$ also depends only
on the difference $\mu=\mu_1-\mu_2$, and let us denote by
$\sC_{i,j}(\mu)$ the set of irreducible components
of the top dimension of any fiber of $\BS_{i,j}^{\fg,\mu_1,\mu_2}$ over
$\Gr^{\mu_1}_{B^-(M_i)}\times \Gr^{\mu_2}_{B(M_j)}$.

As before, the preimage of
$(\Gr^{\lambda_1}_{B(M_i)}\cap \Gr^{\mu_1}_{B^-(M_i)})
\times (\Gr^{\mu_2}_{B(M_j)}\times \Gr^{\lambda_2}_{B^-(M_j)})$
in $\BS_{i,j}^{\fg,\mu_1,\mu_2}$ is naturally a locally closed subset inside
$\overset{\circ}\BS{}^{\lambda_1-\lambda_2}_{\fg,\fb}$ and we obtain a decomposition
\begin{equation} \label{Zastava ij}
\sB_\fg(\lambda)\simeq \underset{\mu_1,\mu_2}\cup\,
\sB_{\fm_i}(\lambda_1-\mu_1)\times \sB_{\fm_j}(\mu_2-\lambda_2)\times
\sC_{i,j}(\mu_1-\mu_2).
\end{equation}
Moreover, this decomposition refines those of \eqref{Zastava decomposition}
and \eqref{Zastava star}.
Hence, the $e_i,f_i, e^*_j,f^*_j$ operations preserve
the decomposition of \eqref{Zastava ij},
and $e_i,f_i$ act along the first factor,
whereas $e^*_j,f^*_j$ act only along the second factor,
and hence, they commute.

\section{The highest weight property} \label{golova}

In this section (with the exception of \secref{finite golova} and
\secref{golova!}), $\fg$ will still be a general Kac-Moody algebra.
We will reduce \conjref{lagrange} to \conjref{golova?}, and verify
the latter in the finite and affine cases.

\ssec{}

Observe that for any $i\in I$, using \eqref{strata of quasi-maps},
by adding a defect at $\bc\in \bC$, we can realize
$\CZ^{\lambda-\alpha_i}_{\fg,\fb}(\bC)$ as a closed subscheme of
$\CZ^\lambda_{\fg,\fb}(\bC)$ of codimension 2.
Similarly, $\fF_{\fg,\fb}^{\lambda-\alpha_i}$ is a closed
subscheme of $\fF^\lambda_{\fg,\fb}$.

\medskip

Let $\bK$ be an irreducible component of
$\overset{\circ}{\fF}{}^\lambda_{\fg,\fb}$, and let $\ol{\bK}$ be its closure
in $\fF^\lambda_{\fg,\fb}$. 
Consider the intersection
$$\ol{\bK}\cap \overset{\circ}{\fF}{}^{\lambda-\alpha_i}_{\fg,\fb}\subset
\fF^\lambda_{\fg,\fb}.$$

\begin{prop}  \label{highest weight of our crystal}
If the above intersection is non-empty, it consists of one irreducible component,
whose dimension is $\dim(\bK)-1$.
\end{prop}

\begin{proof} (of \propref{highest weight of our crystal}.)

Recall the stack $\BS_{\fg,\fb,\fp_i}^{\mu}$ introduced in the
previous section.
Let $\BS_{\fg,\fb,\fp_i}^{\mu,\leq\lambda}$ and $\BS_{\fg,\fb,\fp_i}^{\mu,\lambda}$
be its closed (resp., locally closed) subschemes equal to
$$\overset{\circ}\BS{}_{\fg,\fp_i}^\theta\underset{\Conv_{M_i}}\times
\left(\Gr^0_{B(M_i)}\times (\Gr^\mu_{B(M_i)}\cap\ol{\Gr}^\lambda_{B^-(M_i)})\right),\, 
\overset{\circ}\BS{}_{\fg,\fp_i}^\theta\underset{\Conv_{M_i}}\times
\left(\Gr^0_{B(M_i)}\times (\Gr^\mu_{B(M_i)}\cap\Gr^\lambda_{B^-(M_i)})\right),$$
respectively,
where $\ol{\Gr}^\lambda_{B^-(M_i)}$ is the closure of
$\Gr^\lambda_{B^-(M_i)}$ in $\Gr_{M_i}$. 
(In other words, we impose the condition that the maps
$(\U^\nu)_{\F_{M_i}}\to \L_{\F^0_T}(\langle \lambda,\nu\rangle\cdot \bc)$
are regular (resp., regular bundle maps).) 

We have
$\BS_{\fg,\fb,\fp_i}^{\mu,\leq\lambda}=
\underset{0\leq\lambda'\leq\lambda}\bigcup\,\BS_{\fg,\fb,\fp_i}^{\mu,\lambda'}$, 
and $\BS_{\fg,\fb,\fp_i}^{\mu,\leq\lambda}$ identifies naturally with
a locally closed subset of $\fF^\lambda_{\fg,\fb}$ (of the same
dimension), and 
$$\fF^\lambda_{\fg,\fb}\simeq \underset{\mu}\bigcup\,\,
\BS_{\fg,\fb,\fp_i}^{\mu,\leq \lambda},\,\,\, 
\overset{\circ}{\fF}{}^\lambda_{\fg,\fb}\simeq \underset{\mu}\bigcup\,\,
\BS_{\fg,\fb,\fp_i}^{\mu,\lambda}.$$

We have a
pair of Cartesian squares:
$$
\CD
\BS{}_{\fg,\fb,\fp_i}^{\mu,\lambda} @>>>
\BS_{\fg,\fb,\fp_i}^{\mu,\leq\lambda} @<<<
\BS{}_{\fg,\fb,\fp_i}^{\mu,\lambda-\alpha_i} \\
@VVV   @VVV  @VVV   \\
\overset{\circ}{\fF}{}^\lambda_{\fg,\fb} @>>> 
\fF^\lambda_{\fg,\fb}  @ <{\ol{\iota}_{\alpha_i}}<< 
\overset{\circ}{\fF}{}^{\lambda-\alpha_i}_{\fg,\fb}.
\endCD
$$

Therefore, it is enough to show that for any irreducible
component $\bK$ of $\BS_{\fg,\fb,\fp_i}^{\mu,\lambda}$, the intersection
$$\ol{\bK}\cap \BS_{\fg,\fb,\fp_i}^{\mu,\lambda-\alpha_i}$$
consists of one irreducible component, whose dimension is $\on{dim}(\bK)-1$.

However, since $\BS_{\fg,\fb,\fp_i}^{\mu}\to \Gr^\mu_{B(M_i)}$ splits 
as a product, our assertion follows from the fact that 
$\Gr_{B(M_i)}^\mu\cap \Gr^\lambda_{B^-(M_i)}$ and
$\Gr_{B(M_i)}^\mu\cap \Gr^{\lambda-\alpha_i}_{B^-(M_i)}$ are irreducible
varieties, since $M_i$ is a reductive group of rank $1$.

\end{proof}

In order to be able to apply the uniqueness theorem
of \cite{ks}, we need the following conjecture to be satisfied
for our $\fg$:

\begin{conj}  \label{golova?}
If $\lambda\neq 0$, for every irreducible component
$\bK$ of $\overset{\circ}\fF{}^\lambda_{\fg,\fb}$, the intersection
$\ol{\bK}\cap \overset{\circ}{\fF}{}^{\lambda-\alpha_i}_{\fg,\fb}\subset
\fF^\lambda_{\fg,\fb}$ is non-empty for at least one $i\in I$.
\end{conj}

Note that \conjref{golova?} implies \conjref{lagrange}:

\begin{proof}
We will argue by induction on the length of $\lambda$.
The dimension estimate is obvious for $\lambda=0$, and let
us suppose that it is verified for all $\lambda'<\lambda$.

Let $\bK$ be an irreducible component of
$\overset{\circ}\fF{}^\lambda_{\fg,\fb}$, and let $i\in I$ be such that
$\bK':=\ol{\bK}\cap \overset{\circ}{\fF}{}^{\lambda-\alpha_i}_{\fg,\fb}$
is non-empty.

We know that $\on{dim}(\bK')\leq \langle \lambda-\alpha_i,\check\rho\rangle$.
Hence, by \propref{highest weight of our crystal},
$\on{dim}(\bK)\leq \langle \lambda,\check\rho\rangle$, which is what
we had to prove.
\end{proof}

Note also that if we assume \conjref{golova?}, and thus obtain
a well-defined crystal structure on $\sB_\fg$, the operation
$$\bK\mapsto
\ol{\bK}\cap\overset{\circ}{\fF}{}^{\lambda-\alpha_i}_{\fg,\fb},$$
as in \propref{highest weight of our crystal},
viewed as a map $\sB_\fg(\lambda)\to \sB_\fg(\lambda-\alpha_i)\cup 0$,
equals $f^*_i$.

Hence, if \conjref{golova?} is verified, we obtain, according to
\secref{highest weight}, that our crystal $\sB_\fg$ can be identified
with Kashiwara's crystal $\sB^\infty_\fg$.

In particular, we obtain the following corollary:

\begin{cor} \label{the set of components}
If \conjref{golova?} is verified, then for
a parabolic subalgebra $\fp\subset \fg$ such that
its Levi is of finite type, the dimension of the fibers of
the map $\BS^\nu_{\fg,\fp}\to \Conv^\nu_M$ is
$\leq \langle\nu,\check\rho-2\check\rho_M\rangle$, and the set
of irreducible components of dimension $\langle\nu,\check\rho-2\check\rho_M\rangle$
of (any) such fiber parametrizes a basis
of the space $\on{\Hom}_\fm(V^\mu_M,U(n(\check\fp)))$.
\end{cor}

\ssec{} \label{finite golova}

Although \conjref{golova?} is well-known for $\fg$ of finite type,
we will give here yet another proof, using \thmref{fd Cartier}.
Then we will modify this argument to prove \conjref{golova?} in the
affine case.

We will argue by contradiction. Let $\lambda$ be a minimal
element in $\Lambda^{pos}_\fg$, for which there exists a component
$\bK$ such that all
$\ol{\bK}\cap\overset{\circ}{\fF}{}^{\lambda-\alpha_i}_{\fg,\fb}$
are empty. As in the proof of \conjref{lagrange} above,
we obtain in particular that for all $\lambda'<\lambda$,
the dimension estimate
$\on{dim}(\fF^{\lambda'}_{\fg,\fb})=|\lambda'|$ holds.

\medskip

Let $\bC$ be a genus $0$ curve. According to \cite{bgfm}, we have a smooth map
$\CZ^\lambda_{\fg,\fb}(\bC)\to \ol{\Bun}_B(\bC)$, and let
$\partial(\CZ^\lambda_{\fg,\fb}(\bC))$ be the preimage of $\partial(\ol{\Bun}_B(\bC))$
in $\CZ^\lambda_{\fg,\fb}(\bC)$. According to \thmref{fd Cartier},
$\partial(\CZ^\lambda_{\fg,\fb}(\bC))$ is a quasi-effective Cartier divisor in
$\CZ^\lambda_{\fg,\fb}(\bC)$.
Let us consider the intersection $\ol{\bK}\cap \partial(\CZ^\lambda_{\fg,\fb}(\bC))$.
This intersection is non-empty, since it contains the
point-scheme, corresponding to $\CZ^0_{\fg,\fb}(\bC)\simeq
\fF^0_{\fg,\fb}\subset \CZ^\lambda_{\fg,\fb}(\bC)$.

On the one hand, from \lemref{dimension drop}, we obtain that
$\on{dim}(\ol{\bK}\cap \partial(\CZ^\lambda_\fg)) \geq
\on{dim}(\bK)-1\geq |\lambda|-1$. On the other hand, this intersection is
contained in
$$\underset{\lambda'<\lambda}\cap\, \ol{\bK}\cap \fF^{\lambda'}_{\fg,\fb}.$$
But, according to our assumption, all the $\lambda'$ appearing in the
above expression have the property that $|\lambda'|\leq |\lambda|-2$,
and hence
$\on{dim}(\ol{\bK}\cap \partial(\CZ^\lambda_{\fg,\fb}(\bC)))$ is of dimension
$\leq |\lambda|-2$, which is a contradiction.

\ssec{}  \label{golova!}

In this subsection we will prove \conjref{golova?} for affine Lie
algebras. We will change our notation, and from now on $\fg$ will
denote a finite-dimensional semi-simple Lie algebra, and $\fg_{aff}$
will be the corresponding affine Kac-Moody algebra.

As above, we will argue by contradiction. Let $\mu$ be a minimal
element of $\widehat\Lambda{}^{pos}_\fg$, for which there exists
a component $\bK$ violating \conjref{golova?}. In particular,
we can assume that the dimension estimate
$\on{dim}(\overset{\circ}{\fF}{}^{\lambda'}_{\fg_{aff},\fb^+_{aff}})=|\lambda'|$ holds for all
$\lambda'<\lambda$.

\medskip

Let us consider the corresponding Uhlenbeck space
$\fU^\lambda_{G;B}$. We can view our component $\bK$ as
a subset in $\Bun^\lambda_{G;B}(\bS,\bD_\infty;\bD_0)$,
and let $\ol{\ol{\bK}}$ be the closure of $\bK$ inside
$\fU^\lambda_{G;B}$.

Consider the intersection
$\ol{\ol{\bK}}\cap \partial(\fU^\lambda_{G;B})$.
As before, it is non-empty, since it contains the
point $a\cdot \bc\in \fU^{\mu;\mu,0}_{G,B}\subset
\fU^\lambda_{G;B}$, and hence is of dimension at least
$|\lambda|-1$. On the other hand we will show that if
$\bK$ violates \conjref{golova?}, this intersection will
be of dimension $\leq |\lambda|-2$.

\medskip

Indeed, consider the fiber of $\varrho^\lambda_{\fb^+_{aff}}:
\fU^\lambda_{G;B}\to \overset{\circ}\bC{}^\lambda$ over
$\lambda\cdot \bc\in \overset{\circ}\bC{}^\lambda$.
According to \thmref{stratification of parabolic Uhlenbeck},
it is the union over decompositions
$\lambda=\lambda_1+\lambda_2+b\cdot \delta$ of varieties
$$\overset{\circ}{\fF}{}^{\lambda_1}_{\fg_{aff},\fb^+_{aff}}\times 
\Sym^b(\overset{\circ}\bX-0_\bX),$$
and $\partial(\overset{\circ}\fU{}^\lambda_{G;B})$ corresponds to the locus
where $\lambda_2\neq 0$.
Since $|\delta|$ is the dual Coxeter number, and is $>1$, the intersection
\begin{equation} \label{intersection}
\ol{\ol{\bK}}\cap  
(\overset{\circ}{\fF}{}^{\lambda_1}_{\fg_{aff},\fb^+_{aff}}\times
\Sym^b(\overset{\circ}\bX-0_\bX))
\end{equation}
is of dimension $\leq |\lambda|-2$ unless $b=0$ and $|\lambda_1|=|\lambda|-1$.

However, under the projection $\fU^\mu_{G;B}\to \CZ^\mu_{\fg_{aff},\fb^+_{aff}}(\bC)$,
the intersection \eqref{intersection} gets mapped to
$\ol{\bK}\cap \overset{\circ}{\fF}{}^{\lambda_1}_{\fg_{aff},\fb^+_{aff}}$, and by 
assumption the
latter is empty for $\lambda_1=\lambda-\alpha_i$ for all $i\in I_{aff}$.

This completes the proof of \conjref{golova?} for untwisted affine Lie
algebras. The proof for a twisted affine Lie algebra follows by realizing
it (as well as all the related moduli spaces) as the fixed-point set of a
finite order automorphism of the corresponding untwisted affine Lie algebra.
We will not use the result for the twisted affine Lie algebras, and we leave
the details to the interested reader.

\bigskip

\centerline{{\bf Part V}: {\Large Computation of the IC sheaf.}}

\bigskip

\section{IC stalks: statements}

\ssec{}
The purpose of this section is to formulate some statements about the
behavior of the IC-sheaves on the parabolic Uhlenbeck spaces
$\tU^\theta_{G,P}$ and $\fU^\theta_{G,P}$. Since the material here
is largely parallel to that of \cite{bgfm}, we will assume
that the reader is familiar with the notation of {\it loc.cit.}

\medskip

Recall that we have a natural map
$$\fr_{\fp_{aff}^+}:\wt{\fU}{}^\theta_{G;P}\to
\fU^\theta_{G;P}$$
and the decompositions into locally closed subschemes
$$
\fU^\theta_{G;P}=\bigcup\limits_{\theta_2,b}\fU^{\theta;\theta_2,b}_{G;P}
\qquad\text{and}\qquad
\wt{\fU}^\theta_{G;P}=\bigcup\limits_{\theta_2,b}\wt{\fU}^{\theta;\theta_2,b}_{G;P}
$$
where $\wt{\fU}^{\theta;\theta_2,b}_{G;P}=
\fr_{\fp_{aff}^+}^{-1}(\fU^{\theta;\theta_2,b}_{G;P})$.
Recall also that we have the natural sections
$$\fs_{\fp^+_{aff}}:\overset{\circ}\bC{}^\theta\to \fU^\theta_{G;P}
\qquad\text{and}\qquad
\wt{\fs}_{\fp^+_{aff}}:
\on{Mod}_{M_{aff}}^{\theta,+}\to \wt{\fU}{}^\theta_{G;P}.$$

\ssec{}

For a standard parabolic $\fp\subset \fg$, consider the corresponding
parabolic $\fp^+_{aff}\subset \fg_{aff}$,
and let $\fV_\fp$ be the unipotent radical of the 
parabolic attached to it in the Langlands dual Lie algebra $\check \fg_{aff}$.
(Recall that by $\fV$ we denote $\fV_\fp$ for $\fp=\fg$.)

The corresponding Levi contains the subgroup $\check M_{aff}\simeq
\check M\times \BG_m$
(we ignore the other $\BG_m$-factor because it is central), and
consider $\fV_\fp$ as a $\check M_{aff}$-module. Note that the lattice
$\widehat\Lambda_{\fg,\fp}$ is naturally the lattice of central
characters of $\check M_{aff}$, and $\widehat\Lambda{}^{pos}_{\fg,\fp}\subset
\widehat\Lambda_{\fg,\fp}$ corresponds to those that appear in
$\on{Sym}(\fV_\fp)$.

\medskip

Let us denote by $\on{Loc}^{M_{aff}}$ the canonical equivalence of
categories between $\on{Rep}(\check M_{aff})$ and the category of spherical
perverse sheaves on the affine Grassmannian of $M_{aff}$.
More generally, for a scheme (or stack $S$) mapping to $\Bun_{M_{aff}}$,
we will denote by $\on{Loc}^{M_{aff}}_{S,\bC}$ the corresponding functor
from $\on{Rep}(\check M_{aff})$ to the category of perverse sheaves
on $S\underset{\Bun_{M_{aff}}}\times \H_{M_{aff},\bC}$,
cf. \cite{bgfm}, Theorem 1.12.

For an element $\theta\in \widehat\Lambda{}^{pos}_{\fg,\fp}$ we will denote by
$\fP(\theta)$ elements of the set of partitions of $\theta$ as a sum
$\theta=\underset{k}\Sigma\, n_k\cdot \theta_k$, $\theta_k\in
\widehat\Lambda{}^{pos}_{\fg,\fp}$; let $|\fP(\theta)|$ be the sum
$\underset{k}\Sigma\, n_k$. For a partition $\fP(\theta)$, let
$\on{Sym}^{\fP(\theta)}(\overset{\circ}\bC)$ denote the
corresponding partially symmetrized power of $\overset{\circ}{\bC}$,
i.e., 
$$\on{Sym}^{\fP(\theta)}(\overset{\circ}\bC):=\underset{k}\Pi\, 
\overset{\circ}{\bC}{}^{\theta_k}.$$
Let $\H^{\fP(\theta),+}_{M_{aff},\bC}$ be
the appropriate version of the Hecke stack over $\on{Sym}^{\fP(\theta)}(\overset{\circ}\bC)$.

Recall (cf. Sect. 4.1 of \cite{bgfm}), we have a natural
finite map
$$\wt{\on{Norm}}_{\fP(\theta)}:\H^{\fP(\theta),+}_{M_{aff},\bC}\to
\H^{\theta,+}_{M_{aff},\bC},$$ covering the map
$\on{Norm}_{\fP(\theta)}:\bC^{\fP(\theta)}\to \bC^\theta$.

As in \cite{bgfm}, Theorem 1.12, we have the functor
$\on{Loc}^{M_{aff},\fP(\theta)}_{S,\bC}$ from
$\on{Rep}(\check M_{aff})$ to the category of perverse sheaves
on $S\underset{\Bun_{M_{aff}}}\times \H^{\fP(\theta),+}_{M_{aff},\bC}$.

\begin{thm}  \label{sechenie na volne}
The object $(\wt{\fs}_{\fp^+_{aff}})^!(\on{IC}_{\wt{\fU}_{G;P}^\theta})$
is isomorphic to the direct sum over all $\fP(\theta)$ of
$$(\wt{\on{Norm}}_{\fP(\theta)})_*
\left(\on{Loc}^{M_{aff},\fP(\theta)}_{pt,\bC}(\fV_\fp)\right)[-|\fP(\theta)|],$$
where $pt\to \Bun_M(\bC)$ is the map corresponding to the trivial bundle.
\end{thm}

\ssec{}

Next we will state the theorem that describes the stalks of
$\on{IC}_{\wt{\fU}_{G;P}^\theta}$.

Let $\theta\in \widehat\Lambda{}^{pos}_{\fg,\fp}$ be decomposed as
$\theta=\theta_1+\theta_2+b\cdot \delta$, and recall that the
corresponding stratum
$\wt{\fU}^{\theta;\theta_2,b}_{G;P}$ is isomorphic to
$$
\left(\Bun_{G;P}^{\theta_1}(\bS',\bD'_\infty;\bD'_0)\underset{\Bun_M(\bC,\infty_\bC)}
\times \H_{M_{aff},\bC}^{\theta_2,+}\right)\times
\on{Sym}^b(\overset{\circ\circ}\bS).
$$

Let us fix partitions
$\fP(b)$ and $\fP(\theta_2)$ of $b$ and $\theta_2$, respectively.
Let $\overset{\circ}{\on{Sym}}{}^{\fP(\theta_2)}(\overset{\circ}\bC)\subset
\on{Sym}^{\fP(\theta_2)}(\overset{\circ}\bC)$ and
$\overset{\circ}{\on{Sym}}{}^{\fP(b)}(\overset{\circ}\bS)\subset
\on{Sym}^{\fP(b)}(\overset{\circ}\bS)$ be the open
subsets obtained by removing all the diagonals.
Let us
denote by $\wt{\fU}_{G;P}^{\theta;\fP(\theta_2),\fP(b)}$ the
corresponding locally closed subvariety in
$\wt{\fU}_{G;P}^{\theta;\theta_2,b}$. Note that it is isomorphic to
$$\left(\Bun_{G;P}^{\theta_1}(\bS',\bD'_\infty;\bD'_0)\underset{\Bun_M(\bC,\infty_\bC)}
\times \left(\H^{\fP(\theta_2),+}_{M_{aff},\bC}
\underset{\on{Sym}^{\fP(\theta_2)}(\bC)}\times
\overset{\circ}{\on{Sym}}{}^{\fP(\theta_2)}(\overset{\circ}\bC)\right)
\right)\times \overset{\circ}{\on{Sym}}{}^{\fP(b)}(\overset{\circ\circ}\bS).$$

\begin{thm}  \label{sloi na volne}
The *-restriction of $\on{IC}_{\wt{\fU}_{G;P}^\theta}$ to the stratum
$\wt{\fU}_{G;P}^{\theta;\fP(\theta_2),\fP(b)}$ is isomorphic to
$$\on{Loc}^{M_{aff},\fP(\theta_2)}_{\Bun_{G;P}^{\theta_1}
(\bS,\bD_\infty;\bD_0),\bC}
\left(\underset{i\geq 0}\oplus\, \on{Sym}^i(\fV_\fp)[2i]\right)
\boxtimes \left(\Sym^{\fP(b)}(\fV^f)\otimes
\on{IC}_{\overset{\circ}{\on{Sym}}{}^{\fP(b)}(\overset{\circ\circ}\bS)}\right).$$
\end{thm}

\ssec{}

We will now formulate theorems parallel to \thmref{sechenie na volne}
and \thmref{sloi na volne} for the schemes $\fU^\theta_{G;P}$.

Let $V$ be a representation of the group $\check M_{aff}$,
and let $\fP(\theta):\theta=\underset{k}\Sigma\, n_k \cdot \theta_k$
be a partition.
We will denote by $\ol{\on{Loc}}{}^{\fP(\theta)}(V)$ the semi-simple complex of
sheaves on $\on{Sym}^{\fP(\theta)}(\bC)$ equal to the IC-sheaf
tensored by the complex of vector spaces equal to
$$\underset{k}\bigotimes\, \left(V^f_{\theta_k}\right)^{\otimes n_k},$$
where $V^f$ denotes the kernel of the principal nilpotent
element $f\in \check \fm$ acting on $V^f$ endowed with a
principal grading, which we declare cohomological.

Note that the same procedure makes sense when $V$ is
actually a semi-simple complex of
$\check M_{aff}$-representations.

\begin{thm} \label{sechenie na cherte}
The object $(\fs_{\fp^+_{aff}})^*(\on{IC}_{\fU_{G;P}^\theta})$
is isomorphic to the direct sum over all $\fP(\theta)$ of
$$(\on{Norm}_{\fP(\theta)})_*(\ol{\on{Loc}}{}^{\fP(\theta)}(\fV_\fp))[|\fP(\theta)|].$$
\end{thm}

And finally, we have:

\begin{thm} \label{sloi na cherte}
The *-restriction of $\on{IC}_{\fU_{G;P}^\theta}$ to the stratum
$$
\fU_{G;P}^{\theta;\fP(\theta_2),\fP(b)}\simeq
\Bun_{G;P}^{\theta_1}(\bS,\bD_\infty;\bD_0)\times
\overset{\circ}{\on{Sym}}{}^{\fP(\theta_2)}(\overset{\circ}\bC)\times
\overset{\circ}{\on{Sym}}{}^{\fP(b)}(\overset{\circ\circ}\bS)
$$
is isomorphic to
$$\IC_{\Bun_{G;P}(\bS,\bD_\infty;\bD_0)}\boxtimes
\ol{\on{Loc}}{}^{\fP(\theta_2)}
\left(\underset{i\geq 0}\oplus\,\on{Sym}^i(\fV_\fp)[2i]\right)\boxtimes
\left(\Sym^{\fP(b)}(\fV^f)\otimes
\on{IC}_{\overset{\circ}{\on{Sym}}{}^{\fP(b)}(\overset{\circ\circ}\bS)}\right).
$$
\end{thm}

\begin{cor}\label{cor ic uhlenbeck}
\thmref{ic uhlenbeck} holds.
\end{cor}

To prove \corref{cor ic uhlenbeck} it is enough to note that it is a
particular case of \thmref{sloi na cherte} when $P=G$.

\section{IC stalks: proofs}

\ssec{Logic of the proof}

To prove all the 5 theorems, i.e. \thmref{sechenie na volne},
\thmref{sloi na volne}, \thmref{sechenie na cherte},
\thmref{sloi na cherte} and \thmref{ic uhlenbeck} we will proceed
by induction on the length $|\theta|$ of the parameter $\theta$.

Thus, let $\theta$ be $(\ol{\theta},a)$, and we assume that
the assertions of  \thmref{sechenie na volne},
\thmref{sloi na volne}, \thmref{sechenie na cherte} and 
\thmref{sloi na cherte} hold for all $\theta'<\theta$ 
and that the assertion of \thmref{ic uhlenbeck} holds for all $a'<a$.

As we shall see, the induction hypothesis will give a description of
the restrictions of $\on{IC}_{\fU_{G;P}^\theta}$ and
(resp., $\on{IC}_{\wt{\fU}_{G;P}^\theta}$) on almost all the strata
$\on{IC}_{\fU_{G;P}^{\theta;\fP(\theta_2),\fP(b)}}$
(resp., $\on{IC}_{\wt{\fU}_{G;P}^{\theta;\fP(\theta_2),\fP(b)}}$).

This will allow us to perform the induction step and prove
\thmref{sechenie na volne} for the parameter equal to $\theta$.
Then we will deduce \thmref{sloi na volne}, \thmref{sechenie na cherte}
and \thmref{sloi na cherte} from \thmref{sechenie na volne}.
Finally, as was noted before, \thmref{ic uhlenbeck} is a particular
case of \thmref{sloi na cherte}.

\ssec{} \label{previous strata}

Consider a decomposition $\theta=\theta_1+\theta_2+b\cdot \delta$,
and let $\wt{\fU}_{G;P}^{\theta;\fP(\theta_2),\fP(b)}$ be the corresponding 
stratum. We claim that if
$$\left((\theta_2<\theta)\bigwedge (b<a)\right) \bigvee
(|\fP(\theta_2)|>1)\bigvee (|\fP(b)|>1)$$
then the restriction of $\on{IC}_{\wt{\fU}_{G;P}^\theta}$
to this stratum is known and given by the expression of
\thmref{sloi na volne}, by the induction hypothesis.

We will first establish the required isomorphism locally
in the {\'e}tale topology on
$\wt{\fU}_{G;P}^{\theta;\fP(\theta_2),\fP(b)}$
and then argue that it holds globally. An absolutely
similar argument establishes the description of the
restriction of $\on{IC}_{\fU_{G;P}^\theta}$
to the corresponding stratum in $\fU_{G;P}^\theta$
under the same assumption on $\fP(\theta_2),\fP(b)$.

We claim that every geometric point belonging to
the stratum $\wt{\fU}_{G;P}^{\theta;\fP(\theta_2),\fP(b)}$ has
an {\'e}tale neighborhood, which is smoothly equivalent to
an {\'e}tale neighborhood of a point in
$\wt{\fU}_{G;P}^{\theta';\theta_2,b}$, where
$\theta'<\theta$.

\medskip

Let us first consider the case when $\theta_1\neq 0$.

Recall the map $\varrho^\theta_{\fp^+_{aff}}:\fU_{G;P}^\theta\to
\overset{\circ}\bC{}^\theta$ and the map
$\fr_{\fp^+_{aff}}:\wt{\fU}_{G;P}^\theta\to\fU_{G;P}^\theta$.

Using \propref{description of enhanced strata},
let us write a point $\sigma$ of $\wt{\fU}_{G;P}^{\theta;\theta_2,b}$ as a
triple
$$\left(\sigma',D^{\theta_2},(\F_{M_{aff}},\beta),(\Sigma b_i\cdot \bs_i)\right),$$
where:

\begin{itemize}

\item
$\sigma'$ is a based map of degree $\bC\to
\Bun_{G;P}(\bX,\infty_\bX;0_\bX)$. Let us denote by $\F'_G$
the corresponding principal $G$-bundle on $\bS'$, endowed
with a trivialization along $\bD'_\infty$ and a reduction to
$P$ along $\bD'_0$; let $\F'_M$ denote the
corresponding $M$-bundle on $\bD'_0\simeq \bC$.

\item 
$D^{\theta_2}$ is an element of $\overset{\circ}\bC{}^{\theta_2}$.

\item
$(\F'_M,D^{\theta_2},\F_{M_{aff}},\beta)$ is a point of
$\H^{\theta_2,+}_{M_{aff},\bC}$. 

\item
$\Sigma b_i\cdot \bs_i$ is a $0$-cycle of degree $b$ on
$\overset{\circ\circ}{\bS}$.

\end{itemize}

Assume for a moment that the support of the divisor
$\varrho^\theta_{\fp^+_{aff}}(\sigma')\in \overset{\circ}\bC{}^\theta$
is disjoint from both $\pi_{v}(\Sigma b_i\cdot \bs_i)$ and $D^{\theta_2}$.
Then, according to the factorization property
\propref{horizontal factorization}, our point of $\wt{\fU}_{G;P}^{\theta;\theta_2,b}$
is smoothly equivalent to a point in the product
$\wt{\fU}_{G;P}^{\theta-\theta_1;\theta_2,b}\times \wt{\fU}_{G;P}^{\theta_1;0,0}$,
and our assertion about smooth equivalence follows since
$\wt{\fU}_{G;P}^{\theta_1;0,0}\simeq \Bun_{G;P}^{\theta_1}(\bS;,\bD'_\infty;\bD'_0)$
and the latter is smooth.

\medskip

It remains to analyze the case when the support of
$\varrho^\theta_{\fp^+_{aff}}(\sigma')$ does intersect the supports of
$D^{\theta_2}$ and $\pi_{v}(\Sigma b_i\cdot \bs_i)$. By applying again
\propref{horizontal factorization}, we can assume that there exists a
point $\bc\in \overset{\circ}\bC$, such that
$D^{\theta_2}=\theta_2\cdot \bc$ and
$\pi_{v}(\Sigma b_i\cdot\bs_i)=b\cdot \bc$.

Consider the fiber product
$$\left(\wt{\fU}_{G;P}^\theta \times \Bun_G^{a_1}(\bS',\bD'_\infty)\right)
\underset{\bX^{(a)}\times \bX^{(a_1)}}\times
\left(\bX^{(a)}\times (\overset{\circ}\bX-0_\bX)^{(a_1)}\right)_{disj},$$
which, according to \propref{vertical factorization}, maps smoothly
to $\wt{\fU}_{G;P}^{\theta+a_1\cdot \delta}$.

We can view this morphism as a convolution of a point
$\sigma\in \wt{\fU}_{G;P}^\theta$ with a point $\sigma_1\in
\Bun_G^{a_1}(\bS',\bD'_\infty)$ such that $\varpi_v(\sigma)$ and
$\varpi_v(\sigma_1)$ have disjoint supports. We will denote the
resulting point of $\wt{\fU}_{G;P}^{\theta+a_1\cdot \delta}$ by
$\sigma\circ \sigma_1$.

It is easy to see that if $\sigma\in
\wt{\fU}_{G;P}^{\theta;\theta_2,b}$,
then $\sigma\circ \sigma_1\in \wt{\fU}_{G;P}^{\theta+a_1\cdot \delta;\theta_2,b}$;
moreover, $\sigma\circ \sigma_1$ corresponds to the triple
$$\left(\sigma'\circ \sigma_1,D^{\theta_2},(\F_{M_{aff}},\beta),(\Sigma b_i\cdot \bs_i)\right),$$
where all but the first piece of data remain unchanged.

Therefore, it is sufficient to show that there exists an integer $a_1$
large enough and $\sigma_1\in \Bun_G^{a_1}(\bS',\bD'_\infty)$, such
that the support of $\varrho^{\theta_1+a_1\cdot \delta}_{\fp^+_{aff}}(\sigma\circ \sigma_1)$
is disjoint from $\bc$. The latter means that the $G$-bundle obtained
from $\sigma'\circ \sigma_1$ on the line $\bc\times \bX$ should be
trivial and its reduction to $P$ at $\bc\times 0_\bX$ in the generic
position with respect to the trivialization at $\bc\times\infty_\bX$.

We can view $\sigma_1$ as a pair consisting of a divisor
$D\in (\overset{\circ}\bX-0_\bX)^{(a_1)}$ and a based map
$\bC\to\Gr^{BD,a_1}_{G,\bX,D}$. Then the restriction
of $\sigma\circ \sigma_1$ to $\bc\times \bX$ is obtained
from the restriction of $\sigma$ to $\bc\times \bX$ by a Hecke
transformation corresponding to the value of the above map
$\bC\to\Gr^{BD,a_1}_{G,\bX,D}$ at $\bc$.

Now, it is clear that by choosing $a_1$ to be large enough we can
always find a based map $\bC\to\Gr^{BD,a_1}_{G,\bX,D}$ which would bring
$\sigma'|_{\bc\times \bX}$ to the desired generic position.

\ssec{}  \label{previous strata, continued}

Now, let us assume that $\theta_1=0$, but
$(0<\theta_2<\theta)\bigwedge (0<b)$. Write $\theta_2$ as
$(\ol{\theta}_2,b_2)$ and let us consider separately the cases
$b_2\neq 0$ and $b_2=0$.

If $b_2\neq 0$, the assertion about the smooth equivalence follows
from \propref{vertical factorization}.

If $b_2=0$, the composition
$\wt{\fU}_{G;P}^{\theta;\theta_2,b}\hookrightarrow
\wt{\fU}_{G;P}^\theta\to \fU^a_G\overset{\varpi^a_h}\to \bX^{(a)}$
maps to $(\overset{\circ}\bX-0_\bX)^{(a)}$, and the assertion once
again follows from \propref{vertical factorization}.

\medskip

Thus, it remains to analyze the cases when either $\theta_2=\theta$, but
$|\fP(\theta_2)|>1$ or $b\cdot \delta=\theta$, but $|\fP(b)|>1$.

In the former case the assertion about the smooth equivalence follows
immediately from \propref{horizontal factorization}.
In the latter case, given a point
$$(\Sigma b_i\cdot \bs_i)\in \on{Sym}^b(\overset{\circ\circ}\bS)\simeq
\wt{\fU}_{G;P}^{b\cdot \delta;0,b},$$
which belongs to the stratum $\overset{\circ}{\on{Sym}}{}^{\fP(b)}
(\overset{\circ\circ}\bS)$ with $|\fP(b)|>1$, its projection with
respect to at least one of the projections $\varpi^b_h$ or
$\varpi^b_v$ will be a divisor supported in more than one point.
I.e., we deduce the smooth equivalence
assertion from either \propref{vertical factorization} or
\propref{horizontal factorization}.

\medskip

To summarize, we obtain that there are only two types of strata
not covered by the induction hypothesis. One is when
$\theta_2=\theta$, and the corresponding stratum is isomorphic to
$\overset{\circ}\bC$, which is contained in the image of $\fs_{\fp^+_{aff}}$.
The restriction is locally constant, because of the equivariance 
with respect to the group $\BA_1$ action on the pair $(\bS',\bD'_0)$
by ``horizontal" shifts (i.e., along the $\bC$-factor).

The other type of strata occurs only when $\theta=a\cdot \delta$;
the stratum itself is isomorphic to $\overset{\circ\circ}\bS$, and
it is contained in an open subset of $\wt{\fU}^\theta_{G;P}$ that
projects isomorphically onto an open subset of $\fU^a_G$,
via \propref{vertical factorization}.

Note that the restriction of $\on{IC}_{\wt{\fU}^\theta_{G;P}}$ to
this last stratum is automatically locally constant. Indeed, it is
enough to show that $\on{IC}_{\fU^a_G}|_{\overset{\circ}\bS}$ is
locally constant, but this follows from the fact that $\on{IC}_{\fU^a_G}$
is equivariant with respect to the group of affine-linear
transformations acting on $\bS$.

\medskip

Thus, we have established that the restriction of
$\on{IC}_{\wt{\fU}_{G;P}^\theta}$ to
$\wt{\fU}_{G;P}^{\theta;\fP(\theta_2),\fP(b)}$ locally has the
required form.

Let us show now that the isomorphism in fact holds globally.
First of all, we claim that the complex
$\on{IC}_{\wt{\fU}_{G;P}^\theta}|_{\wt{\fU}_{G;P}^{\theta;\fP(\theta_2),\fP(b)}}$
is semi-simple. This follows from the fact, that by retracing our
calculation in the category of mixed Hodge modules, we obtain that
$\on{IC}_{\wt{\fU}_{G;P}^\theta}|_{\wt{\fU}_{G;P}^{\theta;\fP(\theta_2),\fP(b)}}$
is pure. Hence, the semi-simplicity assertion follows from the
decomposition theorem.

Therefore, it remains to see that the cohomologically shifted
perverse sheaves that appear as direct summands in
\thmref{sloi na volne} have no monodromy. This is shown in the
same way as the corresponding assertion in \cite{bgfm}, Sect. 5.11.

\ssec{The induction step}

We will now perform the induction step and prove
\thmref{sechenie na volne} for the value of the parameter
equal to $\theta$. The proof essentially mimics the
argument of \cite{bgfm}, Sect. 5.

\medskip

Recall the map $\varrho^\theta_{M_{aff}}:
\wt{\fU}_{G;P}^\theta\to \on{Mod}^{\theta,+}_{M_{aff}}$,
and its section
$\wt{\fs}_{\fp^+_{aff}}:\on{Mod}^{\theta,+}_{M_{aff}}\to
\wt{\fU}_{G;P}^\theta$. Recall also
(cf. \corref{contraction of parabolic Uhlenbeck} and
\propref{description of enhanced strata})
that we have the ``new" $\BG_m$-action on $\wt{\fU}_{G;P}^\theta$
that contracts it onto the image of $\wt{\fs}_{\fp^+_{aff}}$.

Under these circumstances, we have as in \cite{bgfm}, Lemma 5.2 and
Proposition 5.3:
\begin{lem}
There is a canonical isomorphism
$$(\wt{\fs}_{\fp^+_{aff}})^!(\on{IC}_{\wt{\fU}_{G;P}^\theta})\simeq
(\varrho^\theta_{M_{aff}})_!(\on{IC}_{\wt{\fU}_{G;P}^\theta}),$$
and $(\wt{\fs}_{\fp^+_{aff}})^!(\on{IC}_{\wt{\fU}_{G;P}^\theta})$
is a semi-simple complex.
\end{lem}

\medskip

By the induction hypothesis and the factorization property,
we may assume that the restriction of
$(\wt{\fs}_{\fp^+_{aff}})^!(\on{IC}_{\wt{\fU}_{G;P}^\theta})$
to the open subset in $\on{Mod}^{\theta,+}_{M_{aff}}$ equal
to the preimage of the complement of the main diagonal
$\overset{\circ}\bC{}^\theta-\Delta(\overset{\circ}\bC)$
has the desired form. I.e.,
$$(\wt{\fs}_{\fp^+_{aff}})^!(\on{IC}_{\wt{\fU}_{G;P}^\theta})\simeq
\underset{\fP(\theta),|\fP(\theta)|>1}\bigoplus (\wt{\on{Norm}}_{\fP(\theta)})_*
\left(\on{Loc}^{M_{aff},\fP(\theta)}_{pt,\bC}(\fV_\fp)\right)[-|\fP(\theta)|]\bigoplus
\CK^\theta,$$
where $\CK^\theta$ is a semi-simple complex on
$\Gr_{M_{aff},\bC}^{\theta,+}\subset \on{Mod}^{\theta,+}_{M_{aff}}$,
and our goal is to show that
\begin{equation}  \label{identification of the direct summand}
\CK^\theta\simeq \on{Loc}^{M_{aff}}_{pt,\bC}\bigl((\fV_\fp)_\theta\bigr)[-1],
\end{equation}
where the subscript $\theta$ means that we are taking the direct
summand of $\fV_\fp$ corresponding to the central character $\theta$.

We pick a point in $\bC$ that we will call $\bc$, and consider the
fiber of $\wt{\fU}_{G;P}^\theta$ over it, which we denote by
$^\fU\fF^\theta_{G;P}$. By definition, we have a projection
$$\varrho^\theta_{M_{aff}}:{}^\fU\fF^\theta_{G;P}\to
\Gr^{\theta,+}_{M_{aff}}.$$

As in \cite{bgfm}, Sect. 5.12, it is sufficient to prove that
$(\varrho^\theta_{M_{aff}})_!(\on{IC}_{\wt{\fU}{}_{G;P}^\theta}|_{^\fU\fF^\theta_{G;P}})$
is concentrated in the (perverse) cohomomological degrees $\leq 0$ and
that its $0$-th cohomology is isomorphic to
$$\on{Loc}^{M_{aff}}\left((U(\fV_\fp))_\theta\right).$$

\ssec{}

We calculate
$(\varrho^\theta_{M_{aff}})_!(\on{IC}_{\wt{\fU}_{G;P}^\theta}|_{^\fU\fF^\theta_{G;P}})$
by intersecting $^\fU\fF^\theta_{G;P}$ with the various strata
$\wt{\fU}_{G;P}^{\theta;\fP(\theta_2),\fP(b)}$ of $\wt{\fU}_{G;P}^\theta$.

\noindent {\it Case 1: The open stratum.}
Set $\oF{}_{G;P}^\theta={}^\fU\fF_{G;P}^\theta\cap
\Bun_{G;P}^\theta(\bS,\bD_0;\bD_\infty)$. 
(Of course, $\oF{}_{G;P}^\theta\simeq \oF{}^\theta_{\fg_{aff},\fp^+_{aff}}$,
in the terminology of Part IV.)
The projection of the latter scheme onto $\Gr^{\theta,+}_{M_{aff}}$ was
studied in \secref{golova}.

For an $M_{aff}$-dominant coweight $\nu$, let $\oF{}_{G;P}^\nu$ denote the
preimage in $\oF{}_{G;P}^\theta$ of $\Gr^\nu_{M_{aff}}$, and we know from
\corref{the set of components} that $\oF{}_{G;P}^\nu$ is of dimension
$\leq \langle \nu, \check\rho_{aff}\rangle$. Hence,

$$
(\varrho^\theta_{M_{aff}})_!(\underline{\BC}_{\oF{}^\nu_{G;P}})
$$
lies in the cohomological degrees
$\leq \langle \theta, 2(\check\rho_{aff}-\check\rho_{M_{aff}})\rangle$.
Moreover, by \corref{the set of components}
its top (i.e., $\langle \theta, 2(\check\rho_{aff}-\check\rho_{M_{aff}})\rangle$)
perverse cohomology is isomorphic to
$$\on{IC}_{\Gr^\nu_{M_{aff}}}\otimes 
\Hom_{\check M_{aff}}(V^\nu_{M_{aff}},U(\fn(\check\fp)).$$

Since, $\Bun_{G;P}^\theta(\bS',\bD'_\infty;\bD'_0)$ is smooth of dimension
$\langle \theta, 2(\check\rho_{aff}-\check\rho_{M_{aff}})\rangle$, we obtain that
$(\varrho^\theta_{M_{aff}})_!(\on{IC}_{\wt{\fU}{}^\theta_{G;P}}|_{\oF{}^\theta_{G;P}})$
indeed lies in the cohomological degrees $\leq 0$, and its $0$-th
perverse cohomology is $\on{Loc}(U(\fn(\check\fp))_{\theta})$.

\medskip

Next, we will show that all other strata
do not contribute to the cohomological
degrees $\geq 0$.

\medskip

\noindent {\it Case 2: The intermediate strata, when
$\theta_2\neq \theta$, $b\cdot \delta\neq \theta$.}

Note that if the intersection of
$\wt{\fU}^{\theta;\fP(\theta_2),\fP(b)}_{G;P}$ with
$^\fU\fF_{G;P}^\theta$ is non-empty, then necessarily
$|\fP(\theta_2)|=1$. In this case, we will denote this intersection
by $\fF^{\theta;\theta_2,\fP(b)}_{G;P}$. As a scheme it
is isomorphic to
$$(\oF{}^{\theta_1}_{G;P}\underset{\Bun_M(\bC)}\times
\H_{M_{aff}}^{\theta_2,+})\times
\overset{\circ}{\on{Sym}}{}^{\fP(b)}(\overset{\circ}\bX-0_\bX).$$
Note also that
$$\oF{}^{\theta_1}_{G;P}\underset{\Bun_M}\times
\H_{M_{aff}}^{\theta_2,+}\simeq
\oF{}^{\theta_1}_{G;P}\underset{\Gr^{\theta_1,+}_{M_{aff}}}\times
(\Gr^{\theta_1,+}_{M_{aff}}\star \Gr^{\theta_2,+}_{M_{aff}}),$$
where
$\Gr^{\theta_1,+}_{M_{aff}}\star \Gr^{\theta_2,+}_{M_{aff}}
\subset \on{Conv}_{M_{aff}}$ is the
corresponding subscheme in the convolution diagram, cf. \secref{convolution}.
We shall view $\fF^{\theta;\theta_2,\fP(b)}_{G;P}$ as a fibration
over $\oF{}^{\theta_1}_{G;P}\times \overset{\circ}{\on{Sym}}{}^{\fP(b)}(\overset{\circ}\bX-0_\bX)$
with the typical fiber $\Gr^{\theta_2,+}_{M_{aff}}$.

In terms of the above identifications, the projection
$\varrho^\theta_{M_{aff}}:{}^\fU\fF_{G;P}^\theta\to \Gr^{\theta,+}_M$
is equal to the composition
$$
\begin{aligned}
\left(\oF{}^{\theta_1}_{G;P}\underset{\Gr^{\theta_1,+}_{M_{aff}}}\times
(\Gr^{\theta_1,+}_{M_{aff}}\star \Gr^{\theta_2,+}_{M_{aff}})\right)\times
\overset{\circ}{\on{Sym}}{}^{\fP(b)}(\overset{\circ}\bX-0_\bX)
\to\\
\Gr^{\theta_1,+}_{M_{aff}}\star \Gr^{\theta_2,+}_{M_{aff}}\to
\Gr_{M_{aff}}^{\theta_1+\theta_2,+}
\to\Gr^{\theta,+}_{M_{aff}},
\end{aligned}
$$
where the latter map is induced by the central cocharacter $b\cdot\delta$.

\medskip

Now, using the assumption that
$\theta_2\neq \theta$ and $b\cdot \delta\neq \theta$, we can use
\secref{previous strata} to write down the restriction of
$\on{IC}_{\wt{\fU}^\theta_{G;P}}$ to $\fF^{\theta;\theta_2,\fP(b)}_{G;P}$.

We obtain that it is equal to the external product of the complex
$$\left(\underline{\BC}_{\oF{}_{G;P}^{\theta_1}}
[\langle \theta_1,2(\check\rho_{aff}-\check\rho_{M_{aff}})\rangle]\right)\boxtimes
\left(\on{Sym}^{\fP(b)}(\fV^f)[2|\fP(b)|]\otimes
\underline{\BC}_{\overset{\circ}{\on{Sym}}{}^{\fP(b)}(\overset{\circ}\bX-0_\bX)}\right)$$
along the base
$\oF{}_{G;P}^{\theta_1}\times
\overset{\circ}{\on{Sym}}{}^{\fP(b)}(\overset{\circ}\bX-0_\bX)$,
and the perverse complex
$$\on{Loc}^{M_{aff}}\left(\underset{i\geq 0}\oplus\,
(\on{Sym}^i(\fV_\fp))_{\theta_2}[2i]\right)  $$
along the fiber $\Gr^{\theta_2,+}_{M_{aff}}$.

\medskip

By \corref{the set of components}, we obtain that the
!-direct image of the restriction
$\on{IC}_{\wt{\fU}^\theta_{G;P}}|_{\fF^{\theta;\theta_2,\fP(b)}_{G;P}}$
onto $\Gr^{\theta_1,+}_{M_{aff}}\star \Gr^{\theta_2,+}_{M_{aff}}$
is a complex of sheaves lying in
strictly negative cohomological degrees with spherical perverse cohomology.
By the exactness of convolution, its further direct image onto
$\Gr^{\theta,+}_{M_{aff}}$ also lies in strictly negative cohomological degrees, which
is what we had to show.

\medskip

\noindent {\it Case 3: The strata with $\theta=b\cdot \delta$.}

The intersection $\wt{\fU}^{b\cdot \delta;0,b}_{G;P}\cap {}^\fU\fF^\theta_{G;P}$
is isomorphic to $(\overset{\circ}\bX-0_\bX)^{(b)}$, which we further
subdivide according to partitions $\fP(b)$ of $b$.
The map
$\wt{\fU}^{b\cdot \delta;0,b}_{G;P}\cap {}^\fU\fF^\theta_{G;P}\to
\Gr^{\theta,+}_{M_{aff}}$ is the composition of the projection
$\wt{\fU}^{b\cdot \delta;0,b}_{G;P}\cap {}^\fU\fF^\theta_{G;P}\to pt$
and the embedding $pt\to \Gr^{\theta,+}_{M_{aff}}$ corresponding
to the central cocharacter $b\cdot \delta$.

\medskip

We will consider separately two cases: (a) when $|\fP(b)|>1$ and
(b) when $|\fP(b)|=1$.

In case (a), the restriction of $\on{IC}_{\wt{\fU}^\theta_{G;P}}$
to the corresponding stratum $\overset{\circ}{\on{Sym}}{}^{\fP(b)}(\overset{\circ\circ}\bS)$
is known, according to \secref{previous strata, continued}.
In particular, we know that this complex lives in the perverse
cohomological degrees $\leq -|\fP(b)|$. Hence, when we further restrict it
to $\overset{\circ}{\on{Sym}}{}^{\fP(b)}(\overset{\circ}\bX-0_\bX)\subset
\overset{\circ}{\on{Sym}}{}^{\fP(b)}(\overset{\circ\circ}\bS)$, it
lives in the perverse cohomological degrees $\leq -2|\fP(b)|<-1-|\fP(b)|$.
Therefore, when we take its direct image along
$\overset{\circ}{\on{Sym}}{}^{\fP(b)}(\overset{\circ}\bX-0_\bX)$
we obtain a complex in the cohomological degrees $<0$.

In case (b), by the definition of intersection cohomology,
the restriction of $\on{IC}_{\wt{\fU}^\theta_{G;P}}$ to
$\overset{\circ\circ}\bS$ lives in strictly negative cohomological
degrees. According to \secref{previous strata, continued}, this
restriction is locally constant; therefore its further restriction to
$\bX-0_\bX\subset \overset{\circ\circ}\bS$ lives in the cohomological
degrees $<-1$. Hence, its cohomology along $\bX-0_\bX$ lives in
the cohomological degrees $<0$.

\medskip

\noindent {\it Case 4: The stratum $\theta_2=\theta$.}

According to \propref{stratification of parabolic Uhlenbeck},
the intersection $\wt{\fU}^{\theta;\theta,0}_{G;P}\cap {}
^\fU\fF^\theta_{G;P}$ is in fact isomorphic to $\Gr_{M_{aff}}^{\theta,+}$, 
such that the map $\varrho^\theta_{M_{aff}}$ is the identity map.

Hence, the required assertion follows from the fact that the
restriction of $\on{IC}_{\wt{\fU}^\theta_{G;P}}$ to $\Gr_{M_{aff}}^{\theta,+}$
lives in the negative cohomological degrees, by the definition of
the IC sheaf.

\ssec{}

Thus, the induction step has been performed and we have established
\thmref{sechenie na volne} for the parameter equal to $\theta$.
Now \thmref{sechenie na cherte} is deduced from \thmref{sechenie na volne}
is the same way as Theorem 7.2 of \cite{bgfm} is deduced from
Theorem 4.5 of {\it loc. cit.}

\medskip

Let us prove now \thmref{sloi na cherte}. As was explained
in \secref{previous strata}, it remains to identify the restriction
of $\on{IC}_{\fU^\theta_{G;P}}$ to strata of two types.

\medskip

\noindent {\it Type 1}:
$\fU^{\theta;\fP(\theta),0}_{G;P}$, where $\fP(\theta)$ is a
$1$-element partition.

Note that $\fU^{\theta;\fP(\theta),0}_{G;P}$ is a closed subset
of $\fU^{\theta;\theta,0}_{G;P}$, and the latter is exactly the image
of the map $\fs_{\fp^+_{aff}}$. Therefore, we have to calculate
the $*$-restriction of
$$\underset{\fP(\theta)}\bigoplus\,
(\on{Norm}_{\fP(\theta)})_*(\ol{\on{Loc}}{}^{\fP(\theta)}(\fV_\fp))[|\fP(\theta)|]$$
to $\Gr^{\theta,+}_{M_{aff},\bC}\subset \on{Mod}_{M_{aff},\bC}^{\theta,+}$,
and as in the proof of Theorem 1.12 of \cite{bgfm} we obtain the
required answer for \thmref{sloi na cherte}.

\noindent {\it Type 2}:
We can assume that $P=G$ and we are dealing with the stratum
$\overset{\circ}{\bS}\subset \on{Sym}^a(\overset{\circ}{\bS})\subset \fU^a_G$.

By \thmref{sechenie na cherte} we know the restriction of
$\on{IC}_{\fU^a_G}$ to $\overset{\circ}\bC\subset \overset{\circ}{\bS}$,
and we obtain
$$\underset{i}\oplus (\on{Sym}^i(\fV^f))_a[2i]\otimes
\underline{\BC}_{\overset{\circ}\bC}[2].$$

In particular, by passing to the category of mixed Hodge modules
and retracing the proof of \thmref{sechenie na cherte}, we obtain that
$\on{IC}_{\fU^a_G}|_{\overset{\circ}{\bC}}$ is pure.

However, since $\on{IC}_{\fU^a_G}|_{\overset{\circ}{\bS}}$ is equivariant
with respect to the group of affine-linear transformations of
$\overset{\circ}{\bS}$, we obtain that $\on{IC}_{\fU^a_G}|_{\overset{\circ}{\bS}}$
is also pure, and hence semi-simple.
Therefore, it is isomorphic to
$\underset{i}\oplus (\on{Sym}^i(\fV^f))_a[2i]\otimes
\underline{\BC}_{\overset{\circ}\bS}[2]$, as required.

\medskip

Finally, let us prove \thmref{sloi na volne}. In this case there
are also two types of strata, not covered by the induction hypothesis:
$\wt{\fU}^{\theta;\fP(\theta),0}_{G;P}$ and
$\wt{\fU}^{a\cdot \delta;0,a}_{G;P}$.

The assertion for the stratum of the first type follows from
\thmref{sechenie na volne} just as the corresponding assertion for
$\fU^\theta_{G,P}$ follows from \thmref{sechenie na cherte}.
The assertion for the stratum of the second type follows from
the corresponding fact for $\fU^a_G$, which has been
established above.

\section{Appendix}  \label{AppA}

We recall the setup of Proposition ~\ref{locally of finite type}.
Assume that $\fg$ is symmetrizable. Then a stronger statement holds:

\begin{thm}  \label{tony}
The scheme $\on{Maps}^\theta(\bC,\CG_{\fg,\fp})$
is of finite type.
\end{thm}

We reproduce here the proof due to V.~Drinfeld and A.~Joseph.
Recall that according to Proposition ~\ref{locally of finite type},
$\on{Maps}^\theta(\bC,\CG_{\fg,\fp})$ is a union of finite type open
subschemes $\on{Maps}^\theta(\bC,\CG_{\fg,\fp}^w), w\in\CW$.
We have to check that
the above union is actually finite, i.e. for $w$ big enough we have
$\on{Maps}^\theta(\bC,\CG_{\fg,\fp})=
\on{Maps}^\theta(\bC,\CG_{\fg,\fp}^w)$.
In other words, we have to find $w\in\CW$
such that any map
$\sigma\in\on{Maps}^\theta(\bC,\CG_{\fg,\fp})$
lands into the open subscheme $\CG_{\fg,\fp}^w\subset\CG_{\fg,\fp}$.
We will give a proof for $\fp=\fb$, and the general case follows immediately.

More precisely, for $\mu=\underset{i\in I}\sum \, n_i\cdot \alpha_i\in\Lambda_\fg^{pos}$ we
set $|\mu|:=\underset{i\in I}\sum \, n_i$. We will prove that
the image of $\sigma(\bC)$ in $\CG_{\fg,\fb}$ never intersects Schubert subvariety
$\ol{\CG}_{\fg,\fb,w}$ for $\ell(w)>2|\mu|$ (the reduced length).

To this end note that $|\mu|$ is the intersection multiplicity of
the curve $\sigma(\bC)$ and the Schubert divisor
$\underset{i\in I}\bigcup\, \ol{\CG}_{\fg,\fb,s_i}$ in $\CG_{\fg,\fb}$.
If $\sigma(\bC)$ passes through a point of $\ol{\CG}_{\fg,\fb,w}$, this
intersection multiplicity cannot be smaller than the multiplicity of
the Schubert divisor at this point. Thus it suffices to prove that the
multiplicity of the Schubert divisor
$\underset{i\in I}\bigcup\, \ol{\CG}_{\fg,\fb,s_i}$ at a point
of the Schubert variety $\ol{\CG}_{\fg,\fb,w}$ is at least $\frac{1}{2}\ell(w)$.

\medskip

We choose $\rhoch\in\fh^*$ such that
$\langle\alpha_i,\rhoch\rangle = 1, \ \forall\ i\in I$.
Let $\CV_\rhoch$ denote the simple $\fg$-module with
highest weight $\rhoch$.  For each $w\in \CW$, let $v_{w\rhoch}$  denote an
extremal  vector in $\CV_\rhoch$ of weight $w\rhoch$; let $v^*_{w\rhoch}$
be an extremal vector in $\CV^*_\rhoch$ of weight $-w\rhoch$.

Recall that under the projective embedding
$\CG_{\fg,\fb}\hookrightarrow\BP(\CV^*_\rhoch)$ the Schubert divisor
$\bigcup_{i\in I}\ol{\CG}_{\fg,\fb,s_i}$ is cut out by the equation $v_\rhoch$.
There is a transversal slice to $\ol{\CG}_{\fg,\fb,w}$
through the $T$-fixed point $w$
isomorphic to $w(\fn)\cap\fn^-$ (and the isomorphism is given by the action
of the corresponding nilpotent Lie group), and in the coordinates
$f\in w(\fn)\cap\fn^-$ the above equation reads
$$\langle v_\rhoch,\exp(f)v^*_{w\rhoch}\rangle=0$$
where $\exp$ is the isomorphism between the nilpotent Lie algebra
$w(\fn)\cap\fn^-$ and the corresponding Lie group.

Hence the multiplicity $m_w$ of the Schubert divisor at the point
$w\in\ol{\CG}_{\fg,\fb,w}$
is the maximal integer $m$ such that
$\langle f_1\ldots f_m v_\rhoch,v^*_{w\rhoch}\rangle=0$ for any
$f_1,\ldots,f_m\in w(\fn)\cap\fn^-$. Let $F^n$ denote the canonical 
filtration on the universal enveloping algebra $U(\fn^-)$.
Now the desired estimate $m_w\geq\frac{1}{2}\ell(w)$ is a consequence
of the following lemma belonging to A.Joseph. We are grateful to him for the
permission to reproduce it here.

\begin{lem}
\label{joseph}
Suppose $v_{w\rhoch}\in F^m(U(\fn^-))v_\rhoch$.
Then $2m\ge \ell(w)$.
\end{lem}

\begin{proof}
The proof uses the Chevalley-Kostant construction of $\CV_\rhoch$ as extended
in ~\cite{j} to the affine case.  The extension generalizes with no significant
change for $\fg$ as above.  The extension from the semisimple case (of
Chevalley-Kostant) needs a little care as infinite sums are involved.  A
slightly more streamlined analysis is given in ~\cite{gj}.  Some details are
given below.

Recall that $\fg$ admits a non-degenerate symmetric invariant bilinear form
( \ \ , \ \ ).  From this one may construct the Clifford algebra  $C(\fg)$
of $\fg$  defined as a quotient of the tensor algebra  $T(\fg)$  of $\fg$
by the ideal generated  by the elements $x\otimes y - y\otimes x - 2(x,y),
\ \forall\ x,y\in\fg$.

One defines a Lie algebra homomorphism $\varphi$ of $\fg$ into a subspace
of certain infinite sums of quadratic elements of  $C(\fg)$.
(See ~\cite{j}, 4.11
or ~\cite{gj}, 4.7).  In this one checks that only finitely many commutators
contribute to a given term in $[\varphi(x),\varphi(y)]$  which therefore
makes sense and by the construction equals  $\varphi[x,y]$
(cf. ~\cite{gj}, 4.7).

It turns out that we need to reorder the expression for $\varphi(x)$ so that
only finitely many of the negative root vectors lie to the right.  This
process is well-defined for root vectors $x$ of weight $\pm\alphach_i,\
i\in I$.  This is extended to $\fg$ via commutation and the Jacobi identities
(see ~\cite{j}, 4.11).  Notably
\begin{equation}
\label{zvezdochka}
\varphi(h) = \rhoch(h) - \frac{1}{2}\sum_{\alphach\in\Delta^+}\sum_r \alphach
(h) e^r_{-\alphach} e^r_\alphach,\ \forall\ h\in\fh,
\end{equation}
where $\Delta^+$ is the set of positive roots, $e^r_\alphach$ is a basis for
$\fg_\alphach$ and $e^r_{-\alphach}$  a dual basis for  $\fg_{-\alphach}$.

Through the diamond lemma one shows that $C(\fg)$ admits triangular
decomposition, that is to say there is a vector space isomorphism
$$ C(\fn^-)\otimes C(\fh)\otimes C(\fn) \ {\buildrel\sim\over\longrightarrow}
\ C(\fg) $$
given by multiplication.  (The left hand side is defined by restricting
the bilinear form: in particular  $C(\fn^-)=\Lambda(\fn^-)$).

View $\fn$ as a subspace of $C(\fg)$  and let $I(\fg)$ denote the left ideal
it generates.  Through multiplication by  $\varphi(x):x\in\fg$,  and the above
ordering, it follows that  $C(\fg)/I(\fg)$  becomes a $\fg$-module.
By ~\cite{j} 4.12 it is a direct sum of  $\dim C(\fh)$  copies of  $\CV_\rhoch$.  In
particular if $v_\rhoch$ is the image of the identity of  $C(\fg)$ in $C(\fg)
/I(\fg)$,  then  $U(\varphi(\fg))v_\rhoch$  is isomorphic to  $\CV_\rhoch$.
(However it does not lie in  $C(\fn^-)\pmod{I(\fg)}$,  only in  $C(\fn^-)\otimes
C(\fh)\pmod{I(\fg)})$.  (Note by ~(\ref{zvezdochka}) above,
$v_\rhoch$  has weight $\rhoch$).

Given $w\in \CW$, set $S(w^{-1})=\{\alphach\in\Delta^+$ s.t. $w^{-1}\alphach\in
-\Delta^+\}$.  Since $\betach\in S(w^{-1})$  is a real root there is a unique
up to scalars root vector  $e_{-\betach}\in\fn^-$  of weight $-\betach$.  Since
$C(\fg)/I(\fg)\approx \Lambda(\fn^-)\otimes C(\fh)$  as an  $\fh$-module it
follows that the subspace of $C(\fg)/I(\fg)$  of weight $w\rhoch$ is just
$$ \left(\prod_{\betach\in S(w^{-1})} e_{-\betach}\right) C(\fh). $$
\par
With respect to the canonical filtration $F^n$ of $C(\fg)$ such an
element has degree between $\ell(w)$ and $\ell(w)+\dim\fh$.  Yet $\varphi
(\fg)\subset F^2 C(\fg)$ and $v_\rhoch\in F^0 C(\fg)$.  The assertion
of the lemma follows.
\end{proof}

This completes the proof of \thmref{tony}.

\section{Erratum}

\centerline{{\it Added in February 2006}}

\bigskip

\ssec{}

We are grateful to Joel Kamnitzer for pointing out the following mistake in the paper:

\medskip

\propref{highest weight of our crystal} is wrong. 
Namely, the intersection 
\begin{equation} \label{suspicious intersection}
\ol{\bK}\cap \overset{\circ}{\fF}{}^{\lambda-\alpha_i}_{\fg,\fb}\subset
\fF^\lambda_{\fg,\fb}
\end{equation}
may be reducible. Indeed, let $\mu$ be such that $\bK\cap \BS_{\fg,\fb,\fp_i}^{\mu,\lambda}$
is dense in $\bK$. However, it might happen that the intersection \eqref{suspicious intersection}
contains an irreducible component $\bK'$, such that $\bK'\cap \BS_{\fg,\fb,\fp_i}^{\mu',\lambda}$
is dense in $\bK'$ for $\mu'\neq \mu$. (J. Kamnitzer showed us a counter-example for
$\fg={\mathsf {sl}_3}$, $\alpha=\alpha_1+2\cdot \alpha_2$ and $i=2$.)

\medskip

Let us analyze separately the following two cases:

\medskip

\noindent{(1)}
The map $\overset{\circ}\BS{}_{\fg,\fp_i}^\theta\to \Gr^\theta_{M_i}$
collapses $\bK$ to a single point. 

We claim that in this case, it is easy to 
see that the intersection \eqref{suspicious intersection} is empty. Indeed,
the image of $\bK$ is a point that lies in 
$\Gr^\mu_{B(M_i)}\cap \ol\Gr^\lambda_{B^-(M_i)}\subset \Gr_{M_i}$,
and the intersection in question is the pre-image of 
$\ol\Gr^{\lambda-\alpha_i}_{B^-(M_i)}\subset \Gr_{M_i}$.

\medskip

\noindent{(2)}
The map $\overset{\circ}\BS{}_{\fg,\fp_i}^\theta\to \Gr^\theta_{M_i}$
{\it does not} collapse $\bK$ to a single point. 

Then the image
of $\bK\cap \BS_{\fg,\fb,\fp_i}^{\mu,\lambda}$ under this map is the
(unique) irreducible component of $\Gr^\mu_{B(M_i)}\cap \Gr^\lambda_{B^-(M_i)}$,
of non-zero dimension. Hence, the intersection 
$$\ol{\bK}\cap \BS_{\fg,\fb,\fp_i}^{\mu,\lambda-\alpha_i}$$
is non-empty and consists of the single irreducible component of dimension 
$\on{dim}(\bK)-1$, equal to the preimage of the (unique) irreducible
component of $\Gr^\mu_{B(M_i)}\cap \Gr^{\lambda-\alpha_i}_{B^-(M_i)}\subset \Gr_{M_i}$.
The closure of $\ol{\bK}\cap \BS_{\fg,\fb,\fp_i}^{\mu,\lambda-\alpha_i}$
in $\overset{\circ}{\fF}{}^{\lambda-\alpha_i}_{\fg,\fb}$ equals by definition
$f_i^*(\bK)$; in particular $f_i^*(\bK)\neq 0$.

\medskip

Thus, instead of the erroneous \propref{highest weight of our crystal} we have
shown the following:

\begin{prop} \label{highest weight, corrected}
If the intersection \eqref{suspicious intersection} is non-empty, then $f_i^*(\bK)\neq 0$
and it corresponds to an irreducible component of dimension $\dim(\bK)-1$ of 
\eqref{suspicious intersection}.
\end{prop}

The implication \conjref{golova?} $\Rightarrow$ \conjref{lagrange} remains valid.

\bigskip

\centerline{{\it Added in November 2012}}

\bigskip

\ssec{}

We are grateful to Hiraku Nakajima for pointing out the following mistakes in 
Theorems~7.10,~16.7,~16.8. 

First of all, in the 4-th paragraph of 16.2, for a partition 
${\mathfrak P}(\theta):\ \theta=\sum_kn_k\cdot\theta_k$ with pairwise
distinct $\theta_k$, replace the 4-th line with 
$$(\overset{\circ}{\mathbf C})^{{\mathfrak P}(\theta)}=
\operatorname{Sym}^{{\mathfrak P}(\theta)}(\overset{\circ}{\mathbf C})=
\prod_k(\overset{\circ}{\mathbf C})^{(n_k)}.$$
Second, in the second paragraph of 16.6, replace the definition of
$\overline{\operatorname{Loc}}{}^{{\mathfrak P}(\theta)}(V)$ by the IC-extension
from the off-diagonal part of 
$\operatorname{Sym}^{{\mathfrak P}(\theta)}(\overset{\circ}{\mathbf C})$
of $\boxtimes_k(V^f_{\theta_k})^{(n_k)}$: the external product of symmetric
powers of the constant local systems on $\overset{\circ}{\mathbf C}$. 
More generally, for a smooth algebraic variety $Y$ in place of 
$\overset{\circ}{\mathbf C}$, denote by 
$\overline{\operatorname{Loc}}{}^{{\mathfrak P}(\theta)}_Y(V)$ the similarly defined
IC sheaf on the symmetric power stratum $Y^{{\mathfrak P}(\theta)}$.

Now Theorem~16.7 (along with its proof) holds true as stated (with the 
corrected definition of 
$\overline{\operatorname{Loc}}{}^{{\mathfrak P}(\theta)}({\mathfrak V}_{\mathfrak p})$).
Furthermore, in~Theorem~16.8 replace the last line with 
$$\IC_{\Bun_{G;P}(\bS,\bD_\infty;\bD_0)}\boxtimes
\ol{\on{Loc}}{}^{\fP(\theta_2)}
\left(\underset{i\geq 0}\oplus\,\on{Sym}^i(\fV_\fp)[2i]\right)\boxtimes
\ol{\on{Loc}}{}^{\fP(b)}_{\overset{\circ\circ}\bS}\left(\underset{i\geq0}\oplus\,
\on{Sym}^i(\fV)[2i]\right).
$$
The proof of~Theorem~16.8 remains valid.

Finally, replace~Theorem~7.10 with the following:
{\em  The restriction $\on{IC}_{\fU^a_G}|_{\fU^{a;\fP(b)}_G}$ is locally constant
and is isomorphic to $\on{IC}_{\on{Bun}^a_G(\bS,\bD_\infty)}$ tensored by the 
(locally constant) complex 
$\ol{\on{Loc}}{}^{\fP(b)}_{\overset{\circ\circ}\bS}\left(\underset{i\geq0}\oplus\,
\on{Sym}^i(\fV)[2i]\right)$ on 
$\operatorname{Sym}^{{\mathfrak P}(b)}(\overset{\circ}{\mathbf S})$.}
The proof of~Theorem~7.10 remains valid.


\begin{thebibliography}{99}

\bibitem[B]{b} R.~Bezrukavnikov,
{\em Perverse coherent sheaves (after Deligne)},
preprint math.AG/0005152.

\bibitem[BBD]{BBD} A.~Beilinson, J.~Bernstein and P.~Deligne,
{\em Faisceaux Pervers}, Ast{\'e}risque {\bf 100} (1982), 5--171.

\bibitem[BaGi]{BaGi} V.~Baranovsky, V.~Ginzburg,
{Algebraic construction of the Uhlenbeck moduli space}, manuscript (1998).

\bibitem[BFGM]{bgfm} A.~Braverman, M.~Finkelberg, D.~Gaitsgory, I.~Mirkovi\'c,
{\em Intersection Cohomology of Drinfeld's compactifications},
Selecta Math. (N.S.) {\bf 8} (2002), 381--418. {\em Erratum},
Selecta Math. (N.S.) {\bf 10} (2004), 429--430.
 
\bibitem[BG]{bg} A.~Braverman, D.~Gaitsgory, {\em Crystals via affine
Grassmannian}, Duke Math. J. {\bf 107} (2001), 561--575.

\bibitem[BG1]{bg1} A.~Braverman, D.~Gaitsgory, {\em Geometric Eisenstein series},
Invent. Math. {\bf 150} (2002), 287--384.

\bibitem[CB]{CB} W.~Crawley-Boevey, 
{\em Normality of Marsden-Weinstein reductions for representations of quivers},
Math. Ann. {\bf 325} (2003), no. 1, 55--79.

\bibitem[Do]{d} S.~K.~Donaldson,
{\em Connections, cohomology and the intersection forms
of four-manifolds}, Jour. Diff. Geom. {\bf 24} (1986), 275--341.

\bibitem[DK]{dk} S.~K.~Donaldson, P.~Kronheimer, {\em The Geometry of Four-Manifolds},
Oxford Mathematical Monographs (1990).

\bibitem[Fa]{Fa} G.~Faltings, {\em Algebraic Loop Groups and Moduli
Spacesof Bundles}, Jour. Eur. Math. Soc. {\bf 5} (2003), 41--68.

\bibitem[FFKM]{ffkm} B.~Feigin, M.~Finkelberg, A.~Kuznetsov, I.~Mirkovi\'c,
{\em Semi-infinite Flags},
The AMS Translations {\bf 194} (1999), 81-148.

\bibitem[FGK]{fgk} M.~Finkelberg, D.~Gaitsgory, A.~Kuznetsov,
{\em Uhlenbeck spaces for $\BA^2$ and affine Lie algebra $\hsl_n$},
Publ. RIMS, Kyoto Univ. {\bf 39} (2003), no. 4, 721--766.

\bibitem[FKMM]{fkmm} M.~Finkelberg, A.~Kuznetsov, N.~Markarian, I.~Mirkovi\'c,
{\em A note on a symplectic structure on the space of $G$-monopoles},
Comm. Math. Phys. {\bf 201} (1999), 411-421.

\bibitem[Ga]{g} D.~Gaitsgory, {\em Construction of central elements in the affine
Hecke algebra via nearby cycles}, Invent. Math. {\bf 144} (2001),
253--280.

\bibitem[Gi]{Gi} D.~Gieseker,
{\em On the moduli of vector bundles on an algebraic surface},
Ann. of Math. (2) {\bf 106} (1977), 45--60.

\bibitem[GrJ]{gj} J.~Greenstein, A.~Joseph,  
{\em A Chevalley-Kostant presentation of basic modules for $\widehat{sl(2)}$ 
and the associated affine KPRV determinants at $q=1$}, Bull. Sci. Math. {\bf 125} (2001), 85--108.

\bibitem[J]{j} A.~Joseph, {\em On an affine quantum KPRV determinant at $q=1$},
Bull. Sci. Math. {\bf 125} (2001), 23--48.

\bibitem[Ka]{k} M.~Kashiwara, {\em The flag manifold of Kac-Moody Lie algebra},
in Algebraic Analysis, Geometry, and Number Theory, Proceedings of the
JAMI Inaugural Conference, The Johns Hopkins University Press (1989),
161--190.

\bibitem[Ka1]{k1} M.~Kashiwara, private communication.

\bibitem[Ka2]{k2} M.~Kashiwara, {\em Kazhdan-Lusztig conjecture for a symmetrizable
Kac-Moody Lie algebra}, Progr. Math. {\bf 87} (1990), 407--433.

\bibitem[Ka3]{k3} M.~Kashiwara, {\em Crystallizing the q-analogue of
universal enveloping algebras}, Duke Math. J. {\bf 63} (1991), 465--516.

\bibitem[Ka4]{k4} M.~Kashiwara, {\em Crystal base and Littelmann's refined
Demazure character formula}, Duke Math. J. {\bf 71} (1993), 839--858.

\bibitem[Ka5]{k5} M.~Kashiwara, {\em On crystal bases},
CMS Conference Proceedings {\bf 16} (1995), 155--197.

\bibitem[KaS]{ks} M.~Kashiwara, Y.~Saito, {\em Geometric construction
of crystal bases}, Duke Math. J. {\bf 89} (1997), 9--36.

\bibitem[KiP]{kp} B.~Kim, R.~Pandharipande, {\em The connectedness of the moduli
space of maps to homogeneous spaces}, preprint math.AG/0003168.

\bibitem[Li]{Li} J.~Li, {\em Algebraic geometric interpretation of Donaldson's
polynomial invariants}, Journal of Differential Geometry {\bf 37} (1993),
417--466.

\bibitem[Mo]{Mo} J.~W.~Morgan, {\em Comparison of the Donaldson polynomial invariants with their
algebro-geometric analogues}, Topology {\bf 32} (1993), 449--488.

\bibitem[Na]{n1} H.~Nakajima, {\em Lectures on Hilbert schemes of points on
surfaces}, The AMS University Lecture Series {\bf 18} (1999).

\bibitem[Pe]{p} N.~Perrin, {\em Courbes rationelles sur les vari{\'e}t{\'e}s
homog{\`e}nes}, Ann. Inst. Fourier {\bf 52} (2002), 105--132.

\bibitem[Th]{t} J.~F.~Thomsen, {\em Irreducibility of
$\overline{M}_{0,n}(G/P,\beta)$}, Internat. J. Math. {\bf 9} (1998), 367--376.

\bibitem[U]{u} K.~K.~Uhlenbeck,
{\em Connections with $L^p$ bounds on curvature},
Comm. Math. Phys. {\bf 83} (1982), 31--42.

\bibitem[Va]{v} G.~Valli,
{\em Interpolation theory, loop groups and instantons},
J. Reine Angew. Math. {\bf 446} (1994), 137--163.


\end{thebibliography}
\end{document}